%% file: martineta.tex.vcleaned.tex
\title{Ramified Partition Algebras}\date{}
\author{P P Martin and A Elgamal%
\footnote{Sadly Ahmed died before this work was completed. 
  His memory lives on.}
\\ \myaddress}
\documentclass[12pt]{article} 
\usepackage{epic} \usepackage{eepic} \usepackage{latexsym}
\usepackage{amssymb}
\input martinew
\newcommand{\Pbase}{p}
\newcommand{\Ttwo}{\Nset{2}}

\newcommand{\trampart}{$\Nset{2}$--ramified partition}
\newcommand{\draw}{{\bf d}}
\newcommand{\proj}{{\mbox{pr}}}
\newcommand{\direct}[1]{\prod_{#1}}
\newcommand{\dip}[1]{\direct{#1}}
\newcommand{\tp}[1]{\bigotimes_{#1}}

\newcommand{\hap}{\ha^P}
\newcommand{\parts}[1]{\pi(#1)}
\newcommand{\hal}{\ha^{-}}
\newcommand{\hae}{\ha^{env}}
     \newcommand{\emptyindex}{()}
     \newcommand{\spineweight}{()}
\newcommand{\nowithshape}{g}
\newcommand{\Bell}[1]{b_{#1}}
\newcommand{\AxA}[1]{(A^{#1},A^{#1})}
\newcommand{\Ax}[1]{(A^{#1},1)}
\newcommand{\xA}[1]{(1,A^{#1})}
\newcommand{\DT}{D^{(T)}}
     \newcommand{\NB}{NB}

\newcommand{\MI}{\cite[\S 1]{Martin96}}%
\newcommand{\QQ}{{\bf Q}}
\newcommand{\SM}{Statistical Mechanics}
\newcommand{\globe}{e^T}
\newcommand{\Glob}{{\cal G}}
\newcommand{\Floc}{{\cal F}}
\begin{document} \maketitle
\begin{abstract}
  For each natural number $n$, poset $T$, and $|T|$--tuple of scalars
  $\QQ$, we introduce the ramified partition algebra $P_n^{(T)}(\QQ)$,
  which is a physically motivated and natural generalization of the
  partition algebra \cite{Martin91,Martin94} (the partition algebra
  coincides with case $|T|=1$).  
For fixed $n$ and $T$ these algebras, like the partition algebra, 
have a basis independent of $\QQ$.  We investigate their
  representation theory in case $T=\Nset{2}:=(\{1,2\},\leq)$.  We show
  that $P_n^{(\Ttwo)}(\QQ)$ is quasi--hereditary over field $k$ when 
$Q_1 Q_2$ is invertible in $k$ and 
  $k$ is such that certain finite group algebras over $k$ are semisimple (e.g.
  when $k$ is algebraically closed, characteristic zero).  
Under these conditions we determine an index set for
  simple modules of $P_n^{(\Ttwo)}(\QQ)$, and construct standard modules
  with this index set.  We show that there are unboundedly many
  choices of $\QQ$ such that $P_n^{(\Ttwo)}(\QQ)$ is {\em not} semisimple
  for sufficiently large $n$, but that it is generically semisimple
  for all $n$.
  
  We construct tensor space representations of certain non--semisimple
  specializations of $P_n^{(\Ttwo)}(\QQ)$, and show how to use these
  to build clock model transfer matrices \cite{Martin91} in arbitrary
  physical dimensions.
\end{abstract}
\section{Introduction}

For $k$ a ring and $Q \in k$ the partition algebras $P_n(Q)$
($n=1,2,..$) are a tower of finite dimensional unital $k$-algebras. 
These algebras first arose in the context of transfer matrix
algebras in \SM\ \cite{Martin91,Martin94}, but are also known to play roles in
invariant theory 
\cite{Jones94,
MartinRollet98corrigendum} 
analogous to that of the
Brauer algebra \cite{Brauer37,Weyl46,HanlonWales94}, 
and in the study of Schur algebras (see \cite{HenkeKoenig01,MartinWoodcock98}
and cf. \cite[\S4]{Donkin01}). 
They
have a rich and (for $k=\C$) tractable
representation theory \cite{Martin96,Xi00}. 
The ramified partition algebras $P_n^{(T)} =  P_n^{(T)}(\QQ)$
are a generalisation 
depending on a poset  $T$, 
which is again
physically motivated (see below), but also natural as an abstract
extension, as we shall see. 
In this paper we construct these new algebras (in \S\ref{ss2}), 
and describe the representation theory of 
the simplest non--trivial cases (\S\ref{ss3}). 

An important tool in our approach is {\em quasi--heredity}
\cite{DlabRingel89B,Donkin93,Donkin98agc}.
When applicable 
this provides a crucial organizational scheme for an algebra,
via constructions for a set of {\em standard modules} 
(the set of {\em heads} of these modules constituting a complete set of simple
modules
\cite[\S A1]{Donkin98agc}).
When, as here, quasi--heredity `commutes' with a globalisation
functor \cite[\S6.2]{Green80} \cite{ClineParshallScott88} 
on a tower of algebras we have a
singularly powerful representation theoretic tool, even when the
algebras are semisimple 
(indeed in semisimple cases it typically {\em reduces} to a
glorified Jones basic construction 
\cite{GoodmandelaHarpeJones89,HalversonRam95}). 
Let $P_{n}^{(T)}\!-\!\mod$ denote the category of left $P_{n}^{(T)}$--modules. 
We may suppose without loss of generality, as we shall see, that the 
underlying set of the poset $T$ is $\{1,2,..,d \}$. 
Write  $\QQ = (Q_1, Q_2, ..) \in k^d$ and set 
$Q^{\pi} = \prod_{t \in T} Q_t \in k$. We shall show 
that for $Q^{\pi}$ invertible in $k$ 
there exists an idempotent $\globe \in P_n^{(T)}$ such that 
$\globe  P_n^{(T)} \globe  \cong  P_{n-1}^{(T)}$. It follows that
there are functors 
$\Floc: P_n^{(T)}\!-\!\mod\rightarrow P_{n-1}^{(T)}\!-\!\mod$ 
and $\Glob: P_{n-1}^{(T)}\!-\!\mod\rightarrow P_{n}^{(T)}\!-\!\mod$, called
{\em localisation} and {\em globalisation} respectively, which port information
between layers in the tower. 
Taken in conjunction with quasi--heredity 
this provides the framework for an iterative mechanism of
analysis 
(with the base case $n=1$ being amenable to brute force).%
\footnote{Both globalisation and quasi--heredity
apply to the ordinary partition algebra, and are key tools in its
analysis \cite{Martin96}. 
The utility of globalisation is a natural consequence of the physical
context of these algebras (cf.  \cite{Martin2000}), but it {\em is}
intriguing to see 
how often quasi--heredity also occurs in this setting 
\cite{ClineParshallScott99,Martin96,Xi00} (in as much as a physical
interpretation is so far lacking). 
}


This paper is concerned with the {\em representation theory} of the
ramified partition algebras. 
This {\em requires} no reference to, 
nor interest in,
the physical applications, 
but a remark on the authors' guiding motivation 
in constructing and studying these algebras 
is in order 
(see \S\ref{physics} for a more detailed exposition).  
The guiding motivation for
this work is the study of Statistical Mechanical lattice models
\cite{Baxter82,Pathria72} of three dimensional physical systems 
by exact computation \cite{BazhanovBaxter94,ReggeZecchina2000}.  
There has been some progress in this study using the ordinary partition algebra
$P_n=P_n(Q)$ \cite{MartinSaleur94b},
which underlies the transfer matrix algebra formalism 
(see \S\ref{physics}) for 
{\em Potts models} in high dimensions 
\cite{Martin91,Baxter82}.  However, as may be observed from
\cite{Martin94,DasmahapatraMartin96} and references therein, there
is a problem with the ordinary partition algebra approach 
in descending in dimension to a physical lattice capturing 
3d (or indeed any finite $d>2$) spatial geometry. 
We suggest that a better foundation would be provided by  
a more general algebra of {\em clock models} 
\cite{ElitzurPearsonShigemitsu,EinhornSavitRabinovici} 
or similar models \cite{AshkinTeller43,Fan72}. That is, one in which the set of 
{\em spin configurations} is equipped, 
though the form of the interaction Hamiltonian (see \S\ref{physics}), 
with a richer structure than the discrete topology 
(which is the structure effectively provided by the Potts Hamiltonian). 
Our reasoning is as follows. 
The defining basis of the partition algebra
is built from the set $N$ of spatial lattice sites  
(with $|N| = n$) by forming the set of partitions of $N$. 
A partition forms a base for a rather crude topology on $N$ 
(a {\em partition topology}). In the partition algebra this 
constrains the spatial geometry {\em and} the spin interactions to be 
very crude --- and in this sense the two limitations are linked. 
Accordingly we look here at a generalization of the partition algebra
suitable for the interactions in clock models. 
The idea is that the {\em set} partitions of $P_n(Q)$ should be 
glorified with a structure interpolating between the discrete 
topology and the metric topology of 3d Euclidean space. 
The challenge is that this structure not destroy the property of
associative algebra. 


The ordinary partition algebra is the case $T=\{ 1 \}$, 
however 
the representation theory of $P_n^{(T)}$ is vastly richer than that of
$P_n$ in general, as we will show.
The simplest 
interesting new case is $T=\Nset{2} := ( \{1,2 \}, \; 1 \leq 2)$. 
In this paper we concentrate mainly on this case. 
We determine
conditions for it to be quasi--hereditary, and construct its standard modules. 
We show that  $P_n^{(\Ttwo)}(\QQ)$ is generically semisimple, but 
(in \S\ref{ss4}) we also exhibit, 
for certain values of $\QQ$, 
{\em special} representations $R$ which may be used to show 
non--semisimplicity at these values of $\QQ$. 
It is also these representations which 
may be used in clock model transfer matrix computations (cf. \cite{Martin91}). 
Results of these computations for clock models will be presented in a
separate paper, but to underline the motivation we 
include a small sample of such results 
here --- zeros of the partition function for $Z_Q$--symmetric
models (recall that zeros converge on the real axis at a phase
transition point \cite{Martin91,MatveevShrock96b}).


The representation $R$ is also the obvious generalization of the
`Potts' representation of $P_n$ which arises in invariant theory
\cite{Martin91,Jones94,MartinWoodcock98}.  
An analysis of this, and of 
the generic representation theory of $P_n^{(T)}$, 
is the first step in
seeking for $P_n^{(T)}$  a
role in invariant theory analogous to that of $P_n$.

\section{Preliminaries}\label{ss2}
Let $(T,\leq)$ be a poset of degree $d$, 
$k$ a field,  $Q_t \in k$ for each $t \in T$, 
and $k$--algebra $P_n(Q_t)$ the ordinary partition algebra. 
The $T$--ramified partition algebra will be defined to be a certain
subalgebra of the tensor product $\bigotimes_{t \in T} P_n(Q_t)$,
depending on $\leq$.
\subsection{Construction of $P^{(T)}_n$}

First we recall the definition of $P_n(Q)$ \MI.
Given a set $U$ we write $\Ee(U)$ for the set of partitions of $U$ (we
identify this with the set of equivalence relations on $U$, writing
$\sim^a$ for the relation corresponding to $a\in \Ee(U)$).  
Note that any bijection $f:U \rightarrow V$ naturally underlies a
corresponding bijection $\Ee()^f:\Ee(U)\rightarrow\Ee(V)$.

For $a\in\Ee(U)$ we write $q\in a$ for $q\subseteq U$ a part of $a$ as
a partition. For $S\subset U$ write $a|_S$ for the restriction of $a$
to $S$. Given $n \in \N$ we write $\Nset{n}=\{1,2,\ldots,n\}$,
$\Nset{n}'=\{1',2',\ldots,n'\}$, and so on. We set $\Pbase_n =
\Ee(\Nset{n} \cup \Nset{n}')$, and $k\Pbase_n$ the $k$--space with
basis $\Pbase_n$.

Let $f:\Nset{n}\cup\Nset{n}' \rightarrow \Nset{n}'\cup\Nset{n}''$ 
be the map $i \mapsto i'$.
Then given $a,b\in \Pbase_n$ we define 
$\draw_u(ab)\in\Ee(\Nset{n}\cup\Nset{n}'\cup\Nset{n}'')$ 
by identifying with the equivalence relation generated by
$a$ and $\Ee()^f(b)$.
Let $g:\Nset{n}\cup\Nset{n}'' \rightarrow \Nset{n}\cup\Nset{n}'$ 
be the map $i \mapsto i$ for $i \in \Nset{n}$ and $i'' \mapsto
i'$ for $i''\in\Nset{n}''$. Then define 
$\draw(ab)=\Ee()^g(\draw_u(ab)|_{\Nset{n}\cup\Nset{n}''}) \in \Pbase_n$. 
Fixing $Q\in k$ we now define the product $ab \in k\Pbase_n$ 
by
\eql(partition algebra product)
ab = Q^{c(a,b)} \draw(ab) \eq
where $c(a,b)=\ha\{ q\subseteq\Nset{n}' \; : q\in\draw_u(ab) \}$. 
In \MI\ it is shown that this product ($k$--linearly extended) endows
$k\Pbase_n$ with the structure of a unital associative $k$--algebra,
which we call the partition algebra, and denote by $P_n(Q)$.
 

Now fix a poset $(T,\leq)$ with $T=\{1,2,..,d\}$
but $\leq$ arbitrary, 
and an element $\QQ =
(Q_1,Q_2,..,Q_d) \in k^d$. Our aim is to define the $T$--ramified
partition algebra as a subalgebra of $P_n^{(d)} = \bigotimes_{t\in T}
P_n(Q_t)$, depending on the partial order $\leq$. We will often
concentrate on the special case $T=\{1,2\}$ with $1\leq 2$ (and where
a poset is intended, $\Nset{2}$ will denote this one, 
and $\Nset{d}$ the natural ordered generalization), although
our construction is general. 

\newcommand{\refine}{\leq}
Given $a,b\in\Ee(U)$ we say that $a$ is a {\em refinement} of $b$ (see
e.g. \cite[ch.4 \S 7]{Mattson93}), denoted
$a \refine b$, if each part of $b$ is a union of one or more parts
of $a$; i.e. if $i\sim^a j$ implies $i \sim^b j$. 
It is easy to show that
\prl(pr0) 
Let
$a,b\in  \Ee(U)$ and $U'\subseteq U$.  
Then $a \refine b$ implies $a|_{U'} \refine b|_{U'}$. 
\end{pr}

For $S=\Ee(U)$ we define $S^{(T)}$ to be the subset of $\dip{T} S$
given by the elements $a =(a_t \; : t \in T )$ for which $t \leq t'$ implies
$a_t$ is a refinement of $a_{t'}$.  We call any such $a \in S^{(T)}$ a $T$--{\em
  ramified partition}. For example $(\{\{1\},\{2\}\},\{\{1,2\}\}) \in
\Ee(\{1,2\})^{(\Nset{2})}$.  

\prl(pr2)
For any $d$--tuple $\QQ=(Q_1,..,Q_d)\in k^d$, the set 
$\Pbase_n^{(T)}$ is a basis for a subalgebra of $\tp{t\in T}
P_n(Q_t)$. 
\end{pr}
{\em Proof:} It is enough to show that if $a,b \in \Pbase_n^{(T)}$
then $\draw(a_tb_t)$ is a refinement of $\draw(a_{t'}b_{t'})$ if $t \leq t'$. 
But if  $t \leq t'$ then $(a_{t},a_{t'})$ and $(b_{t},b_{t'})$ both appear
in $\Pbase_n^{(\Nset{2})}$ 
(where we identify $\Nset{2}=\{1,2\}$ with $\{t,t'\}$) 
so we are done if we can show for 
$\{ (a,b),(c,d) \} \subseteq \Pbase_n^{(\Nset{2})}$ 
that $\draw(ac)$ a refinement of $\draw(bd)$. 
Now by proposition~\ref{pr0} it is enough to show that 
$\draw_u(ac) \refine \draw_u(bd)$. 
Suppose $i \sim^{\draw_u(ac)} j$ and that the connection is established 
through some sequence of connections alternately within $a$ and $c$. 
Each connection within $a$ is within $b$, each connection within 
$c$ is within $d$, so $ i \sim^{\draw_u(bd)} j$. 
\hfill \Qed

Fixing $\QQ$ as above 
we denote the $k$--algebra with basis $\Pbase_n^{(T)}$ by $P_n^{(T)} =
P_n^{(T)}(\QQ)$. 
This is the {\em $T$--ramified partition algebra}. 


Two posets are isomorphic if there is a bijective homomorphism between
the underlying sets which preserves the partial order. It will be
evident that the construction of $P^{(T)}_n$ depends on the partial
order, but not on the names of the elements of $T$. Thus $P^{(T)}_n$
is defined for poset $T$ with {\em any} underlying set, and 
$T\cong T'$ implies  $P^{(T)}_n \cong P^{(T')}_n$.  


Not all choices of $T$ are {\em necessarily} physically interesting, in the
sense of the motivation discussed in the introduction.  From an
abstract point of view, as we will see later, the choice of $T$ as the
poset of divisors of a natural number gives rise to some interesting
representations. But we will concentrate largely on $T=\Ttwo$, which
is among the main motivating examples, captures much of the flavour of
the general case (consider the proof of proposition~\ref{pr2}), and
already has a very rich structure.

Refinement equips any $\Ee(U)$, and in particular $\Pbase_n$, with the
property of lattice. The bottom element $\{ U \}$ is called the
{\em trivial} partition. 


\subsection{Examples and notations} 

Putting $T=\{1\}$ we recover the partition algebra itself. 
Recall from \MI\ 
the special elements in basis $\Pbase_{n}$ denoted $1$, $A^i$, and
$A^{ij}$ as follows: 
\[ 
1= \{\{1,1'\},\{2,2'\},..,\{i,i'\},..,\{n,n'\}\} 
\]
$A^i$ differs from 1 by $\{i,i'\} \leadsto \{i\},\{i'\}$; 
and $A^{ij}$ differs from 1 by 
$\{i,i'\},\{j,j'\} \leadsto \{i,i',j,j'\}$ 
(see also the inner components of figure~\ref{fig1'} later).
Thus $1\in\Pbase_{1}$ is
$\{\{1,1'\}\}$ and $A^1\in\Pbase_{1}$ is $\{\{1\},\{1'\}\} $. 
See \MI\ also for the partition algebra {\em diagram calculus}. 

The next simplest case is
$T=\Nset{2}$. 
An example of the $\Nset{2}$--ramified partition algebra product 
is illustrated in figure~\ref{fig1} 
(the formal interpretation of this diagram is given shortly). 
The basis $p_1^{(\Nset{2})}$ consists of three elements: 
$(1,1)=(\{\{1,1'\}\},\{\{1,1'\}\})$;  
$\; (A^1,1)=(\{\{1\},\{1'\}\},\{\{1,1'\}\})$;  
and $(A^1,A^1)=(\{\{1\},\{1'\}\},\{\{1\},\{1'\}\})$. 
Note that $\{(A^1,A^1)\}$ is a basis for an ideal, 
as is $\{(A^1,A^1),(A^1,1)\}$ 
--- see later. 


\newcommand{\brak}[1]{{ \{ #1 \} }} 
It is convenient to introduce a `serial' shorthand notation and also 
a diagram calculus for $P_n^{(\Nset{d})}$. 
The shorthand for a $\Nset{2}$--ramified partition is 
to write out the more refined (`inner') 
partition and then to group the parts to indicate the less 
refined (`outer') partition by using another nest of brackets. 
Thus in $\Pbase^{(\Nset{2})}_1$, $(1,1)=\{\brak{\{1,1'\} }\}$, 
$(A^1,1)=\{\brak{\{1\},\{1'\} }\}$, and 
$(A^1,A^1)=\{\brak{\{1\} },\brak{\{1'\} }\}$. 
There is an obvious generalization for $T$ any chain. 


The diagram calculus applies the same principle to partition algebra
diagrams 
(in these the elements of $\Nset{n}$,$\Nset{n}'$ are represented as two
rows of vertices arranged with $i'$ vertically below $i$, so
composition involves vertical juxtaposition). 
Inner parts are indicated by `hard wiring' vertices together
as in ordinary partition algebra diagrams.  
Outer parts are indicated
by drawing `islands' around clusters of connected vertices. 
See figure~\ref{fig1} for some examples.  
Note by equation~(\ref{partition algebra product}) that
inner (outer) clusters {\em isolated from the exterior in composition}
are replaced by a factor $Q_1$ (resp. $Q_2$).  As for $P_n(Q)$, there
will not be a strictly planar representation in general.
Consider also the elements in figure~\ref{fig1'}, and 
note that 
\eql(relations)
\AxA{i} \AxA{ij} \AxA{i} = \AxA{i} \hspace{.51in}
\AxA{ij} \AxA{i} \AxA{ij} = \AxA{ij}. 
\eq

It will be evident from the diagrams that $P_n^{(\Nset{2})}$ is 
isomorphic to its opposite (indeed this is clearly true for any $T$). 
We denote the image of an element $a$ 
under this isomorphism by $a^o$ (as a diagram it is 
reflected top to bottom). 


\begin{figure}
\input{./xfig/prodex3.eepic}
\caption{\label{fig1} The algebra product in $P_3^{(\Nset{2})}$ in an 
example: $ab=Q_1 c$. 
Diagrams for $a$ and $b$ (drawn one above the other 
on the left) are juxtaposed to overlay and 
identify the vertices indicated by dotted lines, and these vertices are
then removed in exchange for an appropriate scalar factor (cf. \MI).}
\end{figure}
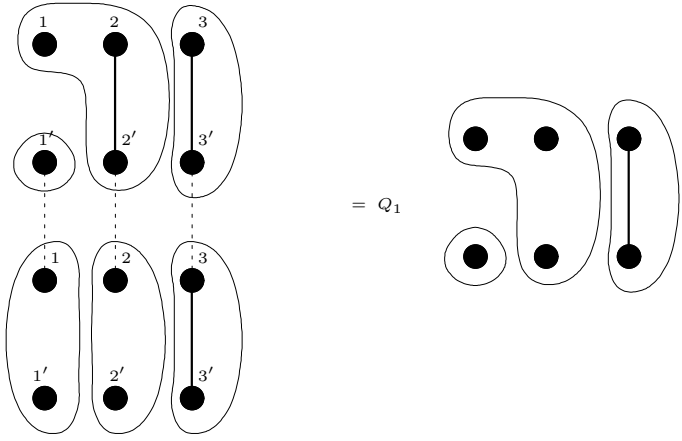
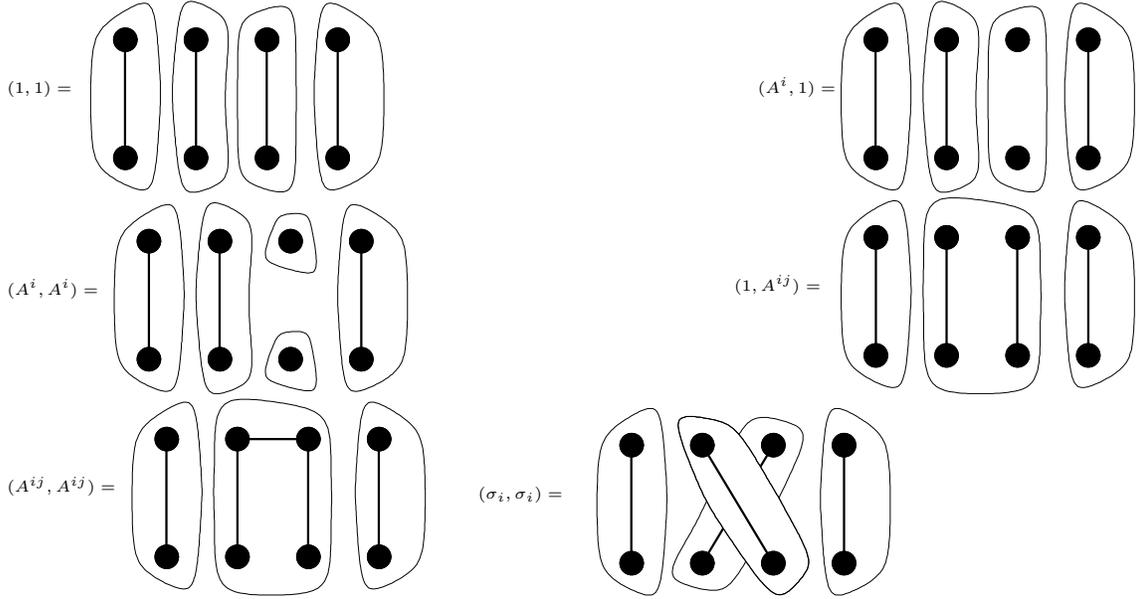
\begin{figure}
\input{./xfig/multunit.eepic}
\input{./xfig/Ai1.eepic}
\input{./xfig/AiAi.eepic}
\input{./xfig/OneAij.eepic}
\input{./xfig/AijAij.eepic}
\input{./xfig/Sigmai.eepic}
\caption{\label{fig1'} The multiplicative unit in 
$P_4^{(\Nset{2})}$ together
with various other elements. In these examples $i=2$ and $j=3$. 
Note that these `diagrams' do not necessarily embed in the plane in
any natural way, and that care must be taken with interpretation of 
`crossing' lines. Over/under information is not relevant, but may 
appear in diagrams as a guide to the eye.}
\end{figure}
\subsection{Canonical inclusions}


Let $S_n$ denote the symmetric group on $\Nset{n}$. 
By 
$ S_n \hookrightarrow P_n $
we denote the natural monomorphism  
(and here, as usual, 
we write $\sigma_i$ for the 2--cycle $(i \; i\!+\!1) \in S_n$, 
and write $\sigma_{ij}$ for the 2--cycle $(ij)$). 

For $t \in T$ define 
\begin{eqnarray*} \proj_t : p_n^{(T)} & \rightarrow & p_n \\
  (a_{1},a_{2},\ldots,a_{d}) & \mapsto & a_t . \end{eqnarray*}
Since $a \in p_n$ is a refinement of itself we have
\prl(diagonal)
The {\em diagonal} inclusion 
\begin{eqnarray*} \DT : p_n & \rightarrow & p_n^{(T)} \\
  a & \mapsto & (a,a,\ldots,a) \end{eqnarray*}
(or simply $D$, where no ambiguity arises) extends to an algebra morphism 
$$ \DT : P_n( Q^{\pi} )
               \hookrightarrow P_n^{(T)}(\QQ). $$
\end{pr} 
For example, set $e^T = \DT(A^n)$, then $e^T e^T = Q^{\pi} e^T$. 

There are various ways to inject $P_{n-1}^{(T)}$ into $P_n^{(T)}$. 
We treat as inclusion that which appends to every basis element a
part which is $\{ n , n' \}$ at every level of $T$. 
For example, 
in case $T=\Nset{2}$ the map on basis elements is given explicitly by 
\[
a \mapsto a \cup \{\{\{ n,n' \}\}\} . 
\]
\newcommand{\eT}{e^T}
\NB, including $P_{n-1}^{(T)}$ in $P_n^{(T)}$ in this way means that 
$\eT \; P_{n-1}^{(T)} = P_{n-1}^{(T)} \; \eT$, 
and hence that  
$\eT P_{n}^{(T)}$ is a left $P_{n-1}^{(T)}$--module. 


By
$\Pbase_n \stackrel{In}{\hookrightarrow} \Pbase_n^{(\Nset{2})}$ 
we denote the {\em inner} inclusion of $\Pbase_n$, 
i.e. where the inner  partition in $In(p)$ is $p$ 
and the outer partition is trivial. 



Note that the serial notation for a \trampart\ is 
unique if  written out so that at each step the numerically lowest possible 
element not yet written is written next 
(for which purpose take $i'=n+i$). 

Let us now fix the case $T=\Ttwo$ {\em until further notice}. 


\section{Structure of $P_n^{(\Ttwo)}$}\label{ss3}
\subsection{Filtration by ideals} 
For a partition $a \in p_n$ (or any restriction thereof) 
let $\parts{a}$ be the set of 
{\em propagating} parts of $a$, i.e. those intersecting both 
$\Nset{n}$ and $\Nset{n}'$. 
Then the {\em propagating number} $\hap(a) := |\parts{a}|$ (see \cite{Martin96}). 
For $a \in \Pbase_n^{(\Nset{2})}$ let 
$\hap(a) \in \N_0 \times \N_0$ denote $(\hap(a_1),\hap(a_2))$, 
and similarly for any $T$. 


For $a \in \Pbase_n^{(\Nset{2})}$ let the {\em propagating index} 
$\hal(a)$ be the list of 
$\hap(a_1|_q)$ over all $q \in \parts{a_2}$ arranged so that 
$\hal(a)_i \geq \hal(a)_{i+1}$. 
For example, the three $\Nset{2}$--ramified partitions illustrated in
figure~\ref{fig1} have $\hal(a)=(1,1)$, $\hal(b)=(1,0,0)$, $\hal(c)=(1,0)$. 
Note that the values taken by the propagating index differ from integer
partitions only in that there may be zero parts. 

A generalization of the propagating index may be formulated for
other choices of $T$, but we will not pursue the point here. 

For $\lambda$ a propagating index 
we denote by $\Pbase_n^{(\Ttwo)}(\lambda)$ the subset of 
elements $a \in \Pbase_n^{(\Ttwo)}$ such that $\hal(a)=\lambda$. 


The {\em propagating envelope} $\hae(a)$ is 
$\sum_i f(\hal(a)_i)$ where $f(0)=1$ and $f(j)=j$ otherwise. 
This number is the minimum $n$ required so that there is some 
$b \in \Pbase_n^{(\Ttwo)}$ with $\hal(b)=\hal(a)$. 
For example $\hal(\{\brak{\{1,2\},\{1',2'\} }\})=(0)$, so 
$\hae(\{\{\brak{1,2\},\{1',2'\} }\})=1$; and $a,b,c$ from
figure~\ref{fig1} have $\hae(a)=2$, $\hae(b)=3$, $\hae(c)=2$. 


\NB, $\hap$ and $\hae$ are determined by $\hal$ (specifically 
$\hap(a_1) = \sum_i \hal(a)_i$ and $\hap(a_2)$ is the number of
(integer) parts of $\hal(a)$). 
The basis $\Pbase_n^{(\Ttwo)}$ is fixed, 
and $\hal$ 
is invariant, under the opposite isomorphism. 

We denote by $\Pbase_n(i)$ the subset of elements $a \in \Pbase_n$ 
such that $\hap(a)=i$. 
Recall
\prl(basis lemma) \mbox{{\rm \cite{Martin96} }} 
Let $a \in \Pbase_n(i)$ and $b \in \Pbase_n(i')$. 
Then there exist $c,d \in \Pbase_n$ 
such that $a=cbd$ iff $i \leq i'$.
\end{pr}
This is a powerful tool in the analysis of $P_n(Q)$. The next few
propositions build up to a $T=\Ttwo$ generalization
(proposition~\ref{ideal filtration}).
\newcommand{\Pgen}{G}
\prl(generators) The subset $\Pgen_n$ of elements of the types illustrated 
in figure~\ref{fig1'} 
(or equivalently $\{ (A^1,1),(A^1,A^1),(1,A^{12}),(A^{12},A^{12}) \} \cup 
\{ (\sigma_i,\sigma_i) ; i=1,2,...,n-1 \}$) 
generates $P_n^{(\Nset{2})}$.
\end{pr}

{\em Proof:} 
Let $a$ be an element of the defining basis. 
We construct it as a product of elements of $\Pgen_n$ as follows.  
\newcommand{\low}{l}%
\newcommand{\lprop}{{\cal X}}%
Step 1 is to construct an element of $\Pbase_n^{(T)}$ with 
`enough' propagating connections 
to act as a base for the construction of $a$. 
For non--empty $q \subseteq \Nset{n}$, $q' \subseteq \Nset{n}'$ let 
$\low{(q \cup q')}$ (resp. $\low'{(q \cup q')}$) 
be the numerically lowest unprimed (resp. primed) element of 
the argument.
For $a \in \Pbase_n^{(\Nset{2})}$ let $\pi_1(a)=\pi(a_1)$ and 
$\pi_2(a) = \{ q \in \pi(a_2) \; | \; q\cap r =\emptyset\;\forall
r\in\pi(a_1) \}$. Let 
$\lprop_t(a)  = \cup_{q \in \pi_t(a)} \{ \low(q), \low'(q) \} $, 
and
 $\lprop^o(a) = \lprop_1(a)\cup\lprop_2(a)$ 
--- the union of the set of low numbered pairs 
for inner propagating parts of $a$ and the set of 
low numbered pairs for outer propagating parts containing 
no inner propagating parts. 
Now consider the set of basis elements of the form 
$w \in S_n \hookrightarrow P_n 
\stackrel{D}{\hookrightarrow} P_n^{(T)}$, 
i.e. those generated by the $(\sigma_{i},\sigma_{i})$s. 
Note that there is always at least one such $w$ which restricts, on 
$\lprop^o(a)$, 
to the same pairings as $a$ does (albeit diagonally). 
Choose any one such $w$.
Step 2: Apply $\AxA{i}$s so as to cut all the connections {\em not} on  
$\lprop^o(a)$.
Step 3: Take the resultant basis element and make all the 
additional outer connections required for $a$, 
i.e. those {\em within} $\Nset{n}$ and those 
within $\Nset{n}'$, 
using $\xA{ij}$s from the left and right as appropriate. 
Step 4: Make the additional inner connections for $a$ using $\AxA{ij}$s as 
appropriate (\NB, these will not make any further 
{\em outer} connections, 
cf. step 3 and the definition of $\Pbase_n^{(\Nset{2})}$). 
Finally make the required inner disconnections, 
i.e. between those pairs in $\lprop^o(a)$
referring to outer but not inner propagating parts, using $\Ax{i}$s. 
\hfill \Qed

\prl(pr3)
For $a,b \in p_n^{(T)}$ \newline
(i)
$
\hap(\draw(ab))_t \leq \min (\hap(a)_t , \hap(b)_t ) \; \forall t \in T.
$
\newline
(ii) If $\hap(a)_t > 0$ and $t \leq t'$ then $\hap(a)_{t'} >0$. 
\newline
And for $T=\Ttwo$ \newline 
(iii)
Let $a,b $ be such that $\hap(\draw(ab))_2=\hap(a)_2$. 
Then $\hal(\draw(ab))_i \leq \hal(a)_i$ for all suitable $i$. 
\newline
(iv)
$\hae(\draw(ab)) \leq \hae(a)$.
\end{pr}
{\em Proof:} The first claim is elementary 
in as much as it holds in $P_n^{(d)}$,
and $P_n^{(d)} \supseteq P_n^{(T)}$. 
The second and third follow from the fact that 
each propagating part of $a_1$ is a subset of a 
propagating part of $a_2$. 
To demonstrate {\em (iv)} 
let $\hae_0(a)$ be the number of trouser legs 
(outer propagating parts) 
containing no 
legs (inner propagating parts), then $\hae(a)=\hap(a_1)+\hae_0(a)$. 
The first summand is non--increasing in any product by {\em (i)}. 
The second summand may increase, but by {\em (i)} case $t=2$ if it 
does so then for each new empty trouser 
(increase of one) there must be at least one  non-empty 
trouser lost, so the first summand must decrease by at least one.  
\hfill \Qed


It follows that $\hal$ provides a means of filtering ideals of 
$P_n^{(\Ttwo)}$. 

\newcommand{\localthing}{{\cal I}_{\lambda}}
\begin{figure}
\input{./xfig/preid01.eepic}
\caption{\label{preid01} The special element $\localthing$ of
$p_{14}^{(\Nset{2})}((3,2^2,1^3,0))$.}
\end{figure}
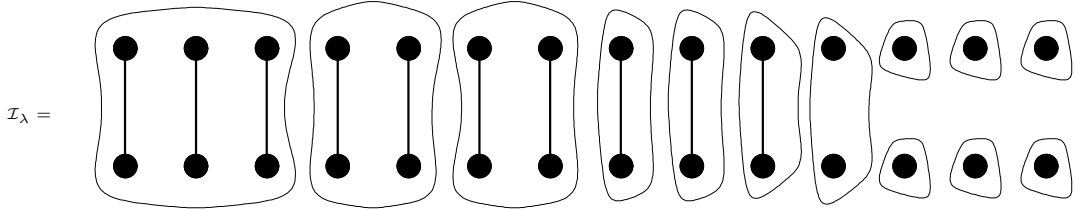

\newcommand{\Xel}{{\cal I}}
\newcommand{\Xell}[1]{{\cal I}_{#1}}
\newcommand{\xel}{x}

Let us pick out a special element $\Xel_{\lambda}$ 
of  $p_n^{(\Ttwo)}(\lambda)$ for each 
$\lambda$ as follows. Define $\xel_0=\{\brak{\{1\},\{1'\}}\}$ and 
\[
\xel_m = \{\brak{\{1,1' \},\{2,2'\}, \ldots , \{ m,m' \}}\} . 
\]
Recall the usual Young subgroup construction whereby 
$S_n \times S_m \hookrightarrow S_{n+m}$. We denote 
by $a \times b$ the analogous combination of elements 
$a \in P_n^{(T)}$ and $b \in P_m^{(T)}$ 
into an element  of $ P_{n+m}^{(T)}$. 
Then 
\[
\Xel_{\lambda} = \left( \times_i \xel_{\lambda_i} \right)
  \times \left( \times_{i=1}^{n-\hae(\lambda)} 
\{\brak{\{1\}},\{\brak{1'\}}\} \right)
\]
 where $\hae(\lambda)$, resp. $\hap(\lambda)$,  
denotes the $\hae$, resp. $\hap$, of any $a$ with 
$\hal(a)=\lambda$. 
An example is given in figure~\ref{preid01}. 


Let $\Lambda$ denote the set of possible propagating indices
$\lambda$. 
We write $\ol(\lam)=
(\ol(\lam)_1^{\ol(\lam)^1},\ol(\lam)_2^{\ol(\lam)^2},...)$  
for
$\lam \in \Lambda$ expressed in the exponent notation
(i.e. a repeated component is replaced by an exponent recording
the number of repeats, thus $(1,1,1)=(1^3)$ and so on). 

Note that for each $\lambda\in\Lambda$ there is an $m \in \N_0$ such 
that $\lambda \in \N_0^m$ ($\N_0^0$ is to be understood as the 
empty list). Each $\N_0^m$ has a natural action of $S_m$ upon it
(permuting the indices), and each $S_m$ orbit in $\N_0^m$ has 
a unique element which is a $\lambda$. 
This is called the {\em dominant} element.  
If we equate a $\lambda$ with an element $\mu$ of $\N_0^m$ we will 
always intend the dominant element in the orbit of $\mu$. 

We define a partial order $\leq$ on $\Lambda$ as follows. 
First write a relation $\lam' \rho \lam$ if \newline  
(0) $\lam'=\lam$; or \newline 
(i) $\lam'=(\lam_1,\lam_2,...,\lam_{i-1},\lam_i -1, \lam_{i+1},...)$; 
or \newline 
(ii) $\lam'=(\lam_1,\lam_2,...,\lam_{i-1},\lam_i +\lam_j,
\lam_{i+1},..., \lam_{j-1},\lam_{j+1},...)$. \newline
Now write $\lam' \leq \lam$ for the transitive completion of this 
relation (it is manifestly antisymmetric). 
Note that in case (i) $\hap(\lam)-\hap(\lam')=(1,0)$; 
and in case (ii)
$\hap(\lam)-\hap(\lam')=(0,1)$. 

\begin{figure}
\input{./xfig/poset1.eepic}
\caption{\label{poset1} The bottom of $(\Lambda,\leq)$.}
\end{figure}
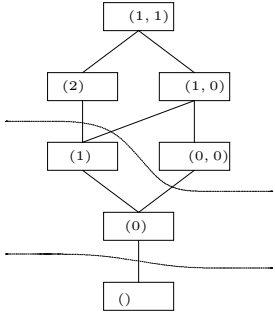

Examples: $() \leq \lam$ for all $\lam$ and, fixing $n$, 
$(1^n) \geq \lam$ for all $\lam$. 
The poset is {\em not} in general a lattice.
The Hasse diagram starts as in figure~\ref{poset1}. 

\prl(filter1') 
For $a,b \in p_n^{(\Ttwo)}$, 
$\hal(\draw(ab)) \leq \hal(a)$. 
\end{pr}
{\em Proof:} It is enough to show the result true 
for $a$ or $b$ a generator (as in proposition~\ref{generators}). 
But multiplication by a generator only implements 
one of (0)--(ii) above. 
\hfill \Qed

\newcommand{\meet}{\wedge}
\prl(ideal filtration) 
Put $P=P_n^{(\Nset{2})}$ and $p=\Pbase_n^{(\Nset{2})}$. 
\newline
(i) 
Let $a \in p(\lam)$ and $b \in p(\lam')$. 
Then there exist $c,d \in p$ 
such that $a=cbd$ iff $\lam \leq \lam'$.
\newline
(ii) $\cup_{\lam' \leq \lam} p(\lam')$ is a basis for $P \Xell{\lam} P$. 
\newline
(iii) 
$\lam' \geq \lam$ if and only if 
$P \Xell{\lam'} P \supseteq P \Xell{\lam} P$. 
\end{pr}
{\em Proof:} 
{\em (i)} The {\em only if} part follows from proposition~\ref{filter1'}.  

For each suitable $i$, there is an outer part of $b$ with  
$\lam_i'$ inner propagating parts. 
For each such $i$, consider acting on $b$ with $\AxA{j}$s 
so as to break up inner parts, until only the 
lowest numbered primed and unprimed elements in each 
inner propagating part remain connected 
(i.e., in a total of $\lam_i'$ primed/unprimed pairs). 
Note by proposition~\ref{basis lemma} that this process is {\em reversible} 
(for example if $j$ was in the same inner part as $k$ 
then $\AxA{jk}$ acting on the left would restore the connection).  
Then act with 
$(\sigma_{ij},\sigma_{ij})$s until these lowest numbered elements 
coincide with the corresponding data for $\Xell{\lam'}$. 
{\em This} step is reversible since $\sigma_{ij}$ is an involution.  
If $\lam = \lam'$ then $a$ may also be taken reversibly to 
$\Xell{\lam'}$. 
To obtain $c$ and $d$ simply compose the series of 
generators taking $b$ to $\Xell{\lam}$ then the reverses 
which would take $\Xell{\lam}$ to $a$. 

Finally, 
by definition, if $\lam < \lam'$ then we can move from $\lam'$ to 
$\lam$ by a sequence of relations $\rho$ of types (i) and (ii). 
Note that $\Xell{\lam'}$ may be taken 
to an element of $\Pbase_n(\lam)$ by a sequence of multiplications 
by generators which tracks this.

It will be evident that 
{\em (i)} implies 
{\em (ii)} implies {\em (iii)}. 
\hfill \Qed


An element $a$ of a $k$--algebra is a {\em pre--idempotent} if there
exists $\alpha\in k$ such that $\alpha a$ is idempotent.

For $\lam\in\Lambda\setminus\{ ()\}$ define 
$\Xell{\lam}'\in p_n^{(\Ttwo)}$
whose inner partition is given by 
$\proj_{inner}(\Xell{\lam}')=\proj_{inner}(\Xell{\lam})$; 
and whose outer partition agrees with $\proj_{outer}(\Xell{\lam})$ 
except on the tail of non--propagating parts (if any), 
whose elements are instead
included in the last (i.e. rightmost as in figure~\ref{preid01}) propagating
part. 
\prl(pre-ids)
Again put $P=P_n^{(\Nset{2})}$ and  
suppose $Q_1$ invertible in $k$. 
Then each $\Xell{\lam}'$ is pre--idempotent, and 
these pre--idempotents commute with each other, and obey
$$ P\Xell{\lam}'P =  P\Xell{\lam}P . $$
\end{pr} 
{\em Proof:} Note that each $\proj_{inner}(\Xell{\lam}')$ may be expressed
in the form $\prod_{j \in J} A^j$ for some $J$, and hence is
pre--idempotent in $P_n(Q_1)$ provided that $Q_1$ invertible. 
Further, the elements of $P_n$ of this form commute. 
Meanwhile 
$\proj_{outer}(\Xell{\lam}')$ may be expressed
in the form $\prod_{(i,j) \in J} A^{ij}$ for some $J$, and hence is
idempotent in $P_n(Q_2)$. 
Again, all such elements of $P_n$ commute with each other. 
Thus the set of objects $\Xell{\lam}'$ are commuting pre--idempotents in
$P_n^{(d)}$, and hence also in $P_n^{(\Ttwo)}$. 

Now  
$\Xell{\lam}'$ differs from $\Xell{\lam}$  only in
the tail (if any) of pairs of parts of the outer partition 
of the form $\{j\}\{j'\}$. If
$\lam \neq ()$ then this is preceded by at least one propagating
part. Let $\{\ldots (j-1) (j-1)' \}$ be the last such. 
Then $A^{(j-1) \; j} \{\ldots (j-1) (j-1)' \} \{j\}\{j'\} A^{(j-1) \; j}
               = \{\ldots (j-1) (j-1)' j j' \}$ 
and $A^{ j} \{\ldots (j-1) (j-1)' j j' \} A^{ j}
               = \{\ldots (j-1) (j-1)' \} \{j\}\{j'\}$ 
establishing the equality of ideals.  
\hfill \Qed


\subsection{On heredity}
Recall that propagating index determines propagating number. 
Define $\Lambda^m \subset \Lambda$ by $\lam\in\Lambda^m$ implies 
$\hap(\proj_{inner}(a))+\hap(\proj_{outer}(a))=m$ for any 
$a\in p_n^{(\Ttwo)}(\lam)$.
For example, $(0,0),(1) \in \Lambda^2$. 
\prl(hered.ids1)
There exist idempotents $\Xell{m}$ such that 
\newline
(i)
$ P_n^{(\Ttwo)} \Xell{m} P_n^{(\Ttwo)} 
      = \sum_{\lam\in\Lambda^m}
      P_n^{(\Ttwo)}\Xell{\lam}'P_n^{(\Ttwo)}$ 
for all $m$;
and \newline
(ii)
$ P_n^{(\Ttwo)} \Xell{m+1} P_n^{(\Ttwo)} \supset 
P_n^{(\Ttwo)} \Xell{m} P_n^{(\Ttwo)} $ for all $m=1,\ldots,2n-1$. 
\end{pr}
{\em Proof:}
Given a pair $(e_x,e_y)$, say, of commuting idempotents 
in an algebra $A$ one may construct an orthogonal
pair  $(e_x (1- e_y),e_y)$ and hence an idempotent $e_x-e_xe_y+e_y$. 
Clearly
$A(e_x-e_xe_y+e_y)A \supseteq
  Ae_x(e_x-e_xe_y+e_y)A + Ae_y(e_x-e_xe_y+e_y)A
   = Ae_xA + Ae_yA$ (saturating the inequality) 
--- \NB, sums inside and outside of brackets have different meanings here. 
Applying this process as many times as necessary to the idempotents 
with $\lam\in\Lambda^m$ 
(noting proposition~\ref{pre-ids})
we obtain our candidate for $\Xell{m}$ 
satisfying {\em (i)}. 
For {\em (ii)}, observe that every $\lam\in\Lambda^m$ ($m<2n$) obeys $\lam<\lam'$
for some $\lam'\in\Lambda^{m+1}$. 
\hfill \Qed

\begin{figure}
\input{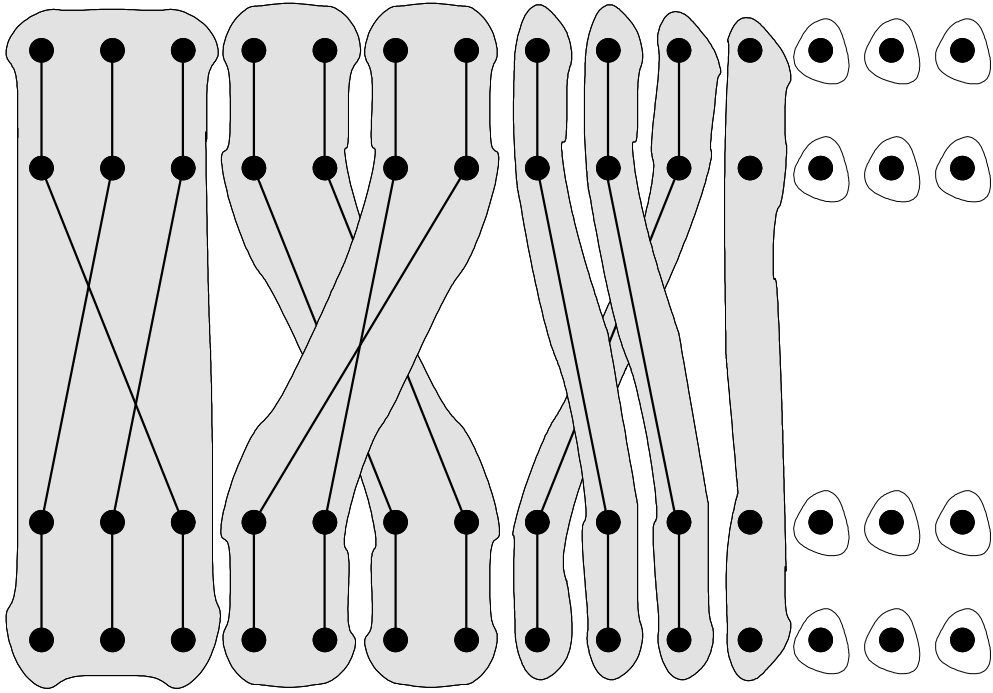}
\caption{\label{IPI02} Constructing an element of $I_{\lambda}PI_{\lambda}$ in
  case $\lambda=(3,2^2,1^3,0)$ (the middle two rows of points are here
  to exhibit membership of $I_{\lambda}PI_{\lambda}$ --- on ignoring
  these the diagram becomes an element of the basis
  $p_{14}^{(\Nset{2})}((3,2^2,1^3,0)) \cap
  I_{\lambda}PI_{\lambda}$). The shading here is a guide to the eye
  {\em only}.}
\end{figure}

Suppose that $Q_1Q_2$ is invertible in $k$. Then
$\Xell{\lam}$ and $\Xell{\lam}'$ 
are pre--idempotent (this being so we
will not further labour the distinction between $\Xell{\lam}$ and its
idempotent form). It is clear in this case that $\Xell{()}$ is a heredity
idempotent in the sense of \cite{DlabRingel89B} (or see
\cite{Martin96}). 
\prl(qh1) 
$P_n^{(\Ttwo)}$ is quasi-hereditary with heredity chain
$(\Xell{2n},\Xell{2n-1},\ldots,\Xell{1},\Xell{()})$. 
\end{pr}
{\em Proof:}
Write $P_i$ for $P_n^{(\Ttwo)}/P_n^{(\Ttwo)}\Xell{i-1}P_n^{(\Ttwo)}$. 
Given proposition~\ref{hered.ids1} it remains to show: 
(i) that subalgebra $\Xell{i}P_n^{(\Ttwo)}\Xell{i}$ is semi--simple in
$P_i$; and (ii) that 
$P_n^{(\Ttwo)}\Xell{i} 
  \otimes_{\Xell{i}P_n^{(\Ttwo)}\Xell{i}}\Xell{i} P_n^{(\Ttwo)}
    \cong P_n^{(\Ttwo)}\Xell{i} P_n^{(\Ttwo)}$ in $P_i$. 
A basis for $P_i$ is $\cup_{\lam\in\Lambda^i} p_n^{(\Ttwo)}(\lam)$. 
It follows from proposition~\ref{filter1'} 
that the product of any two basis elements not
from the same $p_n^{(\Ttwo)}(\lam)$ is zero in $P_i$. 
For (i), one then notes that the subalgebra
$\Xell{i}P_n^{(\Ttwo)}\Xell{i}$ 
will break up as a direct sum with summands of form 
$\Xell{\lam}P_n^{(\Ttwo)}\Xell{\lam}$. 
The basis of {\em this} algebra given by its intersection 
with the given basis of $P_i$ 
(illustrated in figure~\ref{IPI02}) 
is closed under multiplication, with multiplicative 
unit $\Xell{\lam}$, and invertible elements. Thus 
this algebra  
is isomorphic to a group algebra. (We will identify the group in
section~\ref{smodules}.)  

For (ii), one considers the obvious map 
$$  (aI_i,I_ib) \mapsto aI_ib $$ 
in terms of basis elements. 
Again, since the $\{ I_{\lam} \}$ are orthogonal (in $P_i$) it is enough to
consider $  (aI_{\lam},I_{\lam}b) \mapsto aI_{\lam}b $. 
As written, this requires that the kernel coincides with 
elements of the form $(ac,b)-(a,cb)$ where 
$c \in \Xell{i}P_n^{(\Ttwo)}\Xell{i}$. 
The kernel clearly includes this set, so it is enough to show that it 
is no larger. 
The map is realizable at the level of basis elements, from which the 
extent of the kernel is immediate, 
as follows. 

Let $q$ be an element of $\Pbase_n^{(\Ttwo)}(\lam) \cap PI_{\lam}$ (which
set is a basis for $PI_{\lam}$ in $P_i$), and $a$ an element of
$\Pbase_n^{(\Ttwo)}(\lam) \cap I_{\lam}PI_{\lam}$ (similarly a basis for
$I_{\lam}PI_{\lam}$). Up to a scalar $qa$ is in $\Pbase_n^{(\Ttwo)}(\lam)
\cap PI_{\lam}$ and the map
\[
\Pbase_n^{(\Ttwo)}(\lam) \cap I_{\lam}PI_{\lam} 
     \rightarrow \Pbase_n^{(\Ttwo)}(\lam) \cap PI_{\lam}
\]
given in this way by $a \mapsto qa$ is injective for any $q$ (again consider
figure~\ref{IPI02}, for example). Define an equivalence relation on 
$\Pbase_n^{(\Ttwo)}(\lam) \cap PI_{\lam}$ by $q \sim q'$ if $qa=q'$ for some $a$ as
above. Let $q^{\sim}$ be a set of representatives of equivalence
classes. Then each element of $\Pbase_n^{(\Ttwo)}(\lam) \cap PI_{\lam}$ may be
written uniquely in the form $qa$ where $q \in q^{\sim}$ (i.e. $(q,a)
\mapsto qa $ is a bijection). 

The same argument works on the right. 
\hfill \Qed

\subsection{Left ideals, sections and special modules}\label{smodules}

Let $P_n^{(\Ttwo)}[\lam]$ denote the $\lam$ section in the 
filtration given by proposition~\ref{ideal filtration}(ii) 
(i.e. with basis $\Pbase_n^{(\Ttwo)}(\lambda)$). 
Let $P_n^{(\Ttwo)}[\lam]'$ denote the image of $P_n^{(\Ttwo)} \Xell{\lam}$
in this section (cf. \cite[\S 3]{Martin96}). 
The algebra which is the image of $\Xell{\lam} P_n^{(\Ttwo)} \Xell{\lam}$
in this section acts on $P_n^{(\Ttwo)}[\lam]'$ from the right.  
So in addition to being a left  $P_n^{(\Ttwo)}$--module
$P_n^{(\Ttwo)}[\lam]'$ 
is also a right module for the group algebra of a permutation group 
$S[\lam]$ 
whose action is to permute the primed sides of outer 
propagating parts with the same number of inner 
propagating parts; and to permute the primed elements of 
inner propagating parts within each outer propagating part 
--- thus 
$$ S[\lam] = \prod_i S_{\ol(\lam)_i} \wr S_{\ol(\lam)^i} $$ 
(consider the action from below of 
$S_1 \times S_2 \times S_3 \times S_1$ on outer and 
$S_3 \times (S_2 \times S_2) \times (S_1 \times S_1 \times S_1) \times S_0 $ 
on inner parts in figure~\ref{preid01}). 
Indeed, $P_n^{(\Ttwo)}[\lambda]'_{S[\lam]}$ is a sum of copies 
of the right regular module. 
\newcommand{\kk}{k}%
Thus if $L_{\Theta}$ is a simple left $\kk S[\lam]$--module we 
may define  a left $P_n^{(\Ttwo)}$--module 
$$\SS(\Theta) = \SS(\Theta)(n) = 
P_n^{(\Ttwo)}[\lam]' \otimes_{\kk S[\lam]} L_{\Theta} . $$ 
This is an indecomposable summand of $P_n^{(\Ttwo)}[\lambda]'$ and
hence a standard module of $P_n^{(\Ttwo)}$ (in the quasi--hereditary
sense). 


Evidently, every simple $P_n^{(\Ttwo)}$--module appears as a subquotient
of some $P_n^{(\Ttwo)}[\lambda]'$, and each $P_n^{(\Ttwo)}[\lambda]'$ is a sum
of modules of form $\SS(\Theta)$, with 
$\SS(\Theta)\cong \SS(\Theta')$ iff $L_{\Theta}
\cong L_{\Theta'}$. 
Thus every simple $P_n^{(\Ttwo)}$--module appears as a subquotient
of some $\SS(\Theta)$. The following result is well known.
\prl(261000.1)\mbox{\rm \cite[\S I.Appendix B]{Macdonald95} }
Let $G$ be a group with $C_G$ conjugacy classes. Then
equivalence classes of simple $G \wr S_n$--modules are indexed by
$C_G$--tuples of integer partitions of total degree $n$. 
\Qed \end{pr}
It follows that for each $\lambda$ there is a distinct inequivalent
$\SS(\Theta)$ for each {\em tuple} of multipartitions (tuples of integer
partitions are called multipartitions) in which the $i^{th}$
multipartition is as specified in the proposition in case
$S_{\ol(\lam)_i} \wr S_{\ol(\lam)^i} $. 
 
Let us denote by $\Gamma_n$ the set of all such indices $\Theta$ unioned over all
$\lambda$ such that $\hae(\lambda)=n$. 
For example, 
$\Gamma_0$ consists of a single label ($\spineweight$, say) 
corresponding to the nominal irreducible representation of 
a product over an empty set of wreaths 
(it will be evident from the construction that the correct 
interpretation is as a single copy of the ground ring); 
$\Gamma_1$ consists of labels, 
$(0),(1)$, say, 
corresponding to the single irreducible representations of each of $S_1 \wr
S_0$ and $S_1 \wr S_1$ respectively;
$\Gamma_2$ consists of labels, 
$(0^{(2)}),(0^{(1^2)})$ (the irreducibles of $S_2 \wr S_0$), 
$(1,0)$ (the irreducible of $(S_1\wr S_1)\times(S_1\wr S_0)$), 
$(1^{(2)}),(1^{(1^2)})$ (irreducibles of $S_2 \wr S_1$), 
$(1,1)^{(2)}$ and $(1,1)^{(1^2)}$ (irreducibles of $S_1 \wr S_2$), 
say. 
See figure~\ref{fig2} for more examples. 


Note that the $\kk S[\lambda]$ action may be realized as a subalgebra 
of $ P_n^{(T)}$ in a similar way to
$\kk S_n$. 
The extra structure in $\kk S[\lambda]$ is just so that it 
commutes with $\Xell{\lambda}$. 
When we confuse elements of $\kk S[\lambda]$ with elements of 
$ P_n^{(T)}$ it will be through this realization. 
Since $\kk S[\lambda]$ is
a group algebra we may take it semisimple for $\kk$ a suitable field of
characteristic 0 (such as $\C$). Let $e_{\Theta}$ then be a primitive
idempotent of $\kk S[\lambda]$ such that 
$\kk S[\lambda] e_{\Theta}= L_{\Theta}$ and 
$e_{\Theta} \kk S[\lambda] e_{\Theta}= \kk e_{\Theta}$. 
We may define an inner product $<|>$ on $ \SS(\Theta)$ by 
$$ <a|b>  e_{\Theta} \Xell{\lambda} = a^o b $$
(to be understood in the $\lambda$ section, which makes the opposite
right module to $\SS(\Theta)$ also a dual). 

For example, $\Xell{\emptyindex}= \prod_{i=1}^n (A^n,A^n)$ and
$S[\emptyindex] = S_0$ which here means the trivial group (so 
there is only one irreducible representation and 
$e_{\spineweight}=1$). A basis for $P_n^{(T)} \Xell{\emptyindex}$ is 
$\{ wx \; | \; w \in \Ee(\Nset{n})^{(T)}, \; 
x=\{\{\{1' \}\},\{\{ 2' \}\},...,\{\{n' \}\}\} \}$. Thus the gram
matrix $G_{\spineweight}$ of the inner product with respect to this basis
has rows and columns conveniently indexed by $\Ee(\Nset{n})^{(T)}$. 
The $ab^{th}$ entry of $G_{\spineweight}$ is $\prod_i Q_i^{h_i}$ where
$h_i$ is the number of parts in $a_i b_i$ (the normal composition of
partitions). In particular for $n=2$ 
the $\lambda=\emptyset$ module has gram matrix computed as follows:
\[
\input{./xfig/Gram1.eepic}
\]
whereupon 
$$|G_{\spineweight}^{n=2}| = Q_1^3 Q_2^4 (Q_1 -1)(Q_2-1) . $$ 
Similarly we have (for the dimension twelve module at $n=3$) 
$$|G_{\spineweight}^{n=3}| = 
Q_1^{12} Q_2^{20} (Q_1 -2)(Q_2 -2)(Q_1 -1)^{7}(Q_2 -1)^{7}  $$  
and (for the dimension 60 module at $n=4$)   
$$|G_{\spineweight}^{n=4}| = 
Q_1^{60} Q_2^{119} (Q_1 -3)(Q_2 -3)(Q_1 -2)^{13}(Q_2 -2)^{11}
(Q_1 -1)^{45}(Q_2 -1)^{48} . $$   
What is striking about these results is that the determinants
factorise over $\Z[Q_1]\cup\Z[Q_2] $ rather than $\C[Q_1,Q_2]$ 
\footnote{Indeed, computing the determinant of a 
common or garden 60x60 matrix with entries
polynomial in 2 variables is probably beyond the reach of modern
computers. We used the GAP package\cite{GAP99}, some lucky guesses and some
saturated bounds to reduce to a few tens of minutes computation on a
powerful laptop.}. 
\prl(spine generically irreducible)
For all $n$ the module $\SS(\spineweight)$ is generically simple. 
\end{pr}
{\em Proof:}
Consider $G_{\spineweight}$ for $n$ arbitrary. 
It is enough to show that this is non--singular on an open subset of
parameter space. 
Recall that the composition of partitions required to compute the inner
product involves a transitive closure. It follows that 
the diagonal entry in each row of  $G_{\spineweight}$ 
is always of (possibly equal) highest total
degree. Consider the matrix obtained by 
deleting the rows (and columns) where it is uniquely of highest
degree. Iterating this procedure (cf. \cite{MartinSaleur94b}) 
we arrive at the empty matrix. It
follows that the determinant of the gram matrix is a finite polynomial
in $Q_1,Q_2$. \Qed

\subsection{Modules $ \SS(\Theta)(n)$ with $\hae(\Theta)=n$}
Let us extend the definition of $\hae$ so 
$\hae(\Theta) = \hae(\lambda)$, where $\lambda$ is the propagating index 
underlying $\Theta$. 
Then `modules $\SS(\Theta)$ with $\hae(\Theta)=n$' are  
the $P_n^{(\Ttwo)}$--modules $\SS(\Theta)$ whose canonical basis 
elements have propagating index which cannot be realized in
$P_{n-1}^{(\Ttwo)}$. They also obey 
$$ \DT(A^i) \SS(\Theta) = \DT(A^{ij}) \SS(\Theta) =0 $$
for all $i,j$, whereas no $\SS(\Theta)$ with $\hae(\Theta)<n$ has this property 
(for, in this case we may construct a basis element from a basis element for the 
equivalent module at level $n-1$ by 
$a \mapsto a \cup \{\{\{n\}\},\{\{n'\}\}\}$ --- 
and this is not killed by $\DT(A^n)$). 
\prl(simple1)
If $\hae(\Theta)=n$ and $Q_{inner}$ invertible in $k$ 
then $\SS(\Theta)(n)$ is simple. 
\end{pr}
{\em Proof:} 
Consider the gram matrix as constructed in the previous section. 
Because of the quotient (we are in the $\lambda$ section) 
we may apply the pigeonhole principle to see that 
the inner
product $<a|b>$ between any two basis elements is zero unless each
outer part of $a$ meets only one outer part of $b^o$ in the
composition.
Thus the gram matrix breaks up into blocks indexed by the possible 
arrangements of outer parts in the restriction of $a$ to $\Nset{n}'$ 
(or equivalently the restriction of $b^o$ to $\Nset{n}$), 
and this arrangement characterizes all the products in the block. 
Since $\hae(\Theta)=n$ this implies that, 
in a block, each 
outer part from $a$ meeting one from $b^o$ carries the same number
of inner propagating lines. Within the block, then, up to overall factors of
$Q_{inner}$, the inner product coincides with that of the
corresponding simple $\kk S[\lambda ]$--module. Each block determinant 
is therefore the same, and non--zero. 
Hence the whole gram matrix is non--singular. 
\Qed

\subsection{Generic representation theory}

For any $T$ 
recall that $\eT$ commutes with $P_{n-1}^{(T)}$, 
and that $\eT \eT = Q^{\pi} \eT$. 
Thus
\prl(embed) 
For $\prod_{t \in T} Q_t$ invertible in $k$,
\eql(localize)
\eT \; P_n^{(T)} \; \eT \; 
= \; P_{n-1}^{(T)} \; \eT \; \cong P_{n-1}^{(T)}
\eq
is an algebra isomorphism; and 
\begin{eqnarray*} \label{functorloc} 
\Gg: P_{n-1}^{(T)}\mbox{--mod} & \rightarrow & P_{n}^{(T)}\mbox{--mod}
\\
M & \mapsto & \; P_{n}^{(T)} \; \eT \; M
\end{eqnarray*} 
is a full embedding of $P_{n-1}^{(T)}$--mod in $P_{n}^{(T)}$--mod
\cite[\S 6.2]{Green80} (cf. \cite{Martin96}).
\end{pr} 


Note that 
$\Gg(P_{n-1}^{(\Ttwo)} \Xell{\lambda}) = P_{n}^{(\Ttwo)} \Xell{\lambda}$ 
(\NB, strictly speaking $ \Xell{\lambda}$ is a
different object on each side of this equation, i.e. depending on the
level $n$; it is this 
`globalization' identity which justifies omission of the $n$ label). 


By equation~(\ref{localize}), it is the quotient 
\[
P[n] \; := \; P_n^{(T)} / P_n^{(T)}  \eT  P_n^{(T)} 
\]
whose representations will be missed in the embedding 
(or equivalently, those left modules $M$ for which 
$  \eT.  M=0$). 
Thus 
\prl(simpleindex)
If $\Lambda_n^{(T)}$ is an index set for the equivalence classes of 
simple modules of 
$ P_n^{(T)}$ then 
\newline
(i) $|\Lambda_0^{(T)}|=1$; \newline
and if  $ Q_1 Q_2$ invertible in $k$ then \newline 
(ii) $\Lambda_n^{(\Ttwo)}=\Lambda_{n-1}^{(\Ttwo)} \cup \Gamma_n 
 = \cup_{0\leq i \leq n} \Gamma_i$, and \newline 
(iii) $\{ \head \SS(\Theta) \; | \; \Theta \in \Lambda_n^{(\Ttwo)} \}$ is a
complete set of representatives of the equivalence classes of 
simple modules of $ P_n^{(\Ttwo)}$. 
\end{pr}


In other words, the set of equivalence classes of 
irreducible representations is indexed by the possible values of 
$\hal$ together with 
the corresponding set of tuples of multipartitions described in 
section~\ref{smodules}. 


Recall \cite[\S2.2]{Martin2000} that $P_n^m$ is the subalgebra of 
$P_{n+m}$ in which each of the nodes labelled $n+1,\ldots,n+m$ is {\em identified} 
with its primed partner.  
The algebras $P_n^m$ generalize in a natural way to algebras
$P_n^{(T)m}$ (any $T$), 
whereupon $P_n^{(T)} \subset P_n^{(T)1} \subset P_{n+1}^{(T)}$. 
The chain of inclusions from $P_0^{(T)}$ upward obtained from this 
sequence is denoted $P_*^{(T)}$. 

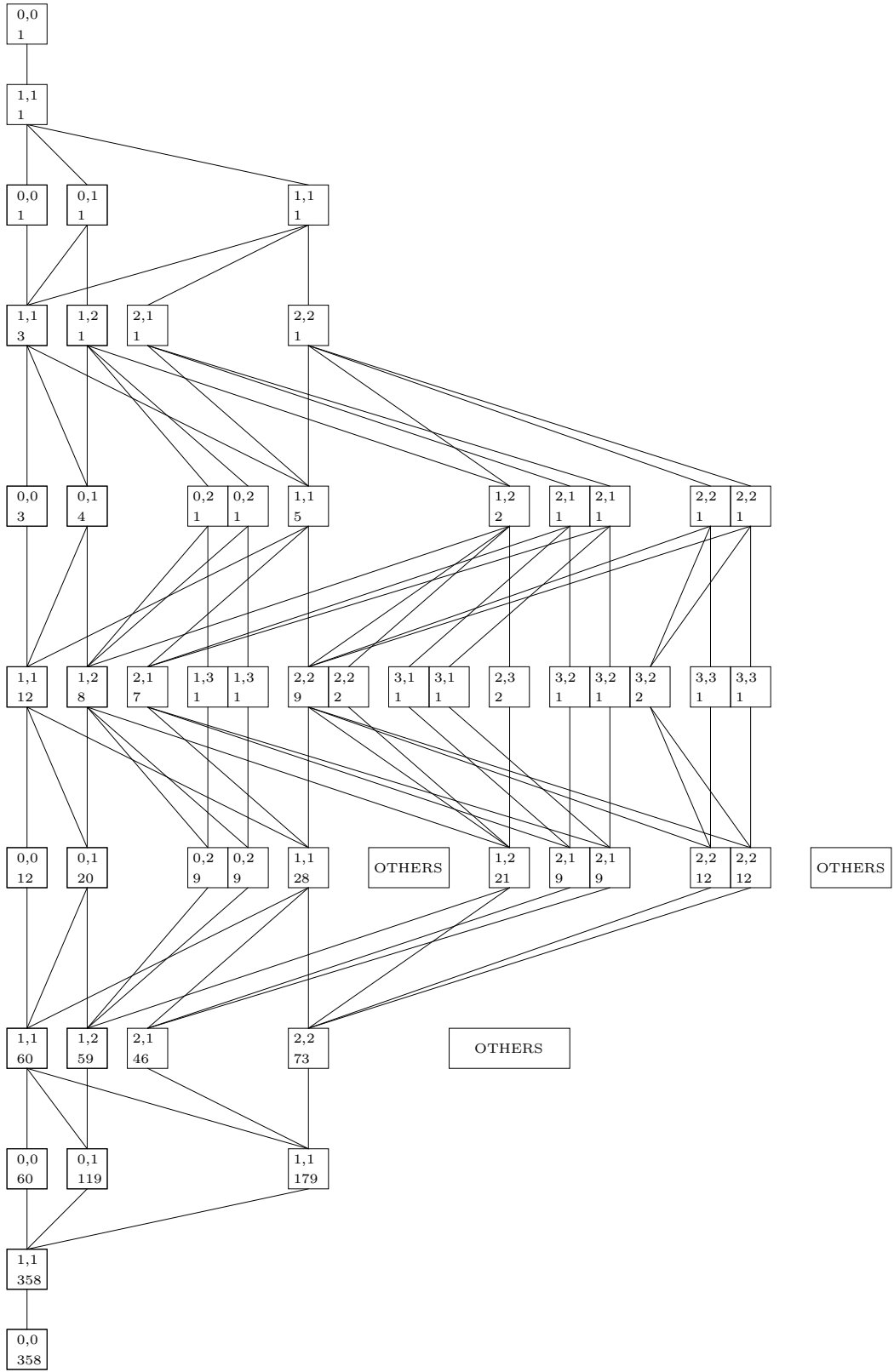
\begin{figure}
\input{./xfig/bratteli.eepic}
\caption{\label{fig2} 
Generic restriction rules for the generic 
irreducible $\hap=(0,0)$ representation
of $P_5^{(\Nset{2})}$ through the chain of restrictions $P_*^{(\Nset{2})}$. 
The top layer is $P_0^{(\Nset{2})}$ and the bottom is $P_5^{(\Nset{2})}$. 
In each layer up to $P_2^{(\Nset{2})1}$ a complete list of irreducible
representations is given. 
Thereafter only irreducibles occurring in the restriction are 
recorded. The irreducibles are labelled by their $\hap$ index
(i.e. ambiguously) and by their dimension.}
\end{figure}

As for restriction rules, it will be clear that the module $\SS({\spineweight})$
in $P_{2n}^{(\Ttwo)}$ restricts to the regular representation of
$P_n^{(\Ttwo)}$. En route we have the restriction 
$$ \Res(P_{n}^{(\Ttwo)},P_{n-1}^{(\Ttwo)}, {\SS({\spineweight})} ) 
     \cong \SS({\spineweight}) \oplus \SS({(0)}) \oplus \SS({(1)}) $$
and so on. The first few layers are depicted in figure~\ref{fig2}. 
We do not need the general result here, but see
proposition~\ref{res metric}.  



Given a module $M \in P_{n-1}^{(T)}\mbox{--mod}$ there is a globalization
functor $\Gg$ and, 
via the injection of $P_{n-1}^{(T)}$ into $P_{n}^{(T)}$, 
an induction functor $\; \mbox{Ind} \;$ acting on it,  
so that 
$\Gg M $ and $ \mbox{Ind} M $ are in $ P_{n}^{(T)}\mbox{--mod}  \;$ 
{\em  for any} $n$. 
In writing down such modules, 
we will not generally give an explicit indication of 
the level $n$ unless it is seems helpful to do so.

We now  assemble the tools for proving generic semisimplicity. 
By construction (or using \cite{Donkin98agc}) we have
\prl(res standard)
$$ \Gg \SS({\lam})(n) \cong \SS({\lam})(n+1) $$
\end{pr}
\newcommand{\qh}{quasi--hereditary}
Recall \cite{Donkin93} 
\prl(proj ff)
Projective modules of \qh\ algebras have filtrations by standard modules.
\end{pr}
We denote the multiplicity in this filtration by
$[P_{\lam}:\SS(\mu)]$. 
The following is readily proved (cf. \cite{Martin96}):
\prl(indres)
For $M \in P_n^{(T)}-\Mod$, 
$$ \Res(n+1,n+2, \Gg \Gg M ) \cong \Ind(n,n+1, M )  $$
\eql(ind G commute) \Indi \Gg M \cong \Gg \Indi M . \eq
\end{pr}
Thus if $\SS(\lam)(n)$ is projective for all $n$ then so is
$\Resi\SS(\lam)(n)$. Whereupon 
\prl(res metric)
$[\Resi\SS({\lam}):\SS(\mu)] \neq 0$ implies 
$\hae(\lam)-\hae(\mu) \in \{-1,0,+1 \}$.
\end{pr}
{\em Proof:}
Consider the restriction 
$\Res(P_{n-1}^{(\Ttwo)},P_{n}^{(\Ttwo)},P_{n}^{(\Ttwo)}[\hat{\lam}]')$ 
where $\hat{\lam}$ is such that $ \SS({\lam})$ is a summand. 
Since $ P_{n}^{(\Ttwo)}[\hat{\lam}]'$ is a direct sum of $
\SS({\lam})$s all with the same $\hae$, it is enough to show an
equivalent result for this module. 

In our basis of $ P_{n}^{(\Ttwo)}[\hat{\lam}]'$
the elements conform to one of the following
possibilities for the connectivity of node $n$ (the distinguished
node, as it were, on restriction).
\newline
(1a) Node $n$ is the sole unprimed element in a propagating part (inner
    not outer) --- this contributes to a basis for a module isomorphic
    to $ P_{n-1}^{(\Ttwo)}[{\lam'}]'$ with
    $\lam' \rho \hat{\lam}$ via a $\rho$--step of type (i) (i.e. with
    $\hae(\lam') = \hae(\lam) -1$).  
\newline
(1b) Node $n$ is the sole unprimed element in a propagating part (inner
    and outer) --- this means that $\hap(\hat{\lam})$ includes a
    singleton, and the situation is as above except that this
    singleton has disappeared in $\lam'$. Again $\hae(\lam') =
    \hae(\lam) -1$. 
\newline
(1c) Node $n$ is the sole unprimed element in a propagating part
(outer not inner) --- this time $\hap(\hat{\lam})$ includes a
    zero element and the situation is as above except that this
    has disappeared in $\lam'$. Again $\hae(\lam') =
    \hae(\lam) -1$. 
\newline 
(2a) Node $n$ is a singleton both with respect to inner
and outer partitions --- the set of such objects forms a basis for a
module isomorphic to $ P_{n-1}^{(\Ttwo)}[\hat{\lam}]'$. 
\newline 
(2b) Node $n$ is connected, but in such a way that regarding it as primed
    (for the purpose of determining the propagating index) does not
    change the propagating index (i.e. both the inner and outer part
    in which it lies have other nodes of both primed and unprimed
    type). Again sets of such objects form bases for modules
    isomorphic to $ P_{n-1}^{(\Ttwo)}[\hat{\lam}]'$, provided we
    quotient by the modules produced at (1abc) above. 
\newline
(3a) Node $n$ is in a part (inner not outer) containing only unprimed
    elements --- in the restriction there is no action on this node,
    so, quotienting by all the modules discussed so far, it may be
    regarded much as if it were primed. Thus we can build a basis for
    a copy of a module   $ P_{n-1}^{(\Ttwo)}[{\lam'}]'$  where $\lam'$
    is such that one element of $\hap(\lam')$ is incremented by
    1 compared to $\hap(\lam)$, and hence  $\hae(\lam') = \hae(\lam) +1$. 
\newline
(3b) Node $n$ is in a part (inner and outer) containing only unprimed
    elements --- as above, except that $\hap(\lam')$ has a new
    singleton element. Again $\hae(\lam') = \hae(\lam) +1$. 
\newline
(3c) Node $n$ is in a part (outer not inner) containing only unprimed
    elements --- as above, except that $\hap(\lam')$ has a new
    zero element. Again $\hae(\lam') = \hae(\lam) +1$. 
\hfill \Qed

 

\begin{theo}
$P_n^{(\Ttwo)}(\QQ)$ is generically semisimple.
\end{theo}
{\em Proof:}
Let Prop$_{\lam}$ be that $\SS(\lam) = P_{\lam} = S_{\lam}$, the simple
module, for all $n$ (i.e. this module lies in a singleton block). 
Then we are done if Prop$_{\lam}$ true for all $\lam$. 
Suppose it is true for all $\{ \lam \; | \; \hae(\lam) \leq m \}$. Work by
induction on $m$. The base case is $m=0$, which is true generically
by proposition~\ref{spine generically irreducible} using
quasi-heredity and Brauer--Humphreys reciprocity \cite{Donkin93}. 

Applying equation~(\ref{ind G commute}) to case $M=1$, $n=0$ we see
that every standard module with $\hae(\nu) \leq m$ appears in the
$m^{th}$ induction (and hence also the $m^{th}$ restriction) 
of $\SS({\spineweight})$. 
Thus for each $\nu$ with $\hae(\nu)=m+1$ there is a $\lam$ with
$\hae(\lam)=m$ such that $[\Resi\SS({\lam}):\SS(\nu)] \neq 0$. 
We now use this $\lam$ and the inductive assumptions to establish Prop$_{\nu}$. 

Consider the Frobenius reciprocity
\eql(FR)
\Hom_{n+1} (\Indi \SS(\kappa) , \SS(\lam) ) \cong 
\Hom_{n} (\SS(\kappa) , \Resi \SS(\lam) )
\eq
By proposition~\ref{res metric}, if $\hae(\kappa)\geq m+2$ then this space is empty
(since it is empty on the left by the inductive assumption). On the other
hand, by assumption, all the indecomposable summands $P_{\chi}$ of 
$ \Resi\SS({\lam})$ are simple except possibly those with
$\hae(\chi)=m+1$. But suppose some $\SS(\phi)$ a submodule of
$P_{\chi}$. Then it is simple by equation~(\ref{FR}) and
quasi--heredity (the only possible maps into it would show up in
equation~(\ref{FR}) for some $\kappa$), and so $\phi=\chi$, and
indeed $P_{\chi}$ is simple by Brauer--Humphreys reciprocity and
quasi--heredity. Thus in particular our $\SS(\nu)=P_{\nu}=S_{\nu}$. 
\hfill \Qed


\section{Exceptional representations}\label{ss4}
\newcommand{\edge}{{\mathfrak E}_}%
\newcommand{\vertex}{{\mathfrak V}_}%
\newcommand{\coxeterA}{A_}%
\newcommand{\coxeterAh}{\hat{A}_}%
\newcommand{\nnset}[1]{\{1,2,..,#1\}}%

\newcommand{\sse}{{\mbox{e}}}%
Fix $Q \in \N$ and let $V_Q$ denote the $k$--space with basis 
$\{ \sse_1 , \ldots , \sse_{Q} \}$.
Then for each ordinary partition algebra $P_n(Q)$ there a
$Q^n$--dimensional {\em Potts representation} $R$ acting on
$V=V_Q^{\otimes n}$.  This is described, for example, in
\cite{MartinRollet98corrigendum} (see also below).  The underlying
idea is that each basis element of $V$ is a colouring of the $n$
lattice sites ($\{1,...,Q\}$ being the colours), and that such a
colouring is `consistent' with a given partition of the sites only if
the elements within each part in the partition have the same colour. 
Fixing $\QQ$ a function from $T$ to $\N$, the tensor product algebra
$\bigotimes_{t} P_n(Q_t)$ may be equipped with a corresponding product
representation.  It will be evident that this representation is
$\left( \prod_t Q_t \right)^n$--dimensional.  This is to say, it acts
on a space isomorphic to that of the Potts representation of $P_n(\prod_t Q_t )$
(recall the `diagonal' isomorphism of $P_n(\prod_t Q_t )$ to a
subalgebra of our algebra).  
We will again write $R$ both for this representation of
$\bigotimes_{t} P_n(Q_t)$ and for its
restriction to $P_n^{(T)}$. 
\begin{theo}
For any $\QQ \in \N^d$ there exists an $m\in \N$ such that
$P_n^{(T)}(\QQ)$ is non--semisimple for all $n>m$. 
\end{theo}
{\em Proof:} The proof is a straightforward generalization of the
corresponding result for the ordinary partition algebra as in
\cite{MartinSaleur94b}. In short, one compares the growth rate (as $n$
increases) of
dimension of standard module $\SS({\spineweight})$ with that of the Potts
representation $R$, having established that, if semisimple, the former is
a summand of the latter (consider $R(\Xel_{\spineweight})$).
\hfill\Qed 
\newline

\subsection{Application to physical models}\label{physics}

For $H$ a graph let $\edge{H}$ (resp. $\vertex{H}$) denote its edge
(resp. vertex) set. 
Define graph $\coxeterA{l}$ by $\vertex{\coxeterA{l}}=\nnset{l}$
and $\edge{\coxeterA{l}}= \{ \{i,j \} | i-j=\pm1 \}$, and 
graph $\coxeterAh{l}$ by $\vertex{\coxeterAh{l}}=\nnset{l}$
and $\edge{\coxeterAh{l}}= \{ \{i,j \} | i-j\equiv \pm1 \mbox{ mod. $l$} \}$,
 as usual. 

Fix $Q \in \N$. 
A {\em spin configuration} on graph $H$ is an element of
$S_H = \hom(\vertex{H},\nnset{Q})$ 
(i.e. a colouring of the vertices of $H$ from $Q$ colours). 
For $s \in S_H$ let $s_i$ denote the `colour' of vertex $i$.

Given a graph $H$ 
let graph $G=H \times \coxeterA{l}$ be given by 
$\vertex{H \times \coxeterA{l}}=\vertex{H}\times\vertex{\coxeterA{l}}$ and 
$\{ (a,b),(c,d) \} \in \edge{H \times \coxeterA{l}}$ 
if either $a=c$ and $\{b,d \} \in \edge{ \coxeterA{l}}$ or $b=d$ and 
$\{ a,c\} \in \edge{H}$. 
We may break graph $G$ up in to $l$ `lateral' layers $H_k$ 
($k=1,2,..,l$) each isomorphic to $H$, and
$l-1$ `transverse' layers $H_{k,k+1}$ each consisting of the edges between a
certain pair of lateral layers. 


A {\em $Q$--state model} on $G$ may be constructed as follows 
\cite{Baxter82,Martin91}. 
An {\em edge Hamiltonian} $ {\cal H}$ is a map from 
$\hom(\vertex{\coxeterA{2}},\nnset{Q})$ to the real numbers 
(here $\coxeterA{2}$ stands for any of the 2 vertex subgraphs of $G$
containing precisely one edge). 
From this we form {\em graph Hamiltonian} 
$$ {\cal H}_G = \sum_{\{i,j\}\in \edge{ G}}
{\cal H}(s_i,s_j) . $$
Given such a function and a complex variable $\beta$ 
define the {\em partition function} 
$$ Z = Z(\beta) = \sum_{s \in S_G} \exp(\beta {\cal H}_G(s)) . $$
In \SM\ one  wants to know the analytic structure of $Z(\beta)$.
\footnote{Actually $Z$ will only be physically interesting, 
  within the framework described, for certain carefully
  chosen $G$ and ${\cal H}$; and
  interesting constructions are possible outside of this framework. 
  Neither point need concern us here.} 
If the Hamiltonian 
${\cal H}(s_i,s_j) $ is integer valued then, regarding $Z$ as a
function of $\exp(\beta)$, this structure consists of roots. 
(More generally of limits of distributions of roots as $l$ and $H$
grow large in some way.) 


For $s,t \in S_H$ let $T_{st}^{(l)}$ denote the partition function in case the
colours of vertices in $H_1$ (resp. $H_l$) have been {\em fixed} to $s$
(resp. $t$). 
Note that, given a bijection $\vertex{H} \leftrightarrow \nnset{n}$, 
$S_H$
is a basis for $V$ in a natural way. 
Thus $T^{(l)}$ is a matrix in $\End(V)$. 
In case $l=1$ the
first and last layers are identified, and $T^{(1)}$ is a diagonal matrix. 
We may similarly introduce a matrix $T^{(-)}$ 
of partition functions associated
to a single transverse layer $H_{k,k+1}$. Here $s$ determines the
colour of vertices in layer $k$, and $t$ in layer $k+1$. 
These layer matrices are called {\em transfer matrices}. 
Note that we may use them to build up $T^{(l)}$ one layer at a time:
$$T^{(l)} = \left( T^{(1)} T^{(-)} \right)^l  T^{(1)} . $$  
Similarly, we can introduce matrices which will build
up $T^{(1)}$ and $T^{(-)}$ one {\em edge} at a time.


For certain choices of ${\cal H}(s_i,s_j) $ the representation $R$ can
be used to build these transfer matrices. 
The representation $R$ of the full tensor product algebra 
may be used to build 
transfer matrices for models with a significantly greater
variety of Hamiltonian than can the restriction of $R$ to $P_n(\prod_t Q_t )$.  
The point about $P_n^{(T)}$ is that $R$ may be restricted to it {\em
  without} losing the greater variety.  
A suitable concrete example will
illustrate these points. 
Let $\QQ=(Q_1,Q_2)=(2,3)$. 
Let $D_{j}$ be the $j \times j$ matrix with all entries 1, 
and $I_{j}$ be the $j \times j$ unit matrix. 
Then
\[
R((A^i,A^i))=I_{6^{i-1}} \otimes D_{Q_1} \otimes D_{Q_2} \otimes I_{6^{n-i-2}}
\]
\[
R((A^i,1))=I_{6^{i-1}} \otimes D_{Q_1} \otimes I_{Q_2} \otimes I_{6^{n-i-2}}
\]
(cf. \cite[\S12.4]{Baxter82}). 


The colour--valued variables $s_i$ associated to the vertices $i$ of $G$ are
called {\em spins}. 
Note that $R$ above acts on a space appropriate for 
a $Q_1Q_2=$6--state model.  
Considering spins $s_i$ as taking values from the integers mod. 6, 
then a $Z_6$--symmetric model 
\cite{Baxter82,EinhornSavitRabinovici} is one in which 
the edge Hamiltonian takes the form
${\cal H}(s_i,s_j)=f(s_i - s_j)$, 
where $f$ may be any real valued even function 
\footnote{Directed graphs and non--even functions produce {\em chiral} models, 
which have also been studied, but need not concern us here.}.
Suppose we want to build the transfer matrix for such a model 
with $f(0)=2$, $f(1)=f(2)=1$, $f(3)=0$ 
(such as the model examined in \cite[p.305]{Martin91}). 
Then the transfer matrix factor building in 
the contribution of that edge to the Hamiltonian is
\[
R\left( \left(  1+(e^{\gamma} -1)(A^{ij},A^{ij}) \right) 
        \left(  1+(e^{\beta} -1))(1,A^{ij}) \right) \right)
\] \[
=
R\left( \left(  
    1
    +(e^{\beta} -1))(1,A^{ij}) 
    +e^{\beta}(e^{\gamma} -1)(A^{ij},A^{ij}) 
  \right) \right)
\]
if lattice edge $(i,j)$ is in the transfer matrix layer 
(i.e. it is in some $H_k$), and 
\[
R\left( \left(  
    e^{\beta}(e^{\gamma} -1)1
    +(e^{\beta} -1)(A^i,1) 
    +(A^i,A^i) 
  \right) \right)
\]
if $(i,j)$ is perpendicular to the transfer matrix layer; 
where $\gamma = 2$, $\beta = -1$. 

\begin{figure}
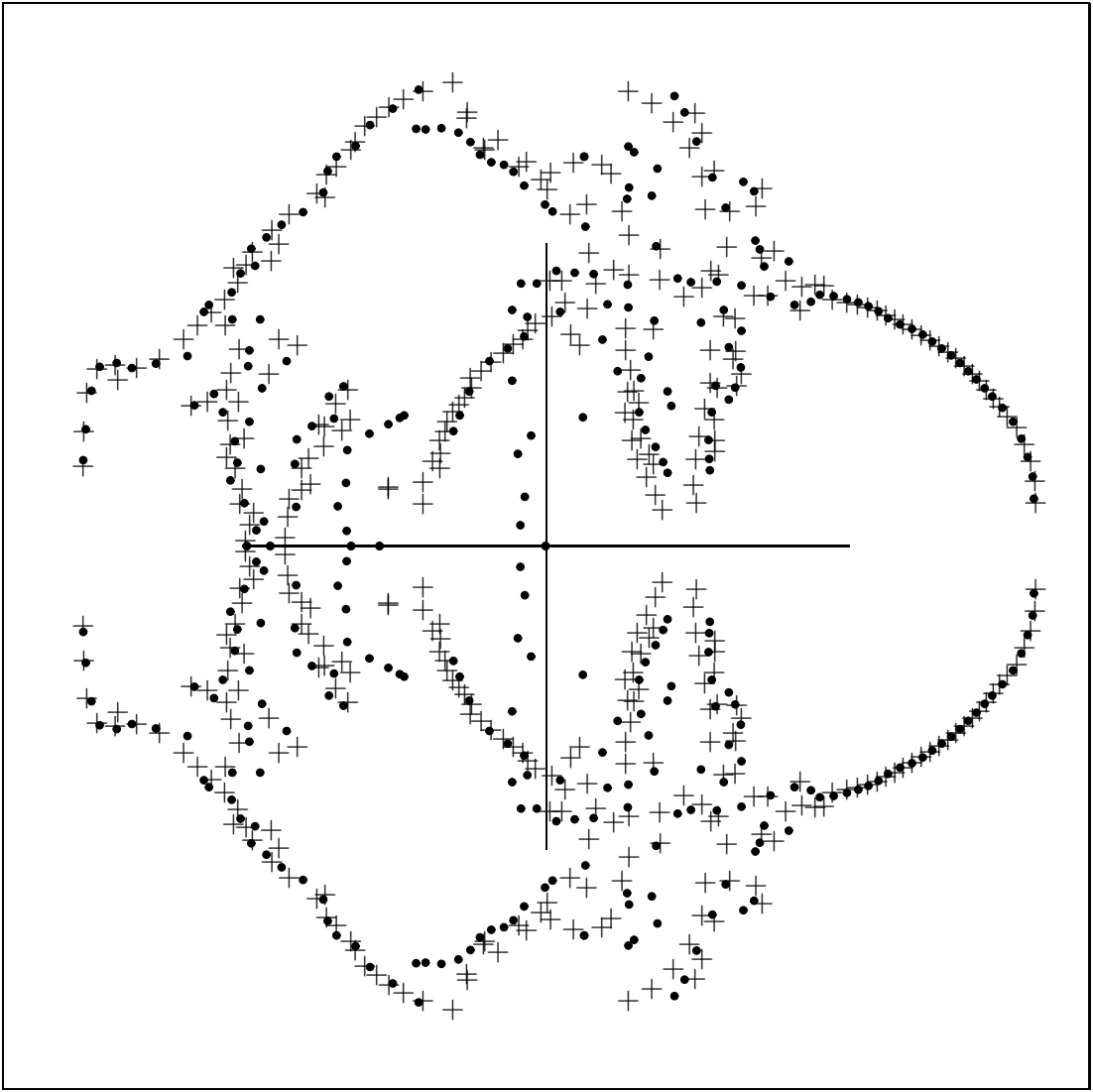

\input ./pfn_mp19/gnuploTeX/q4n7n8l2sq310.tex
\caption{\label{zeros01} Zeros of the partition function in
  $e^{\beta}$ for $7\times 9$ (dots) and $8\times 9$ (crosses)
  square lattice systems with $Q_1=Q_2=2$. 
  }
\end{figure}
\begin{figure}
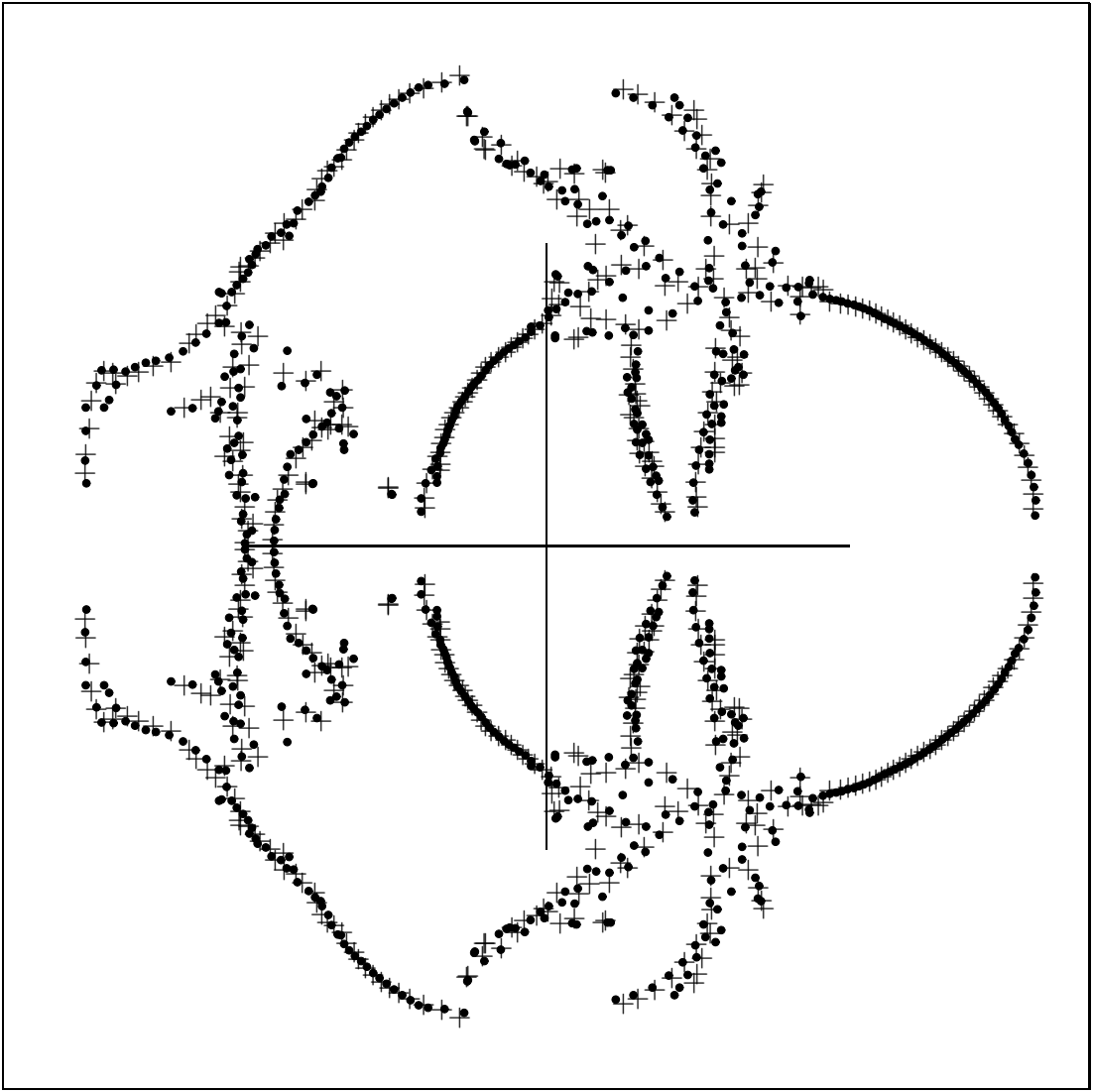

\input ./pfn_mp19/gnuploTeX/q4n8n10l3sq310.tex
\caption{\label{zeros02} Zeros of the partition function in
  $e^{\beta}$ for $8\times 13$ (crosses) and $10\times 13$ (dots)
  square lattice systems with $Q_1=Q_2=2$. 
  }
\end{figure}

The analytic structure of the partition function $Z$ for such models
reveals physical {\em phase transitions} 
as points on the positive $\exp(\beta)$--axis which are
approached by distributions of roots as the system size grows large
(see \cite{Martin91} for details). 
In computing this structure we are, at base, using the algebra as nothing more
than a formal aid to computation 
(the globalization functor $\Gg$ gives 
explicit control of the thermodynamic, i.e. large $n$, limit 
as a global algebra limit,
but there is no way yet to use this in specific calculations). 
None the less, the results are of
considerable interest both in the way they model phase transitions and
in what they tell us about specific physical systems. 
We will postpone such interpretation to a separate paper but, 
in the renaissance spirit, include a couple of 
results here (in the form of distributions of roots), for illustration. 
(There has been renewed interest in results of this kind recently, 
although mostly obtained by brute force calculation,
see \cite{MatveevShrock96b}.) 
Figures~\ref{zeros01} and~\ref{zeros02} 
show the roots of $Z$ in $\exp(\beta)$ for 
$Q_1=Q_2=2$ and $f(0)=3$, $f(1)=1$ and $f(2)=0$, for various graphs. 
The results for various different graphs
are shown, as a rough indication of the sense in which the distribution of
roots may approach a limit as $G$ becomes large. 
\footnote{A better indication is perhaps given by the {\em animated}
  combined version of these figures available online from
  {\tt
    www.city.ac.uk/mathematics/algebra/statinfo/StatMechIntro2.html} 
  or from PPM's homepage. }
The graphs concerned are given by 
$(H=\coxeterAh{7},l=9)$ and $(H=\coxeterAh{8},l=9)$ and 
$(H=\coxeterAh{8},l=13)$ and $(H=\coxeterAh{10},l=13)$ 
(i.e. these examples are models of two--dimensional physical systems). 



\section{Discussion}
There remain a number of interesting open problems on the
representation theory side. 
We have not determined the exceptional structure of
$P_n^{(\Nset{2})}$ (cf. $P_n$, where complete results are known), 
nor addressed its representation theory in positive characteristic 
(nothing is known here for $P_n$ either, although it {\em is} known to be a
hard problem \cite{MartinWoodcock98}). 
As for more exotic posets $T$, 
we do not know even the generic structure of $P_n^{(T)}$ 
(indeed we do not even know if it {\em has} a generic structure). 

As a quasi--hereditary algebra, it is natural to ask if
$P_n^{(\Nset{2})}$ has tilting modules \cite{Donkin93}. 
These have yet to be constructed even for $P_n$.

The assumptions of this paper restrict us to finite $T$, however there
is no serious obstruction to relaxing this constraint. Infinite $T$
lacks, for the moment, physical motivation, however it offers some
intriguing possibilities cf. \cite{VershikYakubovich98} 
(and references therein) for example.

\[ \]
{\bf Acknowledgements.}
PPM would like to thank Anton Cox for invaluable discussions, and
EPSRC  for partial financial support. 

\appendix
\section*{Appendix: Rank of $P_n^{(\Nset{2})}$}

The degree of $\Pbase_n^{(\Ttwo)}$, 
or the set 
$\Ee(N)^{(\Ttwo)}$ 
of $\Ttwo$--ramified partitions of any finite set $N$, 
may be computed as follows. 
The {\em shape} of a partition (of $N$ say) 
is the list of degrees of its parts, expressed as an integer partition. 
For $\mu$ such a shape, let 
$\nowithshape_{\mu}$ be the number of partitions with that shape. 
Let $\Bell{n}$ be the number of partitions of $\Nset{n}$. 
Then 
\[
|p_n^{(\Ttwo)}| = \sum_{\mu} \nowithshape_{\mu} \prod_i \Bell{\mu_i}
\]
The first few of these are
\[
p_1^{(T)} : \hspace{1in}
\begin{array}{cccc}
\nowithshape_{\mu} \times \prod && \mu
\\ \hline
1 \times 1 &&(1,1) \\
1 \times 2 &&(2)
\\ \hline
\sum = 3
\end{array}
\]
\[
\begin{array}{cccc}
\nowithshape_{\mu} \times \prod && \mu
\\ \hline
1 \times 1 &&(1,1,1) \\
3 \times 2 &&(2,1) \\
1 \times 5 &&(3)
\\ \hline
\sum = 12
\end{array}
\]
\[
p_2^{(T)} : \hspace{1in}
\begin{array}{cccc}
\nowithshape_{\mu} \times \prod && \mu
\\ \hline
1 \times 1 &&(1,1,1,1) \\
6 \times 2 &&(2,1,1) \\
4 \times 5 &&(3,1) \\
3 \times 2 \times 2 &&(2,2) \\
1 \times 15 &&(4)
\\ \hline
\sum = 60
\end{array}
\]
\[
\begin{array}{cccc}
\nowithshape_{\mu} \times \prod && \mu
\\ \hline
1 \times 1 &&(1,1,1,1,1) \\
10 \times 2 &&(2,1,1,1) \\
10 \times 5 &&(3,1,1) \\
15 \times 2 \times 2 &&(2,2,1) \\
5 \times 15 &&(4,1) \\
10 \times 2 \times 5 &&(3,2) \\
1 \times 52 &&(5)
\\ \hline
\sum = 358
\end{array}
\]
These figures may be compared with the summed squares of dimensions
of generic irreducible representations of the corresponding algebras
in figure~\ref{fig2}.

\bibliographystyle{amsplain}
\bibliography{new31,main,emma,citePPM}
\end{document}

%% file: martinew.tex
\providecommand{\noglossaryignore}[1]{}
\newcommand{\globalglossaryentry}[3]{\makebox[1.5in][l]{\tt $\backslash${#1}} 
\makebox[1.1in][l]{{$#2$}} \makebox[2.5in][l]{{#3}}\newline} 
\newcommand{\newcommandabbreviation}[3]{\newcommand{#1}{#2}%
\noglossaryignore{\globalglossaryentry{#3}{#2}{}}}
\newcommand{\renewcommandabbreviation}[3]{\renewcommand{#1}{#2}%
\noglossaryignore{\globalglossaryentry{#3}{#2}{}}}
\newcommand{\newcommandmacro}[4]{\newcommand{#1}{#2}%
\noglossaryignore{\globalglossaryentry{#3}{#2}{#4}}}
\newcommand{\gge}[3]{\noglossaryignore{\globalglossaryentry{#1}{#2}{#3}}}
\newcommand{\myaddress}%
{\parbox{3in}{\footnotesize \begin{center} 
Mathematics Department, City University, \\  
Northampton Square, London EC1V 0HB, UK.\end{center}}}
    
\newcounter{minidef}[section]

\newcounter{minicapt}

\newtheorem{theo}{Theorem}         
\newtheorem{de}{Definition}     \newtheorem{pr}{Proposition}

\noglossaryignore{GREEK ETC.\newline}
\newcommandabbreviation{\e}{\epsilon}{e}        
\newcommandabbreviation{\lam}{\lambda}{lam}  
\newcommandabbreviation{\la}{\langle}{la}        
\newcommandabbreviation{\ran}{\rangle}{ran}
\newcommandabbreviation{\ha}{\#}{ha}             
\newcommandabbreviation{\rmap}{\rightarrow}{rmap}
\newcommandabbreviation{\aaa}{\alpha}{aaa}        
\newcommandabbreviation{\ab}{\alpha,\beta}{ab}
\newcommandabbreviation{\aab}{a(\ab )}{aab}       
\noglossaryignore{\newline RINGS\newline}
\newcommandabbreviation{\HH}{H \!\!\! I}{HH}               
\newcommandabbreviation{\C}{\mathbb C}{C}
\newcommandabbreviation{\N}{\mathbb N}{N}   
\newcommandabbreviation{\Z}{\mathbb Z}{Z}      
\renewcommandabbreviation{\Re}{\mathbb R}{Re}
\newcommandabbreviation{\R}{{\mathbb R}}{R}
\newcommandabbreviation{\Q}{\mathbb Q }{Q}
\renewcommandabbreviation{\H}{\mathbb H }{H}
\noglossaryignore{\newline SYMMETRIC GROUP\newline}
\def\Sym(#1){\Sigma(#1)}                   
\gge{Sym(-)}{\Sym(-)}{symmetric group on - objects}
\def\Sy(#1){\Sigma_{#1}}                   
\gge{Sy(-)}{\Sy(-)}{symmetric group irreducible -}
\def\sym(#1){\mbox{\LARGE s}(#1)}        
\gge{sym(-)}{\sym(-)}{symmetric group on - objects (variant)}
\def\sy(#1){\mbox{\LARGE s}({#1})}        
\gge{sy(-)}{\sy(-)}{symmetric group irreducible - (variant)}
\newcommandmacro{\cs}{\C \, \sy(n)}{cs}{symmetric group algebra over $\C$}
\noglossaryignore{\newline PARTITIONS/SETS\newline}
\newcommand{\Nset}[1]{\underline{#1}}
\gge{Nset\{-\}}{\Nset{-}}{set of natural numbers to -} 
\def\nset(#1){ \{ #1 \}_{ \underline{n} }} 
\gge{nset(-)}{\nset(-)}{a set $-\times\Nset$}
\def\ul(#1){_{\underline{#1}}}             
\gge{ul(-)}{{}\ul(-)}{subscript underline -}
\def\Ee(#1){{\bf E}_{#1}}                  
\gge{Ee(-)}{\Ee(-)}{set of equivalence relations on set -}
\def\Eee(#1){{\bf E}_{\{ #1 \}_{\underline{n}}}}   
\gge{Eee(-)}{\Eee(-)}{ditto for nset}
\def\Een(#1,#2){{\bf E}_{\{ #1 \}_{\underline{#2}}}}   
\def\Ssn(#1,#2){{\bf S}_{\{ #1 \}_{\underline{#2}}}}   
\def\Ss(#1){{\bf S}_{#1}}                  
\def\Sss(#1){{\bf S}_{\{ #1 \}_{\underline{n}}}}   
\def\bbc(#1){((\beta_1)(\beta_2)...(\beta_{#1}))}      
\newcommandmacro{\Ln}{{\Gamma}^{n}}{Ln}{large index set}
\newcommandmacro{\LnQ}{{\Gamma}^{n}_Q}{LnQ}{index set}
\newcommandmacro{\Zz}{\zeta}{Zz}{`shape' function}
\noglossaryignore{\newline PARTITION ALGEBRA\newline}
\def\ka(#1){\kappa_{#1}}                   
\def\Sm(#1){\Sigma_{#1}}                   
\newcommandmacro{\com}{\bullet}{com}{bullet composition}
\newcommandmacro{\enm}{\; e^n(\! m\! ) \;}{enm}{product of idempotents}
\def\Ai(#1){ A^{ #1 \cdot } }              
\def\Aij(#1,#2){ A^{ #1  #2 } }            
\newcommandmacro{\One}{\mbox{\bf $1 \!\!\! 1$}}{One}{algebra unit 1}
\newcommandmacro{\Bp}{B_p}{Bp}{partition basis}
\def\Bb(#1){B_p[#1]}                       
\def\Pp(#1){P_n[#1]}                       
\def\Ps(#1){P_n[#1] \! /}                  
\newcommandmacro{\Ph}{\hat{P}}{Ph}{P hat  algebra}
\def\Is(#1){\sim^{#1}}                     
\noglossaryignore{\newline MODULES\newline}
\def\Wm(#1){{\cal S}_{#1}}                 
\gge{Wm(-)}{\Wm(-)}{Weyl module with index -}
\def\wm(#1,#2){{}_{#1}{\cal S}_{#2}}       
\gge{wm(-1,-)}{\wm(-1,-)}{Weyl module with index -}
\def\Ind(#1,#2,#3){\mbox{Ind}_{#1}^{#2}#3} 
\gge{Ind(-1,-2,-)}{\Ind(-1,-2,-)}{induction}
\def\Res(#1,#2,#3){\mbox{Res}_{#1}^{#2}#3} 
\gge{Res(-1,-2,-)}{\Res(-1,-2,-)}{restriction}
\newcommandabbreviation{\weyl}{standard}{weyl}
\newcommandabbreviation{\mod}{\mbox{mod}}{mod}
\newcommandabbreviation{\head}{\mbox{head }}{head}
\newcommandabbreviation{\Weyl}{Weyl}{Weyl}
\def\SS(#1){{\cal S}_{#1}}                 
\gge{SS(-)}{\SS(-)}{Specht/Weyl module index -}
\def\LL(#1){{\cal L}_{#1}}                 
\gge{LL(-)}{\LL(-)}{simple module index -}
\noglossaryignore{\newline FUNCTORS/MAPS\newline}
\newcommandmacro{\Gg}{{\cal G}}{Gg}{G Functor}
\newcommandmacro{\Fg}{{\cal F}}{Fg}{F Functor}
\newcommandmacro{\ra}{\rightarrow}{ra}{}
\def\ses(#1,#2,#3){0\ra #1 \ra #2 \ra #3 \ra 0}   
\gge{ses(1,2,-)}{\ses(1,2,-)}{\hspace{.5in} short exact sequence}
\def\starr(#1){ \stackrel{ #1 }{\longrightarrow} }
\gge{starr(-)}{\starr(-)}{}
\newcommandmacro{\doublerightarrow}{\; -\!\!\! -\!\!\!\!\!\! \gg \;}
{doublerightarrow}{}
\noglossaryignore{\newline PARTITION ALGEBRA MAPS\newline}
\newcommandmacro{\smap}{s}{smap}{`inclusion' map}
\newcommandmacro{\tmap}{t}{tmap}{$ P_n -> S_n$}
\newcommandmacro{\pmap}{\psi}{pmap}{$ S_n -> P_n $}
\noglossaryignore{\newline MISC.\newline}
\def\Amap(#1){{\cal A}_{#1}}               
\gge{Amap(-)}{\Amap(-)}{}
\def\Rr(#1){R_{#1}}                        
\gge{Rr(-)}{\Rr(-)}{restriction of E}
\def\Cr(#1){C_{#1}}                        
\gge{Cr(-)}{\Cr(-)}{restriction of E to N}
\newcommandmacro{\Tm}{{\cal T}}{Tm}{Transfer Matrix}
\def\On(#1){{\cal I}_{#1}}
\gge{On(-)}{\On(-)}{}
\newcommandmacro{\UU}{\underline{\sqcup}}{UU}{}  
\newcommandmacro{\UUU}{\sqcup}{UUU}{}  
\newcommandmacro{\Vq}{V_Q^{\otimes n}}{Vq}{Potts config. space}
\def\bs(#1,#2){\mbox{{\Large $\ast$}}^{#1}_{#2}}  
\gge{bs(-,-)}{\bs(-,-)}{general plumbing multiplier}
\newcommand{\ignore}[1]{}
\gge{ignore\{-\}}{\ignore{-}}{ignore argument!}
\def\choo(#1,#2){ \left( \begin{array}{c} #1 \\ #2 \end{array} \right) } 
\gge{choo(-1,-)}{\choo(-1,-)}{choose}
\newcommand{\Qed}{$\Box$}
\gge{Qed}{\mbox{\Qed}}{QED}
\def\staq(#1){\stackrel{#1}{=}}            
\gge{staq(-)}{\staq(-)}{}
\def\stam(#1){\stackrel{#1}{\rightarrow}}  
\gge{stam(-)}{\stam(-)}{}
\def\mat{ \left( \begin{array} }    
\def\tam{ \end{array}  \right) }
\gge{mat/tam}{...}{matrix delimiters}
\newcommand{\beq}{\begin{equation} }
\def\eql(#1){ \begin{equation} \label{#1} 
}
\newcommand{\eq}{\end{equation} }
\def\eqal(#1){\begin{eqnarray} \label{#1} }
\def\eqa{\end{eqnarray} }
\def\lab(#1){\label{#1}
}
\def\prl(#1){ \begin{pr} \label{#1} 
}
\def\del(#1){ \begin{de} \label{#1} 
}
\gge{smeq\{-\}}{...}{small equation}
\gge{fneq\{-\}}{...}{very small equation}
\noglossaryignore{\newline HECKE/BLOB\newline}
\newcommandmacro{\Hnq}{H_n(q)}{Hnq}{ * freestanding symbol}
\newcommandmacro{\Hn}{H_n}{Hn}{      *-mod etc.}
\newcommandmacro{\A}{{\cal A}}{A}{}
\newcommandmacro{\Cwts}{C}{Cwts}{}
\newcommandmacro{\CA}{{\cal A}}{CA}{}

\newcommandmacro{\calA}{{\cal A}}{calA}{}
\newcommandmacro{\modi}{\mbox{Mod} }{modi}{was mod not modi!}
\newcommandmacro{\Wgen}{{\Bbb S}}{Wgen}{}
\def\ol(#1){\overline{#1}}
\newcommandmacro{\st}{\mbox{St}}{st}{}
\def\CMult(#1,#2){(#1:#2)}
\def\CM(#1,#2){( #1 : #2 )}
\def\FMult#1,#2{(#1:#2)}
\def\CF#1,#2{(#1:#2)}

\newcommandmacro{\Top}{\mbox{Top}}{Top}{}
\newcommandmacro{\Soc}{\mbox{Soc}}{Soc}{}
\newcommandmacro{\Head}{\mbox{Head}}{Head}{}
\newcommandmacro{\Filt}{{\cal F}}{Filt}{}
\newcommandmacro{\Mod}{\mbox{mod}}{Mod}{}
\newcommandmacro{\Resi}{\mbox{Res }}{Resi}{was without i!}
\newcommandmacro{\Indi}{\mbox{Ind }}{Indi}{was without i!}
\def\RR(#1,#2){R^{#1}_{#2}}   
\def\TT(#1,#2){T^{#1}_{#2}}   

\def\Hom{\mbox{Hom}}

\newcommandmacro{\Ann}{\mbox{Ann}}{Ann}{}
\newcommandmacro{\Cen}{\mbox{Cen}}{Cen}{}
\newcommandmacro{\End}{\mbox{End}}{End}{}
\newcommandabbreviation{\semisimple}{semisimple}{semisimple}
\newcommandabbreviation{\Bratteli}{Bratteli}{Bratteli}
\newcommandabbreviation{\JBC}{Jones Basic Construction}{JBC}
\newcommandabbreviation{\pa}{partition algebra}{pa}
\newcommandabbreviation{\TM}{transfer matrix}{TM}
\newcommandabbreviation{\PM}{Potts model}{PM}
\newcommandabbreviation{\QSC}{quantum spin chain}{QSC}
\newcommandabbreviation{\Hamiltonian}{Hamiltonian}{Hamiltonian}
\newcommandabbreviation{\YS}{Young symmetrizer}{YS}



%% file: xfig/prodex3.eepic
\setlength{\unitlength}{0.00041667in}
\begingroup\makeatletter\ifx\SetFigFont\undefined%
\gdef\SetFigFont#1#2#3#4#5{%
  \reset@font\fontsize{#1}{#2pt}%
  \fontfamily{#3}\fontseries{#4}\fontshape{#5}%
  \selectfont}%
\fi\endgroup%
{\renewcommand{\dashlinestretch}{30}
\begin{picture}(8574,5510)(0,-10)
\put(1399,1958){\blacken\ellipse{300}{300}}
\put(1399,1958){\ellipse{300}{300}}
\put(1399,458){\blacken\ellipse{300}{300}}
\put(1399,458){\ellipse{300}{300}}
\path(1924,1958)(1941,1905)(1956,1851)
	(1969,1797)(1981,1743)(1991,1691)
	(1999,1641)(2006,1593)(2013,1548)
	(2018,1506)(2022,1466)(2026,1429)
	(2029,1395)(2031,1363)(2033,1333)
	(2034,1304)(2035,1277)(2036,1252)
	(2036,1226)(2037,1202)(2036,1177)
	(2036,1152)(2035,1126)(2034,1099)
	(2033,1071)(2031,1041)(2029,1010)
	(2026,976)(2022,939)(2018,900)
	(2013,858)(2006,814)(1999,767)
	(1991,718)(1981,666)(1969,614)
	(1956,561)(1941,509)(1924,458)
	(1903,405)(1881,356)(1858,312)
	(1834,272)(1811,237)(1788,206)
	(1767,178)(1746,155)(1726,134)
	(1708,117)(1690,102)(1673,89)
	(1658,77)(1642,68)(1628,59)
	(1613,52)(1599,45)(1585,40)
	(1570,34)(1555,30)(1540,25)
	(1524,22)(1506,18)(1488,15)
	(1469,13)(1448,12)(1426,12)
	(1403,14)(1378,17)(1353,23)
	(1326,32)(1300,45)(1274,61)
	(1249,83)(1226,110)(1204,140)
	(1186,173)(1169,207)(1154,242)
	(1142,276)(1132,308)(1124,340)
	(1118,369)(1113,397)(1109,423)
	(1106,448)(1105,472)(1104,494)
	(1104,516)(1105,537)(1105,558)
	(1106,580)(1107,602)(1108,626)
	(1109,651)(1110,678)(1110,709)
	(1110,742)(1110,778)(1110,818)
	(1109,863)(1107,912)(1106,965)
	(1104,1022)(1102,1083)(1101,1148)
	(1099,1215)(1099,1283)(1100,1355)
	(1101,1425)(1103,1492)(1104,1555)
	(1106,1613)(1108,1666)(1109,1714)
	(1110,1758)(1110,1796)(1110,1831)
	(1110,1862)(1109,1891)(1108,1916)
	(1107,1940)(1106,1962)(1105,1983)
	(1104,2004)(1104,2024)(1104,2045)
	(1105,2067)(1107,2090)(1110,2115)
	(1115,2141)(1121,2168)(1128,2198)
	(1138,2228)(1150,2260)(1165,2293)
	(1182,2325)(1202,2356)(1224,2384)
	(1249,2408)(1274,2425)(1300,2438)
	(1326,2446)(1353,2451)(1378,2454)
	(1403,2454)(1426,2452)(1448,2450)
	(1469,2446)(1488,2442)(1506,2437)
	(1524,2432)(1540,2427)(1555,2421)
	(1570,2415)(1585,2409)(1599,2402)
	(1613,2394)(1628,2386)(1642,2376)
	(1658,2365)(1673,2353)(1690,2338)
	(1708,2322)(1726,2302)(1746,2280)
	(1767,2254)(1788,2225)(1811,2192)
	(1834,2154)(1858,2112)(1881,2065)
	(1903,2013)(1924,1958)
\put(1474,2183){\makebox(0,0)[lb]{\smash{{{\SetFigFont{6}{7.2}{\rmdefault}{\mddefault}{\updefault}$2$}}}}}
\put(1324,683){\makebox(0,0)[lb]{\smash{{{\SetFigFont{6}{7.2}{\rmdefault}{\mddefault}{\updefault}$2'$}}}}}
\put(2374,4958){\blacken\ellipse{300}{300}}
\put(2374,4958){\ellipse{300}{300}}
\put(2374,3458){\blacken\ellipse{300}{300}}
\put(2374,3458){\ellipse{300}{300}}
\thicklines
\path(2374,4958)(2374,3458)
\thinlines
\path(2899,4958)(2916,4905)(2931,4851)
	(2944,4797)(2956,4743)(2966,4691)
	(2974,4641)(2981,4593)(2988,4548)
	(2993,4506)(2997,4466)(3001,4429)
	(3004,4395)(3006,4363)(3008,4333)
	(3009,4304)(3010,4277)(3011,4252)
	(3011,4226)(3012,4202)(3011,4177)
	(3011,4152)(3010,4126)(3009,4099)
	(3008,4071)(3006,4041)(3004,4010)
	(3001,3976)(2997,3939)(2993,3900)
	(2988,3858)(2981,3814)(2974,3767)
	(2966,3718)(2956,3666)(2944,3614)
	(2931,3561)(2916,3509)(2899,3458)
	(2878,3405)(2856,3356)(2832,3312)
	(2809,3272)(2786,3237)(2763,3206)
	(2741,3178)(2720,3155)(2700,3134)
	(2681,3117)(2662,3102)(2645,3089)
	(2629,3077)(2613,3068)(2598,3059)
	(2583,3052)(2568,3045)(2553,3040)
	(2538,3034)(2522,3030)(2506,3025)
	(2490,3022)(2472,3018)(2454,3015)
	(2434,3013)(2414,3012)(2392,3012)
	(2369,3014)(2345,3017)(2321,3023)
	(2296,3032)(2271,3045)(2247,3061)
	(2224,3083)(2203,3110)(2184,3140)
	(2168,3173)(2154,3207)(2142,3242)
	(2133,3276)(2126,3308)(2120,3340)
	(2117,3369)(2114,3397)(2113,3423)
	(2113,3448)(2114,3472)(2116,3494)
	(2118,3516)(2121,3537)(2124,3558)
	(2127,3580)(2131,3602)(2134,3626)
	(2137,3651)(2140,3678)(2143,3709)
	(2145,3742)(2147,3778)(2149,3818)
	(2150,3863)(2151,3912)(2151,3965)
	(2151,4022)(2150,4083)(2150,4148)
	(2149,4215)(2149,4283)(2149,4355)
	(2150,4425)(2150,4492)(2151,4555)
	(2151,4613)(2150,4666)(2149,4714)
	(2148,4758)(2146,4796)(2144,4831)
	(2141,4862)(2138,4891)(2135,4916)
	(2131,4940)(2127,4962)(2124,4983)
	(2121,5004)(2118,5024)(2115,5045)
	(2114,5067)(2113,5090)(2113,5115)
	(2115,5141)(2118,5168)(2123,5198)
	(2130,5228)(2139,5260)(2151,5293)
	(2165,5325)(2182,5356)(2202,5384)
	(2224,5408)(2247,5425)(2271,5438)
	(2296,5446)(2321,5451)(2345,5454)
	(2369,5454)(2392,5452)(2414,5450)
	(2434,5446)(2454,5442)(2472,5437)
	(2490,5432)(2506,5427)(2522,5421)
	(2538,5415)(2553,5409)(2568,5402)
	(2583,5394)(2598,5386)(2613,5376)
	(2629,5365)(2645,5353)(2662,5338)
	(2681,5322)(2700,5302)(2720,5280)
	(2741,5254)(2763,5225)(2786,5192)
	(2809,5154)(2832,5112)(2856,5065)
	(2878,5013)(2899,4958)
\put(2449,5183){\makebox(0,0)[lb]{\smash{{{\SetFigFont{6}{7.2}{\rmdefault}{\mddefault}{\updefault}$3$}}}}}
\put(2449,3683){\makebox(0,0)[lb]{\smash{{{\SetFigFont{6}{7.2}{\rmdefault}{\mddefault}{\updefault}$3'$}}}}}
\put(499,1958){\blacken\ellipse{300}{300}}
\put(499,1958){\ellipse{300}{300}}
\put(499,458){\blacken\ellipse{300}{300}}
\put(499,458){\ellipse{300}{300}}
\put(499,4958){\blacken\ellipse{300}{300}}
\put(499,4958){\ellipse{300}{300}}
\put(499,3458){\blacken\ellipse{300}{300}}
\put(499,3458){\ellipse{300}{300}}
\put(1399,3458){\blacken\ellipse{300}{300}}
\put(1399,3458){\ellipse{300}{300}}
\put(1399,4958){\blacken\ellipse{300}{300}}
\put(1399,4958){\ellipse{300}{300}}
\put(2374,1958){\blacken\ellipse{300}{300}}
\put(2374,1958){\ellipse{300}{300}}
\put(2374,458){\blacken\ellipse{300}{300}}
\put(2374,458){\ellipse{300}{300}}
\put(7924,3758){\blacken\ellipse{300}{300}}
\put(7924,3758){\ellipse{300}{300}}
\put(7924,2258){\blacken\ellipse{300}{300}}
\put(7924,2258){\ellipse{300}{300}}
\put(5974,3758){\blacken\ellipse{300}{300}}
\put(5974,3758){\ellipse{300}{300}}
\put(5974,2258){\blacken\ellipse{300}{300}}
\put(5974,2258){\ellipse{300}{300}}
\put(6874,2258){\blacken\ellipse{300}{300}}
\put(6874,2258){\ellipse{300}{300}}
\put(6874,3758){\blacken\ellipse{300}{300}}
\put(6874,3758){\ellipse{300}{300}}
\dashline{60.000}(499,3458)(499,1958)
\dashline{60.000}(1399,3458)(1399,1958)
\thicklines
\path(1399,4958)(1399,3458)
\path(7924,3758)(7924,2258)
\thinlines
\dashline{60.000}(2374,3458)(2374,1958)
\thicklines
\path(2374,1958)(2374,458)
\thinlines
\path(124,1958)(107,1905)(92,1851)
	(79,1797)(67,1743)(57,1691)
	(49,1641)(42,1593)(35,1548)
	(30,1506)(26,1466)(22,1429)
	(19,1395)(17,1363)(15,1333)
	(14,1304)(13,1277)(12,1252)
	(12,1226)(12,1202)(12,1177)
	(12,1152)(13,1126)(14,1099)
	(15,1071)(17,1041)(19,1010)
	(22,976)(26,939)(30,900)
	(35,858)(42,814)(49,767)
	(57,718)(67,666)(79,614)
	(92,561)(107,509)(124,458)
	(145,405)(167,356)(190,312)
	(214,272)(237,237)(260,206)
	(281,178)(302,155)(322,134)
	(340,117)(358,102)(375,89)
	(390,77)(406,68)(420,59)
	(435,52)(449,45)(463,40)
	(478,34)(493,30)(508,25)
	(524,22)(542,18)(560,15)
	(579,13)(600,12)(622,12)
	(645,14)(670,17)(695,23)
	(722,32)(748,45)(774,61)
	(799,83)(822,110)(844,140)
	(862,173)(879,207)(894,242)
	(906,276)(916,308)(924,340)
	(930,369)(935,397)(939,423)
	(942,448)(943,472)(944,494)
	(944,516)(943,537)(943,558)
	(942,580)(941,602)(940,626)
	(939,651)(938,678)(938,709)
	(938,742)(938,778)(938,818)
	(939,863)(941,912)(942,965)
	(944,1022)(946,1083)(947,1148)
	(949,1215)(949,1283)(948,1355)
	(947,1425)(945,1492)(944,1555)
	(942,1613)(940,1666)(939,1714)
	(938,1758)(938,1796)(938,1831)
	(938,1862)(939,1891)(940,1916)
	(941,1940)(942,1962)(943,1983)
	(944,2004)(944,2024)(944,2045)
	(943,2067)(941,2090)(938,2115)
	(933,2141)(927,2168)(920,2198)
	(910,2228)(898,2260)(883,2293)
	(866,2325)(846,2356)(824,2384)
	(799,2408)(774,2425)(748,2438)
	(722,2446)(695,2451)(670,2454)
	(645,2454)(622,2452)(600,2450)
	(579,2446)(560,2442)(542,2437)
	(524,2432)(508,2427)(493,2421)
	(478,2415)(463,2409)(449,2402)
	(435,2394)(420,2386)(406,2376)
	(390,2365)(375,2353)(358,2338)
	(340,2322)(322,2302)(302,2280)
	(281,2254)(260,2225)(237,2192)
	(214,2154)(190,2112)(167,2065)
	(145,2013)(124,1958)
\path(237,5308)(226,5290)(216,5271)
	(207,5251)(198,5230)(191,5207)
	(184,5183)(178,5158)(173,5133)
	(169,5106)(166,5079)(164,5051)
	(164,5023)(164,4995)(166,4967)
	(168,4939)(172,4911)(178,4884)
	(184,4858)(192,4833)(201,4809)
	(211,4787)(222,4766)(234,4746)
	(247,4728)(262,4711)(277,4696)
	(294,4682)(312,4670)(332,4659)
	(354,4650)(378,4642)(403,4636)
	(429,4630)(457,4625)(486,4622)
	(516,4618)(547,4615)(579,4612)
	(611,4610)(644,4607)(676,4603)
	(708,4599)(739,4594)(770,4589)
	(799,4582)(827,4574)(854,4564)
	(878,4553)(901,4540)(922,4526)
	(942,4509)(959,4490)(974,4469)
	(987,4445)(996,4424)(1003,4401)
	(1010,4376)(1015,4349)(1020,4320)
	(1024,4289)(1026,4256)(1028,4222)
	(1030,4186)(1031,4148)(1031,4109)
	(1031,4069)(1031,4027)(1030,3985)
	(1030,3943)(1029,3899)(1029,3856)
	(1029,3813)(1029,3770)(1029,3727)
	(1030,3685)(1032,3644)(1034,3604)
	(1037,3565)(1041,3528)(1046,3492)
	(1052,3457)(1059,3424)(1067,3393)
	(1076,3363)(1087,3335)(1099,3308)
	(1110,3287)(1123,3267)(1137,3248)
	(1151,3230)(1167,3213)(1184,3197)
	(1202,3182)(1221,3169)(1241,3156)
	(1262,3145)(1284,3135)(1307,3127)
	(1331,3119)(1355,3114)(1380,3109)
	(1405,3106)(1431,3105)(1458,3105)
	(1484,3107)(1511,3110)(1538,3115)
	(1565,3121)(1592,3129)(1618,3139)
	(1644,3150)(1670,3162)(1695,3176)
	(1720,3191)(1744,3208)(1767,3226)
	(1789,3246)(1811,3267)(1832,3289)
	(1852,3313)(1871,3337)(1890,3364)
	(1907,3391)(1924,3420)(1939,3450)
	(1953,3480)(1967,3512)(1980,3545)
	(1992,3580)(2004,3616)(2015,3654)
	(2025,3693)(2035,3734)(2043,3775)
	(2051,3819)(2059,3863)(2065,3908)
	(2071,3954)(2076,4002)(2080,4049)
	(2083,4098)(2085,4147)(2086,4196)
	(2087,4246)(2086,4295)(2085,4344)
	(2083,4393)(2080,4442)(2076,4489)
	(2071,4537)(2065,4583)(2059,4628)
	(2051,4672)(2043,4716)(2035,4757)
	(2025,4798)(2015,4837)(2004,4875)
	(1992,4911)(1980,4946)(1967,4979)
	(1953,5011)(1939,5041)(1924,5070)
	(1904,5105)(1883,5137)(1861,5168)
	(1837,5197)(1812,5224)(1786,5250)
	(1759,5274)(1730,5296)(1701,5318)
	(1670,5337)(1639,5356)(1607,5372)
	(1574,5388)(1541,5402)(1508,5414)
	(1474,5426)(1440,5436)(1407,5445)
	(1374,5452)(1341,5459)(1309,5465)
	(1278,5469)(1247,5473)(1218,5476)
	(1189,5478)(1162,5480)(1136,5481)
	(1111,5482)(1087,5483)(1065,5483)
	(1044,5483)(1024,5483)(991,5483)
	(962,5483)(935,5483)(910,5483)
	(888,5483)(867,5483)(848,5483)
	(831,5483)(815,5483)(800,5483)
	(785,5483)(772,5483)(758,5483)
	(745,5483)(732,5483)(718,5483)
	(703,5483)(687,5483)(672,5483)
	(656,5483)(639,5483)(620,5482)
	(600,5481)(579,5479)(557,5477)
	(534,5474)(510,5470)(486,5465)
	(461,5459)(437,5452)(413,5443)
	(389,5434)(366,5423)(344,5410)
	(323,5397)(303,5382)(284,5366)
	(267,5348)(251,5329)(237,5308)
\path(874,3533)(862,3563)(847,3592)
	(832,3617)(816,3640)(801,3661)
	(787,3678)(774,3694)(762,3707)
	(751,3719)(741,3729)(730,3739)
	(720,3749)(708,3758)(695,3768)
	(681,3777)(664,3787)(644,3797)
	(621,3807)(595,3817)(565,3825)
	(533,3831)(499,3833)(465,3831)
	(433,3825)(403,3817)(377,3807)
	(354,3797)(334,3787)(317,3777)
	(303,3768)(290,3758)(278,3749)
	(268,3739)(257,3729)(247,3719)
	(236,3707)(224,3694)(211,3678)
	(197,3661)(182,3640)(166,3617)
	(151,3592)(136,3563)(124,3533)
	(116,3505)(112,3478)(109,3452)
	(109,3429)(110,3408)(112,3390)
	(114,3374)(117,3361)(120,3349)
	(123,3339)(126,3329)(130,3320)
	(134,3312)(139,3303)(145,3293)
	(152,3283)(159,3271)(169,3258)
	(180,3244)(194,3228)(210,3211)
	(229,3193)(250,3175)(274,3158)
	(302,3142)(331,3130)(358,3120)
	(384,3113)(406,3107)(427,3103)
	(444,3100)(460,3098)(474,3096)
	(487,3096)(499,3096)(511,3096)
	(524,3096)(538,3098)(554,3100)
	(571,3103)(592,3107)(614,3113)
	(640,3120)(667,3130)(696,3142)
	(724,3158)(748,3175)(769,3193)
	(788,3211)(804,3228)(818,3244)
	(829,3258)(839,3271)(846,3283)
	(853,3293)(859,3303)(864,3312)
	(868,3321)(872,3329)(875,3339)
	(878,3349)(881,3361)(884,3374)
	(886,3390)(888,3408)(889,3429)
	(889,3452)(886,3478)(882,3505)(874,3533)
\path(2899,1958)(2916,1905)(2931,1851)
	(2944,1797)(2956,1743)(2966,1691)
	(2974,1641)(2981,1593)(2988,1548)
	(2993,1506)(2997,1466)(3001,1429)
	(3004,1395)(3006,1363)(3008,1333)
	(3009,1304)(3010,1277)(3011,1252)
	(3011,1226)(3012,1202)(3011,1177)
	(3011,1152)(3010,1126)(3009,1099)
	(3008,1071)(3006,1041)(3004,1010)
	(3001,976)(2997,939)(2993,900)
	(2988,858)(2981,814)(2974,767)
	(2966,718)(2956,666)(2944,614)
	(2931,561)(2916,509)(2899,458)
	(2878,405)(2856,356)(2832,312)
	(2809,272)(2786,237)(2763,206)
	(2741,178)(2720,155)(2700,134)
	(2681,117)(2662,102)(2645,89)
	(2629,77)(2613,68)(2598,59)
	(2583,52)(2568,45)(2553,40)
	(2538,34)(2522,30)(2506,25)
	(2490,22)(2472,18)(2454,15)
	(2434,13)(2414,12)(2392,12)
	(2369,14)(2345,17)(2321,23)
	(2296,32)(2271,45)(2247,61)
	(2224,83)(2203,110)(2184,140)
	(2168,173)(2154,207)(2142,242)
	(2133,276)(2126,308)(2120,340)
	(2117,369)(2114,397)(2113,423)
	(2113,448)(2114,472)(2116,494)
	(2118,516)(2121,537)(2124,558)
	(2127,580)(2131,602)(2134,626)
	(2137,651)(2140,678)(2143,709)
	(2145,742)(2147,778)(2149,818)
	(2150,863)(2151,912)(2151,965)
	(2151,1022)(2150,1083)(2150,1148)
	(2149,1215)(2149,1283)(2149,1355)
	(2150,1425)(2150,1492)(2151,1555)
	(2151,1613)(2150,1666)(2149,1714)
	(2148,1758)(2146,1796)(2144,1831)
	(2141,1862)(2138,1891)(2135,1916)
	(2131,1940)(2127,1962)(2124,1983)
	(2121,2004)(2118,2024)(2115,2045)
	(2114,2067)(2113,2090)(2113,2115)
	(2115,2141)(2118,2168)(2123,2198)
	(2130,2228)(2139,2260)(2151,2293)
	(2165,2325)(2182,2356)(2202,2384)
	(2224,2408)(2247,2425)(2271,2438)
	(2296,2446)(2321,2451)(2345,2454)
	(2369,2454)(2392,2452)(2414,2450)
	(2434,2446)(2454,2442)(2472,2437)
	(2490,2432)(2506,2427)(2522,2421)
	(2538,2415)(2553,2409)(2568,2402)
	(2583,2394)(2598,2386)(2613,2376)
	(2629,2365)(2645,2353)(2662,2338)
	(2681,2322)(2700,2302)(2720,2280)
	(2741,2254)(2763,2225)(2786,2192)
	(2809,2154)(2832,2112)(2856,2065)
	(2878,2013)(2899,1958)
\path(5712,4108)(5701,4090)(5691,4071)
	(5682,4051)(5673,4030)(5666,4007)
	(5659,3983)(5653,3958)(5648,3933)
	(5644,3906)(5641,3879)(5639,3851)
	(5639,3823)(5639,3795)(5641,3767)
	(5643,3739)(5647,3711)(5653,3684)
	(5659,3658)(5667,3633)(5676,3609)
	(5686,3587)(5697,3566)(5709,3546)
	(5722,3528)(5737,3511)(5752,3496)
	(5769,3482)(5787,3470)(5807,3459)
	(5829,3450)(5853,3442)(5878,3436)
	(5904,3430)(5932,3425)(5961,3422)
	(5991,3418)(6022,3415)(6054,3412)
	(6086,3410)(6119,3407)(6151,3403)
	(6183,3399)(6214,3394)(6245,3389)
	(6274,3382)(6302,3374)(6329,3364)
	(6353,3353)(6376,3340)(6397,3326)
	(6417,3309)(6434,3290)(6449,3269)
	(6462,3245)(6471,3224)(6478,3201)
	(6485,3176)(6490,3149)(6495,3120)
	(6499,3089)(6501,3056)(6503,3022)
	(6505,2986)(6506,2948)(6506,2909)
	(6506,2869)(6506,2827)(6505,2785)
	(6505,2743)(6504,2699)(6504,2656)
	(6504,2613)(6504,2570)(6504,2527)
	(6505,2485)(6507,2444)(6509,2404)
	(6512,2365)(6516,2328)(6521,2292)
	(6527,2257)(6534,2224)(6542,2193)
	(6551,2163)(6562,2135)(6574,2108)
	(6585,2087)(6598,2067)(6612,2048)
	(6626,2030)(6642,2013)(6659,1997)
	(6677,1982)(6696,1969)(6716,1956)
	(6737,1945)(6759,1935)(6782,1927)
	(6806,1919)(6830,1914)(6855,1909)
	(6880,1906)(6906,1905)(6933,1905)
	(6959,1907)(6986,1910)(7013,1915)
	(7040,1921)(7067,1929)(7093,1939)
	(7119,1950)(7145,1962)(7170,1976)
	(7195,1991)(7219,2008)(7242,2026)
	(7264,2046)(7286,2067)(7307,2089)
	(7327,2113)(7346,2137)(7365,2164)
	(7382,2191)(7399,2220)(7414,2250)
	(7428,2280)(7442,2312)(7455,2345)
	(7467,2380)(7479,2416)(7490,2454)
	(7500,2493)(7510,2534)(7518,2575)
	(7526,2619)(7534,2663)(7540,2708)
	(7546,2754)(7551,2802)(7555,2849)
	(7558,2898)(7560,2947)(7561,2996)
	(7562,3046)(7561,3095)(7560,3144)
	(7558,3193)(7555,3242)(7551,3289)
	(7546,3337)(7540,3383)(7534,3428)
	(7526,3472)(7518,3516)(7510,3557)
	(7500,3598)(7490,3637)(7479,3675)
	(7467,3711)(7455,3746)(7442,3779)
	(7428,3811)(7414,3841)(7399,3870)
	(7379,3905)(7358,3937)(7336,3968)
	(7312,3997)(7287,4024)(7261,4050)
	(7234,4074)(7205,4096)(7176,4118)
	(7145,4137)(7114,4156)(7082,4172)
	(7049,4188)(7016,4202)(6983,4214)
	(6949,4226)(6915,4236)(6882,4245)
	(6849,4252)(6816,4259)(6784,4265)
	(6753,4269)(6722,4273)(6693,4276)
	(6664,4278)(6637,4280)(6611,4281)
	(6586,4282)(6562,4283)(6540,4283)
	(6519,4283)(6499,4283)(6466,4283)
	(6437,4283)(6410,4283)(6385,4283)
	(6363,4283)(6342,4283)(6323,4283)
	(6306,4283)(6290,4283)(6275,4283)
	(6260,4283)(6247,4283)(6233,4283)
	(6220,4283)(6207,4283)(6193,4283)
	(6178,4283)(6162,4283)(6147,4283)
	(6131,4283)(6114,4283)(6095,4282)
	(6075,4281)(6054,4279)(6032,4277)
	(6009,4274)(5985,4270)(5961,4265)
	(5936,4259)(5912,4252)(5888,4243)
	(5864,4234)(5841,4223)(5819,4210)
	(5798,4197)(5778,4182)(5759,4166)
	(5742,4148)(5726,4129)(5712,4108)
\path(6349,2333)(6337,2363)(6322,2392)
	(6307,2417)(6291,2440)(6276,2461)
	(6262,2478)(6249,2494)(6237,2507)
	(6226,2519)(6216,2529)(6205,2539)
	(6195,2549)(6183,2558)(6170,2568)
	(6156,2577)(6139,2587)(6119,2597)
	(6096,2607)(6070,2617)(6040,2625)
	(6008,2631)(5974,2633)(5940,2631)
	(5908,2625)(5878,2617)(5852,2607)
	(5829,2597)(5809,2587)(5792,2577)
	(5778,2568)(5765,2558)(5753,2549)
	(5743,2539)(5732,2529)(5722,2519)
	(5711,2507)(5699,2494)(5686,2478)
	(5672,2461)(5657,2440)(5641,2417)
	(5626,2392)(5611,2363)(5599,2333)
	(5591,2305)(5587,2278)(5584,2252)
	(5584,2229)(5585,2208)(5587,2190)
	(5589,2174)(5592,2161)(5595,2149)
	(5598,2139)(5601,2129)(5605,2120)
	(5609,2112)(5614,2103)(5620,2093)
	(5627,2083)(5634,2071)(5644,2058)
	(5655,2044)(5669,2028)(5685,2011)
	(5704,1993)(5725,1975)(5749,1958)
	(5777,1942)(5806,1930)(5833,1920)
	(5859,1913)(5881,1907)(5902,1903)
	(5919,1900)(5935,1898)(5949,1896)
	(5962,1896)(5974,1896)(5986,1896)
	(5999,1896)(6013,1898)(6029,1900)
	(6046,1903)(6067,1907)(6089,1913)
	(6115,1920)(6142,1930)(6171,1942)
	(6199,1958)(6223,1975)(6244,1993)
	(6263,2011)(6279,2028)(6293,2044)
	(6304,2058)(6314,2071)(6321,2083)
	(6328,2093)(6334,2103)(6339,2112)
	(6343,2121)(6347,2129)(6350,2139)
	(6353,2149)(6356,2161)(6359,2174)
	(6361,2190)(6363,2208)(6364,2229)
	(6364,2252)(6361,2278)(6357,2305)(6349,2333)
\path(8449,3758)(8466,3705)(8481,3651)
	(8494,3597)(8506,3543)(8516,3491)
	(8524,3441)(8531,3393)(8538,3348)
	(8543,3306)(8547,3266)(8551,3229)
	(8554,3195)(8556,3163)(8558,3133)
	(8559,3104)(8560,3077)(8561,3052)
	(8561,3026)(8562,3002)(8561,2977)
	(8561,2952)(8560,2926)(8559,2899)
	(8558,2871)(8556,2841)(8554,2810)
	(8551,2776)(8547,2739)(8543,2700)
	(8538,2658)(8531,2614)(8524,2567)
	(8516,2518)(8506,2466)(8494,2414)
	(8481,2361)(8466,2309)(8449,2258)
	(8428,2205)(8406,2156)(8382,2112)
	(8359,2072)(8336,2037)(8313,2006)
	(8291,1978)(8270,1955)(8250,1934)
	(8231,1917)(8212,1902)(8195,1889)
	(8179,1877)(8163,1868)(8148,1859)
	(8133,1852)(8118,1845)(8103,1840)
	(8088,1834)(8072,1830)(8056,1825)
	(8040,1822)(8022,1818)(8004,1815)
	(7984,1813)(7964,1812)(7942,1812)
	(7919,1814)(7895,1817)(7871,1823)
	(7846,1832)(7821,1845)(7797,1861)
	(7774,1883)(7753,1910)(7734,1940)
	(7718,1973)(7704,2007)(7692,2042)
	(7683,2076)(7676,2108)(7670,2140)
	(7667,2169)(7664,2197)(7663,2223)
	(7663,2248)(7664,2272)(7666,2294)
	(7668,2316)(7671,2337)(7674,2358)
	(7677,2380)(7681,2402)(7684,2426)
	(7687,2451)(7690,2478)(7693,2509)
	(7695,2542)(7697,2578)(7699,2618)
	(7700,2663)(7701,2712)(7701,2765)
	(7701,2822)(7700,2883)(7700,2948)
	(7699,3015)(7699,3083)(7699,3155)
	(7700,3225)(7700,3292)(7701,3355)
	(7701,3413)(7700,3466)(7699,3514)
	(7698,3558)(7696,3596)(7694,3631)
	(7691,3662)(7688,3691)(7685,3716)
	(7681,3740)(7677,3762)(7674,3783)
	(7671,3804)(7668,3824)(7665,3845)
	(7664,3867)(7663,3890)(7663,3915)
	(7665,3941)(7668,3968)(7673,3998)
	(7680,4028)(7689,4060)(7701,4093)
	(7715,4125)(7732,4156)(7752,4184)
	(7774,4208)(7797,4225)(7821,4238)
	(7846,4246)(7871,4251)(7895,4254)
	(7919,4254)(7942,4252)(7964,4250)
	(7984,4246)(8004,4242)(8022,4237)
	(8040,4232)(8056,4227)(8072,4221)
	(8088,4215)(8103,4209)(8118,4202)
	(8133,4194)(8148,4186)(8163,4176)
	(8179,4165)(8195,4153)(8212,4138)
	(8231,4122)(8250,4102)(8270,4080)
	(8291,4054)(8313,4025)(8336,3992)
	(8359,3954)(8382,3912)(8406,3865)
	(8428,3813)(8449,3758)
\put(424,5183){\makebox(0,0)[lb]{\smash{{{\SetFigFont{6}{7.2}{\rmdefault}{\mddefault}{\updefault}$1$}}}}}
\put(1324,5183){\makebox(0,0)[lb]{\smash{{{\SetFigFont{6}{7.2}{\rmdefault}{\mddefault}{\updefault}$2$}}}}}
\put(574,2183){\makebox(0,0)[lb]{\smash{{{\SetFigFont{6}{7.2}{\rmdefault}{\mddefault}{\updefault}$1$}}}}}
\put(424,3683){\makebox(0,0)[lb]{\smash{{{\SetFigFont{6}{7.2}{\rmdefault}{\mddefault}{\updefault}$1'$}}}}}
\put(349,683){\makebox(0,0)[lb]{\smash{{{\SetFigFont{6}{7.2}{\rmdefault}{\mddefault}{\updefault}$1'$}}}}}
\put(1474,3683){\makebox(0,0)[lb]{\smash{{{\SetFigFont{6}{7.2}{\rmdefault}{\mddefault}{\updefault}$2'$}}}}}
\put(2449,2183){\makebox(0,0)[lb]{\smash{{{\SetFigFont{6}{7.2}{\rmdefault}{\mddefault}{\updefault}$3$}}}}}
\put(2449,683){\makebox(0,0)[lb]{\smash{{{\SetFigFont{6}{7.2}{\rmdefault}{\mddefault}{\updefault}$3'$}}}}}
\put(4399,2858){\makebox(0,0)[lb]{\smash{{{\SetFigFont{6}{7.2}{\rmdefault}{\mddefault}{\updefault}$=\; Q_1$}}}}}
\end{picture}
}

%% file: xfig/multunit.eepic
\setlength{\unitlength}{0.00041667in}
\begingroup\makeatletter\ifx\SetFigFont\undefined%
\gdef\SetFigFont#1#2#3#4#5{%
  \reset@font\fontsize{#1}{#2pt}%
  \fontfamily{#3}\fontseries{#4}\fontshape{#5}%
  \selectfont}%
\fi\endgroup%
{\renewcommand{\dashlinestretch}{30}
\begin{picture}(4800,2466)(0,-10)
\put(1500,1947){\blacken\ellipse{300}{300}}
\put(1500,1947){\ellipse{300}{300}}
\put(2400,1947){\blacken\ellipse{300}{300}}
\put(2400,1947){\ellipse{300}{300}}
\put(2400,447){\blacken\ellipse{300}{300}}
\put(2400,447){\ellipse{300}{300}}
\put(1500,447){\blacken\ellipse{300}{300}}
\put(1500,447){\ellipse{300}{300}}
\put(3300,1947){\blacken\ellipse{300}{300}}
\put(3300,1947){\ellipse{300}{300}}
\put(4200,1947){\blacken\ellipse{300}{300}}
\put(4200,1947){\ellipse{300}{300}}
\put(4200,447){\blacken\ellipse{300}{300}}
\put(4200,447){\ellipse{300}{300}}
\put(3300,447){\blacken\ellipse{300}{300}}
\put(3300,447){\ellipse{300}{300}}
\thicklines
\path(2400,1947)(2400,447)
\path(1500,1947)(1500,447)
\path(4200,1947)(4200,447)
\path(3300,1947)(3300,447)
\thinlines
\path(1125,1947)(1117,1921)(1109,1893)
	(1102,1863)(1096,1832)(1091,1799)
	(1086,1764)(1082,1728)(1078,1690)
	(1075,1650)(1072,1609)(1070,1567)
	(1068,1523)(1066,1478)(1065,1432)
	(1064,1385)(1063,1338)(1062,1290)
	(1062,1242)(1062,1193)(1062,1145)
	(1062,1097)(1063,1049)(1064,1002)
	(1065,955)(1066,910)(1068,865)
	(1070,821)(1072,779)(1075,738)
	(1078,699)(1082,662)(1086,626)
	(1091,591)(1096,559)(1102,528)
	(1109,500)(1117,472)(1125,447)
	(1138,415)(1153,385)(1169,358)
	(1188,332)(1208,308)(1231,285)
	(1255,263)(1280,242)(1307,222)
	(1335,203)(1364,184)(1394,166)
	(1425,149)(1455,133)(1486,117)
	(1517,103)(1547,90)(1576,77)
	(1605,67)(1633,58)(1659,50)
	(1684,45)(1707,42)(1729,42)
	(1749,45)(1768,50)(1785,59)
	(1800,72)(1812,85)(1823,102)
	(1833,121)(1842,142)(1850,166)
	(1858,193)(1865,222)(1872,253)
	(1878,286)(1884,321)(1889,357)
	(1894,395)(1898,435)(1903,475)
	(1907,517)(1910,559)(1914,602)
	(1917,645)(1920,688)(1923,731)
	(1926,774)(1929,816)(1931,858)
	(1934,900)(1936,941)(1939,981)
	(1941,1020)(1943,1058)(1945,1095)
	(1947,1132)(1948,1168)(1949,1203)
	(1950,1238)(1950,1272)(1950,1308)
	(1949,1345)(1948,1381)(1946,1419)
	(1944,1457)(1942,1495)(1940,1534)
	(1938,1574)(1935,1615)(1932,1656)
	(1930,1698)(1927,1741)(1924,1783)
	(1920,1826)(1917,1869)(1914,1912)
	(1910,1954)(1906,1996)(1902,2037)
	(1897,2076)(1892,2115)(1887,2152)
	(1881,2187)(1875,2220)(1868,2252)
	(1861,2280)(1853,2307)(1844,2331)
	(1834,2352)(1824,2370)(1812,2385)
	(1800,2397)(1786,2407)(1770,2413)
	(1753,2416)(1735,2416)(1715,2414)
	(1693,2410)(1671,2403)(1647,2395)
	(1622,2384)(1596,2372)(1569,2359)
	(1541,2345)(1513,2329)(1484,2312)
	(1455,2295)(1427,2277)(1398,2258)
	(1370,2239)(1343,2219)(1316,2199)
	(1291,2178)(1266,2156)(1243,2134)
	(1222,2111)(1201,2087)(1183,2062)
	(1166,2035)(1151,2008)(1137,1978)(1125,1947)
\path(2775,2172)(2783,2151)(2790,2129)
	(2796,2105)(2801,2079)(2805,2052)
	(2807,2023)(2809,1993)(2810,1962)
	(2810,1930)(2810,1896)(2808,1862)
	(2807,1827)(2805,1792)(2802,1756)
	(2800,1720)(2797,1683)(2794,1647)
	(2791,1610)(2788,1574)(2786,1537)
	(2783,1502)(2781,1466)(2780,1431)
	(2778,1396)(2777,1362)(2776,1328)
	(2775,1295)(2775,1262)(2775,1230)
	(2775,1197)(2775,1164)(2775,1132)
	(2775,1099)(2776,1065)(2777,1032)
	(2778,997)(2780,962)(2781,927)
	(2783,891)(2786,855)(2788,818)
	(2791,782)(2794,745)(2797,708)
	(2800,671)(2802,634)(2805,598)
	(2807,562)(2808,527)(2810,493)
	(2810,459)(2810,427)(2809,396)
	(2807,366)(2805,338)(2801,311)
	(2796,286)(2790,263)(2783,242)
	(2775,222)(2764,202)(2752,183)
	(2737,166)(2721,151)(2704,136)
	(2685,123)(2664,110)(2642,98)
	(2619,86)(2595,75)(2570,64)
	(2545,54)(2520,45)(2494,37)
	(2469,29)(2443,23)(2419,18)
	(2395,14)(2373,12)(2351,12)
	(2331,15)(2312,20)(2295,28)
	(2278,39)(2264,54)(2250,72)
	(2240,89)(2231,108)(2223,129)
	(2215,152)(2208,178)(2201,206)
	(2194,236)(2188,267)(2182,301)
	(2176,336)(2171,372)(2166,410)
	(2161,449)(2156,489)(2152,529)
	(2148,571)(2143,612)(2139,654)
	(2136,696)(2132,738)(2128,780)
	(2125,822)(2121,863)(2118,903)
	(2115,943)(2113,983)(2110,1021)
	(2108,1059)(2105,1096)(2104,1132)
	(2102,1168)(2101,1203)(2100,1238)
	(2100,1272)(2100,1308)(2101,1345)
	(2102,1381)(2104,1418)(2106,1456)
	(2108,1494)(2111,1533)(2114,1572)
	(2117,1612)(2120,1652)(2124,1693)
	(2127,1735)(2131,1776)(2135,1818)
	(2139,1860)(2143,1901)(2148,1942)
	(2152,1983)(2157,2023)(2162,2062)
	(2167,2100)(2173,2137)(2179,2172)
	(2185,2206)(2191,2237)(2198,2267)
	(2205,2295)(2213,2320)(2221,2343)
	(2230,2364)(2240,2382)(2250,2397)
	(2265,2413)(2281,2426)(2299,2434)
	(2318,2438)(2339,2439)(2362,2437)
	(2386,2433)(2411,2426)(2437,2418)
	(2464,2407)(2492,2396)(2520,2383)
	(2547,2369)(2574,2355)(2601,2340)
	(2627,2324)(2651,2308)(2675,2292)
	(2696,2274)(2716,2256)(2734,2238)
	(2749,2218)(2763,2196)(2775,2172)
\path(4725,1947)(4733,1921)(4741,1893)
	(4748,1863)(4754,1832)(4759,1799)
	(4764,1764)(4768,1728)(4772,1690)
	(4775,1650)(4778,1609)(4780,1567)
	(4782,1523)(4784,1478)(4785,1432)
	(4786,1385)(4787,1338)(4788,1290)
	(4788,1242)(4788,1193)(4788,1145)
	(4788,1097)(4787,1049)(4786,1002)
	(4785,955)(4784,910)(4782,865)
	(4780,821)(4778,779)(4775,738)
	(4772,699)(4768,662)(4764,626)
	(4759,591)(4754,559)(4748,528)
	(4741,500)(4733,472)(4725,447)
	(4712,415)(4697,385)(4681,358)
	(4662,332)(4642,308)(4619,285)
	(4595,263)(4570,242)(4543,222)
	(4515,203)(4486,184)(4456,166)
	(4425,149)(4395,133)(4364,117)
	(4333,103)(4303,90)(4274,77)
	(4245,67)(4217,58)(4191,50)
	(4166,45)(4143,42)(4121,42)
	(4101,45)(4082,50)(4065,59)
	(4050,72)(4038,85)(4027,102)
	(4017,121)(4008,142)(4000,166)
	(3992,193)(3985,222)(3978,253)
	(3972,286)(3966,321)(3961,357)
	(3956,395)(3952,435)(3947,475)
	(3943,517)(3940,559)(3936,602)
	(3933,645)(3930,688)(3927,731)
	(3924,774)(3921,816)(3919,858)
	(3916,900)(3914,941)(3911,981)
	(3909,1020)(3907,1058)(3905,1095)
	(3903,1132)(3902,1168)(3901,1203)
	(3900,1238)(3900,1272)(3900,1308)
	(3901,1345)(3902,1381)(3904,1419)
	(3906,1457)(3908,1495)(3910,1534)
	(3912,1574)(3915,1615)(3918,1656)
	(3920,1698)(3923,1741)(3926,1783)
	(3930,1826)(3933,1869)(3936,1912)
	(3940,1954)(3944,1996)(3948,2037)
	(3953,2076)(3958,2115)(3963,2152)
	(3969,2187)(3975,2220)(3982,2252)
	(3989,2280)(3997,2307)(4006,2331)
	(4016,2352)(4026,2370)(4038,2385)
	(4050,2397)(4064,2407)(4080,2413)
	(4097,2416)(4115,2416)(4135,2414)
	(4157,2410)(4179,2403)(4203,2395)
	(4228,2384)(4254,2372)(4281,2359)
	(4309,2345)(4337,2329)(4366,2312)
	(4395,2295)(4423,2277)(4452,2258)
	(4480,2239)(4507,2219)(4534,2199)
	(4559,2178)(4584,2156)(4607,2134)
	(4628,2111)(4649,2087)(4667,2062)
	(4684,2035)(4699,2008)(4713,1978)(4725,1947)
\path(3000,2172)(2989,2152)(2980,2131)
	(2971,2108)(2964,2083)(2957,2056)
	(2952,2028)(2947,1998)(2943,1967)
	(2940,1935)(2937,1901)(2935,1867)
	(2934,1832)(2932,1796)(2932,1760)
	(2931,1723)(2930,1686)(2930,1649)
	(2930,1612)(2930,1576)(2929,1539)
	(2929,1503)(2929,1467)(2928,1432)
	(2928,1397)(2927,1362)(2927,1329)
	(2926,1295)(2926,1262)(2925,1230)
	(2925,1197)(2925,1164)(2926,1132)
	(2926,1099)(2927,1065)(2927,1032)
	(2928,997)(2928,962)(2929,927)
	(2929,891)(2929,855)(2930,818)
	(2930,782)(2930,745)(2930,708)
	(2931,671)(2932,634)(2932,598)
	(2934,562)(2935,527)(2937,493)
	(2940,459)(2943,427)(2947,396)
	(2952,366)(2957,338)(2964,311)
	(2971,286)(2980,263)(2989,242)
	(3000,222)(3014,202)(3029,183)
	(3047,166)(3066,151)(3086,137)
	(3109,123)(3133,111)(3159,99)
	(3185,88)(3213,77)(3242,67)
	(3271,57)(3300,49)(3329,40)
	(3358,33)(3387,27)(3415,22)
	(3441,19)(3467,17)(3491,17)
	(3514,20)(3534,24)(3553,31)
	(3571,42)(3586,55)(3600,72)
	(3610,88)(3619,107)(3627,129)
	(3634,153)(3641,179)(3646,207)
	(3651,237)(3655,269)(3658,303)
	(3661,339)(3663,375)(3665,414)
	(3667,453)(3668,493)(3669,534)
	(3669,575)(3670,616)(3670,658)
	(3670,700)(3670,741)(3671,782)
	(3671,823)(3671,863)(3672,903)
	(3672,942)(3672,980)(3673,1017)
	(3674,1054)(3674,1091)(3675,1126)
	(3675,1162)(3675,1197)(3675,1232)
	(3675,1268)(3674,1303)(3674,1340)
	(3673,1377)(3672,1414)(3672,1452)
	(3672,1491)(3671,1531)(3671,1571)
	(3671,1612)(3670,1653)(3670,1694)
	(3670,1736)(3670,1778)(3669,1819)
	(3669,1860)(3668,1901)(3667,1941)
	(3665,1980)(3663,2019)(3661,2055)
	(3658,2091)(3655,2125)(3651,2157)
	(3646,2187)(3641,2215)(3634,2241)
	(3627,2265)(3619,2287)(3610,2306)
	(3600,2322)(3586,2339)(3571,2352)
	(3553,2363)(3534,2370)(3514,2374)
	(3491,2377)(3467,2377)(3441,2375)
	(3415,2372)(3387,2367)(3358,2361)
	(3329,2354)(3300,2345)(3271,2337)
	(3242,2327)(3213,2317)(3185,2306)
	(3159,2295)(3133,2283)(3109,2271)
	(3086,2257)(3066,2243)(3047,2228)
	(3029,2211)(3014,2192)(3000,2172)
\put(0,1272){\makebox(0,0)[lb]{\smash{{{\SetFigFont{6}{7.2}{\rmdefault}{\mddefault}{\updefault}$(1,1)=$}}}}}
\end{picture}
}

%% file: xfig/Ai1.eepic
\setlength{\unitlength}{0.00041667in}
\begingroup\makeatletter\ifx\SetFigFont\undefined%
\gdef\SetFigFont#1#2#3#4#5{%
  \reset@font\fontsize{#1}{#2pt}%
  \fontfamily{#3}\fontseries{#4}\fontshape{#5}%
  \selectfont}%
\fi\endgroup%
{\renewcommand{\dashlinestretch}{30}
\begin{picture}(4800,2466)(0,-10)
\put(1500,1947){\blacken\ellipse{300}{300}}
\put(1500,1947){\ellipse{300}{300}}
\put(2400,1947){\blacken\ellipse{300}{300}}
\put(2400,1947){\ellipse{300}{300}}
\put(2400,447){\blacken\ellipse{300}{300}}
\put(2400,447){\ellipse{300}{300}}
\put(1500,447){\blacken\ellipse{300}{300}}
\put(1500,447){\ellipse{300}{300}}
\put(3300,1947){\blacken\ellipse{300}{300}}
\put(3300,1947){\ellipse{300}{300}}
\put(4200,1947){\blacken\ellipse{300}{300}}
\put(4200,1947){\ellipse{300}{300}}
\put(4200,447){\blacken\ellipse{300}{300}}
\put(4200,447){\ellipse{300}{300}}
\put(3300,447){\blacken\ellipse{300}{300}}
\put(3300,447){\ellipse{300}{300}}
\thicklines
\path(2400,1947)(2400,447)
\path(1500,1947)(1500,447)
\path(4200,1947)(4200,447)
\thinlines
\path(1125,1947)(1117,1921)(1109,1893)
	(1102,1863)(1096,1832)(1091,1799)
	(1086,1764)(1082,1728)(1078,1690)
	(1075,1650)(1072,1609)(1070,1567)
	(1068,1523)(1066,1478)(1065,1432)
	(1064,1385)(1063,1338)(1062,1290)
	(1062,1242)(1062,1193)(1062,1145)
	(1062,1097)(1063,1049)(1064,1002)
	(1065,955)(1066,910)(1068,865)
	(1070,821)(1072,779)(1075,738)
	(1078,699)(1082,662)(1086,626)
	(1091,591)(1096,559)(1102,528)
	(1109,500)(1117,472)(1125,447)
	(1138,415)(1153,385)(1169,358)
	(1188,332)(1208,308)(1231,285)
	(1255,263)(1280,242)(1307,222)
	(1335,203)(1364,184)(1394,166)
	(1425,149)(1455,133)(1486,117)
	(1517,103)(1547,90)(1576,77)
	(1605,67)(1633,58)(1659,50)
	(1684,45)(1707,42)(1729,42)
	(1749,45)(1768,50)(1785,59)
	(1800,72)(1812,85)(1823,102)
	(1833,121)(1842,142)(1850,166)
	(1858,193)(1865,222)(1872,253)
	(1878,286)(1884,321)(1889,357)
	(1894,395)(1898,435)(1903,475)
	(1907,517)(1910,559)(1914,602)
	(1917,645)(1920,688)(1923,731)
	(1926,774)(1929,816)(1931,858)
	(1934,900)(1936,941)(1939,981)
	(1941,1020)(1943,1058)(1945,1095)
	(1947,1132)(1948,1168)(1949,1203)
	(1950,1238)(1950,1272)(1950,1308)
	(1949,1345)(1948,1381)(1946,1419)
	(1944,1457)(1942,1495)(1940,1534)
	(1938,1574)(1935,1615)(1932,1656)
	(1930,1698)(1927,1741)(1924,1783)
	(1920,1826)(1917,1869)(1914,1912)
	(1910,1954)(1906,1996)(1902,2037)
	(1897,2076)(1892,2115)(1887,2152)
	(1881,2187)(1875,2220)(1868,2252)
	(1861,2280)(1853,2307)(1844,2331)
	(1834,2352)(1824,2370)(1812,2385)
	(1800,2397)(1786,2407)(1770,2413)
	(1753,2416)(1735,2416)(1715,2414)
	(1693,2410)(1671,2403)(1647,2395)
	(1622,2384)(1596,2372)(1569,2359)
	(1541,2345)(1513,2329)(1484,2312)
	(1455,2295)(1427,2277)(1398,2258)
	(1370,2239)(1343,2219)(1316,2199)
	(1291,2178)(1266,2156)(1243,2134)
	(1222,2111)(1201,2087)(1183,2062)
	(1166,2035)(1151,2008)(1137,1978)(1125,1947)
\path(2775,2172)(2783,2151)(2790,2129)
	(2796,2105)(2801,2079)(2805,2052)
	(2807,2023)(2809,1993)(2810,1962)
	(2810,1930)(2810,1896)(2808,1862)
	(2807,1827)(2805,1792)(2802,1756)
	(2800,1720)(2797,1683)(2794,1647)
	(2791,1610)(2788,1574)(2786,1537)
	(2783,1502)(2781,1466)(2780,1431)
	(2778,1396)(2777,1362)(2776,1328)
	(2775,1295)(2775,1262)(2775,1230)
	(2775,1197)(2775,1164)(2775,1132)
	(2775,1099)(2776,1065)(2777,1032)
	(2778,997)(2780,962)(2781,927)
	(2783,891)(2786,855)(2788,818)
	(2791,782)(2794,745)(2797,708)
	(2800,671)(2802,634)(2805,598)
	(2807,562)(2808,527)(2810,493)
	(2810,459)(2810,427)(2809,396)
	(2807,366)(2805,338)(2801,311)
	(2796,286)(2790,263)(2783,242)
	(2775,222)(2764,202)(2752,183)
	(2737,166)(2721,151)(2704,136)
	(2685,123)(2664,110)(2642,98)
	(2619,86)(2595,75)(2570,64)
	(2545,54)(2520,45)(2494,37)
	(2469,29)(2443,23)(2419,18)
	(2395,14)(2373,12)(2351,12)
	(2331,15)(2312,20)(2295,28)
	(2278,39)(2264,54)(2250,72)
	(2240,89)(2231,108)(2223,129)
	(2215,152)(2208,178)(2201,206)
	(2194,236)(2188,267)(2182,301)
	(2176,336)(2171,372)(2166,410)
	(2161,449)(2156,489)(2152,529)
	(2148,571)(2143,612)(2139,654)
	(2136,696)(2132,738)(2128,780)
	(2125,822)(2121,863)(2118,903)
	(2115,943)(2113,983)(2110,1021)
	(2108,1059)(2105,1096)(2104,1132)
	(2102,1168)(2101,1203)(2100,1238)
	(2100,1272)(2100,1308)(2101,1345)
	(2102,1381)(2104,1418)(2106,1456)
	(2108,1494)(2111,1533)(2114,1572)
	(2117,1612)(2120,1652)(2124,1693)
	(2127,1735)(2131,1776)(2135,1818)
	(2139,1860)(2143,1901)(2148,1942)
	(2152,1983)(2157,2023)(2162,2062)
	(2167,2100)(2173,2137)(2179,2172)
	(2185,2206)(2191,2237)(2198,2267)
	(2205,2295)(2213,2320)(2221,2343)
	(2230,2364)(2240,2382)(2250,2397)
	(2265,2413)(2281,2426)(2299,2434)
	(2318,2438)(2339,2439)(2362,2437)
	(2386,2433)(2411,2426)(2437,2418)
	(2464,2407)(2492,2396)(2520,2383)
	(2547,2369)(2574,2355)(2601,2340)
	(2627,2324)(2651,2308)(2675,2292)
	(2696,2274)(2716,2256)(2734,2238)
	(2749,2218)(2763,2196)(2775,2172)
\path(4725,1947)(4733,1921)(4741,1893)
	(4748,1863)(4754,1832)(4759,1799)
	(4764,1764)(4768,1728)(4772,1690)
	(4775,1650)(4778,1609)(4780,1567)
	(4782,1523)(4784,1478)(4785,1432)
	(4786,1385)(4787,1338)(4788,1290)
	(4788,1242)(4788,1193)(4788,1145)
	(4788,1097)(4787,1049)(4786,1002)
	(4785,955)(4784,910)(4782,865)
	(4780,821)(4778,779)(4775,738)
	(4772,699)(4768,662)(4764,626)
	(4759,591)(4754,559)(4748,528)
	(4741,500)(4733,472)(4725,447)
	(4712,415)(4697,385)(4681,358)
	(4662,332)(4642,308)(4619,285)
	(4595,263)(4570,242)(4543,222)
	(4515,203)(4486,184)(4456,166)
	(4425,149)(4395,133)(4364,117)
	(4333,103)(4303,90)(4274,77)
	(4245,67)(4217,58)(4191,50)
	(4166,45)(4143,42)(4121,42)
	(4101,45)(4082,50)(4065,59)
	(4050,72)(4038,85)(4027,102)
	(4017,121)(4008,142)(4000,166)
	(3992,193)(3985,222)(3978,253)
	(3972,286)(3966,321)(3961,357)
	(3956,395)(3952,435)(3947,475)
	(3943,517)(3940,559)(3936,602)
	(3933,645)(3930,688)(3927,731)
	(3924,774)(3921,816)(3919,858)
	(3916,900)(3914,941)(3911,981)
	(3909,1020)(3907,1058)(3905,1095)
	(3903,1132)(3902,1168)(3901,1203)
	(3900,1238)(3900,1272)(3900,1308)
	(3901,1345)(3902,1381)(3904,1419)
	(3906,1457)(3908,1495)(3910,1534)
	(3912,1574)(3915,1615)(3918,1656)
	(3920,1698)(3923,1741)(3926,1783)
	(3930,1826)(3933,1869)(3936,1912)
	(3940,1954)(3944,1996)(3948,2037)
	(3953,2076)(3958,2115)(3963,2152)
	(3969,2187)(3975,2220)(3982,2252)
	(3989,2280)(3997,2307)(4006,2331)
	(4016,2352)(4026,2370)(4038,2385)
	(4050,2397)(4064,2407)(4080,2413)
	(4097,2416)(4115,2416)(4135,2414)
	(4157,2410)(4179,2403)(4203,2395)
	(4228,2384)(4254,2372)(4281,2359)
	(4309,2345)(4337,2329)(4366,2312)
	(4395,2295)(4423,2277)(4452,2258)
	(4480,2239)(4507,2219)(4534,2199)
	(4559,2178)(4584,2156)(4607,2134)
	(4628,2111)(4649,2087)(4667,2062)
	(4684,2035)(4699,2008)(4713,1978)(4725,1947)
\path(3000,2172)(2989,2152)(2980,2131)
	(2971,2108)(2964,2083)(2957,2056)
	(2952,2028)(2947,1998)(2943,1967)
	(2940,1935)(2937,1901)(2935,1867)
	(2934,1832)(2932,1796)(2932,1760)
	(2931,1723)(2930,1686)(2930,1649)
	(2930,1612)(2930,1576)(2929,1539)
	(2929,1503)(2929,1467)(2928,1432)
	(2928,1397)(2927,1362)(2927,1329)
	(2926,1295)(2926,1262)(2925,1230)
	(2925,1197)(2925,1164)(2926,1132)
	(2926,1099)(2927,1065)(2927,1032)
	(2928,997)(2928,962)(2929,927)
	(2929,891)(2929,855)(2930,818)
	(2930,782)(2930,745)(2930,708)
	(2931,671)(2932,634)(2932,598)
	(2934,562)(2935,527)(2937,493)
	(2940,459)(2943,427)(2947,396)
	(2952,366)(2957,338)(2964,311)
	(2971,286)(2980,263)(2989,242)
	(3000,222)(3014,202)(3029,183)
	(3047,166)(3066,151)(3086,137)
	(3109,123)(3133,111)(3159,99)
	(3185,88)(3213,77)(3242,67)
	(3271,57)(3300,49)(3329,40)
	(3358,33)(3387,27)(3415,22)
	(3441,19)(3467,17)(3491,17)
	(3514,20)(3534,24)(3553,31)
	(3571,42)(3586,55)(3600,72)
	(3610,88)(3619,107)(3627,129)
	(3634,153)(3641,179)(3646,207)
	(3651,237)(3655,269)(3658,303)
	(3661,339)(3663,375)(3665,414)
	(3667,453)(3668,493)(3669,534)
	(3669,575)(3670,616)(3670,658)
	(3670,700)(3670,741)(3671,782)
	(3671,823)(3671,863)(3672,903)
	(3672,942)(3672,980)(3673,1017)
	(3674,1054)(3674,1091)(3675,1126)
	(3675,1162)(3675,1197)(3675,1232)
	(3675,1268)(3674,1303)(3674,1340)
	(3673,1377)(3672,1414)(3672,1452)
	(3672,1491)(3671,1531)(3671,1571)
	(3671,1612)(3670,1653)(3670,1694)
	(3670,1736)(3670,1778)(3669,1819)
	(3669,1860)(3668,1901)(3667,1941)
	(3665,1980)(3663,2019)(3661,2055)
	(3658,2091)(3655,2125)(3651,2157)
	(3646,2187)(3641,2215)(3634,2241)
	(3627,2265)(3619,2287)(3610,2306)
	(3600,2322)(3586,2339)(3571,2352)
	(3553,2363)(3534,2370)(3514,2374)
	(3491,2377)(3467,2377)(3441,2375)
	(3415,2372)(3387,2367)(3358,2361)
	(3329,2354)(3300,2345)(3271,2337)
	(3242,2327)(3213,2317)(3185,2306)
	(3159,2295)(3133,2283)(3109,2271)
	(3086,2257)(3066,2243)(3047,2228)
	(3029,2211)(3014,2192)(3000,2172)
\put(0,1272){\makebox(0,0)[lb]{\smash{{{\SetFigFont{6}{7.2}{\rmdefault}{\mddefault}{\updefault}$(A^i,1)=$}}}}}
\end{picture}
}

%% file: xfig/AiAi.eepic
\setlength{\unitlength}{0.00041667in}
\begingroup\makeatletter\ifx\SetFigFont\undefined%
\gdef\SetFigFont#1#2#3#4#5{%
  \reset@font\fontsize{#1}{#2pt}%
  \fontfamily{#3}\fontseries{#4}\fontshape{#5}%
  \selectfont}%
\fi\endgroup%
{\renewcommand{\dashlinestretch}{30}
\begin{picture}(5100,2466)(0,-10)
\put(1800,1947){\blacken\ellipse{300}{300}}
\put(1800,1947){\ellipse{300}{300}}
\put(2700,1947){\blacken\ellipse{300}{300}}
\put(2700,1947){\ellipse{300}{300}}
\put(2700,447){\blacken\ellipse{300}{300}}
\put(2700,447){\ellipse{300}{300}}
\put(1800,447){\blacken\ellipse{300}{300}}
\put(1800,447){\ellipse{300}{300}}
\put(3600,1947){\blacken\ellipse{300}{300}}
\put(3600,1947){\ellipse{300}{300}}
\put(4500,1947){\blacken\ellipse{300}{300}}
\put(4500,1947){\ellipse{300}{300}}
\put(4500,447){\blacken\ellipse{300}{300}}
\put(4500,447){\ellipse{300}{300}}
\put(3600,447){\blacken\ellipse{300}{300}}
\put(3600,447){\ellipse{300}{300}}
\thicklines
\path(2700,1947)(2700,447)
\path(1800,1947)(1800,447)
\path(4500,1947)(4500,447)
\thinlines
\path(1425,1947)(1417,1921)(1409,1893)
	(1402,1863)(1396,1832)(1391,1799)
	(1386,1764)(1382,1728)(1378,1690)
	(1375,1650)(1372,1609)(1370,1567)
	(1368,1523)(1366,1478)(1365,1432)
	(1364,1385)(1363,1338)(1362,1290)
	(1362,1242)(1362,1193)(1362,1145)
	(1362,1097)(1363,1049)(1364,1002)
	(1365,955)(1366,910)(1368,865)
	(1370,821)(1372,779)(1375,738)
	(1378,699)(1382,662)(1386,626)
	(1391,591)(1396,559)(1402,528)
	(1409,500)(1417,472)(1425,447)
	(1438,415)(1453,385)(1469,358)
	(1488,332)(1508,308)(1531,285)
	(1555,263)(1580,242)(1607,222)
	(1635,203)(1664,184)(1694,166)
	(1725,149)(1755,133)(1786,117)
	(1817,103)(1847,90)(1876,77)
	(1905,67)(1933,58)(1959,50)
	(1984,45)(2007,42)(2029,42)
	(2049,45)(2068,50)(2085,59)
	(2100,72)(2112,85)(2123,102)
	(2133,121)(2142,142)(2150,166)
	(2158,193)(2165,222)(2172,253)
	(2178,286)(2184,321)(2189,357)
	(2194,395)(2198,435)(2203,475)
	(2207,517)(2210,559)(2214,602)
	(2217,645)(2220,688)(2223,731)
	(2226,774)(2229,816)(2231,858)
	(2234,900)(2236,941)(2239,981)
	(2241,1020)(2243,1058)(2245,1095)
	(2247,1132)(2248,1168)(2249,1203)
	(2250,1238)(2250,1272)(2250,1308)
	(2249,1345)(2248,1381)(2246,1419)
	(2244,1457)(2242,1495)(2240,1534)
	(2238,1574)(2235,1615)(2232,1656)
	(2230,1698)(2227,1741)(2224,1783)
	(2220,1826)(2217,1869)(2214,1912)
	(2210,1954)(2206,1996)(2202,2037)
	(2197,2076)(2192,2115)(2187,2152)
	(2181,2187)(2175,2220)(2168,2252)
	(2161,2280)(2153,2307)(2144,2331)
	(2134,2352)(2124,2370)(2112,2385)
	(2100,2397)(2086,2407)(2070,2413)
	(2053,2416)(2035,2416)(2015,2414)
	(1993,2410)(1971,2403)(1947,2395)
	(1922,2384)(1896,2372)(1869,2359)
	(1841,2345)(1813,2329)(1784,2312)
	(1755,2295)(1727,2277)(1698,2258)
	(1670,2239)(1643,2219)(1616,2199)
	(1591,2178)(1566,2156)(1543,2134)
	(1522,2111)(1501,2087)(1483,2062)
	(1466,2035)(1451,2008)(1437,1978)(1425,1947)
\path(3075,2172)(3083,2151)(3090,2129)
	(3096,2105)(3101,2079)(3105,2052)
	(3107,2023)(3109,1993)(3110,1962)
	(3110,1930)(3110,1896)(3108,1862)
	(3107,1827)(3105,1792)(3102,1756)
	(3100,1720)(3097,1683)(3094,1647)
	(3091,1610)(3088,1574)(3086,1537)
	(3083,1502)(3081,1466)(3080,1431)
	(3078,1396)(3077,1362)(3076,1328)
	(3075,1295)(3075,1262)(3075,1230)
	(3075,1197)(3075,1164)(3075,1132)
	(3075,1099)(3076,1065)(3077,1032)
	(3078,997)(3080,962)(3081,927)
	(3083,891)(3086,855)(3088,818)
	(3091,782)(3094,745)(3097,708)
	(3100,671)(3102,634)(3105,598)
	(3107,562)(3108,527)(3110,493)
	(3110,459)(3110,427)(3109,396)
	(3107,366)(3105,338)(3101,311)
	(3096,286)(3090,263)(3083,242)
	(3075,222)(3064,202)(3052,183)
	(3037,166)(3021,151)(3004,136)
	(2985,123)(2964,110)(2942,98)
	(2919,86)(2895,75)(2870,64)
	(2845,54)(2820,45)(2794,37)
	(2769,29)(2743,23)(2719,18)
	(2695,14)(2673,12)(2651,12)
	(2631,15)(2612,20)(2595,28)
	(2578,39)(2564,54)(2550,72)
	(2540,89)(2531,108)(2523,129)
	(2515,152)(2508,178)(2501,206)
	(2494,236)(2488,267)(2482,301)
	(2476,336)(2471,372)(2466,410)
	(2461,449)(2456,489)(2452,529)
	(2448,571)(2443,612)(2439,654)
	(2436,696)(2432,738)(2428,780)
	(2425,822)(2421,863)(2418,903)
	(2415,943)(2413,983)(2410,1021)
	(2408,1059)(2405,1096)(2404,1132)
	(2402,1168)(2401,1203)(2400,1238)
	(2400,1272)(2400,1308)(2401,1345)
	(2402,1381)(2404,1418)(2406,1456)
	(2408,1494)(2411,1533)(2414,1572)
	(2417,1612)(2420,1652)(2424,1693)
	(2427,1735)(2431,1776)(2435,1818)
	(2439,1860)(2443,1901)(2448,1942)
	(2452,1983)(2457,2023)(2462,2062)
	(2467,2100)(2473,2137)(2479,2172)
	(2485,2206)(2491,2237)(2498,2267)
	(2505,2295)(2513,2320)(2521,2343)
	(2530,2364)(2540,2382)(2550,2397)
	(2565,2413)(2581,2426)(2599,2434)
	(2618,2438)(2639,2439)(2662,2437)
	(2686,2433)(2711,2426)(2737,2418)
	(2764,2407)(2792,2396)(2820,2383)
	(2847,2369)(2874,2355)(2901,2340)
	(2927,2324)(2951,2308)(2975,2292)
	(2996,2274)(3016,2256)(3034,2238)
	(3049,2218)(3063,2196)(3075,2172)
\path(5025,1947)(5033,1921)(5041,1893)
	(5048,1863)(5054,1832)(5059,1799)
	(5064,1764)(5068,1728)(5072,1690)
	(5075,1650)(5078,1609)(5080,1567)
	(5082,1523)(5084,1478)(5085,1432)
	(5086,1385)(5087,1338)(5088,1290)
	(5088,1242)(5088,1193)(5088,1145)
	(5088,1097)(5087,1049)(5086,1002)
	(5085,955)(5084,910)(5082,865)
	(5080,821)(5078,779)(5075,738)
	(5072,699)(5068,662)(5064,626)
	(5059,591)(5054,559)(5048,528)
	(5041,500)(5033,472)(5025,447)
	(5012,415)(4997,385)(4981,358)
	(4962,332)(4942,308)(4919,285)
	(4895,263)(4870,242)(4843,222)
	(4815,203)(4786,184)(4756,166)
	(4725,149)(4695,133)(4664,117)
	(4633,103)(4603,90)(4574,77)
	(4545,67)(4517,58)(4491,50)
	(4466,45)(4443,42)(4421,42)
	(4401,45)(4382,50)(4365,59)
	(4350,72)(4338,85)(4327,102)
	(4317,121)(4308,142)(4300,166)
	(4292,193)(4285,222)(4278,253)
	(4272,286)(4266,321)(4261,357)
	(4256,395)(4252,435)(4247,475)
	(4243,517)(4240,559)(4236,602)
	(4233,645)(4230,688)(4227,731)
	(4224,774)(4221,816)(4219,858)
	(4216,900)(4214,941)(4211,981)
	(4209,1020)(4207,1058)(4205,1095)
	(4203,1132)(4202,1168)(4201,1203)
	(4200,1238)(4200,1272)(4200,1308)
	(4201,1345)(4202,1381)(4204,1419)
	(4206,1457)(4208,1495)(4210,1534)
	(4212,1574)(4215,1615)(4218,1656)
	(4220,1698)(4223,1741)(4226,1783)
	(4230,1826)(4233,1869)(4236,1912)
	(4240,1954)(4244,1996)(4248,2037)
	(4253,2076)(4258,2115)(4263,2152)
	(4269,2187)(4275,2220)(4282,2252)
	(4289,2280)(4297,2307)(4306,2331)
	(4316,2352)(4326,2370)(4338,2385)
	(4350,2397)(4364,2407)(4380,2413)
	(4397,2416)(4415,2416)(4435,2414)
	(4457,2410)(4479,2403)(4503,2395)
	(4528,2384)(4554,2372)(4581,2359)
	(4609,2345)(4637,2329)(4666,2312)
	(4695,2295)(4723,2277)(4752,2258)
	(4780,2239)(4807,2219)(4834,2199)
	(4859,2178)(4884,2156)(4907,2134)
	(4928,2111)(4949,2087)(4967,2062)
	(4984,2035)(4999,2008)(5013,1978)(5025,1947)
\path(3300,222)(3311,209)(3324,196)
	(3341,184)(3360,173)(3382,162)
	(3407,151)(3435,141)(3465,130)
	(3496,120)(3529,110)(3562,100)
	(3596,91)(3631,82)(3664,74)
	(3697,66)(3729,60)(3759,55)
	(3787,51)(3812,49)(3835,49)
	(3856,51)(3873,56)(3888,62)
	(3900,72)(3909,83)(3916,97)
	(3922,114)(3926,133)(3928,155)
	(3929,180)(3929,207)(3927,236)
	(3925,267)(3922,300)(3918,333)
	(3913,368)(3908,402)(3903,437)
	(3897,472)(3892,506)(3886,539)
	(3880,570)(3874,600)(3867,628)
	(3861,654)(3855,677)(3848,698)
	(3841,717)(3833,733)(3825,747)
	(3813,762)(3799,775)(3784,784)
	(3767,791)(3749,797)(3729,801)
	(3708,803)(3686,804)(3664,804)
	(3641,803)(3618,802)(3596,799)
	(3574,797)(3552,793)(3532,789)
	(3513,783)(3495,777)(3479,769)
	(3464,759)(3450,747)(3438,734)
	(3427,719)(3417,701)(3406,681)
	(3395,658)(3385,634)(3373,608)
	(3362,580)(3351,551)(3340,522)
	(3329,492)(3319,461)(3310,431)
	(3301,402)(3294,374)(3288,347)
	(3284,321)(3282,298)(3283,276)
	(3286,256)(3291,238)(3300,222)
\path(3300,1722)(3311,1709)(3324,1696)
	(3341,1684)(3360,1673)(3382,1662)
	(3407,1651)(3435,1641)(3465,1630)
	(3496,1620)(3529,1610)(3562,1600)
	(3596,1591)(3631,1582)(3664,1574)
	(3697,1566)(3729,1560)(3759,1555)
	(3787,1551)(3812,1549)(3835,1549)
	(3856,1551)(3873,1556)(3888,1562)
	(3900,1572)(3909,1583)(3916,1597)
	(3922,1614)(3926,1633)(3928,1655)
	(3929,1680)(3929,1707)(3927,1736)
	(3925,1767)(3922,1800)(3918,1833)
	(3913,1868)(3908,1902)(3903,1937)
	(3897,1972)(3892,2006)(3886,2039)
	(3880,2070)(3874,2100)(3867,2128)
	(3861,2154)(3855,2177)(3848,2198)
	(3841,2217)(3833,2233)(3825,2247)
	(3813,2262)(3799,2275)(3784,2284)
	(3767,2291)(3749,2297)(3729,2301)
	(3708,2303)(3686,2304)(3664,2304)
	(3641,2303)(3618,2302)(3596,2299)
	(3574,2297)(3552,2293)(3532,2289)
	(3513,2283)(3495,2277)(3479,2269)
	(3464,2259)(3450,2247)(3438,2234)
	(3427,2219)(3417,2201)(3406,2181)
	(3395,2158)(3385,2134)(3373,2108)
	(3362,2080)(3351,2051)(3340,2022)
	(3329,1992)(3319,1961)(3310,1931)
	(3301,1902)(3294,1874)(3288,1847)
	(3284,1821)(3282,1798)(3283,1776)
	(3286,1756)(3291,1738)(3300,1722)
\put(0,1272){\makebox(0,0)[lb]{\smash{{{\SetFigFont{6}{7.2}{\rmdefault}{\mddefault}{\updefault}$(A^i,A^i)=$}}}}}
\end{picture}
}

%% file: xfig/OneAij.eepic
\setlength{\unitlength}{0.00041667in}
\begingroup\makeatletter\ifx\SetFigFont\undefined%
\gdef\SetFigFont#1#2#3#4#5{%
  \reset@font\fontsize{#1}{#2pt}%
  \fontfamily{#3}\fontseries{#4}\fontshape{#5}%
  \selectfont}%
\fi\endgroup%
{\renewcommand{\dashlinestretch}{30}
\begin{picture}(5100,2529)(0,-10)
\put(1800,1996){\blacken\ellipse{300}{300}}
\put(1800,1996){\ellipse{300}{300}}
\put(2700,1996){\blacken\ellipse{300}{300}}
\put(2700,1996){\ellipse{300}{300}}
\put(2700,496){\blacken\ellipse{300}{300}}
\put(2700,496){\ellipse{300}{300}}
\put(1800,496){\blacken\ellipse{300}{300}}
\put(1800,496){\ellipse{300}{300}}
\put(3600,1996){\blacken\ellipse{300}{300}}
\put(3600,1996){\ellipse{300}{300}}
\put(4500,1996){\blacken\ellipse{300}{300}}
\put(4500,1996){\ellipse{300}{300}}
\put(4500,496){\blacken\ellipse{300}{300}}
\put(4500,496){\ellipse{300}{300}}
\put(3600,496){\blacken\ellipse{300}{300}}
\put(3600,496){\ellipse{300}{300}}
\thicklines
\path(2700,1996)(2700,496)
\path(1800,1996)(1800,496)
\path(4500,1996)(4500,496)
\path(3600,1996)(3600,496)
\thinlines
\path(1425,1996)(1417,1970)(1409,1942)
	(1402,1912)(1396,1881)(1391,1848)
	(1386,1813)(1382,1777)(1378,1739)
	(1375,1699)(1372,1658)(1370,1616)
	(1368,1572)(1366,1527)(1365,1481)
	(1364,1434)(1363,1387)(1362,1339)
	(1362,1291)(1362,1242)(1362,1194)
	(1362,1146)(1363,1098)(1364,1051)
	(1365,1004)(1366,959)(1368,914)
	(1370,870)(1372,828)(1375,787)
	(1378,748)(1382,711)(1386,675)
	(1391,640)(1396,608)(1402,577)
	(1409,549)(1417,521)(1425,496)
	(1438,464)(1453,434)(1469,407)
	(1488,381)(1508,357)(1531,334)
	(1555,312)(1580,291)(1607,271)
	(1635,252)(1664,233)(1694,215)
	(1725,198)(1755,182)(1786,166)
	(1817,152)(1847,139)(1876,126)
	(1905,116)(1933,107)(1959,99)
	(1984,94)(2007,91)(2029,91)
	(2049,94)(2068,99)(2085,108)
	(2100,121)(2112,134)(2123,151)
	(2133,170)(2142,191)(2150,215)
	(2158,242)(2165,271)(2172,302)
	(2178,335)(2184,370)(2189,406)
	(2194,444)(2198,484)(2203,524)
	(2207,566)(2210,608)(2214,651)
	(2217,694)(2220,737)(2223,780)
	(2226,823)(2229,865)(2231,907)
	(2234,949)(2236,990)(2239,1030)
	(2241,1069)(2243,1107)(2245,1144)
	(2247,1181)(2248,1217)(2249,1252)
	(2250,1287)(2250,1321)(2250,1357)
	(2249,1394)(2248,1430)(2246,1468)
	(2244,1506)(2242,1544)(2240,1583)
	(2238,1623)(2235,1664)(2232,1705)
	(2230,1747)(2227,1790)(2224,1832)
	(2220,1875)(2217,1918)(2214,1961)
	(2210,2003)(2206,2045)(2202,2086)
	(2197,2125)(2192,2164)(2187,2201)
	(2181,2236)(2175,2269)(2168,2301)
	(2161,2329)(2153,2356)(2144,2380)
	(2134,2401)(2124,2419)(2112,2434)
	(2100,2446)(2086,2456)(2070,2462)
	(2053,2465)(2035,2465)(2015,2463)
	(1993,2459)(1971,2452)(1947,2444)
	(1922,2433)(1896,2421)(1869,2408)
	(1841,2394)(1813,2378)(1784,2361)
	(1755,2344)(1727,2326)(1698,2307)
	(1670,2288)(1643,2268)(1616,2248)
	(1591,2227)(1566,2205)(1543,2183)
	(1522,2160)(1501,2136)(1483,2111)
	(1466,2084)(1451,2057)(1437,2027)(1425,1996)
\path(3750,2296)(3770,2276)(3788,2255)
	(3804,2232)(3819,2207)(3831,2181)
	(3842,2153)(3852,2124)(3860,2094)
	(3867,2062)(3872,2029)(3877,1995)
	(3880,1961)(3883,1926)(3885,1890)
	(3887,1853)(3888,1817)(3889,1780)
	(3890,1743)(3890,1706)(3891,1669)
	(3891,1632)(3892,1595)(3892,1559)
	(3893,1523)(3894,1487)(3895,1452)
	(3896,1417)(3897,1382)(3898,1348)
	(3899,1314)(3900,1280)(3900,1246)
	(3900,1212)(3899,1177)(3898,1143)
	(3897,1107)(3896,1071)(3895,1034)
	(3894,996)(3893,958)(3892,919)
	(3892,880)(3891,839)(3891,799)
	(3890,758)(3890,717)(3889,675)
	(3888,634)(3887,593)(3885,553)
	(3883,513)(3880,474)(3877,436)
	(3872,399)(3867,363)(3860,329)
	(3852,296)(3842,265)(3831,236)
	(3819,209)(3804,184)(3788,161)
	(3770,140)(3750,121)(3729,105)
	(3706,92)(3682,79)(3655,69)
	(3626,59)(3595,51)(3562,44)
	(3527,37)(3490,32)(3451,27)
	(3411,23)(3370,20)(3327,17)
	(3284,15)(3240,14)(3195,13)
	(3150,12)(3105,12)(3060,12)
	(3016,13)(2973,15)(2930,17)
	(2889,20)(2849,23)(2810,28)
	(2773,33)(2738,39)(2705,47)
	(2674,55)(2645,65)(2618,77)
	(2594,90)(2571,104)(2550,121)
	(2531,140)(2514,160)(2498,183)
	(2485,208)(2472,235)(2461,264)
	(2452,294)(2444,326)(2437,360)
	(2431,395)(2426,431)(2422,469)
	(2419,507)(2416,547)(2414,587)
	(2413,627)(2412,668)(2411,710)
	(2410,751)(2410,792)(2409,833)
	(2409,874)(2408,915)(2408,955)
	(2407,994)(2407,1033)(2406,1071)
	(2405,1109)(2404,1145)(2403,1182)
	(2402,1217)(2401,1252)(2400,1287)
	(2400,1321)(2400,1357)(2401,1394)
	(2402,1430)(2403,1467)(2404,1505)
	(2405,1543)(2406,1581)(2407,1620)
	(2408,1660)(2408,1700)(2409,1741)
	(2409,1782)(2410,1823)(2410,1864)
	(2411,1906)(2412,1947)(2413,1988)
	(2415,2028)(2417,2068)(2420,2106)
	(2423,2144)(2428,2181)(2433,2216)
	(2440,2250)(2448,2281)(2458,2312)
	(2469,2340)(2481,2366)(2496,2389)
	(2512,2411)(2530,2430)(2550,2446)
	(2571,2459)(2594,2470)(2618,2480)
	(2645,2487)(2674,2493)(2705,2497)
	(2738,2500)(2773,2501)(2810,2502)
	(2849,2501)(2889,2499)(2930,2497)
	(2973,2493)(3016,2489)(3060,2484)
	(3105,2479)(3150,2473)(3195,2467)
	(3240,2460)(3284,2453)(3327,2445)
	(3370,2437)(3411,2429)(3451,2420)
	(3490,2411)(3527,2401)(3562,2391)
	(3595,2380)(3626,2368)(3655,2356)
	(3682,2343)(3706,2328)(3729,2313)(3750,2296)
\path(5025,1996)(5033,1970)(5041,1942)
	(5048,1912)(5054,1881)(5059,1848)
	(5064,1813)(5068,1777)(5072,1739)
	(5075,1699)(5078,1658)(5080,1616)
	(5082,1572)(5084,1527)(5085,1481)
	(5086,1434)(5087,1387)(5088,1339)
	(5088,1291)(5088,1242)(5088,1194)
	(5088,1146)(5087,1098)(5086,1051)
	(5085,1004)(5084,959)(5082,914)
	(5080,870)(5078,828)(5075,787)
	(5072,748)(5068,711)(5064,675)
	(5059,640)(5054,608)(5048,577)
	(5041,549)(5033,521)(5025,496)
	(5012,464)(4997,434)(4981,407)
	(4962,381)(4942,357)(4919,334)
	(4895,312)(4870,291)(4843,271)
	(4815,252)(4786,233)(4756,215)
	(4725,198)(4695,182)(4664,166)
	(4633,152)(4603,139)(4574,126)
	(4545,116)(4517,107)(4491,99)
	(4466,94)(4443,91)(4421,91)
	(4401,94)(4382,99)(4365,108)
	(4350,121)(4338,134)(4327,151)
	(4317,170)(4308,191)(4300,215)
	(4292,242)(4285,271)(4278,302)
	(4272,335)(4266,370)(4261,406)
	(4256,444)(4252,484)(4247,524)
	(4243,566)(4240,608)(4236,651)
	(4233,694)(4230,737)(4227,780)
	(4224,823)(4221,865)(4219,907)
	(4216,949)(4214,990)(4211,1030)
	(4209,1069)(4207,1107)(4205,1144)
	(4203,1181)(4202,1217)(4201,1252)
	(4200,1287)(4200,1321)(4200,1357)
	(4201,1394)(4202,1430)(4204,1468)
	(4206,1506)(4208,1544)(4210,1583)
	(4212,1623)(4215,1664)(4218,1705)
	(4220,1747)(4223,1790)(4226,1832)
	(4230,1875)(4233,1918)(4236,1961)
	(4240,2003)(4244,2045)(4248,2086)
	(4253,2125)(4258,2164)(4263,2201)
	(4269,2236)(4275,2269)(4282,2301)
	(4289,2329)(4297,2356)(4306,2380)
	(4316,2401)(4326,2419)(4338,2434)
	(4350,2446)(4364,2456)(4380,2462)
	(4397,2465)(4415,2465)(4435,2463)
	(4457,2459)(4479,2452)(4503,2444)
	(4528,2433)(4554,2421)(4581,2408)
	(4609,2394)(4637,2378)(4666,2361)
	(4695,2344)(4723,2326)(4752,2307)
	(4780,2288)(4807,2268)(4834,2248)
	(4859,2227)(4884,2205)(4907,2183)
	(4928,2160)(4949,2136)(4967,2111)
	(4984,2084)(4999,2057)(5013,2027)(5025,1996)
\put(0,1321){\makebox(0,0)[lb]{\smash{{{\SetFigFont{6}{7.2}{\rmdefault}{\mddefault}{\updefault}$(1,A^{ij})=$}}}}}
\end{picture}
}

%% file: xfig/AijAij.eepic
\setlength{\unitlength}{0.00041667in}
\begingroup\makeatletter\ifx\SetFigFont\undefined%
\gdef\SetFigFont#1#2#3#4#5{%
  \reset@font\fontsize{#1}{#2pt}%
  \fontfamily{#3}\fontseries{#4}\fontshape{#5}%
  \selectfont}%
\fi\endgroup%
{\renewcommand{\dashlinestretch}{30}
\begin{picture}(5325,2529)(0,-10)
\put(2025,1996){\blacken\ellipse{300}{300}}
\put(2025,1996){\ellipse{300}{300}}
\put(2925,1996){\blacken\ellipse{300}{300}}
\put(2925,1996){\ellipse{300}{300}}
\put(2925,496){\blacken\ellipse{300}{300}}
\put(2925,496){\ellipse{300}{300}}
\put(2025,496){\blacken\ellipse{300}{300}}
\put(2025,496){\ellipse{300}{300}}
\put(3825,1996){\blacken\ellipse{300}{300}}
\put(3825,1996){\ellipse{300}{300}}
\put(4725,1996){\blacken\ellipse{300}{300}}
\put(4725,1996){\ellipse{300}{300}}
\put(4725,496){\blacken\ellipse{300}{300}}
\put(4725,496){\ellipse{300}{300}}
\put(3825,496){\blacken\ellipse{300}{300}}
\put(3825,496){\ellipse{300}{300}}
\thicklines
\path(2925,1996)(2925,496)
\path(2025,1996)(2025,496)
\path(4725,1996)(4725,496)
\path(3825,1996)(3825,496)
\path(2925,1996)(3825,1996)
\thinlines
\path(1650,1996)(1642,1970)(1634,1942)
	(1627,1912)(1621,1881)(1616,1848)
	(1611,1813)(1607,1777)(1603,1739)
	(1600,1699)(1597,1658)(1595,1616)
	(1593,1572)(1591,1527)(1590,1481)
	(1589,1434)(1588,1387)(1587,1339)
	(1587,1291)(1587,1242)(1587,1194)
	(1587,1146)(1588,1098)(1589,1051)
	(1590,1004)(1591,959)(1593,914)
	(1595,870)(1597,828)(1600,787)
	(1603,748)(1607,711)(1611,675)
	(1616,640)(1621,608)(1627,577)
	(1634,549)(1642,521)(1650,496)
	(1663,464)(1678,434)(1694,407)
	(1713,381)(1733,357)(1756,334)
	(1780,312)(1805,291)(1832,271)
	(1860,252)(1889,233)(1919,215)
	(1950,198)(1980,182)(2011,166)
	(2042,152)(2072,139)(2101,126)
	(2130,116)(2158,107)(2184,99)
	(2209,94)(2232,91)(2254,91)
	(2274,94)(2293,99)(2310,108)
	(2325,121)(2337,134)(2348,151)
	(2358,170)(2367,191)(2375,215)
	(2383,242)(2390,271)(2397,302)
	(2403,335)(2409,370)(2414,406)
	(2419,444)(2423,484)(2428,524)
	(2432,566)(2435,608)(2439,651)
	(2442,694)(2445,737)(2448,780)
	(2451,823)(2454,865)(2456,907)
	(2459,949)(2461,990)(2464,1030)
	(2466,1069)(2468,1107)(2470,1144)
	(2472,1181)(2473,1217)(2474,1252)
	(2475,1287)(2475,1321)(2475,1357)
	(2474,1394)(2473,1430)(2471,1468)
	(2469,1506)(2467,1544)(2465,1583)
	(2463,1623)(2460,1664)(2457,1705)
	(2455,1747)(2452,1790)(2449,1832)
	(2445,1875)(2442,1918)(2439,1961)
	(2435,2003)(2431,2045)(2427,2086)
	(2422,2125)(2417,2164)(2412,2201)
	(2406,2236)(2400,2269)(2393,2301)
	(2386,2329)(2378,2356)(2369,2380)
	(2359,2401)(2349,2419)(2337,2434)
	(2325,2446)(2311,2456)(2295,2462)
	(2278,2465)(2260,2465)(2240,2463)
	(2218,2459)(2196,2452)(2172,2444)
	(2147,2433)(2121,2421)(2094,2408)
	(2066,2394)(2038,2378)(2009,2361)
	(1980,2344)(1952,2326)(1923,2307)
	(1895,2288)(1868,2268)(1841,2248)
	(1816,2227)(1791,2205)(1768,2183)
	(1747,2160)(1726,2136)(1708,2111)
	(1691,2084)(1676,2057)(1662,2027)(1650,1996)
\path(3975,2296)(3995,2276)(4013,2255)
	(4029,2232)(4044,2207)(4056,2181)
	(4067,2153)(4077,2124)(4085,2094)
	(4092,2062)(4097,2029)(4102,1995)
	(4105,1961)(4108,1926)(4110,1890)
	(4112,1853)(4113,1817)(4114,1780)
	(4115,1743)(4115,1706)(4116,1669)
	(4116,1632)(4117,1595)(4117,1559)
	(4118,1523)(4119,1487)(4120,1452)
	(4121,1417)(4122,1382)(4123,1348)
	(4124,1314)(4125,1280)(4125,1246)
	(4125,1212)(4124,1177)(4123,1143)
	(4122,1107)(4121,1071)(4120,1034)
	(4119,996)(4118,958)(4117,919)
	(4117,880)(4116,839)(4116,799)
	(4115,758)(4115,717)(4114,675)
	(4113,634)(4112,593)(4110,553)
	(4108,513)(4105,474)(4102,436)
	(4097,399)(4092,363)(4085,329)
	(4077,296)(4067,265)(4056,236)
	(4044,209)(4029,184)(4013,161)
	(3995,140)(3975,121)(3954,105)
	(3931,92)(3907,79)(3880,69)
	(3851,59)(3820,51)(3787,44)
	(3752,37)(3715,32)(3676,27)
	(3636,23)(3595,20)(3552,17)
	(3509,15)(3465,14)(3420,13)
	(3375,12)(3330,12)(3285,12)
	(3241,13)(3198,15)(3155,17)
	(3114,20)(3074,23)(3035,28)
	(2998,33)(2963,39)(2930,47)
	(2899,55)(2870,65)(2843,77)
	(2819,90)(2796,104)(2775,121)
	(2756,140)(2739,160)(2723,183)
	(2710,208)(2697,235)(2686,264)
	(2677,294)(2669,326)(2662,360)
	(2656,395)(2651,431)(2647,469)
	(2644,507)(2641,547)(2639,587)
	(2638,627)(2637,668)(2636,710)
	(2635,751)(2635,792)(2634,833)
	(2634,874)(2633,915)(2633,955)
	(2632,994)(2632,1033)(2631,1071)
	(2630,1109)(2629,1145)(2628,1182)
	(2627,1217)(2626,1252)(2625,1287)
	(2625,1321)(2625,1357)(2626,1394)
	(2627,1430)(2628,1467)(2629,1505)
	(2630,1543)(2631,1581)(2632,1620)
	(2633,1660)(2633,1700)(2634,1741)
	(2634,1782)(2635,1823)(2635,1864)
	(2636,1906)(2637,1947)(2638,1988)
	(2640,2028)(2642,2068)(2645,2106)
	(2648,2144)(2653,2181)(2658,2216)
	(2665,2250)(2673,2281)(2683,2312)
	(2694,2340)(2706,2366)(2721,2389)
	(2737,2411)(2755,2430)(2775,2446)
	(2796,2459)(2819,2470)(2843,2480)
	(2870,2487)(2899,2493)(2930,2497)
	(2963,2500)(2998,2501)(3035,2502)
	(3074,2501)(3114,2499)(3155,2497)
	(3198,2493)(3241,2489)(3285,2484)
	(3330,2479)(3375,2473)(3420,2467)
	(3465,2460)(3509,2453)(3552,2445)
	(3595,2437)(3636,2429)(3676,2420)
	(3715,2411)(3752,2401)(3787,2391)
	(3820,2380)(3851,2368)(3880,2356)
	(3907,2343)(3931,2328)(3954,2313)(3975,2296)
\path(5250,1996)(5258,1970)(5266,1942)
	(5273,1912)(5279,1881)(5284,1848)
	(5289,1813)(5293,1777)(5297,1739)
	(5300,1699)(5303,1658)(5305,1616)
	(5307,1572)(5309,1527)(5310,1481)
	(5311,1434)(5312,1387)(5313,1339)
	(5313,1291)(5313,1242)(5313,1194)
	(5313,1146)(5312,1098)(5311,1051)
	(5310,1004)(5309,959)(5307,914)
	(5305,870)(5303,828)(5300,787)
	(5297,748)(5293,711)(5289,675)
	(5284,640)(5279,608)(5273,577)
	(5266,549)(5258,521)(5250,496)
	(5237,464)(5222,434)(5206,407)
	(5187,381)(5167,357)(5144,334)
	(5120,312)(5095,291)(5068,271)
	(5040,252)(5011,233)(4981,215)
	(4950,198)(4920,182)(4889,166)
	(4858,152)(4828,139)(4799,126)
	(4770,116)(4742,107)(4716,99)
	(4691,94)(4668,91)(4646,91)
	(4626,94)(4607,99)(4590,108)
	(4575,121)(4563,134)(4552,151)
	(4542,170)(4533,191)(4525,215)
	(4517,242)(4510,271)(4503,302)
	(4497,335)(4491,370)(4486,406)
	(4481,444)(4477,484)(4472,524)
	(4468,566)(4465,608)(4461,651)
	(4458,694)(4455,737)(4452,780)
	(4449,823)(4446,865)(4444,907)
	(4441,949)(4439,990)(4436,1030)
	(4434,1069)(4432,1107)(4430,1144)
	(4428,1181)(4427,1217)(4426,1252)
	(4425,1287)(4425,1321)(4425,1357)
	(4426,1394)(4427,1430)(4429,1468)
	(4431,1506)(4433,1544)(4435,1583)
	(4437,1623)(4440,1664)(4443,1705)
	(4445,1747)(4448,1790)(4451,1832)
	(4455,1875)(4458,1918)(4461,1961)
	(4465,2003)(4469,2045)(4473,2086)
	(4478,2125)(4483,2164)(4488,2201)
	(4494,2236)(4500,2269)(4507,2301)
	(4514,2329)(4522,2356)(4531,2380)
	(4541,2401)(4551,2419)(4563,2434)
	(4575,2446)(4589,2456)(4605,2462)
	(4622,2465)(4640,2465)(4660,2463)
	(4682,2459)(4704,2452)(4728,2444)
	(4753,2433)(4779,2421)(4806,2408)
	(4834,2394)(4862,2378)(4891,2361)
	(4920,2344)(4948,2326)(4977,2307)
	(5005,2288)(5032,2268)(5059,2248)
	(5084,2227)(5109,2205)(5132,2183)
	(5153,2160)(5174,2136)(5192,2111)
	(5209,2084)(5224,2057)(5238,2027)(5250,1996)
\put(0,1321){\makebox(0,0)[lb]{\smash{{{\SetFigFont{6}{7.2}{\rmdefault}{\mddefault}{\updefault}$(A^{ij},A^{ij})=$}}}}}
\end{picture}
}

%% file: xfig/Sigmai.eepic
\setlength{\unitlength}{0.00041667in}
\begingroup\makeatletter\ifx\SetFigFont\undefined%
\gdef\SetFigFont#1#2#3#4#5{%
  \reset@font\fontsize{#1}{#2pt}%
  \fontfamily{#3}\fontseries{#4}\fontshape{#5}%
  \selectfont}%
\fi\endgroup%
{\renewcommand{\dashlinestretch}{30}
\begin{picture}(5250,2413)(0,-10)
\put(1950,1917){\blacken\ellipse{300}{300}}
\put(1950,1917){\ellipse{300}{300}}
\put(1950,417){\blacken\ellipse{300}{300}}
\put(1950,417){\ellipse{300}{300}}
\put(3750,1917){\blacken\ellipse{300}{300}}
\put(3750,1917){\ellipse{300}{300}}
\put(4650,1917){\blacken\ellipse{300}{300}}
\put(4650,1917){\ellipse{300}{300}}
\put(4650,417){\blacken\ellipse{300}{300}}
\put(4650,417){\ellipse{300}{300}}
\put(2850,417){\blacken\ellipse{300}{300}}
\put(2850,417){\ellipse{300}{300}}
\thicklines
\path(1950,1917)(1950,417)
\path(4650,1917)(4650,417)
\path(3750,1917)(2850,417)
\thinlines
\path(1575,1917)(1567,1891)(1559,1863)
	(1552,1833)(1546,1802)(1541,1769)
	(1536,1734)(1532,1698)(1528,1660)
	(1525,1620)(1522,1579)(1520,1537)
	(1518,1493)(1516,1448)(1515,1402)
	(1514,1355)(1513,1308)(1512,1260)
	(1512,1212)(1512,1163)(1512,1115)
	(1512,1067)(1513,1019)(1514,972)
	(1515,925)(1516,880)(1518,835)
	(1520,791)(1522,749)(1525,708)
	(1528,669)(1532,632)(1536,596)
	(1541,561)(1546,529)(1552,498)
	(1559,470)(1567,442)(1575,417)
	(1588,385)(1603,355)(1619,328)
	(1638,302)(1658,278)(1681,255)
	(1705,233)(1730,212)(1757,192)
	(1785,173)(1814,154)(1844,136)
	(1875,119)(1905,103)(1936,87)
	(1967,73)(1997,60)(2026,47)
	(2055,37)(2083,28)(2109,20)
	(2134,15)(2157,12)(2179,12)
	(2199,15)(2218,20)(2235,29)
	(2250,42)(2262,55)(2273,72)
	(2283,91)(2292,112)(2300,136)
	(2308,163)(2315,192)(2322,223)
	(2328,256)(2334,291)(2339,327)
	(2344,365)(2348,405)(2353,445)
	(2357,487)(2360,529)(2364,572)
	(2367,615)(2370,658)(2373,701)
	(2376,744)(2379,786)(2381,828)
	(2384,870)(2386,911)(2389,951)
	(2391,990)(2393,1028)(2395,1065)
	(2397,1102)(2398,1138)(2399,1173)
	(2400,1208)(2400,1242)(2400,1278)
	(2399,1315)(2398,1351)(2396,1389)
	(2394,1427)(2392,1465)(2390,1504)
	(2388,1544)(2385,1585)(2382,1626)
	(2380,1668)(2377,1711)(2374,1753)
	(2370,1796)(2367,1839)(2364,1882)
	(2360,1924)(2356,1966)(2352,2007)
	(2347,2046)(2342,2085)(2337,2122)
	(2331,2157)(2325,2190)(2318,2222)
	(2311,2250)(2303,2277)(2294,2301)
	(2284,2322)(2274,2340)(2262,2355)
	(2250,2367)(2236,2377)(2220,2383)
	(2203,2386)(2185,2386)(2165,2384)
	(2143,2380)(2121,2373)(2097,2365)
	(2072,2354)(2046,2342)(2019,2329)
	(1991,2315)(1963,2299)(1934,2282)
	(1905,2265)(1877,2247)(1848,2228)
	(1820,2209)(1793,2189)(1766,2169)
	(1741,2148)(1716,2126)(1693,2104)
	(1672,2081)(1651,2057)(1633,2032)
	(1616,2005)(1601,1978)(1587,1948)(1575,1917)
\path(5175,1917)(5183,1891)(5191,1863)
	(5198,1833)(5204,1802)(5209,1769)
	(5214,1734)(5218,1698)(5222,1660)
	(5225,1620)(5228,1579)(5230,1537)
	(5232,1493)(5234,1448)(5235,1402)
	(5236,1355)(5237,1308)(5238,1260)
	(5238,1212)(5238,1163)(5238,1115)
	(5238,1067)(5237,1019)(5236,972)
	(5235,925)(5234,880)(5232,835)
	(5230,791)(5228,749)(5225,708)
	(5222,669)(5218,632)(5214,596)
	(5209,561)(5204,529)(5198,498)
	(5191,470)(5183,442)(5175,417)
	(5162,385)(5147,355)(5131,328)
	(5112,302)(5092,278)(5069,255)
	(5045,233)(5020,212)(4993,192)
	(4965,173)(4936,154)(4906,136)
	(4875,119)(4845,103)(4814,87)
	(4783,73)(4753,60)(4724,47)
	(4695,37)(4667,28)(4641,20)
	(4616,15)(4593,12)(4571,12)
	(4551,15)(4532,20)(4515,29)
	(4500,42)(4488,55)(4477,72)
	(4467,91)(4458,112)(4450,136)
	(4442,163)(4435,192)(4428,223)
	(4422,256)(4416,291)(4411,327)
	(4406,365)(4402,405)(4397,445)
	(4393,487)(4390,529)(4386,572)
	(4383,615)(4380,658)(4377,701)
	(4374,744)(4371,786)(4369,828)
	(4366,870)(4364,911)(4361,951)
	(4359,990)(4357,1028)(4355,1065)
	(4353,1102)(4352,1138)(4351,1173)
	(4350,1208)(4350,1242)(4350,1278)
	(4351,1315)(4352,1351)(4354,1389)
	(4356,1427)(4358,1465)(4360,1504)
	(4362,1544)(4365,1585)(4368,1626)
	(4370,1668)(4373,1711)(4376,1753)
	(4380,1796)(4383,1839)(4386,1882)
	(4390,1924)(4394,1966)(4398,2007)
	(4403,2046)(4408,2085)(4413,2122)
	(4419,2157)(4425,2190)(4432,2222)
	(4439,2250)(4447,2277)(4456,2301)
	(4466,2322)(4476,2340)(4488,2355)
	(4500,2367)(4514,2377)(4530,2383)
	(4547,2386)(4565,2386)(4585,2384)
	(4607,2380)(4629,2373)(4653,2365)
	(4678,2354)(4704,2342)(4731,2329)
	(4759,2315)(4787,2299)(4816,2282)
	(4845,2265)(4873,2247)(4902,2228)
	(4930,2209)(4957,2189)(4984,2169)
	(5009,2148)(5034,2126)(5057,2104)
	(5078,2081)(5099,2057)(5117,2032)
	(5134,2005)(5149,1978)(5163,1948)(5175,1917)
\path(3450,2217)(3431,2202)(3413,2184)
	(3394,2164)(3376,2142)(3357,2118)
	(3339,2091)(3320,2063)(3302,2033)
	(3284,2000)(3265,1967)(3247,1932)
	(3228,1895)(3210,1858)(3191,1819)
	(3173,1780)(3154,1740)(3136,1700)
	(3118,1660)(3099,1619)(3081,1579)
	(3063,1539)(3045,1500)(3028,1461)
	(3010,1422)(2993,1385)(2976,1348)
	(2959,1313)(2943,1278)(2926,1245)
	(2910,1213)(2895,1181)(2879,1151)
	(2865,1121)(2850,1092)(2834,1060)
	(2818,1028)(2803,997)(2787,966)
	(2771,935)(2755,905)(2739,874)
	(2722,844)(2706,814)(2689,783)
	(2671,753)(2654,723)(2637,693)
	(2620,663)(2603,634)(2586,604)
	(2570,576)(2555,547)(2540,519)
	(2526,491)(2514,465)(2502,438)
	(2492,413)(2484,389)(2478,365)
	(2473,343)(2470,322)(2470,302)
	(2471,284)(2475,267)(2482,250)
	(2491,235)(2503,221)(2518,208)
	(2535,196)(2554,184)(2575,173)
	(2599,163)(2624,153)(2650,143)
	(2678,133)(2707,124)(2737,116)
	(2767,108)(2797,100)(2828,93)
	(2858,87)(2888,82)(2918,78)
	(2947,75)(2975,73)(3003,73)
	(3029,75)(3055,79)(3080,85)
	(3104,93)(3127,104)(3150,117)
	(3170,131)(3190,147)(3209,165)
	(3229,186)(3248,208)(3268,233)
	(3288,259)(3307,287)(3327,316)
	(3346,347)(3366,379)(3386,412)
	(3405,446)(3425,481)(3445,516)
	(3464,552)(3484,588)(3503,624)
	(3522,660)(3541,696)(3560,732)
	(3579,767)(3597,802)(3616,837)
	(3634,871)(3651,904)(3669,936)
	(3686,968)(3702,1000)(3719,1031)
	(3735,1061)(3750,1092)(3765,1122)
	(3780,1153)(3794,1184)(3809,1215)
	(3824,1247)(3839,1279)(3854,1311)
	(3869,1344)(3885,1378)(3901,1411)
	(3917,1445)(3933,1480)(3949,1514)
	(3965,1549)(3981,1584)(3997,1618)
	(4013,1653)(4028,1687)(4043,1720)
	(4056,1753)(4070,1786)(4082,1817)
	(4093,1848)(4103,1878)(4111,1906)
	(4118,1933)(4124,1959)(4128,1984)
	(4130,2007)(4130,2029)(4129,2049)
	(4125,2067)(4118,2086)(4109,2103)
	(4097,2119)(4082,2134)(4065,2147)
	(4046,2160)(4025,2172)(4001,2183)
	(3976,2194)(3950,2204)(3922,2214)
	(3893,2223)(3863,2232)(3833,2240)
	(3803,2248)(3772,2255)(3742,2261)
	(3712,2265)(3682,2269)(3653,2271)
	(3625,2272)(3597,2271)(3571,2268)
	(3545,2263)(3520,2256)(3496,2246)
	(3473,2233)(3450,2217)
\put(0,1242){\makebox(0,0)[lb]{\smash{{{\SetFigFont{6}{7.2}{\rmdefault}{\mddefault}{\updefault}$(\sigma_i,\sigma_i)=$}}}}}
\whiten\path(3225,2142)(3245,2126)(3265,2108)
	(3284,2089)(3304,2068)(3323,2046)
	(3343,2022)(3363,1997)(3382,1971)
	(3402,1943)(3421,1914)(3441,1885)
	(3461,1854)(3480,1823)(3500,1792)
	(3520,1760)(3539,1727)(3559,1695)
	(3578,1662)(3597,1629)(3616,1597)
	(3635,1565)(3654,1533)(3672,1501)
	(3691,1471)(3709,1440)(3726,1410)
	(3744,1381)(3761,1352)(3777,1324)
	(3794,1296)(3810,1269)(3825,1242)
	(3841,1213)(3857,1185)(3872,1156)
	(3888,1127)(3903,1099)(3919,1070)
	(3934,1041)(3950,1012)(3966,982)
	(3982,953)(3998,923)(4014,894)
	(4030,864)(4046,834)(4062,805)
	(4077,776)(4092,746)(4107,717)
	(4121,689)(4134,661)(4146,633)
	(4158,606)(4168,580)(4177,554)
	(4185,529)(4191,505)(4196,481)
	(4199,459)(4200,437)(4200,417)
	(4198,396)(4193,376)(4187,356)
	(4179,337)(4169,317)(4158,298)
	(4146,278)(4132,258)(4117,238)
	(4101,218)(4085,198)(4068,178)
	(4050,158)(4031,138)(4012,119)
	(3993,101)(3974,84)(3955,69)
	(3935,54)(3915,42)(3896,32)
	(3876,25)(3856,20)(3835,17)
	(3815,18)(3794,23)(3772,30)
	(3750,42)(3732,54)(3713,67)
	(3694,84)(3674,102)(3653,123)
	(3632,145)(3610,170)(3587,197)
	(3564,225)(3540,255)(3515,287)
	(3490,320)(3465,354)(3439,389)
	(3413,425)(3387,462)(3360,500)
	(3334,538)(3308,576)(3281,614)
	(3255,653)(3229,691)(3203,729)
	(3178,766)(3153,803)(3129,840)
	(3105,875)(3082,910)(3060,945)
	(3039,978)(3018,1011)(2998,1044)
	(2978,1075)(2960,1106)(2942,1137)
	(2925,1167)(2907,1201)(2889,1235)
	(2872,1269)(2855,1304)(2839,1338)
	(2823,1374)(2807,1410)(2791,1446)
	(2775,1483)(2759,1520)(2744,1558)
	(2728,1596)(2712,1634)(2697,1672)
	(2682,1710)(2667,1748)(2653,1786)
	(2639,1823)(2625,1860)(2612,1896)
	(2600,1931)(2589,1965)(2579,1998)
	(2570,2029)(2562,2059)(2555,2088)
	(2550,2114)(2547,2139)(2545,2162)
	(2545,2182)(2546,2201)(2550,2217)
	(2557,2233)(2566,2246)(2578,2257)
	(2593,2266)(2610,2273)(2629,2278)
	(2650,2281)(2674,2283)(2699,2283)
	(2725,2282)(2753,2281)(2782,2278)
	(2812,2275)(2842,2271)(2872,2266)
	(2903,2261)(2933,2255)(2963,2249)
	(2993,2242)(3022,2235)(3050,2227)
	(3078,2218)(3104,2208)(3130,2198)
	(3155,2186)(3179,2173)(3202,2158)(3225,2142)
\path(3225,2142)(3245,2126)(3265,2108)
	(3284,2089)(3304,2068)(3323,2046)
	(3343,2022)(3363,1997)(3382,1971)
	(3402,1943)(3421,1914)(3441,1885)
	(3461,1854)(3480,1823)(3500,1792)
	(3520,1760)(3539,1727)(3559,1695)
	(3578,1662)(3597,1629)(3616,1597)
	(3635,1565)(3654,1533)(3672,1501)
	(3691,1471)(3709,1440)(3726,1410)
	(3744,1381)(3761,1352)(3777,1324)
	(3794,1296)(3810,1269)(3825,1242)
	(3841,1213)(3857,1185)(3872,1156)
	(3888,1127)(3903,1099)(3919,1070)
	(3934,1041)(3950,1012)(3966,982)
	(3982,953)(3998,923)(4014,894)
	(4030,864)(4046,834)(4062,805)
	(4077,776)(4092,746)(4107,717)
	(4121,689)(4134,661)(4146,633)
	(4158,606)(4168,580)(4177,554)
	(4185,529)(4191,505)(4196,481)
	(4199,459)(4200,437)(4200,417)
	(4198,396)(4193,376)(4187,356)
	(4179,337)(4169,317)(4158,298)
	(4146,278)(4132,258)(4117,238)
	(4101,218)(4085,198)(4068,178)
	(4050,158)(4031,138)(4012,119)
	(3993,101)(3974,84)(3955,69)
	(3935,54)(3915,42)(3896,32)
	(3876,25)(3856,20)(3835,17)
	(3815,18)(3794,23)(3772,30)
	(3750,42)(3732,54)(3713,67)
	(3694,84)(3674,102)(3653,123)
	(3632,145)(3610,170)(3587,197)
	(3564,225)(3540,255)(3515,287)
	(3490,320)(3465,354)(3439,389)
	(3413,425)(3387,462)(3360,500)
	(3334,538)(3308,576)(3281,614)
	(3255,653)(3229,691)(3203,729)
	(3178,766)(3153,803)(3129,840)
	(3105,875)(3082,910)(3060,945)
	(3039,978)(3018,1011)(2998,1044)
	(2978,1075)(2960,1106)(2942,1137)
	(2925,1167)(2907,1201)(2889,1235)
	(2872,1269)(2855,1304)(2839,1338)
	(2823,1374)(2807,1410)(2791,1446)
	(2775,1483)(2759,1520)(2744,1558)
	(2728,1596)(2712,1634)(2697,1672)
	(2682,1710)(2667,1748)(2653,1786)
	(2639,1823)(2625,1860)(2612,1896)
	(2600,1931)(2589,1965)(2579,1998)
	(2570,2029)(2562,2059)(2555,2088)
	(2550,2114)(2547,2139)(2545,2162)
	(2545,2182)(2546,2201)(2550,2217)
	(2557,2233)(2566,2246)(2578,2257)
	(2593,2266)(2610,2273)(2629,2278)
	(2650,2281)(2674,2283)(2699,2283)
	(2725,2282)(2753,2281)(2782,2278)
	(2812,2275)(2842,2271)(2872,2266)
	(2903,2261)(2933,2255)(2963,2249)
	(2993,2242)(3022,2235)(3050,2227)
	(3078,2218)(3104,2208)(3130,2198)
	(3155,2186)(3179,2173)(3202,2158)(3225,2142)
\put(2850,1917){\blacken\ellipse{300}{300}}
\put(2850,1917){\ellipse{300}{300}}
\put(3750,417){\blacken\ellipse{300}{300}}
\put(3750,417){\ellipse{300}{300}}
\thicklines
\path(2850,1917)(3750,417)
\end{picture}
}

%% file: xfig/preid01.eepic
\setlength{\unitlength}{0.00041667in}
\begingroup\makeatletter\ifx\SetFigFont\undefined%
\gdef\SetFigFont#1#2#3#4#5{%
  \reset@font\fontsize{#1}{#2pt}%
  \fontfamily{#3}\fontseries{#4}\fontshape{#5}%
  \selectfont}%
\fi\endgroup%
{\renewcommand{\dashlinestretch}{30}
\begin{picture}(13541,2664)(0,-10)
\put(3300,2037){\blacken\ellipse{300}{300}}
\put(3300,2037){\ellipse{300}{300}}
\put(3300,537){\blacken\ellipse{300}{300}}
\put(3300,537){\ellipse{300}{300}}
\thicklines
\path(3300,2037)(3300,537)
\thinlines
\put(2400,2037){\blacken\ellipse{300}{300}}
\put(2400,2037){\ellipse{300}{300}}
\put(2400,537){\blacken\ellipse{300}{300}}
\put(2400,537){\ellipse{300}{300}}
\thicklines
\path(2400,2037)(2400,537)
\thinlines
\put(1500,2037){\blacken\ellipse{300}{300}}
\put(1500,2037){\ellipse{300}{300}}
\put(1500,537){\blacken\ellipse{300}{300}}
\put(1500,537){\ellipse{300}{300}}
\thicklines
\path(1500,2037)(1500,537)
\thinlines
\put(6000,2037){\blacken\ellipse{300}{300}}
\put(6000,2037){\ellipse{300}{300}}
\put(6000,537){\blacken\ellipse{300}{300}}
\put(6000,537){\ellipse{300}{300}}
\thicklines
\path(6000,2037)(6000,537)
\thinlines
\put(6900,2037){\blacken\ellipse{300}{300}}
\put(6900,2037){\ellipse{300}{300}}
\put(6900,537){\blacken\ellipse{300}{300}}
\put(6900,537){\ellipse{300}{300}}
\thicklines
\path(6900,2037)(6900,537)
\thinlines
\put(4200,2037){\blacken\ellipse{300}{300}}
\put(4200,2037){\ellipse{300}{300}}
\put(4200,537){\blacken\ellipse{300}{300}}
\put(4200,537){\ellipse{300}{300}}
\thicklines
\path(4200,2037)(4200,537)
\thinlines
\put(5100,2037){\blacken\ellipse{300}{300}}
\put(5100,2037){\ellipse{300}{300}}
\put(5100,537){\blacken\ellipse{300}{300}}
\put(5100,537){\ellipse{300}{300}}
\thicklines
\path(5100,2037)(5100,537)
\thinlines
\put(11400,2037){\blacken\ellipse{300}{300}}
\put(11400,2037){\ellipse{300}{300}}
\put(11400,537){\blacken\ellipse{300}{300}}
\put(11400,537){\ellipse{300}{300}}
\path(11100,312)(11111,299)(11124,286)
	(11141,274)(11160,263)(11182,252)
	(11207,241)(11235,231)(11265,220)
	(11296,210)(11329,200)(11362,190)
	(11396,181)(11431,172)(11464,164)
	(11497,156)(11529,150)(11559,145)
	(11587,141)(11612,139)(11635,139)
	(11656,141)(11673,146)(11688,152)
	(11700,162)(11709,173)(11716,187)
	(11722,204)(11726,223)(11728,245)
	(11729,270)(11729,297)(11727,326)
	(11725,357)(11722,390)(11718,423)
	(11713,458)(11708,492)(11703,527)
	(11697,562)(11692,596)(11686,629)
	(11680,660)(11674,690)(11667,718)
	(11661,744)(11655,767)(11648,788)
	(11641,807)(11633,823)(11625,837)
	(11613,852)(11599,865)(11584,874)
	(11567,881)(11549,887)(11529,891)
	(11508,893)(11486,894)(11464,894)
	(11441,893)(11418,892)(11396,889)
	(11374,887)(11352,883)(11332,879)
	(11313,873)(11295,867)(11279,859)
	(11264,849)(11250,837)(11238,824)
	(11227,809)(11217,791)(11206,771)
	(11195,748)(11185,724)(11173,698)
	(11162,670)(11151,641)(11140,612)
	(11129,582)(11119,551)(11110,521)
	(11101,492)(11094,464)(11088,437)
	(11084,411)(11082,388)(11083,366)
	(11086,346)(11091,328)(11100,312)
\path(11100,1812)(11111,1799)(11124,1786)
	(11141,1774)(11160,1763)(11182,1752)
	(11207,1741)(11235,1731)(11265,1720)
	(11296,1710)(11329,1700)(11362,1690)
	(11396,1681)(11431,1672)(11464,1664)
	(11497,1656)(11529,1650)(11559,1645)
	(11587,1641)(11612,1639)(11635,1639)
	(11656,1641)(11673,1646)(11688,1652)
	(11700,1662)(11709,1673)(11716,1687)
	(11722,1704)(11726,1723)(11728,1745)
	(11729,1770)(11729,1797)(11727,1826)
	(11725,1857)(11722,1890)(11718,1923)
	(11713,1958)(11708,1992)(11703,2027)
	(11697,2062)(11692,2096)(11686,2129)
	(11680,2160)(11674,2190)(11667,2218)
	(11661,2244)(11655,2267)(11648,2288)
	(11641,2307)(11633,2323)(11625,2337)
	(11613,2352)(11599,2365)(11584,2374)
	(11567,2381)(11549,2387)(11529,2391)
	(11508,2393)(11486,2394)(11464,2394)
	(11441,2393)(11418,2392)(11396,2389)
	(11374,2387)(11352,2383)(11332,2379)
	(11313,2373)(11295,2367)(11279,2359)
	(11264,2349)(11250,2337)(11238,2324)
	(11227,2309)(11217,2291)(11206,2271)
	(11195,2248)(11185,2224)(11173,2198)
	(11162,2170)(11151,2141)(11140,2112)
	(11129,2082)(11119,2051)(11110,2021)
	(11101,1992)(11094,1964)(11088,1937)
	(11084,1911)(11082,1888)(11083,1866)
	(11086,1846)(11091,1828)(11100,1812)
\put(12300,2037){\blacken\ellipse{300}{300}}
\put(12300,2037){\ellipse{300}{300}}
\put(12300,537){\blacken\ellipse{300}{300}}
\put(12300,537){\ellipse{300}{300}}
\path(12000,312)(12011,299)(12024,286)
	(12041,274)(12060,263)(12082,252)
	(12107,241)(12135,231)(12165,220)
	(12196,210)(12229,200)(12262,190)
	(12296,181)(12331,172)(12364,164)
	(12397,156)(12429,150)(12459,145)
	(12487,141)(12512,139)(12535,139)
	(12556,141)(12573,146)(12588,152)
	(12600,162)(12609,173)(12616,187)
	(12622,204)(12626,223)(12628,245)
	(12629,270)(12629,297)(12627,326)
	(12625,357)(12622,390)(12618,423)
	(12613,458)(12608,492)(12603,527)
	(12597,562)(12592,596)(12586,629)
	(12580,660)(12574,690)(12567,718)
	(12561,744)(12555,767)(12548,788)
	(12541,807)(12533,823)(12525,837)
	(12513,852)(12499,865)(12484,874)
	(12467,881)(12449,887)(12429,891)
	(12408,893)(12386,894)(12364,894)
	(12341,893)(12318,892)(12296,889)
	(12274,887)(12252,883)(12232,879)
	(12213,873)(12195,867)(12179,859)
	(12164,849)(12150,837)(12138,824)
	(12127,809)(12117,791)(12106,771)
	(12095,748)(12085,724)(12073,698)
	(12062,670)(12051,641)(12040,612)
	(12029,582)(12019,551)(12010,521)
	(12001,492)(11994,464)(11988,437)
	(11984,411)(11982,388)(11983,366)
	(11986,346)(11991,328)(12000,312)
\path(12000,1812)(12011,1799)(12024,1786)
	(12041,1774)(12060,1763)(12082,1752)
	(12107,1741)(12135,1731)(12165,1720)
	(12196,1710)(12229,1700)(12262,1690)
	(12296,1681)(12331,1672)(12364,1664)
	(12397,1656)(12429,1650)(12459,1645)
	(12487,1641)(12512,1639)(12535,1639)
	(12556,1641)(12573,1646)(12588,1652)
	(12600,1662)(12609,1673)(12616,1687)
	(12622,1704)(12626,1723)(12628,1745)
	(12629,1770)(12629,1797)(12627,1826)
	(12625,1857)(12622,1890)(12618,1923)
	(12613,1958)(12608,1992)(12603,2027)
	(12597,2062)(12592,2096)(12586,2129)
	(12580,2160)(12574,2190)(12567,2218)
	(12561,2244)(12555,2267)(12548,2288)
	(12541,2307)(12533,2323)(12525,2337)
	(12513,2352)(12499,2365)(12484,2374)
	(12467,2381)(12449,2387)(12429,2391)
	(12408,2393)(12386,2394)(12364,2394)
	(12341,2393)(12318,2392)(12296,2389)
	(12274,2387)(12252,2383)(12232,2379)
	(12213,2373)(12195,2367)(12179,2359)
	(12164,2349)(12150,2337)(12138,2324)
	(12127,2309)(12117,2291)(12106,2271)
	(12095,2248)(12085,2224)(12073,2198)
	(12062,2170)(12051,2141)(12040,2112)
	(12029,2082)(12019,2051)(12010,2021)
	(12001,1992)(11994,1964)(11988,1937)
	(11984,1911)(11982,1888)(11983,1866)
	(11986,1846)(11991,1828)(12000,1812)
\put(13200,2037){\blacken\ellipse{300}{300}}
\put(13200,2037){\ellipse{300}{300}}
\put(13200,537){\blacken\ellipse{300}{300}}
\put(13200,537){\ellipse{300}{300}}
\path(12900,312)(12911,299)(12924,286)
	(12941,274)(12960,263)(12982,252)
	(13007,241)(13035,231)(13065,220)
	(13096,210)(13129,200)(13162,190)
	(13196,181)(13231,172)(13264,164)
	(13297,156)(13329,150)(13359,145)
	(13387,141)(13412,139)(13435,139)
	(13456,141)(13473,146)(13488,152)
	(13500,162)(13509,173)(13516,187)
	(13522,204)(13526,223)(13528,245)
	(13529,270)(13529,297)(13527,326)
	(13525,357)(13522,390)(13518,423)
	(13513,458)(13508,492)(13503,527)
	(13497,562)(13492,596)(13486,629)
	(13480,660)(13474,690)(13467,718)
	(13461,744)(13455,767)(13448,788)
	(13441,807)(13433,823)(13425,837)
	(13413,852)(13399,865)(13384,874)
	(13367,881)(13349,887)(13329,891)
	(13308,893)(13286,894)(13264,894)
	(13241,893)(13218,892)(13196,889)
	(13174,887)(13152,883)(13132,879)
	(13113,873)(13095,867)(13079,859)
	(13064,849)(13050,837)(13038,824)
	(13027,809)(13017,791)(13006,771)
	(12995,748)(12985,724)(12973,698)
	(12962,670)(12951,641)(12940,612)
	(12929,582)(12919,551)(12910,521)
	(12901,492)(12894,464)(12888,437)
	(12884,411)(12882,388)(12883,366)
	(12886,346)(12891,328)(12900,312)
\path(12900,1812)(12911,1799)(12924,1786)
	(12941,1774)(12960,1763)(12982,1752)
	(13007,1741)(13035,1731)(13065,1720)
	(13096,1710)(13129,1700)(13162,1690)
	(13196,1681)(13231,1672)(13264,1664)
	(13297,1656)(13329,1650)(13359,1645)
	(13387,1641)(13412,1639)(13435,1639)
	(13456,1641)(13473,1646)(13488,1652)
	(13500,1662)(13509,1673)(13516,1687)
	(13522,1704)(13526,1723)(13528,1745)
	(13529,1770)(13529,1797)(13527,1826)
	(13525,1857)(13522,1890)(13518,1923)
	(13513,1958)(13508,1992)(13503,2027)
	(13497,2062)(13492,2096)(13486,2129)
	(13480,2160)(13474,2190)(13467,2218)
	(13461,2244)(13455,2267)(13448,2288)
	(13441,2307)(13433,2323)(13425,2337)
	(13413,2352)(13399,2365)(13384,2374)
	(13367,2381)(13349,2387)(13329,2391)
	(13308,2393)(13286,2394)(13264,2394)
	(13241,2393)(13218,2392)(13196,2389)
	(13174,2387)(13152,2383)(13132,2379)
	(13113,2373)(13095,2367)(13079,2359)
	(13064,2349)(13050,2337)(13038,2324)
	(13027,2309)(13017,2291)(13006,2271)
	(12995,2248)(12985,2224)(12973,2198)
	(12962,2170)(12951,2141)(12940,2112)
	(12929,2082)(12919,2051)(12910,2021)
	(12901,1992)(12894,1964)(12888,1937)
	(12884,1911)(12882,1888)(12883,1866)
	(12886,1846)(12891,1828)(12900,1812)
\put(7800,2037){\blacken\ellipse{300}{300}}
\put(7800,2037){\ellipse{300}{300}}
\put(7800,537){\blacken\ellipse{300}{300}}
\put(7800,537){\ellipse{300}{300}}
\thicklines
\path(7800,2037)(7800,537)
\thinlines
\path(8175,2262)(8183,2241)(8190,2219)
	(8196,2195)(8201,2169)(8205,2142)
	(8207,2113)(8209,2083)(8210,2052)
	(8210,2020)(8210,1986)(8208,1952)
	(8207,1917)(8205,1882)(8202,1846)
	(8200,1810)(8197,1773)(8194,1737)
	(8191,1700)(8188,1664)(8186,1627)
	(8183,1592)(8181,1556)(8180,1521)
	(8178,1486)(8177,1452)(8176,1418)
	(8175,1385)(8175,1352)(8175,1320)
	(8175,1287)(8175,1254)(8175,1222)
	(8175,1189)(8176,1155)(8177,1122)
	(8178,1087)(8180,1052)(8181,1017)
	(8183,981)(8186,945)(8188,908)
	(8191,872)(8194,835)(8197,798)
	(8200,761)(8202,724)(8205,688)
	(8207,652)(8208,617)(8210,583)
	(8210,549)(8210,517)(8209,486)
	(8207,456)(8205,428)(8201,401)
	(8196,376)(8190,353)(8183,332)
	(8175,312)(8164,292)(8152,273)
	(8137,256)(8121,241)(8104,226)
	(8085,213)(8064,200)(8042,188)
	(8019,176)(7995,165)(7970,154)
	(7945,144)(7920,135)(7894,127)
	(7869,119)(7843,113)(7819,108)
	(7795,104)(7773,102)(7751,102)
	(7731,105)(7712,110)(7695,118)
	(7678,129)(7664,144)(7650,162)
	(7640,179)(7631,198)(7623,219)
	(7615,242)(7608,268)(7601,296)
	(7594,326)(7588,357)(7582,391)
	(7576,426)(7571,462)(7566,500)
	(7561,539)(7556,579)(7552,619)
	(7548,661)(7543,702)(7539,744)
	(7536,786)(7532,828)(7528,870)
	(7525,912)(7521,953)(7518,993)
	(7515,1033)(7513,1073)(7510,1111)
	(7508,1149)(7505,1186)(7504,1222)
	(7502,1258)(7501,1293)(7500,1328)
	(7500,1362)(7500,1398)(7501,1435)
	(7502,1471)(7504,1508)(7506,1546)
	(7508,1584)(7511,1623)(7514,1662)
	(7517,1702)(7520,1742)(7524,1783)
	(7527,1825)(7531,1866)(7535,1908)
	(7539,1950)(7543,1991)(7548,2032)
	(7552,2073)(7557,2113)(7562,2152)
	(7567,2190)(7573,2227)(7579,2262)
	(7585,2296)(7591,2327)(7598,2357)
	(7605,2385)(7613,2410)(7621,2433)
	(7630,2454)(7640,2472)(7650,2487)
	(7665,2503)(7681,2516)(7699,2524)
	(7718,2528)(7739,2529)(7762,2527)
	(7786,2523)(7811,2516)(7837,2508)
	(7864,2497)(7892,2486)(7920,2473)
	(7947,2459)(7974,2445)(8001,2430)
	(8027,2414)(8051,2398)(8075,2382)
	(8096,2364)(8116,2346)(8134,2328)
	(8149,2308)(8163,2286)(8175,2262)
\put(8700,2037){\blacken\ellipse{300}{300}}
\put(8700,2037){\ellipse{300}{300}}
\put(8700,537){\blacken\ellipse{300}{300}}
\put(8700,537){\ellipse{300}{300}}
\thicklines
\path(8700,2037)(8700,537)
\thinlines
\path(9075,2262)(9083,2241)(9090,2219)
	(9096,2195)(9101,2169)(9105,2142)
	(9107,2113)(9109,2083)(9110,2052)
	(9110,2020)(9110,1986)(9108,1952)
	(9107,1917)(9105,1882)(9102,1846)
	(9100,1810)(9097,1773)(9094,1737)
	(9091,1700)(9088,1664)(9086,1627)
	(9083,1592)(9081,1556)(9080,1521)
	(9078,1486)(9077,1452)(9076,1418)
	(9075,1385)(9075,1352)(9075,1320)
	(9075,1287)(9075,1254)(9075,1222)
	(9075,1189)(9076,1155)(9077,1122)
	(9078,1087)(9080,1052)(9081,1017)
	(9083,981)(9086,945)(9088,908)
	(9091,872)(9094,835)(9097,798)
	(9100,761)(9102,724)(9105,688)
	(9107,652)(9108,617)(9110,583)
	(9110,549)(9110,517)(9109,486)
	(9107,456)(9105,428)(9101,401)
	(9096,376)(9090,353)(9083,332)
	(9075,312)(9064,292)(9052,273)
	(9037,256)(9021,241)(9004,226)
	(8985,213)(8964,200)(8942,188)
	(8919,176)(8895,165)(8870,154)
	(8845,144)(8820,135)(8794,127)
	(8769,119)(8743,113)(8719,108)
	(8695,104)(8673,102)(8651,102)
	(8631,105)(8612,110)(8595,118)
	(8578,129)(8564,144)(8550,162)
	(8540,179)(8531,198)(8523,219)
	(8515,242)(8508,268)(8501,296)
	(8494,326)(8488,357)(8482,391)
	(8476,426)(8471,462)(8466,500)
	(8461,539)(8456,579)(8452,619)
	(8448,661)(8443,702)(8439,744)
	(8436,786)(8432,828)(8428,870)
	(8425,912)(8421,953)(8418,993)
	(8415,1033)(8413,1073)(8410,1111)
	(8408,1149)(8405,1186)(8404,1222)
	(8402,1258)(8401,1293)(8400,1328)
	(8400,1362)(8400,1398)(8401,1435)
	(8402,1471)(8404,1508)(8406,1546)
	(8408,1584)(8411,1623)(8414,1662)
	(8417,1702)(8420,1742)(8424,1783)
	(8427,1825)(8431,1866)(8435,1908)
	(8439,1950)(8443,1991)(8448,2032)
	(8452,2073)(8457,2113)(8462,2152)
	(8467,2190)(8473,2227)(8479,2262)
	(8485,2296)(8491,2327)(8498,2357)
	(8505,2385)(8513,2410)(8521,2433)
	(8530,2454)(8540,2472)(8550,2487)
	(8565,2503)(8581,2516)(8599,2524)
	(8618,2528)(8639,2529)(8662,2527)
	(8686,2523)(8711,2516)(8737,2508)
	(8764,2497)(8792,2486)(8820,2473)
	(8847,2459)(8874,2445)(8901,2430)
	(8927,2414)(8951,2398)(8975,2382)
	(8996,2364)(9016,2346)(9034,2328)
	(9049,2308)(9063,2286)(9075,2262)
\put(9600,537){\blacken\ellipse{300}{300}}
\put(9600,537){\ellipse{300}{300}}
\put(9600,2037){\blacken\ellipse{300}{300}}
\put(9600,2037){\ellipse{300}{300}}
\put(10500,2037){\blacken\ellipse{300}{300}}
\put(10500,2037){\ellipse{300}{300}}
\put(10500,537){\blacken\ellipse{300}{300}}
\put(10500,537){\ellipse{300}{300}}
\thicklines
\path(9600,2037)(9600,537)
\thinlines
\path(6450,2637)(6477,2636)(6504,2634)
	(6531,2631)(6559,2627)(6587,2622)
	(6616,2616)(6645,2609)(6675,2602)
	(6705,2595)(6736,2587)(6767,2578)
	(6798,2569)(6829,2560)(6860,2550)
	(6891,2540)(6921,2530)(6951,2519)
	(6980,2507)(7009,2495)(7036,2481)
	(7062,2468)(7087,2453)(7110,2437)
	(7131,2420)(7151,2401)(7169,2382)
	(7186,2360)(7200,2337)(7211,2315)
	(7221,2292)(7229,2267)(7235,2241)
	(7241,2214)(7245,2185)(7248,2155)
	(7249,2124)(7250,2092)(7250,2059)
	(7249,2025)(7247,1990)(7245,1955)
	(7242,1920)(7239,1884)(7235,1847)
	(7231,1811)(7228,1774)(7224,1737)
	(7220,1701)(7217,1664)(7213,1628)
	(7210,1593)(7208,1557)(7205,1522)
	(7204,1488)(7202,1453)(7201,1419)
	(7201,1386)(7200,1353)(7200,1320)
	(7200,1287)(7200,1254)(7200,1221)
	(7201,1188)(7201,1154)(7202,1120)
	(7204,1086)(7205,1051)(7208,1016)
	(7210,980)(7213,944)(7217,908)
	(7220,871)(7224,834)(7228,797)
	(7231,760)(7235,723)(7239,687)
	(7242,650)(7245,614)(7247,579)
	(7249,544)(7250,510)(7250,477)
	(7249,445)(7248,414)(7245,385)
	(7241,356)(7235,329)(7229,304)
	(7221,280)(7211,258)(7200,237)
	(7184,214)(7167,193)(7147,174)
	(7125,157)(7101,142)(7075,128)
	(7048,116)(7019,105)(6989,95)
	(6958,87)(6926,79)(6893,71)
	(6860,65)(6827,59)(6793,53)
	(6760,47)(6726,42)(6694,37)
	(6661,32)(6629,28)(6598,24)
	(6568,20)(6537,17)(6508,14)
	(6479,13)(6450,12)(6421,13)
	(6392,14)(6362,17)(6332,20)
	(6301,24)(6270,28)(6238,32)
	(6205,37)(6172,42)(6138,47)
	(6104,53)(6070,59)(6036,65)
	(6003,71)(5970,79)(5937,87)
	(5906,95)(5876,105)(5847,116)
	(5820,128)(5794,142)(5771,157)
	(5750,174)(5731,193)(5714,214)
	(5700,237)(5690,258)(5682,280)
	(5676,304)(5671,329)(5667,356)
	(5666,385)(5665,414)(5666,445)
	(5668,477)(5671,510)(5675,544)
	(5680,579)(5686,614)(5692,650)
	(5699,687)(5706,723)(5713,760)
	(5720,797)(5727,834)(5734,871)
	(5741,908)(5747,944)(5753,980)
	(5758,1016)(5762,1051)(5766,1086)
	(5769,1120)(5771,1154)(5773,1188)
	(5774,1221)(5775,1254)(5775,1287)
	(5775,1320)(5774,1353)(5773,1386)
	(5771,1419)(5769,1453)(5766,1488)
	(5762,1522)(5758,1557)(5753,1593)
	(5747,1628)(5741,1664)(5734,1701)
	(5727,1737)(5720,1774)(5713,1811)
	(5706,1847)(5699,1884)(5692,1920)
	(5686,1955)(5680,1990)(5675,2025)
	(5671,2059)(5668,2092)(5666,2124)
	(5665,2155)(5666,2185)(5667,2214)
	(5671,2241)(5676,2267)(5682,2292)
	(5690,2315)(5700,2337)(5713,2360)
	(5728,2382)(5746,2401)(5765,2420)
	(5786,2437)(5809,2453)(5833,2468)
	(5859,2481)(5886,2495)(5915,2507)
	(5944,2519)(5974,2530)(6005,2540)
	(6036,2550)(6068,2560)(6099,2569)
	(6131,2578)(6162,2587)(6193,2595)
	(6224,2602)(6254,2609)(6284,2616)
	(6313,2622)(6341,2627)(6369,2631)
	(6396,2634)(6423,2636)(6450,2637)
\path(4650,2637)(4678,2636)(4707,2634)
	(4736,2631)(4766,2627)(4797,2622)
	(4828,2616)(4861,2609)(4894,2602)
	(4927,2595)(4961,2587)(4996,2578)
	(5031,2569)(5066,2560)(5101,2550)
	(5136,2540)(5170,2530)(5204,2519)
	(5236,2507)(5268,2495)(5298,2481)
	(5327,2468)(5354,2453)(5380,2437)
	(5403,2420)(5425,2401)(5444,2382)
	(5460,2360)(5475,2337)(5486,2315)
	(5495,2292)(5502,2267)(5508,2241)
	(5512,2214)(5514,2185)(5515,2155)
	(5515,2124)(5513,2092)(5510,2059)
	(5506,2025)(5502,1990)(5496,1955)
	(5490,1920)(5483,1884)(5476,1847)
	(5469,1811)(5462,1774)(5455,1737)
	(5448,1701)(5441,1664)(5434,1628)
	(5428,1593)(5423,1557)(5418,1522)
	(5414,1488)(5410,1453)(5407,1419)
	(5405,1386)(5403,1353)(5401,1320)
	(5400,1287)(5399,1254)(5398,1221)
	(5398,1188)(5398,1154)(5398,1120)
	(5399,1086)(5401,1051)(5403,1016)
	(5405,980)(5408,944)(5412,908)
	(5415,871)(5419,834)(5423,797)
	(5428,760)(5432,723)(5436,687)
	(5439,650)(5442,614)(5445,579)
	(5447,544)(5449,510)(5449,477)
	(5449,445)(5447,414)(5444,385)
	(5440,356)(5435,329)(5429,304)
	(5421,280)(5411,258)(5400,237)
	(5384,214)(5367,193)(5347,174)
	(5325,157)(5301,142)(5275,128)
	(5248,116)(5219,105)(5189,95)
	(5158,87)(5126,79)(5093,71)
	(5060,65)(5027,59)(4993,53)
	(4960,47)(4926,42)(4894,37)
	(4861,32)(4829,28)(4798,24)
	(4768,20)(4737,17)(4708,14)
	(4679,13)(4650,12)(4621,13)
	(4592,14)(4563,17)(4532,20)
	(4502,24)(4471,28)(4439,32)
	(4406,37)(4374,42)(4340,47)
	(4307,53)(4273,59)(4240,65)
	(4207,71)(4174,79)(4142,87)
	(4111,95)(4081,105)(4052,116)
	(4025,128)(3999,142)(3975,157)
	(3953,174)(3933,193)(3916,214)
	(3900,237)(3889,258)(3879,280)
	(3871,304)(3865,329)(3859,356)
	(3855,385)(3852,414)(3851,445)
	(3850,477)(3850,510)(3851,544)
	(3853,579)(3855,614)(3858,650)
	(3861,687)(3865,723)(3869,760)
	(3872,797)(3876,834)(3880,871)
	(3883,908)(3887,944)(3890,980)
	(3892,1016)(3895,1051)(3896,1086)
	(3898,1120)(3899,1154)(3899,1188)
	(3900,1221)(3900,1254)(3900,1287)
	(3900,1320)(3900,1353)(3899,1386)
	(3899,1419)(3898,1453)(3896,1488)
	(3895,1522)(3892,1557)(3890,1593)
	(3887,1628)(3883,1664)(3880,1701)
	(3876,1737)(3872,1774)(3869,1811)
	(3865,1847)(3861,1884)(3858,1920)
	(3855,1955)(3853,1990)(3851,2025)
	(3850,2059)(3850,2092)(3851,2124)
	(3852,2155)(3855,2185)(3859,2214)
	(3865,2241)(3871,2267)(3879,2292)
	(3889,2315)(3900,2337)(3914,2360)
	(3931,2382)(3949,2401)(3968,2420)
	(3990,2437)(4013,2453)(4037,2468)
	(4063,2481)(4090,2495)(4118,2507)
	(4146,2519)(4176,2530)(4206,2540)
	(4236,2550)(4267,2560)(4298,2569)
	(4329,2578)(4359,2587)(4390,2595)
	(4420,2602)(4450,2609)(4479,2616)
	(4509,2622)(4537,2627)(4566,2631)
	(4594,2634)(4622,2636)(4650,2637)
\path(10050,2037)(10061,2013)(10070,1989)
	(10077,1964)(10082,1938)(10086,1911)
	(10089,1883)(10090,1855)(10090,1826)
	(10089,1797)(10087,1767)(10085,1737)
	(10082,1707)(10078,1676)(10074,1645)
	(10071,1615)(10067,1584)(10063,1553)
	(10060,1523)(10057,1493)(10055,1463)
	(10053,1433)(10052,1403)(10051,1374)
	(10050,1345)(10050,1316)(10050,1287)
	(10050,1258)(10050,1229)(10051,1200)
	(10052,1171)(10053,1141)(10055,1111)
	(10057,1080)(10060,1050)(10063,1019)
	(10067,988)(10071,957)(10074,926)
	(10078,894)(10082,863)(10085,833)
	(10087,802)(10089,772)(10090,743)
	(10090,714)(10089,686)(10086,659)
	(10082,632)(10077,607)(10070,582)
	(10061,559)(10050,537)(10039,517)
	(10025,498)(10011,479)(9994,460)
	(9976,441)(9956,422)(9935,402)
	(9912,382)(9888,362)(9864,341)
	(9838,320)(9811,299)(9784,278)
	(9757,258)(9730,238)(9703,219)
	(9676,202)(9650,186)(9624,171)
	(9600,158)(9576,148)(9554,140)
	(9533,136)(9514,134)(9496,135)
	(9479,140)(9464,149)(9450,162)
	(9439,175)(9429,192)(9420,211)
	(9412,232)(9404,256)(9396,283)
	(9389,312)(9383,343)(9377,376)
	(9371,411)(9366,447)(9361,485)
	(9356,525)(9352,565)(9348,607)
	(9344,649)(9340,692)(9336,735)
	(9333,778)(9329,821)(9326,864)
	(9323,906)(9320,948)(9317,990)
	(9314,1031)(9312,1071)(9309,1110)
	(9307,1148)(9305,1185)(9303,1222)
	(9302,1258)(9301,1293)(9300,1328)
	(9300,1362)(9300,1398)(9301,1435)
	(9302,1471)(9304,1509)(9306,1547)
	(9308,1585)(9310,1624)(9313,1664)
	(9316,1705)(9319,1746)(9322,1788)
	(9325,1831)(9329,1873)(9332,1916)
	(9336,1959)(9340,2002)(9344,2044)
	(9348,2086)(9353,2127)(9357,2166)
	(9362,2205)(9368,2242)(9374,2277)
	(9380,2310)(9386,2342)(9394,2370)
	(9401,2397)(9410,2421)(9419,2442)
	(9428,2460)(9439,2475)(9450,2487)
	(9465,2498)(9482,2504)(9500,2506)
	(9520,2503)(9541,2497)(9564,2487)
	(9589,2475)(9615,2459)(9642,2441)
	(9670,2420)(9699,2398)(9728,2374)
	(9757,2350)(9786,2325)(9815,2299)
	(9844,2273)(9871,2248)(9898,2222)
	(9923,2198)(9947,2173)(9969,2150)
	(9989,2127)(10007,2104)(10023,2082)
	(10038,2060)(10050,2037)
\path(10950,1962)(10961,1938)(10970,1914)
	(10977,1889)(10982,1863)(10986,1836)
	(10989,1808)(10990,1780)(10990,1751)
	(10989,1722)(10987,1692)(10985,1662)
	(10982,1632)(10978,1601)(10974,1570)
	(10971,1540)(10967,1509)(10963,1478)
	(10960,1448)(10957,1418)(10955,1388)
	(10953,1358)(10952,1328)(10951,1299)
	(10950,1270)(10950,1241)(10950,1212)
	(10950,1183)(10950,1154)(10951,1125)
	(10952,1096)(10953,1066)(10955,1036)
	(10957,1005)(10960,975)(10963,944)
	(10967,913)(10971,882)(10974,851)
	(10978,819)(10982,788)(10985,758)
	(10987,727)(10989,697)(10990,668)
	(10990,639)(10989,611)(10986,584)
	(10982,557)(10977,532)(10970,507)
	(10961,484)(10950,462)(10939,442)
	(10925,423)(10911,404)(10894,385)
	(10876,366)(10856,347)(10835,327)
	(10812,307)(10788,287)(10764,266)
	(10738,245)(10711,224)(10684,203)
	(10657,183)(10630,163)(10603,144)
	(10576,127)(10550,111)(10524,96)
	(10500,83)(10476,73)(10454,65)
	(10433,61)(10414,59)(10396,60)
	(10379,65)(10364,74)(10350,87)
	(10339,100)(10329,117)(10320,136)
	(10312,157)(10304,181)(10296,208)
	(10289,237)(10283,268)(10277,301)
	(10271,336)(10266,372)(10261,410)
	(10256,450)(10252,490)(10248,532)
	(10244,574)(10240,617)(10236,660)
	(10233,703)(10229,746)(10226,789)
	(10223,831)(10220,873)(10217,915)
	(10214,956)(10212,996)(10209,1035)
	(10207,1073)(10205,1110)(10203,1147)
	(10202,1183)(10201,1218)(10200,1253)
	(10200,1287)(10200,1323)(10201,1360)
	(10202,1396)(10204,1434)(10206,1472)
	(10208,1510)(10210,1549)(10213,1589)
	(10216,1630)(10219,1671)(10222,1713)
	(10225,1756)(10229,1798)(10232,1841)
	(10236,1884)(10240,1927)(10244,1969)
	(10248,2011)(10253,2052)(10257,2091)
	(10262,2130)(10268,2167)(10274,2202)
	(10280,2235)(10286,2267)(10294,2295)
	(10301,2322)(10310,2346)(10319,2367)
	(10328,2385)(10339,2400)(10350,2412)
	(10365,2423)(10382,2429)(10400,2431)
	(10420,2428)(10441,2422)(10464,2412)
	(10489,2400)(10515,2384)(10542,2366)
	(10570,2345)(10599,2323)(10628,2299)
	(10657,2275)(10686,2250)(10715,2224)
	(10744,2198)(10771,2173)(10798,2147)
	(10823,2123)(10847,2098)(10869,2075)
	(10889,2052)(10907,2029)(10923,2007)
	(10938,1985)(10950,1962)
\path(2400,2562)(2365,2562)(2329,2561)
	(2293,2559)(2257,2557)(2221,2555)
	(2183,2552)(2145,2549)(2107,2546)
	(2068,2542)(2029,2539)(1989,2535)
	(1948,2531)(1908,2527)(1867,2523)
	(1826,2519)(1785,2514)(1744,2509)
	(1703,2504)(1663,2499)(1623,2493)
	(1584,2487)(1546,2480)(1509,2473)
	(1473,2465)(1438,2457)(1404,2448)
	(1372,2437)(1342,2426)(1314,2414)
	(1287,2401)(1262,2387)(1239,2372)
	(1219,2355)(1200,2337)(1182,2316)
	(1167,2294)(1154,2270)(1143,2245)
	(1135,2218)(1128,2189)(1124,2160)
	(1121,2129)(1120,2097)(1120,2064)
	(1122,2030)(1124,1995)(1128,1960)
	(1133,1924)(1138,1887)(1144,1851)
	(1150,1814)(1156,1777)(1162,1740)
	(1168,1703)(1174,1666)(1179,1630)
	(1184,1594)(1188,1558)(1191,1523)
	(1194,1488)(1196,1454)(1198,1420)
	(1199,1386)(1200,1353)(1200,1320)
	(1200,1287)(1200,1254)(1200,1221)
	(1199,1188)(1198,1154)(1196,1120)
	(1194,1086)(1191,1051)(1188,1016)
	(1184,980)(1179,944)(1174,908)
	(1168,871)(1162,834)(1156,797)
	(1150,760)(1144,723)(1138,687)
	(1133,650)(1128,614)(1124,579)
	(1122,544)(1120,510)(1120,477)
	(1121,445)(1124,414)(1128,385)
	(1135,356)(1143,329)(1154,304)
	(1167,280)(1182,258)(1200,237)
	(1219,219)(1239,202)(1262,187)
	(1287,173)(1314,160)(1342,148)
	(1372,137)(1404,126)(1438,117)
	(1473,109)(1509,101)(1546,94)
	(1584,87)(1623,81)(1663,75)
	(1703,70)(1744,65)(1785,60)
	(1826,55)(1867,51)(1908,47)
	(1948,43)(1989,39)(2029,35)
	(2068,32)(2107,28)(2145,25)
	(2183,22)(2221,19)(2257,17)
	(2293,15)(2329,13)(2365,12)
	(2400,12)(2435,12)(2471,13)
	(2507,15)(2543,17)(2580,19)
	(2617,22)(2655,25)(2694,28)
	(2733,32)(2772,35)(2813,39)
	(2853,43)(2894,47)(2936,51)
	(2977,55)(3018,60)(3060,65)
	(3101,70)(3141,75)(3181,81)
	(3221,87)(3259,94)(3296,101)
	(3332,109)(3367,117)(3401,126)
	(3432,137)(3462,148)(3490,160)
	(3516,173)(3541,187)(3563,202)
	(3582,219)(3600,237)(3616,258)
	(3630,280)(3642,304)(3650,329)
	(3657,356)(3661,385)(3664,414)
	(3664,445)(3662,477)(3659,510)
	(3654,544)(3648,579)(3641,614)
	(3633,650)(3624,687)(3615,723)
	(3606,760)(3596,797)(3587,834)
	(3578,871)(3569,908)(3561,944)
	(3554,980)(3547,1016)(3541,1051)
	(3536,1086)(3532,1120)(3529,1154)
	(3527,1188)(3526,1221)(3525,1254)
	(3525,1287)(3525,1320)(3526,1353)
	(3527,1386)(3529,1420)(3532,1454)
	(3536,1488)(3541,1523)(3547,1558)
	(3554,1594)(3561,1630)(3569,1666)
	(3578,1703)(3587,1740)(3596,1777)
	(3606,1814)(3615,1851)(3624,1887)
	(3633,1924)(3641,1960)(3648,1995)
	(3654,2030)(3659,2064)(3662,2097)
	(3664,2129)(3664,2160)(3661,2189)
	(3657,2218)(3650,2245)(3642,2270)
	(3630,2294)(3616,2316)(3600,2337)
	(3582,2355)(3563,2372)(3541,2387)
	(3516,2401)(3490,2414)(3462,2426)
	(3432,2437)(3401,2448)(3367,2457)
	(3332,2465)(3296,2473)(3259,2480)
	(3221,2487)(3181,2493)(3141,2499)
	(3101,2504)(3060,2509)(3018,2514)
	(2977,2519)(2936,2523)(2894,2527)
	(2853,2531)(2813,2535)(2772,2539)
	(2733,2542)(2694,2546)(2655,2549)
	(2617,2552)(2580,2555)(2543,2557)
	(2507,2559)(2471,2561)(2435,2562)(2400,2562)
\put(0,1137){\makebox(0,0)[lb]{\smash{{{\SetFigFont{6}{7.2}{\rmdefault}{\mddefault}{\updefault}$\localthing =$}}}}}
\end{picture}
}

%% file: xfig/poset1.eepic
\setlength{\unitlength}{0.00025000in}
\begingroup\makeatletter\ifx\SetFigFont\undefined%
\gdef\SetFigFont#1#2#3#4#5{%
  \reset@font\fontsize{#1}{#2pt}%
  \fontfamily{#3}\fontseries{#4}\fontshape{#5}%
  \selectfont}%
\fi\endgroup%
{\renewcommand{\dashlinestretch}{30}
\begin{picture}(5724,6639)(0,-10)
\path(2112,612)(3612,612)(3612,12)
	(2112,12)(2112,612)
\path(2112,2112)(3612,2112)(3612,1512)
	(2112,1512)(2112,2112)
\path(2862,1512)(2862,612)
\path(912,3612)(2412,3612)(2412,3012)
	(912,3012)(912,3612)
\path(3312,3612)(4812,3612)(4812,3012)
	(3312,3012)(3312,3612)
\path(3312,5112)(4812,5112)(4812,4512)
	(3312,4512)(3312,5112)
\path(912,5112)(2412,5112)(2412,4512)
	(912,4512)(912,5112)
\path(2112,6612)(3612,6612)(3612,6012)
	(2112,6012)(2112,6612)
\path(2862,6012)(1662,5112)
\path(2862,6012)(4062,5112)
\path(1662,4512)(1662,3612)
\path(4062,4512)(4062,3612)
\path(1662,3612)(4062,4512)
\path(1662,3012)(2862,2112)
\path(2862,2112)(4062,3012)
\dashline{60.000}(12,4062)(14,4062)(19,4062)
	(28,4062)(42,4062)(61,4062)
	(86,4062)(116,4062)(151,4062)
	(191,4062)(234,4062)(279,4062)
	(326,4062)(374,4062)(421,4062)
	(468,4062)(514,4062)(559,4062)
	(602,4062)(644,4062)(684,4062)
	(722,4062)(760,4062)(797,4062)
	(833,4062)(869,4062)(904,4062)
	(940,4062)(976,4062)(1012,4062)
	(1043,4062)(1074,4062)(1106,4062)
	(1138,4062)(1172,4061)(1206,4061)
	(1240,4060)(1276,4059)(1312,4058)
	(1349,4056)(1386,4054)(1424,4052)
	(1463,4049)(1502,4046)(1542,4042)
	(1582,4038)(1622,4033)(1662,4027)
	(1702,4021)(1742,4015)(1782,4007)
	(1822,3999)(1861,3991)(1900,3981)
	(1938,3971)(1975,3960)(2012,3949)
	(2048,3937)(2084,3924)(2118,3910)
	(2152,3896)(2186,3881)(2218,3865)
	(2250,3848)(2281,3831)(2312,3812)
	(2342,3792)(2372,3772)(2402,3750)
	(2432,3727)(2462,3703)(2491,3678)
	(2521,3652)(2551,3625)(2581,3597)
	(2610,3568)(2640,3538)(2670,3507)
	(2700,3476)(2730,3444)(2760,3411)
	(2789,3378)(2819,3345)(2848,3312)
	(2877,3279)(2906,3246)(2934,3213)
	(2963,3180)(2990,3148)(3018,3117)
	(3044,3086)(3071,3056)(3097,3027)
	(3122,2999)(3147,2972)(3171,2946)
	(3196,2921)(3219,2897)(3243,2874)
	(3266,2852)(3289,2832)(3312,2812)
	(3338,2791)(3364,2771)(3390,2753)
	(3416,2735)(3443,2719)(3471,2703)
	(3499,2689)(3527,2675)(3557,2662)
	(3586,2650)(3617,2639)(3648,2629)
	(3679,2620)(3711,2611)(3743,2604)
	(3776,2597)(3809,2591)(3842,2585)
	(3876,2580)(3909,2576)(3943,2573)
	(3976,2570)(4010,2568)(4043,2566)
	(4077,2565)(4110,2564)(4144,2563)
	(4177,2563)(4211,2562)(4244,2562)
	(4278,2562)(4312,2562)(4343,2562)
	(4374,2562)(4406,2562)(4439,2562)
	(4473,2562)(4508,2562)(4544,2562)
	(4582,2562)(4622,2562)(4664,2562)
	(4708,2562)(4753,2562)(4802,2562)
	(4852,2562)(4905,2562)(4960,2562)
	(5017,2562)(5076,2562)(5135,2562)
	(5196,2562)(5256,2562)(5316,2562)
	(5374,2562)(5430,2562)(5482,2562)
	(5529,2562)(5572,2562)(5609,2562)
	(5640,2562)(5665,2562)(5684,2562)
	(5697,2562)(5705,2562)(5710,2562)(5712,2562)
\dashline{60.000}(12,1212)(15,1212)(22,1212)
	(34,1212)(53,1212)(79,1212)
	(111,1212)(148,1212)(191,1212)
	(236,1212)(284,1212)(331,1212)
	(378,1212)(424,1212)(466,1212)
	(507,1212)(544,1212)(578,1212)
	(609,1212)(637,1212)(663,1212)
	(687,1212)(708,1212)(728,1212)
	(746,1212)(762,1212)(784,1212)
	(804,1212)(822,1212)(839,1212)
	(856,1212)(872,1212)(887,1212)
	(903,1212)(919,1212)(935,1212)
	(953,1212)(971,1212)(990,1212)
	(1010,1212)(1033,1212)(1057,1212)
	(1083,1212)(1112,1212)(1130,1212)
	(1149,1212)(1169,1212)(1190,1212)
	(1212,1212)(1236,1212)(1262,1211)
	(1288,1211)(1316,1211)(1346,1210)
	(1376,1210)(1408,1209)(1441,1208)
	(1475,1207)(1511,1206)(1547,1205)
	(1583,1204)(1621,1202)(1659,1200)
	(1697,1199)(1736,1197)(1775,1194)
	(1814,1192)(1853,1189)(1892,1187)
	(1931,1184)(1969,1181)(2008,1177)
	(2047,1174)(2085,1170)(2123,1166)
	(2162,1162)(2196,1158)(2231,1154)
	(2267,1150)(2302,1145)(2339,1140)
	(2375,1135)(2413,1130)(2451,1125)
	(2490,1119)(2529,1113)(2570,1107)
	(2610,1101)(2651,1095)(2693,1088)
	(2735,1082)(2777,1075)(2819,1069)
	(2862,1062)(2905,1055)(2947,1049)
	(2989,1042)(3031,1036)(3073,1029)
	(3114,1023)(3154,1017)(3195,1011)
	(3234,1005)(3273,999)(3311,994)
	(3349,989)(3385,984)(3422,979)
	(3457,974)(3493,970)(3528,966)
	(3562,962)(3598,958)(3635,954)
	(3671,951)(3707,947)(3744,944)
	(3780,941)(3817,939)(3855,936)
	(3892,933)(3930,931)(3968,929)
	(4007,927)(4045,925)(4083,923)
	(4122,922)(4160,920)(4198,919)
	(4236,918)(4273,917)(4310,916)
	(4347,915)(4383,914)(4418,914)
	(4453,913)(4487,913)(4520,913)
	(4553,912)(4584,912)(4615,912)
	(4646,912)(4676,912)(4705,912)
	(4734,912)(4762,912)(4794,912)
	(4826,912)(4857,912)(4889,912)
	(4922,912)(4955,912)(4989,912)
	(5024,912)(5060,912)(5098,912)
	(5138,912)(5179,912)(5221,912)
	(5266,912)(5311,912)(5357,912)
	(5404,912)(5450,912)(5494,912)
	(5536,912)(5575,912)(5610,912)
	(5639,912)(5664,912)(5683,912)
	(5696,912)(5705,912)(5710,912)(5712,912)
\put(2412,162){\makebox(0,0)[lb]{\smash{{{\SetFigFont{5}{6.0}{\rmdefault}{\mddefault}{\updefault}$()$}}}}}
\put(2562,1737){\makebox(0,0)[lb]{\smash{{{\SetFigFont{5}{6.0}{\rmdefault}{\mddefault}{\updefault}$(0)$}}}}}
\put(1362,3237){\makebox(0,0)[lb]{\smash{{{\SetFigFont{5}{6.0}{\rmdefault}{\mddefault}{\updefault}$(1)$}}}}}
\put(3912,3237){\makebox(0,0)[lb]{\smash{{{\SetFigFont{5}{6.0}{\rmdefault}{\mddefault}{\updefault}$(0,0)$}}}}}
\put(3837,4737){\makebox(0,0)[lb]{\smash{{{\SetFigFont{5}{6.0}{\rmdefault}{\mddefault}{\updefault}$(1,0)$}}}}}
\put(1212,4737){\makebox(0,0)[lb]{\smash{{{\SetFigFont{5}{6.0}{\rmdefault}{\mddefault}{\updefault}$(2)$}}}}}
\put(2637,6237){\makebox(0,0)[lb]{\smash{{{\SetFigFont{5}{6.0}{\rmdefault}{\mddefault}{\updefault}$(1,1)$}}}}}
\end{picture}
}

%% file: xfig/Gram1.eepic
\setlength{\unitlength}{0.00025000in}
\begingroup\makeatletter\ifx\SetFigFont\undefined%
\gdef\SetFigFont#1#2#3#4#5{%
  \reset@font\fontsize{#1}{#2pt}%
  \fontfamily{#3}\fontseries{#4}\fontshape{#5}%
  \selectfont}%
\fi\endgroup%
{\renewcommand{\dashlinestretch}{30}
\begin{picture}(12324,13539)(0,-10)
\put(3912,11412){\blacken\ellipse{300}{300}}
\put(3912,11412){\ellipse{300}{300}}
\path(3612,11187)(3623,11174)(3636,11161)
	(3653,11149)(3672,11138)(3694,11127)
	(3719,11116)(3747,11106)(3777,11095)
	(3808,11085)(3841,11075)(3874,11065)
	(3908,11056)(3943,11047)(3976,11039)
	(4009,11031)(4041,11025)(4071,11020)
	(4099,11016)(4124,11014)(4147,11014)
	(4168,11016)(4185,11021)(4200,11027)
	(4212,11037)(4221,11048)(4228,11062)
	(4234,11079)(4238,11098)(4240,11120)
	(4241,11145)(4241,11172)(4239,11201)
	(4237,11232)(4234,11265)(4230,11298)
	(4225,11333)(4220,11367)(4215,11402)
	(4209,11437)(4204,11471)(4198,11504)
	(4192,11535)(4186,11565)(4179,11593)
	(4173,11619)(4167,11642)(4160,11663)
	(4153,11682)(4145,11698)(4137,11712)
	(4125,11727)(4111,11740)(4096,11749)
	(4079,11756)(4061,11762)(4041,11766)
	(4020,11768)(3998,11769)(3976,11769)
	(3953,11768)(3930,11767)(3908,11764)
	(3886,11762)(3864,11758)(3844,11754)
	(3825,11748)(3807,11742)(3791,11734)
	(3776,11724)(3762,11712)(3750,11699)
	(3739,11684)(3729,11666)(3718,11646)
	(3707,11623)(3697,11599)(3685,11573)
	(3674,11545)(3663,11516)(3652,11487)
	(3641,11457)(3631,11426)(3622,11396)
	(3613,11367)(3606,11339)(3600,11312)
	(3596,11286)(3594,11263)(3595,11241)
	(3598,11221)(3603,11203)(3612,11187)
\put(4962,11412){\blacken\ellipse{300}{300}}
\put(4962,11412){\ellipse{300}{300}}
\path(4662,11187)(4673,11174)(4686,11161)
	(4703,11149)(4722,11138)(4744,11127)
	(4769,11116)(4797,11106)(4827,11095)
	(4858,11085)(4891,11075)(4924,11065)
	(4958,11056)(4993,11047)(5026,11039)
	(5059,11031)(5091,11025)(5121,11020)
	(5149,11016)(5174,11014)(5197,11014)
	(5218,11016)(5235,11021)(5250,11027)
	(5262,11037)(5271,11048)(5278,11062)
	(5284,11079)(5288,11098)(5290,11120)
	(5291,11145)(5291,11172)(5289,11201)
	(5287,11232)(5284,11265)(5280,11298)
	(5275,11333)(5270,11367)(5265,11402)
	(5259,11437)(5254,11471)(5248,11504)
	(5242,11535)(5236,11565)(5229,11593)
	(5223,11619)(5217,11642)(5210,11663)
	(5203,11682)(5195,11698)(5187,11712)
	(5175,11727)(5161,11740)(5146,11749)
	(5129,11756)(5111,11762)(5091,11766)
	(5070,11768)(5048,11769)(5026,11769)
	(5003,11768)(4980,11767)(4958,11764)
	(4936,11762)(4914,11758)(4894,11754)
	(4875,11748)(4857,11742)(4841,11734)
	(4826,11724)(4812,11712)(4800,11699)
	(4789,11684)(4779,11666)(4768,11646)
	(4757,11623)(4747,11599)(4735,11573)
	(4724,11545)(4713,11516)(4702,11487)
	(4691,11457)(4681,11426)(4672,11396)
	(4663,11367)(4656,11339)(4650,11312)
	(4646,11286)(4644,11263)(4645,11241)
	(4648,11221)(4653,11203)(4662,11187)
\put(3912,12912){\blacken\ellipse{300}{300}}
\put(3912,12912){\ellipse{300}{300}}
\path(3612,12687)(3623,12674)(3636,12661)
	(3653,12649)(3672,12638)(3694,12627)
	(3719,12616)(3747,12606)(3777,12595)
	(3808,12585)(3841,12575)(3874,12565)
	(3908,12556)(3943,12547)(3976,12539)
	(4009,12531)(4041,12525)(4071,12520)
	(4099,12516)(4124,12514)(4147,12514)
	(4168,12516)(4185,12521)(4200,12527)
	(4212,12537)(4221,12548)(4228,12562)
	(4234,12579)(4238,12598)(4240,12620)
	(4241,12645)(4241,12672)(4239,12701)
	(4237,12732)(4234,12765)(4230,12798)
	(4225,12833)(4220,12867)(4215,12902)
	(4209,12937)(4204,12971)(4198,13004)
	(4192,13035)(4186,13065)(4179,13093)
	(4173,13119)(4167,13142)(4160,13163)
	(4153,13182)(4145,13198)(4137,13212)
	(4125,13227)(4111,13240)(4096,13249)
	(4079,13256)(4061,13262)(4041,13266)
	(4020,13268)(3998,13269)(3976,13269)
	(3953,13268)(3930,13267)(3908,13264)
	(3886,13262)(3864,13258)(3844,13254)
	(3825,13248)(3807,13242)(3791,13234)
	(3776,13224)(3762,13212)(3750,13199)
	(3739,13184)(3729,13166)(3718,13146)
	(3707,13123)(3697,13099)(3685,13073)
	(3674,13045)(3663,13016)(3652,12987)
	(3641,12957)(3631,12926)(3622,12896)
	(3613,12867)(3606,12839)(3600,12812)
	(3596,12786)(3594,12763)(3595,12741)
	(3598,12721)(3603,12703)(3612,12687)
\put(4962,12912){\blacken\ellipse{300}{300}}
\put(4962,12912){\ellipse{300}{300}}
\path(4662,12687)(4673,12674)(4686,12661)
	(4703,12649)(4722,12638)(4744,12627)
	(4769,12616)(4797,12606)(4827,12595)
	(4858,12585)(4891,12575)(4924,12565)
	(4958,12556)(4993,12547)(5026,12539)
	(5059,12531)(5091,12525)(5121,12520)
	(5149,12516)(5174,12514)(5197,12514)
	(5218,12516)(5235,12521)(5250,12527)
	(5262,12537)(5271,12548)(5278,12562)
	(5284,12579)(5288,12598)(5290,12620)
	(5291,12645)(5291,12672)(5289,12701)
	(5287,12732)(5284,12765)(5280,12798)
	(5275,12833)(5270,12867)(5265,12902)
	(5259,12937)(5254,12971)(5248,13004)
	(5242,13035)(5236,13065)(5229,13093)
	(5223,13119)(5217,13142)(5210,13163)
	(5203,13182)(5195,13198)(5187,13212)
	(5175,13227)(5161,13240)(5146,13249)
	(5129,13256)(5111,13262)(5091,13266)
	(5070,13268)(5048,13269)(5026,13269)
	(5003,13268)(4980,13267)(4958,13264)
	(4936,13262)(4914,13258)(4894,13254)
	(4875,13248)(4857,13242)(4841,13234)
	(4826,13224)(4812,13212)(4800,13199)
	(4789,13184)(4779,13166)(4768,13146)
	(4757,13123)(4747,13099)(4735,13073)
	(4724,13045)(4713,13016)(4702,12987)
	(4691,12957)(4681,12926)(4672,12896)
	(4663,12867)(4656,12839)(4650,12812)
	(4646,12786)(4644,12763)(4645,12741)
	(4648,12721)(4653,12703)(4662,12687)
\put(762,8187){\blacken\ellipse{300}{300}}
\put(762,8187){\ellipse{300}{300}}
\path(462,7962)(473,7949)(486,7936)
	(503,7924)(522,7913)(544,7902)
	(569,7891)(597,7881)(627,7870)
	(658,7860)(691,7850)(724,7840)
	(758,7831)(793,7822)(826,7814)
	(859,7806)(891,7800)(921,7795)
	(949,7791)(974,7789)(997,7789)
	(1018,7791)(1035,7796)(1050,7802)
	(1062,7812)(1071,7823)(1078,7837)
	(1084,7854)(1088,7873)(1090,7895)
	(1091,7920)(1091,7947)(1089,7976)
	(1087,8007)(1084,8040)(1080,8073)
	(1075,8108)(1070,8142)(1065,8177)
	(1059,8212)(1054,8246)(1048,8279)
	(1042,8310)(1036,8340)(1029,8368)
	(1023,8394)(1017,8417)(1010,8438)
	(1003,8457)(995,8473)(987,8487)
	(975,8502)(961,8515)(946,8524)
	(929,8531)(911,8537)(891,8541)
	(870,8543)(848,8544)(826,8544)
	(803,8543)(780,8542)(758,8539)
	(736,8537)(714,8533)(694,8529)
	(675,8523)(657,8517)(641,8509)
	(626,8499)(612,8487)(600,8474)
	(589,8459)(579,8441)(568,8421)
	(557,8398)(547,8374)(535,8348)
	(524,8320)(513,8291)(502,8262)
	(491,8232)(481,8201)(472,8171)
	(463,8142)(456,8114)(450,8087)
	(446,8061)(444,8038)(445,8016)
	(448,7996)(453,7978)(462,7962)
\put(1812,8187){\blacken\ellipse{300}{300}}
\put(1812,8187){\ellipse{300}{300}}
\path(1512,7962)(1523,7949)(1536,7936)
	(1553,7924)(1572,7913)(1594,7902)
	(1619,7891)(1647,7881)(1677,7870)
	(1708,7860)(1741,7850)(1774,7840)
	(1808,7831)(1843,7822)(1876,7814)
	(1909,7806)(1941,7800)(1971,7795)
	(1999,7791)(2024,7789)(2047,7789)
	(2068,7791)(2085,7796)(2100,7802)
	(2112,7812)(2121,7823)(2128,7837)
	(2134,7854)(2138,7873)(2140,7895)
	(2141,7920)(2141,7947)(2139,7976)
	(2137,8007)(2134,8040)(2130,8073)
	(2125,8108)(2120,8142)(2115,8177)
	(2109,8212)(2104,8246)(2098,8279)
	(2092,8310)(2086,8340)(2079,8368)
	(2073,8394)(2067,8417)(2060,8438)
	(2053,8457)(2045,8473)(2037,8487)
	(2025,8502)(2011,8515)(1996,8524)
	(1979,8531)(1961,8537)(1941,8541)
	(1920,8543)(1898,8544)(1876,8544)
	(1853,8543)(1830,8542)(1808,8539)
	(1786,8537)(1764,8533)(1744,8529)
	(1725,8523)(1707,8517)(1691,8509)
	(1676,8499)(1662,8487)(1650,8474)
	(1639,8459)(1629,8441)(1618,8421)
	(1607,8398)(1597,8374)(1585,8348)
	(1574,8320)(1563,8291)(1552,8262)
	(1541,8232)(1531,8201)(1522,8171)
	(1513,8142)(1506,8114)(1500,8087)
	(1496,8061)(1494,8038)(1495,8016)
	(1498,7996)(1503,7978)(1512,7962)
\put(762,9687){\blacken\ellipse{300}{300}}
\put(762,9687){\ellipse{300}{300}}
\path(462,9462)(473,9449)(486,9436)
	(503,9424)(522,9413)(544,9402)
	(569,9391)(597,9381)(627,9370)
	(658,9360)(691,9350)(724,9340)
	(758,9331)(793,9322)(826,9314)
	(859,9306)(891,9300)(921,9295)
	(949,9291)(974,9289)(997,9289)
	(1018,9291)(1035,9296)(1050,9302)
	(1062,9312)(1071,9323)(1078,9337)
	(1084,9354)(1088,9373)(1090,9395)
	(1091,9420)(1091,9447)(1089,9476)
	(1087,9507)(1084,9540)(1080,9573)
	(1075,9608)(1070,9642)(1065,9677)
	(1059,9712)(1054,9746)(1048,9779)
	(1042,9810)(1036,9840)(1029,9868)
	(1023,9894)(1017,9917)(1010,9938)
	(1003,9957)(995,9973)(987,9987)
	(975,10002)(961,10015)(946,10024)
	(929,10031)(911,10037)(891,10041)
	(870,10043)(848,10044)(826,10044)
	(803,10043)(780,10042)(758,10039)
	(736,10037)(714,10033)(694,10029)
	(675,10023)(657,10017)(641,10009)
	(626,9999)(612,9987)(600,9974)
	(589,9959)(579,9941)(568,9921)
	(557,9898)(547,9874)(535,9848)
	(524,9820)(513,9791)(502,9762)
	(491,9732)(481,9701)(472,9671)
	(463,9642)(456,9614)(450,9587)
	(446,9561)(444,9538)(445,9516)
	(448,9496)(453,9478)(462,9462)
\put(1812,9687){\blacken\ellipse{300}{300}}
\put(1812,9687){\ellipse{300}{300}}
\path(1512,9462)(1523,9449)(1536,9436)
	(1553,9424)(1572,9413)(1594,9402)
	(1619,9391)(1647,9381)(1677,9370)
	(1708,9360)(1741,9350)(1774,9340)
	(1808,9331)(1843,9322)(1876,9314)
	(1909,9306)(1941,9300)(1971,9295)
	(1999,9291)(2024,9289)(2047,9289)
	(2068,9291)(2085,9296)(2100,9302)
	(2112,9312)(2121,9323)(2128,9337)
	(2134,9354)(2138,9373)(2140,9395)
	(2141,9420)(2141,9447)(2139,9476)
	(2137,9507)(2134,9540)(2130,9573)
	(2125,9608)(2120,9642)(2115,9677)
	(2109,9712)(2104,9746)(2098,9779)
	(2092,9810)(2086,9840)(2079,9868)
	(2073,9894)(2067,9917)(2060,9938)
	(2053,9957)(2045,9973)(2037,9987)
	(2025,10002)(2011,10015)(1996,10024)
	(1979,10031)(1961,10037)(1941,10041)
	(1920,10043)(1898,10044)(1876,10044)
	(1853,10043)(1830,10042)(1808,10039)
	(1786,10037)(1764,10033)(1744,10029)
	(1725,10023)(1707,10017)(1691,10009)
	(1676,9999)(1662,9987)(1650,9974)
	(1639,9959)(1629,9941)(1618,9921)
	(1607,9898)(1597,9874)(1585,9848)
	(1574,9820)(1563,9791)(1552,9762)
	(1541,9732)(1531,9701)(1522,9671)
	(1513,9642)(1506,9614)(1500,9587)
	(1496,9561)(1494,9538)(1495,9516)
	(1498,9496)(1503,9478)(1512,9462)
\put(7137,12912){\blacken\ellipse{300}{300}}
\put(7137,12912){\ellipse{300}{300}}
\path(6837,12687)(6848,12674)(6861,12661)
	(6878,12649)(6897,12638)(6919,12627)
	(6944,12616)(6972,12606)(7002,12595)
	(7033,12585)(7066,12575)(7099,12565)
	(7133,12556)(7168,12547)(7201,12539)
	(7234,12531)(7266,12525)(7296,12520)
	(7324,12516)(7349,12514)(7372,12514)
	(7393,12516)(7410,12521)(7425,12527)
	(7437,12537)(7446,12548)(7453,12562)
	(7459,12579)(7463,12598)(7465,12620)
	(7466,12645)(7466,12672)(7464,12701)
	(7462,12732)(7459,12765)(7455,12798)
	(7450,12833)(7445,12867)(7440,12902)
	(7434,12937)(7429,12971)(7423,13004)
	(7417,13035)(7411,13065)(7404,13093)
	(7398,13119)(7392,13142)(7385,13163)
	(7378,13182)(7370,13198)(7362,13212)
	(7350,13227)(7336,13240)(7321,13249)
	(7304,13256)(7286,13262)(7266,13266)
	(7245,13268)(7223,13269)(7201,13269)
	(7178,13268)(7155,13267)(7133,13264)
	(7111,13262)(7089,13258)(7069,13254)
	(7050,13248)(7032,13242)(7016,13234)
	(7001,13224)(6987,13212)(6975,13199)
	(6964,13184)(6954,13166)(6943,13146)
	(6932,13123)(6922,13099)(6910,13073)
	(6899,13045)(6888,13016)(6877,12987)
	(6866,12957)(6856,12926)(6847,12896)
	(6838,12867)(6831,12839)(6825,12812)
	(6821,12786)(6819,12763)(6820,12741)
	(6823,12721)(6828,12703)(6837,12687)
\put(8187,12912){\blacken\ellipse{300}{300}}
\put(8187,12912){\ellipse{300}{300}}
\path(7887,12687)(7898,12674)(7911,12661)
	(7928,12649)(7947,12638)(7969,12627)
	(7994,12616)(8022,12606)(8052,12595)
	(8083,12585)(8116,12575)(8149,12565)
	(8183,12556)(8218,12547)(8251,12539)
	(8284,12531)(8316,12525)(8346,12520)
	(8374,12516)(8399,12514)(8422,12514)
	(8443,12516)(8460,12521)(8475,12527)
	(8487,12537)(8496,12548)(8503,12562)
	(8509,12579)(8513,12598)(8515,12620)
	(8516,12645)(8516,12672)(8514,12701)
	(8512,12732)(8509,12765)(8505,12798)
	(8500,12833)(8495,12867)(8490,12902)
	(8484,12937)(8479,12971)(8473,13004)
	(8467,13035)(8461,13065)(8454,13093)
	(8448,13119)(8442,13142)(8435,13163)
	(8428,13182)(8420,13198)(8412,13212)
	(8400,13227)(8386,13240)(8371,13249)
	(8354,13256)(8336,13262)(8316,13266)
	(8295,13268)(8273,13269)(8251,13269)
	(8228,13268)(8205,13267)(8183,13264)
	(8161,13262)(8139,13258)(8119,13254)
	(8100,13248)(8082,13242)(8066,13234)
	(8051,13224)(8037,13212)(8025,13199)
	(8014,13184)(8004,13166)(7993,13146)
	(7982,13123)(7972,13099)(7960,13073)
	(7949,13045)(7938,13016)(7927,12987)
	(7916,12957)(7906,12926)(7897,12896)
	(7888,12867)(7881,12839)(7875,12812)
	(7871,12786)(7869,12763)(7870,12741)
	(7873,12721)(7878,12703)(7887,12687)
\put(7137,11412){\blacken\ellipse{300}{300}}
\put(7137,11412){\ellipse{300}{300}}
\put(8187,11412){\blacken\ellipse{300}{300}}
\put(8187,11412){\ellipse{300}{300}}
\path(6762,11712)(6778,11721)(6796,11729)
	(6815,11736)(6837,11742)(6861,11748)
	(6886,11752)(6912,11756)(6940,11760)
	(6970,11762)(7000,11765)(7031,11767)
	(7063,11768)(7096,11770)(7129,11771)
	(7162,11772)(7196,11772)(7229,11773)
	(7263,11774)(7296,11775)(7329,11776)
	(7362,11777)(7395,11778)(7427,11780)
	(7460,11781)(7491,11783)(7523,11784)
	(7555,11786)(7587,11787)(7617,11788)
	(7648,11789)(7679,11790)(7711,11792)
	(7744,11793)(7778,11795)(7812,11796)
	(7847,11798)(7883,11800)(7919,11803)
	(7956,11805)(7993,11808)(8031,11810)
	(8068,11813)(8106,11815)(8144,11817)
	(8181,11819)(8217,11821)(8253,11822)
	(8288,11823)(8322,11823)(8355,11823)
	(8387,11822)(8417,11820)(8445,11817)
	(8472,11813)(8497,11808)(8521,11802)
	(8542,11795)(8562,11787)(8582,11776)
	(8601,11764)(8618,11749)(8633,11733)
	(8648,11716)(8661,11697)(8674,11676)
	(8686,11654)(8698,11631)(8709,11607)
	(8720,11582)(8730,11557)(8739,11532)
	(8747,11506)(8755,11481)(8761,11455)
	(8766,11431)(8770,11407)(8772,11385)
	(8772,11363)(8769,11343)(8764,11324)
	(8756,11307)(8745,11290)(8730,11276)
	(8712,11262)(8695,11252)(8676,11243)
	(8655,11235)(8631,11227)(8605,11220)
	(8577,11213)(8546,11206)(8513,11200)
	(8479,11194)(8443,11188)(8405,11183)
	(8365,11178)(8325,11173)(8283,11168)
	(8241,11164)(8198,11160)(8154,11155)
	(8111,11151)(8067,11148)(8024,11144)
	(7981,11140)(7938,11137)(7896,11133)
	(7855,11130)(7816,11127)(7777,11125)
	(7739,11122)(7703,11120)(7668,11117)
	(7635,11116)(7602,11114)(7571,11113)
	(7541,11112)(7512,11112)(7475,11112)
	(7439,11114)(7404,11116)(7370,11118)
	(7337,11121)(7304,11125)(7271,11129)
	(7238,11133)(7206,11138)(7174,11143)
	(7143,11148)(7112,11153)(7081,11159)
	(7050,11165)(7020,11171)(6991,11177)
	(6963,11184)(6935,11191)(6908,11198)
	(6883,11205)(6859,11213)(6836,11222)
	(6815,11231)(6795,11241)(6778,11251)
	(6762,11262)(6746,11276)(6734,11292)
	(6723,11309)(6714,11328)(6707,11348)
	(6701,11369)(6697,11392)(6693,11416)
	(6691,11440)(6689,11465)(6688,11491)
	(6688,11516)(6689,11541)(6691,11565)
	(6693,11588)(6697,11610)(6703,11631)
	(6710,11651)(6719,11669)(6731,11685)
	(6745,11699)(6762,11712)
\put(10362,12912){\blacken\ellipse{300}{300}}
\put(10362,12912){\ellipse{300}{300}}
\path(10062,12687)(10073,12674)(10086,12661)
	(10103,12649)(10122,12638)(10144,12627)
	(10169,12616)(10197,12606)(10227,12595)
	(10258,12585)(10291,12575)(10324,12565)
	(10358,12556)(10393,12547)(10426,12539)
	(10459,12531)(10491,12525)(10521,12520)
	(10549,12516)(10574,12514)(10597,12514)
	(10618,12516)(10635,12521)(10650,12527)
	(10662,12537)(10671,12548)(10678,12562)
	(10684,12579)(10688,12598)(10690,12620)
	(10691,12645)(10691,12672)(10689,12701)
	(10687,12732)(10684,12765)(10680,12798)
	(10675,12833)(10670,12867)(10665,12902)
	(10659,12937)(10654,12971)(10648,13004)
	(10642,13035)(10636,13065)(10629,13093)
	(10623,13119)(10617,13142)(10610,13163)
	(10603,13182)(10595,13198)(10587,13212)
	(10575,13227)(10561,13240)(10546,13249)
	(10529,13256)(10511,13262)(10491,13266)
	(10470,13268)(10448,13269)(10426,13269)
	(10403,13268)(10380,13267)(10358,13264)
	(10336,13262)(10314,13258)(10294,13254)
	(10275,13248)(10257,13242)(10241,13234)
	(10226,13224)(10212,13212)(10200,13199)
	(10189,13184)(10179,13166)(10168,13146)
	(10157,13123)(10147,13099)(10135,13073)
	(10124,13045)(10113,13016)(10102,12987)
	(10091,12957)(10081,12926)(10072,12896)
	(10063,12867)(10056,12839)(10050,12812)
	(10046,12786)(10044,12763)(10045,12741)
	(10048,12721)(10053,12703)(10062,12687)
\put(11412,12912){\blacken\ellipse{300}{300}}
\put(11412,12912){\ellipse{300}{300}}
\path(11112,12687)(11123,12674)(11136,12661)
	(11153,12649)(11172,12638)(11194,12627)
	(11219,12616)(11247,12606)(11277,12595)
	(11308,12585)(11341,12575)(11374,12565)
	(11408,12556)(11443,12547)(11476,12539)
	(11509,12531)(11541,12525)(11571,12520)
	(11599,12516)(11624,12514)(11647,12514)
	(11668,12516)(11685,12521)(11700,12527)
	(11712,12537)(11721,12548)(11728,12562)
	(11734,12579)(11738,12598)(11740,12620)
	(11741,12645)(11741,12672)(11739,12701)
	(11737,12732)(11734,12765)(11730,12798)
	(11725,12833)(11720,12867)(11715,12902)
	(11709,12937)(11704,12971)(11698,13004)
	(11692,13035)(11686,13065)(11679,13093)
	(11673,13119)(11667,13142)(11660,13163)
	(11653,13182)(11645,13198)(11637,13212)
	(11625,13227)(11611,13240)(11596,13249)
	(11579,13256)(11561,13262)(11541,13266)
	(11520,13268)(11498,13269)(11476,13269)
	(11453,13268)(11430,13267)(11408,13264)
	(11386,13262)(11364,13258)(11344,13254)
	(11325,13248)(11307,13242)(11291,13234)
	(11276,13224)(11262,13212)(11250,13199)
	(11239,13184)(11229,13166)(11218,13146)
	(11207,13123)(11197,13099)(11185,13073)
	(11174,13045)(11163,13016)(11152,12987)
	(11141,12957)(11131,12926)(11122,12896)
	(11113,12867)(11106,12839)(11100,12812)
	(11096,12786)(11094,12763)(11095,12741)
	(11098,12721)(11103,12703)(11112,12687)
\put(10362,11412){\blacken\ellipse{300}{300}}
\put(10362,11412){\ellipse{300}{300}}
\put(11412,11412){\blacken\ellipse{300}{300}}
\put(11412,11412){\ellipse{300}{300}}
\thicklines
\path(10362,11412)(11562,11412)
\thinlines
\path(9987,11712)(10003,11721)(10021,11729)
	(10040,11736)(10062,11742)(10086,11748)
	(10111,11752)(10137,11756)(10165,11760)
	(10195,11762)(10225,11765)(10256,11767)
	(10288,11768)(10321,11770)(10354,11771)
	(10387,11772)(10421,11772)(10454,11773)
	(10488,11774)(10521,11775)(10554,11776)
	(10587,11777)(10620,11778)(10652,11780)
	(10685,11781)(10716,11783)(10748,11784)
	(10780,11786)(10812,11787)(10842,11788)
	(10873,11789)(10904,11790)(10936,11792)
	(10969,11793)(11003,11795)(11037,11796)
	(11072,11798)(11108,11800)(11144,11803)
	(11181,11805)(11218,11808)(11256,11810)
	(11293,11813)(11331,11815)(11369,11817)
	(11406,11819)(11442,11821)(11478,11822)
	(11513,11823)(11547,11823)(11580,11823)
	(11612,11822)(11642,11820)(11670,11817)
	(11697,11813)(11722,11808)(11746,11802)
	(11767,11795)(11787,11787)(11807,11776)
	(11826,11764)(11843,11749)(11858,11733)
	(11873,11716)(11886,11697)(11899,11676)
	(11911,11654)(11923,11631)(11934,11607)
	(11945,11582)(11955,11557)(11964,11532)
	(11972,11506)(11980,11481)(11986,11455)
	(11991,11431)(11995,11407)(11997,11385)
	(11997,11363)(11994,11343)(11989,11324)
	(11981,11307)(11970,11290)(11955,11276)
	(11937,11262)(11920,11252)(11901,11243)
	(11880,11235)(11856,11227)(11830,11220)
	(11802,11213)(11771,11206)(11738,11200)
	(11704,11194)(11668,11188)(11630,11183)
	(11590,11178)(11550,11173)(11508,11168)
	(11466,11164)(11423,11160)(11379,11155)
	(11336,11151)(11292,11148)(11249,11144)
	(11206,11140)(11163,11137)(11121,11133)
	(11080,11130)(11041,11127)(11002,11125)
	(10964,11122)(10928,11120)(10893,11117)
	(10860,11116)(10827,11114)(10796,11113)
	(10766,11112)(10737,11112)(10700,11112)
	(10664,11114)(10629,11116)(10595,11118)
	(10562,11121)(10529,11125)(10496,11129)
	(10463,11133)(10431,11138)(10399,11143)
	(10368,11148)(10337,11153)(10306,11159)
	(10275,11165)(10245,11171)(10216,11177)
	(10188,11184)(10160,11191)(10133,11198)
	(10108,11205)(10084,11213)(10061,11222)
	(10040,11231)(10020,11241)(10003,11251)
	(9987,11262)(9971,11276)(9959,11292)
	(9948,11309)(9939,11328)(9932,11348)
	(9926,11369)(9922,11392)(9918,11416)
	(9916,11440)(9914,11465)(9913,11491)
	(9913,11516)(9914,11541)(9916,11565)
	(9918,11588)(9922,11610)(9928,11631)
	(9935,11651)(9944,11669)(9956,11685)
	(9970,11699)(9987,11712)
\put(612,537){\blacken\ellipse{300}{300}}
\put(612,537){\ellipse{300}{300}}
\path(312,762)(323,775)(336,788)
	(353,800)(372,811)(394,822)
	(419,833)(447,843)(477,854)
	(508,864)(541,874)(574,884)
	(608,893)(643,902)(676,910)
	(709,918)(741,924)(771,929)
	(799,933)(824,935)(847,935)
	(868,933)(885,928)(900,922)
	(912,912)(921,901)(928,887)
	(934,870)(938,851)(940,829)
	(941,804)(941,777)(939,748)
	(937,717)(934,684)(930,651)
	(925,616)(920,582)(915,547)
	(909,512)(904,478)(898,445)
	(892,414)(886,384)(879,356)
	(873,330)(867,307)(860,286)
	(853,267)(845,251)(837,237)
	(825,222)(811,209)(796,200)
	(779,193)(761,187)(741,183)
	(720,181)(698,180)(676,180)
	(653,181)(630,182)(608,185)
	(586,187)(564,191)(544,195)
	(525,201)(507,207)(491,215)
	(476,225)(462,237)(450,250)
	(439,265)(429,283)(418,303)
	(407,326)(397,350)(385,376)
	(374,404)(363,433)(352,462)
	(341,492)(331,523)(322,553)
	(313,582)(306,610)(300,637)
	(296,663)(294,686)(295,708)
	(298,728)(303,746)(312,762)
\put(1662,537){\blacken\ellipse{300}{300}}
\put(1662,537){\ellipse{300}{300}}
\path(1362,762)(1373,775)(1386,788)
	(1403,800)(1422,811)(1444,822)
	(1469,833)(1497,843)(1527,854)
	(1558,864)(1591,874)(1624,884)
	(1658,893)(1693,902)(1726,910)
	(1759,918)(1791,924)(1821,929)
	(1849,933)(1874,935)(1897,935)
	(1918,933)(1935,928)(1950,922)
	(1962,912)(1971,901)(1978,887)
	(1984,870)(1988,851)(1990,829)
	(1991,804)(1991,777)(1989,748)
	(1987,717)(1984,684)(1980,651)
	(1975,616)(1970,582)(1965,547)
	(1959,512)(1954,478)(1948,445)
	(1942,414)(1936,384)(1929,356)
	(1923,330)(1917,307)(1910,286)
	(1903,267)(1895,251)(1887,237)
	(1875,222)(1861,209)(1846,200)
	(1829,193)(1811,187)(1791,183)
	(1770,181)(1748,180)(1726,180)
	(1703,181)(1680,182)(1658,185)
	(1636,187)(1614,191)(1594,195)
	(1575,201)(1557,207)(1541,215)
	(1526,225)(1512,237)(1500,250)
	(1489,265)(1479,283)(1468,303)
	(1457,326)(1447,350)(1435,376)
	(1424,404)(1413,433)(1402,462)
	(1391,492)(1381,523)(1372,553)
	(1363,582)(1356,610)(1350,637)
	(1346,663)(1344,686)(1345,708)
	(1348,728)(1353,746)(1362,762)
\put(612,2037){\blacken\ellipse{300}{300}}
\put(612,2037){\ellipse{300}{300}}
\put(1662,2037){\blacken\ellipse{300}{300}}
\put(1662,2037){\ellipse{300}{300}}
\thicklines
\path(612,2037)(1812,2037)
\thinlines
\path(237,1737)(253,1728)(271,1720)
	(290,1713)(312,1707)(336,1701)
	(361,1697)(387,1693)(415,1689)
	(445,1687)(475,1684)(506,1682)
	(538,1681)(571,1679)(604,1678)
	(637,1677)(671,1677)(704,1676)
	(738,1675)(771,1674)(804,1673)
	(837,1672)(870,1671)(902,1669)
	(935,1668)(966,1666)(998,1665)
	(1030,1663)(1062,1662)(1092,1661)
	(1123,1660)(1154,1659)(1186,1657)
	(1219,1656)(1253,1654)(1287,1653)
	(1322,1651)(1358,1649)(1394,1646)
	(1431,1644)(1468,1641)(1506,1639)
	(1543,1636)(1581,1634)(1619,1632)
	(1656,1630)(1692,1628)(1728,1627)
	(1763,1626)(1797,1626)(1830,1626)
	(1862,1627)(1892,1629)(1920,1632)
	(1947,1636)(1972,1641)(1996,1647)
	(2017,1654)(2037,1662)(2057,1673)
	(2076,1685)(2093,1700)(2108,1716)
	(2123,1733)(2136,1752)(2149,1773)
	(2161,1795)(2173,1818)(2184,1842)
	(2195,1867)(2205,1892)(2214,1917)
	(2222,1943)(2230,1968)(2236,1994)
	(2241,2018)(2245,2042)(2247,2064)
	(2247,2086)(2244,2106)(2239,2125)
	(2231,2142)(2220,2159)(2205,2173)
	(2187,2187)(2170,2197)(2151,2206)
	(2130,2214)(2106,2222)(2080,2229)
	(2052,2236)(2021,2243)(1988,2249)
	(1954,2255)(1918,2261)(1880,2266)
	(1840,2271)(1800,2276)(1758,2281)
	(1716,2285)(1673,2289)(1629,2294)
	(1586,2298)(1542,2301)(1499,2305)
	(1456,2309)(1413,2312)(1371,2316)
	(1330,2319)(1291,2322)(1252,2324)
	(1214,2327)(1178,2329)(1143,2332)
	(1110,2333)(1077,2335)(1046,2336)
	(1016,2337)(987,2337)(950,2337)
	(914,2335)(879,2333)(845,2331)
	(812,2328)(779,2324)(746,2320)
	(713,2316)(681,2311)(649,2306)
	(618,2301)(587,2296)(556,2290)
	(525,2284)(495,2278)(466,2272)
	(438,2265)(410,2258)(383,2251)
	(358,2244)(334,2236)(311,2227)
	(290,2218)(270,2208)(253,2198)
	(237,2187)(221,2173)(209,2157)
	(198,2140)(189,2121)(182,2101)
	(176,2080)(172,2057)(168,2033)
	(166,2009)(164,1984)(163,1958)
	(163,1933)(164,1908)(166,1884)
	(168,1861)(172,1839)(178,1818)
	(185,1798)(194,1780)(206,1764)
	(220,1750)(237,1737)
\put(612,4362){\blacken\ellipse{300}{300}}
\put(612,4362){\ellipse{300}{300}}
\path(312,4587)(323,4600)(336,4613)
	(353,4625)(372,4636)(394,4647)
	(419,4658)(447,4668)(477,4679)
	(508,4689)(541,4699)(574,4709)
	(608,4718)(643,4727)(676,4735)
	(709,4743)(741,4749)(771,4754)
	(799,4758)(824,4760)(847,4760)
	(868,4758)(885,4753)(900,4747)
	(912,4737)(921,4726)(928,4712)
	(934,4695)(938,4676)(940,4654)
	(941,4629)(941,4602)(939,4573)
	(937,4542)(934,4509)(930,4476)
	(925,4441)(920,4407)(915,4372)
	(909,4337)(904,4303)(898,4270)
	(892,4239)(886,4209)(879,4181)
	(873,4155)(867,4132)(860,4111)
	(853,4092)(845,4076)(837,4062)
	(825,4047)(811,4034)(796,4025)
	(779,4018)(761,4012)(741,4008)
	(720,4006)(698,4005)(676,4005)
	(653,4006)(630,4007)(608,4010)
	(586,4012)(564,4016)(544,4020)
	(525,4026)(507,4032)(491,4040)
	(476,4050)(462,4062)(450,4075)
	(439,4090)(429,4108)(418,4128)
	(407,4151)(397,4175)(385,4201)
	(374,4229)(363,4258)(352,4287)
	(341,4317)(331,4348)(322,4378)
	(313,4407)(306,4435)(300,4462)
	(296,4488)(294,4511)(295,4533)
	(298,4553)(303,4571)(312,4587)
\put(1662,4362){\blacken\ellipse{300}{300}}
\put(1662,4362){\ellipse{300}{300}}
\path(1362,4587)(1373,4600)(1386,4613)
	(1403,4625)(1422,4636)(1444,4647)
	(1469,4658)(1497,4668)(1527,4679)
	(1558,4689)(1591,4699)(1624,4709)
	(1658,4718)(1693,4727)(1726,4735)
	(1759,4743)(1791,4749)(1821,4754)
	(1849,4758)(1874,4760)(1897,4760)
	(1918,4758)(1935,4753)(1950,4747)
	(1962,4737)(1971,4726)(1978,4712)
	(1984,4695)(1988,4676)(1990,4654)
	(1991,4629)(1991,4602)(1989,4573)
	(1987,4542)(1984,4509)(1980,4476)
	(1975,4441)(1970,4407)(1965,4372)
	(1959,4337)(1954,4303)(1948,4270)
	(1942,4239)(1936,4209)(1929,4181)
	(1923,4155)(1917,4132)(1910,4111)
	(1903,4092)(1895,4076)(1887,4062)
	(1875,4047)(1861,4034)(1846,4025)
	(1829,4018)(1811,4012)(1791,4008)
	(1770,4006)(1748,4005)(1726,4005)
	(1703,4006)(1680,4007)(1658,4010)
	(1636,4012)(1614,4016)(1594,4020)
	(1575,4026)(1557,4032)(1541,4040)
	(1526,4050)(1512,4062)(1500,4075)
	(1489,4090)(1479,4108)(1468,4128)
	(1457,4151)(1447,4175)(1435,4201)
	(1424,4229)(1413,4258)(1402,4287)
	(1391,4317)(1381,4348)(1372,4378)
	(1363,4407)(1356,4435)(1350,4462)
	(1346,4488)(1344,4511)(1345,4533)
	(1348,4553)(1353,4571)(1362,4587)
\put(612,5862){\blacken\ellipse{300}{300}}
\put(612,5862){\ellipse{300}{300}}
\put(1662,5862){\blacken\ellipse{300}{300}}
\put(1662,5862){\ellipse{300}{300}}
\path(237,5562)(253,5553)(271,5545)
	(290,5538)(312,5532)(336,5526)
	(361,5522)(387,5518)(415,5514)
	(445,5512)(475,5509)(506,5507)
	(538,5506)(571,5504)(604,5503)
	(637,5502)(671,5502)(704,5501)
	(738,5500)(771,5499)(804,5498)
	(837,5497)(870,5496)(902,5494)
	(935,5493)(966,5491)(998,5490)
	(1030,5488)(1062,5487)(1092,5486)
	(1123,5485)(1154,5484)(1186,5482)
	(1219,5481)(1253,5479)(1287,5478)
	(1322,5476)(1358,5474)(1394,5471)
	(1431,5469)(1468,5466)(1506,5464)
	(1543,5461)(1581,5459)(1619,5457)
	(1656,5455)(1692,5453)(1728,5452)
	(1763,5451)(1797,5451)(1830,5451)
	(1862,5452)(1892,5454)(1920,5457)
	(1947,5461)(1972,5466)(1996,5472)
	(2017,5479)(2037,5487)(2057,5498)
	(2076,5510)(2093,5525)(2108,5541)
	(2123,5558)(2136,5577)(2149,5598)
	(2161,5620)(2173,5643)(2184,5667)
	(2195,5692)(2205,5717)(2214,5742)
	(2222,5768)(2230,5793)(2236,5819)
	(2241,5843)(2245,5867)(2247,5889)
	(2247,5911)(2244,5931)(2239,5950)
	(2231,5967)(2220,5984)(2205,5998)
	(2187,6012)(2170,6022)(2151,6031)
	(2130,6039)(2106,6047)(2080,6054)
	(2052,6061)(2021,6068)(1988,6074)
	(1954,6080)(1918,6086)(1880,6091)
	(1840,6096)(1800,6101)(1758,6106)
	(1716,6110)(1673,6114)(1629,6119)
	(1586,6123)(1542,6126)(1499,6130)
	(1456,6134)(1413,6137)(1371,6141)
	(1330,6144)(1291,6147)(1252,6149)
	(1214,6152)(1178,6154)(1143,6157)
	(1110,6158)(1077,6160)(1046,6161)
	(1016,6162)(987,6162)(950,6162)
	(914,6160)(879,6158)(845,6156)
	(812,6153)(779,6149)(746,6145)
	(713,6141)(681,6136)(649,6131)
	(618,6126)(587,6121)(556,6115)
	(525,6109)(495,6103)(466,6097)
	(438,6090)(410,6083)(383,6076)
	(358,6069)(334,6061)(311,6052)
	(290,6043)(270,6033)(253,6023)
	(237,6012)(221,5998)(209,5982)
	(198,5965)(189,5946)(182,5926)
	(176,5905)(172,5882)(168,5858)
	(166,5834)(164,5809)(163,5783)
	(163,5758)(164,5733)(166,5709)
	(168,5686)(172,5664)(178,5643)
	(185,5623)(194,5605)(206,5589)
	(220,5575)(237,5562)
\path(2712,13512)(2712,12)
\path(12,10512)(12312,10512)
\put(3387,8712){\makebox(0,0)[lb]{\smash{{{\SetFigFont{6}{7.2}{\rmdefault}{\mddefault}{\updefault}$Q_1^2 Q_2^2$}}}}}
\put(6762,8787){\makebox(0,0)[lb]{\smash{{{\SetFigFont{6}{7.2}{\rmdefault}{\mddefault}{\updefault}$Q_1 Q_2^2$}}}}}
\put(9987,8712){\makebox(0,0)[lb]{\smash{{{\SetFigFont{6}{7.2}{\rmdefault}{\mddefault}{\updefault}$Q_1 Q_2$}}}}}
\put(3387,4812){\makebox(0,0)[lb]{\smash{{{\SetFigFont{6}{7.2}{\rmdefault}{\mddefault}{\updefault}$Q_1 Q_2^2$}}}}}
\put(6762,4812){\makebox(0,0)[lb]{\smash{{{\SetFigFont{6}{7.2}{\rmdefault}{\mddefault}{\updefault}$Q_1 Q_2^2$}}}}}
\put(9987,4887){\makebox(0,0)[lb]{\smash{{{\SetFigFont{6}{7.2}{\rmdefault}{\mddefault}{\updefault}$Q_1 Q_2$}}}}}
\put(3387,987){\makebox(0,0)[lb]{\smash{{{\SetFigFont{6}{7.2}{\rmdefault}{\mddefault}{\updefault}$Q_1 Q_2$}}}}}
\put(6762,987){\makebox(0,0)[lb]{\smash{{{\SetFigFont{6}{7.2}{\rmdefault}{\mddefault}{\updefault}$Q_1 Q_2$}}}}}
\put(9987,987){\makebox(0,0)[lb]{\smash{{{\SetFigFont{6}{7.2}{\rmdefault}{\mddefault}{\updefault}$Q_1 Q_2$}}}}}
\end{picture}
}

%% file: xfig/bratteli.eepic
\setlength{\unitlength}{0.00041667in}
\begingroup\makeatletter\ifx\SetFigFont\undefined%
\gdef\SetFigFont#1#2#3#4#5{%
  \reset@font\fontsize{#1}{#2pt}%
  \fontfamily{#3}\fontseries{#4}\fontshape{#5}%
  \selectfont}%
\fi\endgroup%
{\renewcommand{\dashlinestretch}{30}
\begin{picture}(13224,20439)(0,-10)
\path(12,20412)(612,20412)(612,19812)
	(12,19812)(12,20412)
\put(162,20187){\makebox(0,0)[lb]{\smash{{{\SetFigFont{6}{7.2}{\rmdefault}{\mddefault}{\updefault}0,0}}}}}
\path(12,13212)(612,13212)(612,12612)
	(12,12612)(12,13212)
\path(12,12612)(612,12612)(612,13212)
	(12,13212)(12,12612)
\put(162,12987){\makebox(0,0)[lb]{\smash{{{\SetFigFont{6}{7.2}{\rmdefault}{\mddefault}{\updefault}0,0}}}}}
\path(12,19212)(612,19212)(612,18612)
	(12,18612)(12,19212)
\put(162,18987){\makebox(0,0)[lb]{\smash{{{\SetFigFont{6}{7.2}{\rmdefault}{\mddefault}{\updefault}1,1}}}}}
\path(12,17112)(612,17112)(612,17712)
	(12,17712)(12,17112)
\path(12,17712)(612,17712)(612,17112)
	(12,17112)(12,17712)
\put(162,17487){\makebox(0,0)[lb]{\smash{{{\SetFigFont{6}{7.2}{\rmdefault}{\mddefault}{\updefault}0,0}}}}}
\path(12,15912)(612,15912)(612,15312)
	(12,15312)(12,15912)
\path(12,15312)(612,15312)(612,15912)
	(12,15912)(12,15312)
\put(162,15687){\makebox(0,0)[lb]{\smash{{{\SetFigFont{6}{7.2}{\rmdefault}{\mddefault}{\updefault}1,1}}}}}
\path(12,10512)(612,10512)(612,9912)
	(12,9912)(12,10512)
\path(12,9912)(612,9912)(612,10512)
	(12,10512)(12,9912)
\put(162,10287){\makebox(0,0)[lb]{\smash{{{\SetFigFont{6}{7.2}{\rmdefault}{\mddefault}{\updefault}1,1}}}}}
\path(912,17112)(1512,17112)(1512,17712)
	(912,17712)(912,17112)
\path(912,17712)(1512,17712)(1512,17112)
	(912,17112)(912,17712)
\put(1062,17487){\makebox(0,0)[lb]{\smash{{{\SetFigFont{6}{7.2}{\rmdefault}{\mddefault}{\updefault}0,1}}}}}
\path(912,15312)(1512,15312)(1512,15912)
	(912,15912)(912,15312)
\path(912,15912)(1512,15912)(1512,15312)
	(912,15312)(912,15912)
\put(1062,15687){\makebox(0,0)[lb]{\smash{{{\SetFigFont{6}{7.2}{\rmdefault}{\mddefault}{\updefault}1,2}}}}}
\path(912,13212)(1512,13212)(1512,12612)
	(912,12612)(912,13212)
\path(912,12612)(1512,12612)(1512,13212)
	(912,13212)(912,12612)
\put(1062,12987){\makebox(0,0)[lb]{\smash{{{\SetFigFont{6}{7.2}{\rmdefault}{\mddefault}{\updefault}0,1}}}}}
\path(912,9912)(1512,9912)(1512,10512)
	(912,10512)(912,9912)
\path(912,10512)(1512,10512)(1512,9912)
	(912,9912)(912,10512)
\put(1062,10287){\makebox(0,0)[lb]{\smash{{{\SetFigFont{6}{7.2}{\rmdefault}{\mddefault}{\updefault}1,2}}}}}
\path(1812,10512)(2412,10512)(2412,9912)
	(1812,9912)(1812,10512)
\put(1887,10287){\makebox(0,0)[lb]{\smash{{{\SetFigFont{6}{7.2}{\rmdefault}{\mddefault}{\updefault}2,1}}}}}
\path(1812,15912)(2412,15912)(2412,15312)
	(1812,15312)(1812,15912)
\put(1887,15687){\makebox(0,0)[lb]{\smash{{{\SetFigFont{6}{7.2}{\rmdefault}{\mddefault}{\updefault}2,1}}}}}
\path(2712,13212)(3312,13212)(3312,12612)
	(2712,12612)(2712,13212)
\path(3312,13212)(3912,13212)(3912,12612)
	(3312,12612)(3312,13212)
\put(2787,12987){\makebox(0,0)[lb]{\smash{{{\SetFigFont{6}{7.2}{\rmdefault}{\mddefault}{\updefault}0,2}}}}}
\put(3387,12987){\makebox(0,0)[lb]{\smash{{{\SetFigFont{6}{7.2}{\rmdefault}{\mddefault}{\updefault}0,2}}}}}
\path(4212,12612)(4812,12612)(4812,13212)
	(4212,13212)(4212,12612)
\put(4287,12987){\makebox(0,0)[lb]{\smash{{{\SetFigFont{6}{7.2}{\rmdefault}{\mddefault}{\updefault}1,1}}}}}
\path(4212,15912)(4812,15912)(4812,15312)
	(4212,15312)(4212,15912)
\put(4287,15687){\makebox(0,0)[lb]{\smash{{{\SetFigFont{6}{7.2}{\rmdefault}{\mddefault}{\updefault}2,2}}}}}
\path(4212,17712)(4812,17712)(4812,17112)
	(4212,17112)(4212,17712)
\put(4287,17487){\makebox(0,0)[lb]{\smash{{{\SetFigFont{6}{7.2}{\rmdefault}{\mddefault}{\updefault}1,1}}}}}
\path(2712,10512)(3312,10512)(3312,9912)
	(2712,9912)(2712,10512)
\path(3312,10512)(3912,10512)(3912,9912)
	(3312,9912)(3312,10512)
\put(2787,10287){\makebox(0,0)[lb]{\smash{{{\SetFigFont{6}{7.2}{\rmdefault}{\mddefault}{\updefault}1,3}}}}}
\put(3387,10287){\makebox(0,0)[lb]{\smash{{{\SetFigFont{6}{7.2}{\rmdefault}{\mddefault}{\updefault}1,3}}}}}
\path(4212,10512)(4812,10512)(4812,9912)
	(4212,9912)(4212,10512)
\put(4287,10287){\makebox(0,0)[lb]{\smash{{{\SetFigFont{6}{7.2}{\rmdefault}{\mddefault}{\updefault}2,2}}}}}
\path(4812,10512)(5412,10512)(5412,9912)
	(4812,9912)(4812,10512)
\put(4887,10287){\makebox(0,0)[lb]{\smash{{{\SetFigFont{6}{7.2}{\rmdefault}{\mddefault}{\updefault}2,2}}}}}
\path(5712,10512)(6312,10512)(6312,9912)
	(5712,9912)(5712,10512)
\path(6312,10512)(6912,10512)(6912,9912)
	(6312,9912)(6312,10512)
\put(5787,10287){\makebox(0,0)[lb]{\smash{{{\SetFigFont{6}{7.2}{\rmdefault}{\mddefault}{\updefault}3,1}}}}}
\put(6387,10287){\makebox(0,0)[lb]{\smash{{{\SetFigFont{6}{7.2}{\rmdefault}{\mddefault}{\updefault}3,1}}}}}
\path(12,7812)(612,7812)(612,7212)
	(12,7212)(12,7812)
\path(12,7212)(612,7212)(612,7812)
	(12,7812)(12,7212)
\put(162,7587){\makebox(0,0)[lb]{\smash{{{\SetFigFont{6}{7.2}{\rmdefault}{\mddefault}{\updefault}0,0}}}}}
\path(912,7812)(1512,7812)(1512,7212)
	(912,7212)(912,7812)
\path(912,7212)(1512,7212)(1512,7812)
	(912,7812)(912,7212)
\put(1062,7587){\makebox(0,0)[lb]{\smash{{{\SetFigFont{6}{7.2}{\rmdefault}{\mddefault}{\updefault}0,1}}}}}
\path(2712,7812)(3312,7812)(3312,7212)
	(2712,7212)(2712,7812)
\path(3312,7812)(3912,7812)(3912,7212)
	(3312,7212)(3312,7812)
\put(2787,7587){\makebox(0,0)[lb]{\smash{{{\SetFigFont{6}{7.2}{\rmdefault}{\mddefault}{\updefault}0,2}}}}}
\put(3387,7587){\makebox(0,0)[lb]{\smash{{{\SetFigFont{6}{7.2}{\rmdefault}{\mddefault}{\updefault}0,2}}}}}
\path(4212,7212)(4812,7212)(4812,7812)
	(4212,7812)(4212,7212)
\put(4287,7587){\makebox(0,0)[lb]{\smash{{{\SetFigFont{6}{7.2}{\rmdefault}{\mddefault}{\updefault}1,1}}}}}
\path(12,5112)(612,5112)(612,4512)
	(12,4512)(12,5112)
\path(12,4512)(612,4512)(612,5112)
	(12,5112)(12,4512)
\put(162,4887){\makebox(0,0)[lb]{\smash{{{\SetFigFont{6}{7.2}{\rmdefault}{\mddefault}{\updefault}1,1}}}}}
\path(912,4512)(1512,4512)(1512,5112)
	(912,5112)(912,4512)
\path(912,5112)(1512,5112)(1512,4512)
	(912,4512)(912,5112)
\put(1062,4887){\makebox(0,0)[lb]{\smash{{{\SetFigFont{6}{7.2}{\rmdefault}{\mddefault}{\updefault}1,2}}}}}
\path(1812,5112)(2412,5112)(2412,4512)
	(1812,4512)(1812,5112)
\put(1887,4887){\makebox(0,0)[lb]{\smash{{{\SetFigFont{6}{7.2}{\rmdefault}{\mddefault}{\updefault}2,1}}}}}
\path(4212,5112)(4812,5112)(4812,4512)
	(4212,4512)(4212,5112)
\put(4287,4887){\makebox(0,0)[lb]{\smash{{{\SetFigFont{6}{7.2}{\rmdefault}{\mddefault}{\updefault}2,2}}}}}
\path(12,3312)(612,3312)(612,2712)
	(12,2712)(12,3312)
\path(12,2712)(612,2712)(612,3312)
	(12,3312)(12,2712)
\put(162,3087){\makebox(0,0)[lb]{\smash{{{\SetFigFont{6}{7.2}{\rmdefault}{\mddefault}{\updefault}0,0}}}}}
\path(912,3312)(1512,3312)(1512,2712)
	(912,2712)(912,3312)
\path(912,2712)(1512,2712)(1512,3312)
	(912,3312)(912,2712)
\put(1062,3087){\makebox(0,0)[lb]{\smash{{{\SetFigFont{6}{7.2}{\rmdefault}{\mddefault}{\updefault}0,1}}}}}
\path(4212,2712)(4812,2712)(4812,3312)
	(4212,3312)(4212,2712)
\put(4287,3087){\makebox(0,0)[lb]{\smash{{{\SetFigFont{6}{7.2}{\rmdefault}{\mddefault}{\updefault}1,1}}}}}
\path(12,1812)(612,1812)(612,1212)
	(12,1212)(12,1812)
\path(12,1212)(612,1212)(612,1812)
	(12,1812)(12,1212)
\put(162,1587){\makebox(0,0)[lb]{\smash{{{\SetFigFont{6}{7.2}{\rmdefault}{\mddefault}{\updefault}1,1}}}}}
\path(12,612)(612,612)(612,12)
	(12,12)(12,612)
\path(12,12)(612,12)(612,612)
	(12,612)(12,12)
\put(162,387){\makebox(0,0)[lb]{\smash{{{\SetFigFont{6}{7.2}{\rmdefault}{\mddefault}{\updefault}0,0}}}}}
\path(8112,7812)(8712,7812)(8712,7212)
	(8112,7212)(8112,7812)
\path(8712,7812)(9312,7812)(9312,7212)
	(8712,7212)(8712,7812)
\put(8187,7587){\makebox(0,0)[lb]{\smash{{{\SetFigFont{6}{7.2}{\rmdefault}{\mddefault}{\updefault}2,1}}}}}
\put(8787,7587){\makebox(0,0)[lb]{\smash{{{\SetFigFont{6}{7.2}{\rmdefault}{\mddefault}{\updefault}2,1}}}}}
\put(8187,7287){\makebox(0,0)[lb]{\smash{{{\SetFigFont{6}{7.2}{\rmdefault}{\mddefault}{\updefault}9}}}}}
\put(8787,7287){\makebox(0,0)[lb]{\smash{{{\SetFigFont{6}{7.2}{\rmdefault}{\mddefault}{\updefault}9}}}}}
\path(8112,10512)(8712,10512)(8712,9912)
	(8112,9912)(8112,10512)
\path(8712,10512)(9312,10512)(9312,9912)
	(8712,9912)(8712,10512)
\path(9312,10512)(9912,10512)(9912,9912)
	(9312,9912)(9312,10512)
\put(8187,10287){\makebox(0,0)[lb]{\smash{{{\SetFigFont{6}{7.2}{\rmdefault}{\mddefault}{\updefault}3,2}}}}}
\put(8787,10287){\makebox(0,0)[lb]{\smash{{{\SetFigFont{6}{7.2}{\rmdefault}{\mddefault}{\updefault}3,2}}}}}
\put(9387,10287){\makebox(0,0)[lb]{\smash{{{\SetFigFont{6}{7.2}{\rmdefault}{\mddefault}{\updefault}3,2}}}}}
\put(8187,9987){\makebox(0,0)[lb]{\smash{{{\SetFigFont{6}{7.2}{\rmdefault}{\mddefault}{\updefault}1}}}}}
\put(8787,9987){\makebox(0,0)[lb]{\smash{{{\SetFigFont{6}{7.2}{\rmdefault}{\mddefault}{\updefault}1}}}}}
\put(9387,9987){\makebox(0,0)[lb]{\smash{{{\SetFigFont{6}{7.2}{\rmdefault}{\mddefault}{\updefault}2}}}}}
\path(8112,13212)(8712,13212)(8712,12612)
	(8112,12612)(8112,13212)
\path(8712,13212)(9312,13212)(9312,12612)
	(8712,12612)(8712,13212)
\put(8187,12987){\makebox(0,0)[lb]{\smash{{{\SetFigFont{6}{7.2}{\rmdefault}{\mddefault}{\updefault}2,1}}}}}
\put(8787,12987){\makebox(0,0)[lb]{\smash{{{\SetFigFont{6}{7.2}{\rmdefault}{\mddefault}{\updefault}2,1}}}}}
\put(8187,12687){\makebox(0,0)[lb]{\smash{{{\SetFigFont{6}{7.2}{\rmdefault}{\mddefault}{\updefault}1}}}}}
\put(8787,12687){\makebox(0,0)[lb]{\smash{{{\SetFigFont{6}{7.2}{\rmdefault}{\mddefault}{\updefault}1}}}}}
\path(10212,10512)(10812,10512)(10812,9912)
	(10212,9912)(10212,10512)
\path(10812,10512)(11412,10512)(11412,9912)
	(10812,9912)(10812,10512)
\put(10287,10287){\makebox(0,0)[lb]{\smash{{{\SetFigFont{6}{7.2}{\rmdefault}{\mddefault}{\updefault}3,3}}}}}
\put(10887,10287){\makebox(0,0)[lb]{\smash{{{\SetFigFont{6}{7.2}{\rmdefault}{\mddefault}{\updefault}3,3}}}}}
\put(10287,9987){\makebox(0,0)[lb]{\smash{{{\SetFigFont{6}{7.2}{\rmdefault}{\mddefault}{\updefault}1}}}}}
\put(10887,9987){\makebox(0,0)[lb]{\smash{{{\SetFigFont{6}{7.2}{\rmdefault}{\mddefault}{\updefault}1}}}}}
\path(10212,13212)(10812,13212)(10812,12612)
	(10212,12612)(10212,13212)
\put(10287,12987){\makebox(0,0)[lb]{\smash{{{\SetFigFont{6}{7.2}{\rmdefault}{\mddefault}{\updefault}2,2}}}}}
\path(10812,13212)(11412,13212)(11412,12612)
	(10812,12612)(10812,13212)
\put(10887,12987){\makebox(0,0)[lb]{\smash{{{\SetFigFont{6}{7.2}{\rmdefault}{\mddefault}{\updefault}2,2}}}}}
\put(10287,12687){\makebox(0,0)[lb]{\smash{{{\SetFigFont{6}{7.2}{\rmdefault}{\mddefault}{\updefault}1}}}}}
\put(10887,12687){\makebox(0,0)[lb]{\smash{{{\SetFigFont{6}{7.2}{\rmdefault}{\mddefault}{\updefault}1}}}}}
\path(10212,7812)(10812,7812)(10812,7212)
	(10212,7212)(10212,7812)
\put(10287,7587){\makebox(0,0)[lb]{\smash{{{\SetFigFont{6}{7.2}{\rmdefault}{\mddefault}{\updefault}2,2}}}}}
\path(10812,7812)(11412,7812)(11412,7212)
	(10812,7212)(10812,7812)
\put(10887,7587){\makebox(0,0)[lb]{\smash{{{\SetFigFont{6}{7.2}{\rmdefault}{\mddefault}{\updefault}2,2}}}}}
\put(10287,7287){\makebox(0,0)[lb]{\smash{{{\SetFigFont{6}{7.2}{\rmdefault}{\mddefault}{\updefault}12}}}}}
\put(10887,7287){\makebox(0,0)[lb]{\smash{{{\SetFigFont{6}{7.2}{\rmdefault}{\mddefault}{\updefault}12}}}}}
\path(12012,7812)(13212,7812)(13212,7212)
	(12012,7212)(12012,7812)
\put(12087,7437){\makebox(0,0)[lb]{\smash{{{\SetFigFont{6}{7.2}{\rmdefault}{\mddefault}{\updefault}OTHERS}}}}}
\path(312,19812)(312,19212)
\path(312,15912)(312,17112)
\path(312,15912)(1212,17112)
\path(312,15912)(4512,17112)
\path(312,18612)(312,17712)
\path(312,18612)(1212,17712)
\path(312,18612)(4512,17712)
\path(1212,17112)(1212,15912)
\path(4512,17112)(2112,15912)
\path(4512,17112)(4512,15912)
\path(312,9912)(312,7812)
\path(312,9912)(1212,7812)
\path(1212,7812)(1212,9912)
\path(4512,7812)(2112,9912)
\path(4512,7812)(4512,9912)
\path(1212,9912)(3012,7812)
\path(1212,9912)(3612,7812)
\path(1212,9912)(7512,7812)
\path(2112,9912)(8412,7812)
\path(2112,9912)(9012,7812)
\path(4512,9912)(7512,7812)
\path(4512,9912)(10512,7812)
\path(4512,9912)(11112,7812)
\path(312,9912)(4512,7812)
\path(7512,7812)(7512,9912)
\path(8412,7812)(8412,9912)
\path(9012,7812)(9012,9912)
\path(10512,7812)(10512,9912)
\path(11112,7812)(11112,9912)
\path(10512,7812)(9612,9912)
\path(11112,7812)(9612,9912)
\path(6012,9912)(8412,7812)
\path(6612,9912)(9012,7812)
\path(5112,9912)(7512,7812)
\path(312,10512)(312,12612)
\path(312,10512)(1212,12612)
\path(1212,12612)(1212,10512)
\path(4512,12612)(2112,10512)
\path(4512,12612)(4512,10512)
\path(1212,10512)(3012,12612)
\path(1212,10512)(3612,12612)
\path(1212,10512)(7512,12612)
\path(2112,10512)(8412,12612)
\path(2112,10512)(9012,12612)
\path(4512,10512)(7512,12612)
\path(4512,10512)(10512,12612)
\path(4512,10512)(11112,12612)
\path(312,10512)(4512,12612)
\path(7512,12612)(7512,10512)
\path(8412,12612)(8412,10512)
\path(9012,12612)(9012,10512)
\path(10512,12612)(10512,10512)
\path(11112,12612)(11112,10512)
\path(10512,12612)(9612,10512)
\path(11112,12612)(9612,10512)
\path(6012,10512)(8412,12612)
\path(6612,10512)(9012,12612)
\path(5112,10512)(7512,12612)
\path(312,15312)(312,13212)
\path(312,15312)(1212,13212)
\path(1212,13212)(1212,15312)
\path(4512,13212)(2112,15312)
\path(4512,13212)(4512,15312)
\path(1212,15312)(3012,13212)
\path(1212,15312)(3612,13212)
\path(1212,15312)(7512,13212)
\path(2112,15312)(8412,13212)
\path(2112,15312)(9012,13212)
\path(4512,15312)(7512,13212)
\path(4512,15312)(10512,13212)
\path(4512,15312)(11112,13212)
\path(312,15312)(4512,13212)
\path(312,5112)(312,7212)
\path(312,5112)(1212,7212)
\path(1212,7212)(1212,5112)
\path(4512,7212)(2112,5112)
\path(4512,7212)(4512,5112)
\path(1212,5112)(3012,7212)
\path(1212,5112)(3612,7212)
\path(1212,5112)(7512,7212)
\path(2112,5112)(8412,7212)
\path(2112,5112)(9012,7212)
\path(4512,5112)(7512,7212)
\path(4512,5112)(10512,7212)
\path(4512,5112)(11112,7212)
\path(312,5112)(4512,7212)
\path(312,4512)(312,3312)
\path(312,4512)(1212,3312)
\path(312,4512)(4512,3312)
\path(1212,4512)(1212,3312)
\path(2112,4512)(4512,3312)
\path(4512,4512)(4512,3312)
\path(312,2712)(312,1812)
\path(1212,2712)(312,1812)
\path(4512,2712)(312,1812)
\path(312,1212)(312,612)
\path(5412,7812)(6612,7812)(6612,7212)
	(5412,7212)(5412,7812)
\path(6612,5112)(8412,5112)(8412,4512)
	(6612,4512)(6612,5112)
\path(3012,12612)(3012,10512)
\path(3612,12612)(3612,10512)
\path(3012,9912)(3012,7812)
\path(3612,9912)(3612,7812)
\path(7212,13212)(7812,13212)(7812,12612)
	(7212,12612)(7212,13212)
\path(7212,7812)(7812,7812)(7812,7212)
	(7212,7212)(7212,7812)
\path(7212,10512)(7812,10512)(7812,9912)
	(7212,9912)(7212,10512)
\put(162,19887){\makebox(0,0)[lb]{\smash{{{\SetFigFont{6}{7.2}{\rmdefault}{\mddefault}{\updefault}1}}}}}
\put(162,18687){\makebox(0,0)[lb]{\smash{{{\SetFigFont{6}{7.2}{\rmdefault}{\mddefault}{\updefault}1}}}}}
\put(162,17187){\makebox(0,0)[lb]{\smash{{{\SetFigFont{6}{7.2}{\rmdefault}{\mddefault}{\updefault}1}}}}}
\put(1062,17187){\makebox(0,0)[lb]{\smash{{{\SetFigFont{6}{7.2}{\rmdefault}{\mddefault}{\updefault}1}}}}}
\put(4287,17187){\makebox(0,0)[lb]{\smash{{{\SetFigFont{6}{7.2}{\rmdefault}{\mddefault}{\updefault}1}}}}}
\put(162,15387){\makebox(0,0)[lb]{\smash{{{\SetFigFont{6}{7.2}{\rmdefault}{\mddefault}{\updefault}3}}}}}
\put(1062,15387){\makebox(0,0)[lb]{\smash{{{\SetFigFont{6}{7.2}{\rmdefault}{\mddefault}{\updefault}1}}}}}
\put(1887,15387){\makebox(0,0)[lb]{\smash{{{\SetFigFont{6}{7.2}{\rmdefault}{\mddefault}{\updefault}1}}}}}
\put(4287,15387){\makebox(0,0)[lb]{\smash{{{\SetFigFont{6}{7.2}{\rmdefault}{\mddefault}{\updefault}1}}}}}
\put(162,12687){\makebox(0,0)[lb]{\smash{{{\SetFigFont{6}{7.2}{\rmdefault}{\mddefault}{\updefault}3}}}}}
\put(1062,12687){\makebox(0,0)[lb]{\smash{{{\SetFigFont{6}{7.2}{\rmdefault}{\mddefault}{\updefault}4}}}}}
\put(2787,12687){\makebox(0,0)[lb]{\smash{{{\SetFigFont{6}{7.2}{\rmdefault}{\mddefault}{\updefault}1}}}}}
\put(3387,12687){\makebox(0,0)[lb]{\smash{{{\SetFigFont{6}{7.2}{\rmdefault}{\mddefault}{\updefault}1}}}}}
\put(4287,12687){\makebox(0,0)[lb]{\smash{{{\SetFigFont{6}{7.2}{\rmdefault}{\mddefault}{\updefault}5}}}}}
\put(162,9987){\makebox(0,0)[lb]{\smash{{{\SetFigFont{6}{7.2}{\rmdefault}{\mddefault}{\updefault}12}}}}}
\put(1062,9987){\makebox(0,0)[lb]{\smash{{{\SetFigFont{6}{7.2}{\rmdefault}{\mddefault}{\updefault}8}}}}}
\put(1887,9987){\makebox(0,0)[lb]{\smash{{{\SetFigFont{6}{7.2}{\rmdefault}{\mddefault}{\updefault}7}}}}}
\put(2787,9987){\makebox(0,0)[lb]{\smash{{{\SetFigFont{6}{7.2}{\rmdefault}{\mddefault}{\updefault}1}}}}}
\put(3387,9987){\makebox(0,0)[lb]{\smash{{{\SetFigFont{6}{7.2}{\rmdefault}{\mddefault}{\updefault}1}}}}}
\put(4287,9987){\makebox(0,0)[lb]{\smash{{{\SetFigFont{6}{7.2}{\rmdefault}{\mddefault}{\updefault}9}}}}}
\put(4887,9987){\makebox(0,0)[lb]{\smash{{{\SetFigFont{6}{7.2}{\rmdefault}{\mddefault}{\updefault}2}}}}}
\put(5787,9987){\makebox(0,0)[lb]{\smash{{{\SetFigFont{6}{7.2}{\rmdefault}{\mddefault}{\updefault}1}}}}}
\put(6387,9987){\makebox(0,0)[lb]{\smash{{{\SetFigFont{6}{7.2}{\rmdefault}{\mddefault}{\updefault}1}}}}}
\put(162,7287){\makebox(0,0)[lb]{\smash{{{\SetFigFont{6}{7.2}{\rmdefault}{\mddefault}{\updefault}12}}}}}
\put(1062,7287){\makebox(0,0)[lb]{\smash{{{\SetFigFont{6}{7.2}{\rmdefault}{\mddefault}{\updefault}20}}}}}
\put(2787,7287){\makebox(0,0)[lb]{\smash{{{\SetFigFont{6}{7.2}{\rmdefault}{\mddefault}{\updefault}9}}}}}
\put(3387,7287){\makebox(0,0)[lb]{\smash{{{\SetFigFont{6}{7.2}{\rmdefault}{\mddefault}{\updefault}9}}}}}
\put(4287,7287){\makebox(0,0)[lb]{\smash{{{\SetFigFont{6}{7.2}{\rmdefault}{\mddefault}{\updefault}28}}}}}
\put(162,4587){\makebox(0,0)[lb]{\smash{{{\SetFigFont{6}{7.2}{\rmdefault}{\mddefault}{\updefault}60}}}}}
\put(1887,4587){\makebox(0,0)[lb]{\smash{{{\SetFigFont{6}{7.2}{\rmdefault}{\mddefault}{\updefault}46}}}}}
\put(162,2787){\makebox(0,0)[lb]{\smash{{{\SetFigFont{6}{7.2}{\rmdefault}{\mddefault}{\updefault}60}}}}}
\put(162,1287){\makebox(0,0)[lb]{\smash{{{\SetFigFont{6}{7.2}{\rmdefault}{\mddefault}{\updefault}358}}}}}
\put(162,87){\makebox(0,0)[lb]{\smash{{{\SetFigFont{6}{7.2}{\rmdefault}{\mddefault}{\updefault}358}}}}}
\put(7287,12987){\makebox(0,0)[lb]{\smash{{{\SetFigFont{6}{7.2}{\rmdefault}{\mddefault}{\updefault}1,2}}}}}
\put(7287,7587){\makebox(0,0)[lb]{\smash{{{\SetFigFont{6}{7.2}{\rmdefault}{\mddefault}{\updefault}1,2}}}}}
\put(7287,7287){\makebox(0,0)[lb]{\smash{{{\SetFigFont{6}{7.2}{\rmdefault}{\mddefault}{\updefault}21}}}}}
\put(7287,9987){\makebox(0,0)[lb]{\smash{{{\SetFigFont{6}{7.2}{\rmdefault}{\mddefault}{\updefault}2}}}}}
\put(7287,12687){\makebox(0,0)[lb]{\smash{{{\SetFigFont{6}{7.2}{\rmdefault}{\mddefault}{\updefault}2}}}}}
\put(4287,4587){\makebox(0,0)[lb]{\smash{{{\SetFigFont{6}{7.2}{\rmdefault}{\mddefault}{\updefault}73}}}}}
\put(1062,4587){\makebox(0,0)[lb]{\smash{{{\SetFigFont{6}{7.2}{\rmdefault}{\mddefault}{\updefault}59}}}}}
\put(5487,7437){\makebox(0,0)[lb]{\smash{{{\SetFigFont{6}{7.2}{\rmdefault}{\mddefault}{\updefault}OTHERS}}}}}
\put(6987,4737){\makebox(0,0)[lb]{\smash{{{\SetFigFont{6}{7.2}{\rmdefault}{\mddefault}{\updefault}OTHERS}}}}}
\put(1062,2787){\makebox(0,0)[lb]{\smash{{{\SetFigFont{6}{7.2}{\rmdefault}{\mddefault}{\updefault}119}}}}}
\put(4287,2787){\makebox(0,0)[lb]{\smash{{{\SetFigFont{6}{7.2}{\rmdefault}{\mddefault}{\updefault}179}}}}}
\put(7287,10287){\makebox(0,0)[lb]{\smash{{{\SetFigFont{6}{7.2}{\rmdefault}{\mddefault}{\updefault}2,3}}}}}
\end{picture}
}

%% file: pfn_mp19/gnuploTeX/q4n7n8l2sq310.tex
\setlength{\unitlength}{0.240900pt}
\ifx\plotpoint\undefined\newsavebox{\plotpoint}\fi
\sbox{\plotpoint}{\rule[-0.200pt]{0.400pt}{0.400pt}}%
\begin{picture}(3000,1800)(0,0)
\font\gnuplot=cmr10 at 10pt
\gnuplot
\sbox{\plotpoint}{\rule[-0.200pt]{0.400pt}{0.400pt}}%
\put(440.0,900.0){\rule[-0.200pt]{115.632pt}{0.400pt}}
\put(920.0,900.0){\rule[-0.200pt]{115.632pt}{0.400pt}}
\put(920.0,420.0){\rule[-0.200pt]{0.400pt}{115.632pt}}
\put(920.0,900.0){\rule[-0.200pt]{0.400pt}{115.632pt}}
\put(60.0,40.0){\rule[-0.200pt]{414.348pt}{0.400pt}}
\put(1780.0,40.0){\rule[-0.200pt]{0.400pt}{414.348pt}}
\put(60.0,1760.0){\rule[-0.200pt]{414.348pt}{0.400pt}}
\put(60.0,40.0){\rule[-0.200pt]{0.400pt}{414.348pt}}
\put(1181,1210){\makebox(0,0){$+$}}
\put(1181,590){\makebox(0,0){$+$}}
\put(1055,1179){\makebox(0,0){$+$}}
\put(1055,621){\makebox(0,0){$+$}}
\put(1059,1100){\makebox(0,0){$+$}}
\put(1059,700){\makebox(0,0){$+$}}
\put(1046,1112){\makebox(0,0){$+$}}
\put(1046,688){\makebox(0,0){$+$}}
\put(974,1218){\makebox(0,0){$+$}}
\put(974,582){\makebox(0,0){$+$}}
\put(986,1276){\makebox(0,0){$+$}}
\put(986,524){\makebox(0,0){$+$}}
\put(989,1365){\makebox(0,0){$+$}}
\put(989,435){\makebox(0,0){$+$}}
\put(1065,1037){\makebox(0,0){$+$}}
\put(1065,763){\makebox(0,0){$+$}}
\put(1052,1329){\makebox(0,0){$+$}}
\put(1052,471){\makebox(0,0){$+$}}
\put(1057,1068){\makebox(0,0){$+$}}
\put(1057,732){\makebox(0,0){$+$}}
\put(1091,1029){\makebox(0,0){$+$}}
\put(1091,771){\makebox(0,0){$+$}}
\put(1047,1210){\makebox(0,0){$+$}}
\put(1047,590){\makebox(0,0){$+$}}
\put(741,1035){\makebox(0,0){$+$}}
\put(741,765){\makebox(0,0){$+$}}
\put(792,1135){\makebox(0,0){$+$}}
\put(792,665){\makebox(0,0){$+$}}
\put(869,1219){\makebox(0,0){$+$}}
\put(869,581){\makebox(0,0){$+$}}
\put(754,1047){\makebox(0,0){$+$}}
\put(754,753){\makebox(0,0){$+$}}
\put(834,1192){\makebox(0,0){$+$}}
\put(834,608){\makebox(0,0){$+$}}
\put(892,1241){\makebox(0,0){$+$}}
\put(892,559){\makebox(0,0){$+$}}
\put(726,966){\makebox(0,0){$+$}}
\put(726,834){\makebox(0,0){$+$}}
\put(803,1151){\makebox(0,0){$+$}}
\put(803,649){\makebox(0,0){$+$}}
\put(884,1228){\makebox(0,0){$+$}}
\put(884,572){\makebox(0,0){$+$}}
\put(1071,1070){\makebox(0,0){$+$}}
\put(1071,730){\makebox(0,0){$+$}}
\put(1047,1245){\makebox(0,0){$+$}}
\put(1047,555){\makebox(0,0){$+$}}
\put(1080,1009){\makebox(0,0){$+$}}
\put(1080,791){\makebox(0,0){$+$}}
\put(1050,1145){\makebox(0,0){$+$}}
\put(1050,655){\makebox(0,0){$+$}}
\put(726,1001){\makebox(0,0){$+$}}
\put(726,799){\makebox(0,0){$+$}}
\put(764,1097){\makebox(0,0){$+$}}
\put(764,703){\makebox(0,0){$+$}}
\put(752,1067){\makebox(0,0){$+$}}
\put(752,733){\makebox(0,0){$+$}}
\put(753,1024){\makebox(0,0){$+$}}
\put(753,776){\makebox(0,0){$+$}}
\put(761,1082){\makebox(0,0){$+$}}
\put(761,718){\makebox(0,0){$+$}}
\put(853,1206){\makebox(0,0){$+$}}
\put(853,594){\makebox(0,0){$+$}}
\put(930,1264){\makebox(0,0){$+$}}
\put(930,536){\makebox(0,0){$+$}}
\put(946,1321){\makebox(0,0){$+$}}
\put(946,479){\makebox(0,0){$+$}}
\put(960,1236){\makebox(0,0){$+$}}
\put(960,564){\makebox(0,0){$+$}}
\put(1068,1127){\makebox(0,0){$+$}}
\put(1068,673){\makebox(0,0){$+$}}
\put(1091,1244){\makebox(0,0){$+$}}
\put(1091,556){\makebox(0,0){$+$}}
\put(1101,1322){\makebox(0,0){$+$}}
\put(1101,478){\makebox(0,0){$+$}}
\put(1102,1371){\makebox(0,0){$+$}}
\put(1102,429){\makebox(0,0){$+$}}
\put(1094,981){\makebox(0,0){$+$}}
\put(1094,819){\makebox(0,0){$+$}}
\put(1084,1045){\makebox(0,0){$+$}}
\put(1084,755){\makebox(0,0){$+$}}
\put(1060,1146){\makebox(0,0){$+$}}
\put(1060,654){\makebox(0,0){$+$}}
\put(512,946){\makebox(0,0){$+$}}
\put(512,854){\makebox(0,0){$+$}}
\put(569,1058){\makebox(0,0){$+$}}
\put(569,742){\makebox(0,0){$+$}}
\put(607,1147){\makebox(0,0){$+$}}
\put(607,653){\makebox(0,0){$+$}}
\put(548,998){\makebox(0,0){$+$}}
\put(548,802){\makebox(0,0){$+$}}
\put(598,1083){\makebox(0,0){$+$}}
\put(598,717){\makebox(0,0){$+$}}
\put(671,994){\makebox(0,0){$+$}}
\put(671,806){\makebox(0,0){$+$}}
\put(783,1124){\makebox(0,0){$+$}}
\put(783,676){\makebox(0,0){$+$}}
\put(1105,958){\makebox(0,0){$+$}}
\put(1105,842){\makebox(0,0){$+$}}
\put(1181,1159){\makebox(0,0){$+$}}
\put(1181,641){\makebox(0,0){$+$}}
\put(1182,1336){\makebox(0,0){$+$}}
\put(1182,464){\makebox(0,0){$+$}}
\put(1192,1150){\makebox(0,0){$+$}}
\put(1192,650){\makebox(0,0){$+$}}
\put(1194,1329){\makebox(0,0){$+$}}
\put(1194,471){\makebox(0,0){$+$}}
\put(1172,1117){\makebox(0,0){$+$}}
\put(1172,683){\makebox(0,0){$+$}}
\put(1168,1309){\makebox(0,0){$+$}}
\put(1168,491){\makebox(0,0){$+$}}
\put(1139,1295){\makebox(0,0){$+$}}
\put(1139,505){\makebox(0,0){$+$}}
\put(1158,1037){\makebox(0,0){$+$}}
\put(1158,763){\makebox(0,0){$+$}}
\put(1163,1074){\makebox(0,0){$+$}}
\put(1163,726){\makebox(0,0){$+$}}
\put(1188,1067){\makebox(0,0){$+$}}
\put(1188,733){\makebox(0,0){$+$}}
\put(1202,1263){\makebox(0,0){$+$}}
\put(1202,537){\makebox(0,0){$+$}}
\put(1207,1373){\makebox(0,0){$+$}}
\put(1207,427){\makebox(0,0){$+$}}
\put(1217,1196){\makebox(0,0){$+$}}
\put(1217,604){\makebox(0,0){$+$}}
\put(1187,1100){\makebox(0,0){$+$}}
\put(1187,700){\makebox(0,0){$+$}}
\put(534,989){\makebox(0,0){$+$}}
\put(534,811){\makebox(0,0){$+$}}
\put(561,1093){\makebox(0,0){$+$}}
\put(561,707){\makebox(0,0){$+$}}
\put(570,1089){\makebox(0,0){$+$}}
\put(570,711){\makebox(0,0){$+$}}
\put(588,1126){\makebox(0,0){$+$}}
\put(588,674){\makebox(0,0){$+$}}
\put(611,1101){\makebox(0,0){$+$}}
\put(611,699){\makebox(0,0){$+$}}
\put(514,974){\makebox(0,0){$+$}}
\put(514,826){\makebox(0,0){$+$}}
\put(527,1218){\makebox(0,0){$+$}}
\put(527,582){\makebox(0,0){$+$}}
\put(534,1024){\makebox(0,0){$+$}}
\put(534,776){\makebox(0,0){$+$}}
\put(545,1039){\makebox(0,0){$+$}}
\put(545,761){\makebox(0,0){$+$}}
\put(508,914){\makebox(0,0){$+$}}
\put(508,886){\makebox(0,0){$+$}}
\put(482,1173){\makebox(0,0){$+$}}
\put(482,627){\makebox(0,0){$+$}}
\put(671,991){\makebox(0,0){$+$}}
\put(671,809){\makebox(0,0){$+$}}
\put(801,1167){\makebox(0,0){$+$}}
\put(801,633){\makebox(0,0){$+$}}
\put(951,1286){\makebox(0,0){$+$}}
\put(951,514){\makebox(0,0){$+$}}
\put(1000,1315){\makebox(0,0){$+$}}
\put(1000,485){\makebox(0,0){$+$}}
\put(1028,1337){\makebox(0,0){$+$}}
\put(1028,463){\makebox(0,0){$+$}}
\put(1154,997){\makebox(0,0){$+$}}
\put(1154,803){\makebox(0,0){$+$}}
\put(1221,1209){\makebox(0,0){$+$}}
\put(1221,591){\makebox(0,0){$+$}}
\put(1250,1297){\makebox(0,0){$+$}}
\put(1250,503){\makebox(0,0){$+$}}
\put(1272,1297){\makebox(0,0){$+$}}
\put(1272,503){\makebox(0,0){$+$}}
\put(1262,1356){\makebox(0,0){$+$}}
\put(1262,444){\makebox(0,0){$+$}}
\put(1220,1260){\makebox(0,0){$+$}}
\put(1220,540){\makebox(0,0){$+$}}
\put(773,1113){\makebox(0,0){$+$}}
\put(773,687){\makebox(0,0){$+$}}
\put(904,1253){\makebox(0,0){$+$}}
\put(904,547){\makebox(0,0){$+$}}
\put(927,1321){\makebox(0,0){$+$}}
\put(927,479){\makebox(0,0){$+$}}
\put(959,1426){\makebox(0,0){$+$}}
\put(959,374){\makebox(0,0){$+$}}
\put(452,933){\makebox(0,0){$+$}}
\put(452,867){\makebox(0,0){$+$}}
\put(456,1366){\makebox(0,0){$+$}}
\put(456,434){\makebox(0,0){$+$}}
\put(486,1351){\makebox(0,0){$+$}}
\put(486,449){\makebox(0,0){$+$}}
\put(498,1379){\makebox(0,0){$+$}}
\put(498,421){\makebox(0,0){$+$}}
\put(434,1128){\makebox(0,0){$+$}}
\put(434,672){\makebox(0,0){$+$}}
\put(446,1346){\makebox(0,0){$+$}}
\put(446,454){\makebox(0,0){$+$}}
\put(498,1228){\makebox(0,0){$+$}}
\put(498,572){\makebox(0,0){$+$}}
\put(440,991){\makebox(0,0){$+$}}
\put(440,809){\makebox(0,0){$+$}}
\put(417,1098){\makebox(0,0){$+$}}
\put(417,702){\makebox(0,0){$+$}}
\put(433,1317){\makebox(0,0){$+$}}
\put(433,483){\makebox(0,0){$+$}}
\put(443,1070){\makebox(0,0){$+$}}
\put(443,730){\makebox(0,0){$+$}}
\put(445,891){\makebox(0,0){$+$}}
\put(445,909){\makebox(0,0){$+$}}
\put(1159,968){\makebox(0,0){$+$}}
\put(1159,832){\makebox(0,0){$+$}}
\put(1230,1173){\makebox(0,0){$+$}}
\put(1230,627){\makebox(0,0){$+$}}
\put(1323,1273){\makebox(0,0){$+$}}
\put(1323,527){\makebox(0,0){$+$}}
\put(1300,1320){\makebox(0,0){$+$}}
\put(1300,480){\makebox(0,0){$+$}}
\put(1282,1367){\makebox(0,0){$+$}}
\put(1282,433){\makebox(0,0){$+$}}
\put(1212,1431){\makebox(0,0){$+$}}
\put(1212,369){\makebox(0,0){$+$}}
\put(1188,1050){\makebox(0,0){$+$}}
\put(1188,750){\makebox(0,0){$+$}}
\put(1052,1393){\makebox(0,0){$+$}}
\put(1052,407){\makebox(0,0){$+$}}
\put(985,1441){\makebox(0,0){$+$}}
\put(985,359){\makebox(0,0){$+$}}
\put(923,1465){\makebox(0,0){$+$}}
\put(923,335){\makebox(0,0){$+$}}
\put(1041,1431){\makebox(0,0){$+$}}
\put(1041,369){\makebox(0,0){$+$}}
\put(458,953){\makebox(0,0){$+$}}
\put(458,847){\makebox(0,0){$+$}}
\put(415,1147){\makebox(0,0){$+$}}
\put(415,653){\makebox(0,0){$+$}}
\put(435,1212){\makebox(0,0){$+$}}
\put(435,588){\makebox(0,0){$+$}}
\put(415,1041){\makebox(0,0){$+$}}
\put(415,759){\makebox(0,0){$+$}}
\put(422,1174){\makebox(0,0){$+$}}
\put(422,626){\makebox(0,0){$+$}}
\put(1223,1153){\makebox(0,0){$+$}}
\put(1223,647){\makebox(0,0){$+$}}
\put(1374,1291){\makebox(0,0){$+$}}
\put(1374,509){\makebox(0,0){$+$}}
\put(1361,1312){\makebox(0,0){$+$}}
\put(1361,488){\makebox(0,0){$+$}}
\put(1347,1314){\makebox(0,0){$+$}}
\put(1347,486){\makebox(0,0){$+$}}
\put(1326,1311){\makebox(0,0){$+$}}
\put(1326,489){\makebox(0,0){$+$}}
\put(1459,1266){\makebox(0,0){$+$}}
\put(1459,534){\makebox(0,0){$+$}}
\put(1430,1278){\makebox(0,0){$+$}}
\put(1430,522){\makebox(0,0){$+$}}
\put(1414,1282){\makebox(0,0){$+$}}
\put(1414,518){\makebox(0,0){$+$}}
\put(1397,1286){\makebox(0,0){$+$}}
\put(1397,514){\makebox(0,0){$+$}}
\put(1529,1226){\makebox(0,0){$+$}}
\put(1529,574){\makebox(0,0){$+$}}
\put(1473,1259){\makebox(0,0){$+$}}
\put(1473,541){\makebox(0,0){$+$}}
\put(1445,1273){\makebox(0,0){$+$}}
\put(1445,527){\makebox(0,0){$+$}}
\put(818,1178){\makebox(0,0){$+$}}
\put(818,622){\makebox(0,0){$+$}}
\put(878,1501){\makebox(0,0){$+$}}
\put(878,299){\makebox(0,0){$+$}}
\put(913,1481){\makebox(0,0){$+$}}
\put(913,319){\makebox(0,0){$+$}}
\put(928,1492){\makebox(0,0){$+$}}
\put(928,308){\makebox(0,0){$+$}}
\put(824,1527){\makebox(0,0){$+$}}
\put(824,273){\makebox(0,0){$+$}}
\put(571,1453){\makebox(0,0){$+$}}
\put(571,347){\makebox(0,0){$+$}}
\put(421,1061){\makebox(0,0){$+$}}
\put(421,739){\makebox(0,0){$+$}}
\put(413,1249){\makebox(0,0){$+$}}
\put(413,551){\makebox(0,0){$+$}}
\put(426,1341){\makebox(0,0){$+$}}
\put(426,459){\makebox(0,0){$+$}}
\put(436,967){\makebox(0,0){$+$}}
\put(436,833){\makebox(0,0){$+$}}
\put(890,1509){\makebox(0,0){$+$}}
\put(890,291){\makebox(0,0){$+$}}
\put(1009,1504){\makebox(0,0){$+$}}
\put(1009,296){\makebox(0,0){$+$}}
\put(1051,1620){\makebox(0,0){$+$}}
\put(1051,180){\makebox(0,0){$+$}}
\put(1024,1490){\makebox(0,0){$+$}}
\put(1024,310){\makebox(0,0){$+$}}
\put(385,671){\makebox(0,0){$+$}}
\put(385,1129){\makebox(0,0){$+$}}
\put(391,1270){\makebox(0,0){$+$}}
\put(391,530){\makebox(0,0){$+$}}
\put(1607,1161){\makebox(0,0){$+$}}
\put(1607,639){\makebox(0,0){$+$}}
\put(1515,1234){\makebox(0,0){$+$}}
\put(1515,566){\makebox(0,0){$+$}}
\put(1486,1251){\makebox(0,0){$+$}}
\put(1486,549){\makebox(0,0){$+$}}
\put(1570,1197){\makebox(0,0){$+$}}
\put(1570,603){\makebox(0,0){$+$}}
\put(1500,1243){\makebox(0,0){$+$}}
\put(1500,557){\makebox(0,0){$+$}}
\put(1639,1120){\makebox(0,0){$+$}}
\put(1639,680){\makebox(0,0){$+$}}
\put(1543,1217){\makebox(0,0){$+$}}
\put(1543,583){\makebox(0,0){$+$}}
\put(359,1123){\makebox(0,0){$+$}}
\put(359,677){\makebox(0,0){$+$}}
\put(412,1290){\makebox(0,0){$+$}}
\put(412,510){\makebox(0,0){$+$}}
\put(573,1488){\makebox(0,0){$+$}}
\put(573,312){\makebox(0,0){$+$}}
\put(612,1527){\makebox(0,0){$+$}}
\put(612,273){\makebox(0,0){$+$}}
\put(634,1565){\makebox(0,0){$+$}}
\put(634,235){\makebox(0,0){$+$}}
\put(822,1531){\makebox(0,0){$+$}}
\put(822,269){\makebox(0,0){$+$}}
\put(796,1587){\makebox(0,0){$+$}}
\put(796,213){\makebox(0,0){$+$}}
\put(1148,1531){\makebox(0,0){$+$}}
\put(1148,269){\makebox(0,0){$+$}}
\put(1173,1433){\makebox(0,0){$+$}}
\put(1173,367){\makebox(0,0){$+$}}
\put(1168,1486){\makebox(0,0){$+$}}
\put(1168,314){\makebox(0,0){$+$}}
\put(558,1458){\makebox(0,0){$+$}}
\put(558,342){\makebox(0,0){$+$}}
\put(429,1023){\makebox(0,0){$+$}}
\put(429,777){\makebox(0,0){$+$}}
\put(695,1608){\makebox(0,0){$+$}}
\put(695,192){\makebox(0,0){$+$}}
\put(1678,739){\makebox(0,0){$+$}}
\put(1678,1061){\makebox(0,0){$+$}}
\put(1583,1185){\makebox(0,0){$+$}}
\put(1583,615){\makebox(0,0){$+$}}
\put(1618,1148){\makebox(0,0){$+$}}
\put(1618,652){\makebox(0,0){$+$}}
\put(1651,1105){\makebox(0,0){$+$}}
\put(1651,695){\makebox(0,0){$+$}}
\put(1557,1207){\makebox(0,0){$+$}}
\put(1557,593){\makebox(0,0){$+$}}
\put(514,1426){\makebox(0,0){$+$}}
\put(514,374){\makebox(0,0){$+$}}
\put(369,1249){\makebox(0,0){$+$}}
\put(369,551){\makebox(0,0){$+$}}
\put(309,1196){\makebox(0,0){$+$}}
\put(309,604){\makebox(0,0){$+$}}
\put(347,1227){\makebox(0,0){$+$}}
\put(347,573){\makebox(0,0){$+$}}
\put(672,1595){\makebox(0,0){$+$}}
\put(672,205){\makebox(0,0){$+$}}
\put(845,1543){\makebox(0,0){$+$}}
\put(845,257){\makebox(0,0){$+$}}
\put(1168,1554){\makebox(0,0){$+$}}
\put(1168,246){\makebox(0,0){$+$}}
\put(1187,1495){\makebox(0,0){$+$}}
\put(1187,305){\makebox(0,0){$+$}}
\put(964,1507){\makebox(0,0){$+$}}
\put(964,293){\makebox(0,0){$+$}}
\put(243,1163){\makebox(0,0){$+$}}
\put(243,637){\makebox(0,0){$+$}}
\put(273,1182){\makebox(0,0){$+$}}
\put(273,618){\makebox(0,0){$+$}}
\put(189,1081){\makebox(0,0){$+$}}
\put(189,719){\makebox(0,0){$+$}}
\put(239,1186){\makebox(0,0){$+$}}
\put(239,614){\makebox(0,0){$+$}}
\put(726,1621){\makebox(0,0){$+$}}
\put(726,179){\makebox(0,0){$+$}}
\put(795,1578){\makebox(0,0){$+$}}
\put(795,222){\makebox(0,0){$+$}}
\put(1122,1571){\makebox(0,0){$+$}}
\put(1122,229){\makebox(0,0){$+$}}
\put(1263,1466){\makebox(0,0){$+$}}
\put(1263,334){\makebox(0,0){$+$}}
\put(1695,1003){\makebox(0,0){$+$}}
\put(1695,797){\makebox(0,0){$+$}}
\put(1596,1173){\makebox(0,0){$+$}}
\put(1596,627){\makebox(0,0){$+$}}
\put(1688,1034){\makebox(0,0){$+$}}
\put(1688,766){\makebox(0,0){$+$}}
\put(1666,1088){\makebox(0,0){$+$}}
\put(1666,712){\makebox(0,0){$+$}}
\put(589,1501){\makebox(0,0){$+$}}
\put(589,299){\makebox(0,0){$+$}}
\put(487,1400){\makebox(0,0){$+$}}
\put(487,400){\makebox(0,0){$+$}}
\put(653,1579){\makebox(0,0){$+$}}
\put(653,221){\makebox(0,0){$+$}}
\put(1157,1586){\makebox(0,0){$+$}}
\put(1157,214){\makebox(0,0){$+$}}
\put(1253,1438){\makebox(0,0){$+$}}
\put(1253,362){\makebox(0,0){$+$}}
\put(188,1026){\makebox(0,0){$+$}}
\put(188,774){\makebox(0,0){$+$}}
\put(210,1181){\makebox(0,0){$+$}}
\put(210,619){\makebox(0,0){$+$}}
\put(1696,968){\makebox(0,0){$+$}}
\put(1696,832){\makebox(0,0){$+$}}
\put(1629,1134){\makebox(0,0){$+$}}
\put(1629,666){\makebox(0,0){$+$}}
\put(773,1635){\makebox(0,0){$+$}}
\put(773,165){\makebox(0,0){$+$}}
\put(619,1540){\makebox(0,0){$+$}}
\put(619,260){\makebox(0,0){$+$}}
\put(194,658){\makebox(0,0){$+$}}
\put(194,1142){\makebox(0,0){$+$}}
\put(1088,1602){\makebox(0,0){$+$}}
\put(1088,198){\makebox(0,0){$+$}}
\put(1314,1282){\circle*{12}}
\put(1314,518){\circle*{12}}
\put(1202,1274){\circle*{12}}
\put(1202,526){\circle*{12}}
\put(1071,1166){\circle*{12}}
\put(1071,634){\circle*{12}}
\put(1050,1314){\circle*{12}}
\put(1050,486){\circle*{12}}
\put(1018,1283){\circle*{12}}
\put(1018,517){\circle*{12}}
\put(979,1104){\circle*{12}}
\put(979,696){\circle*{12}}
\put(943,1271){\circle*{12}}
\put(943,529){\circle*{12}}
\put(876,1046){\circle*{12}}
\put(876,754){\circle*{12}}
\put(886,1232){\circle*{12}}
\put(886,568){\circle*{12}}
\put(906,1316){\circle*{12}}
\put(906,484){\circle*{12}}
\put(897,1075){\circle*{12}}
\put(897,725){\circle*{12}}
\put(880,933){\circle*{12}}
\put(880,867){\circle*{12}}
\put(891,1263){\circle*{12}}
\put(891,537){\circle*{12}}
\put(887,978){\circle*{12}}
\put(887,822){\circle*{12}}
\put(657,900){\circle*{12}}
\put(657,900){\circle*{12}}
\put(774,1082){\circle*{12}}
\put(774,718){\circle*{12}}
\put(860,1213){\circle*{12}}
\put(860,587){\circle*{12}}
\put(881,1316){\circle*{12}}
\put(881,484){\circle*{12}}
\put(1068,1112){\circle*{12}}
\put(1068,688){\circle*{12}}
\put(1051,1278){\circle*{12}}
\put(1051,522){\circle*{12}}
\put(1034,1177){\circle*{12}}
\put(1034,623){\circle*{12}}
\put(1010,1227){\circle*{12}}
\put(1010,573){\circle*{12}}
\put(996,1331){\circle*{12}}
\put(996,469){\circle*{12}}
\put(1094,1057){\circle*{12}}
\put(1094,743){\circle*{12}}
\put(1095,1375){\circle*{12}}
\put(1095,425){\circle*{12}}
\put(1092,1257){\circle*{12}}
\put(1092,543){\circle*{12}}
\put(1083,1200){\circle*{12}}
\put(1083,600){\circle*{12}}
\put(1078,1084){\circle*{12}}
\put(1078,716){\circle*{12}}
\put(1106,1033){\circle*{12}}
\put(1106,767){\circle*{12}}
\put(1113,1145){\circle*{12}}
\put(1113,655){\circle*{12}}
\put(1119,1122){\circle*{12}}
\put(1119,678){\circle*{12}}
\put(1113,1016){\circle*{12}}
\put(1113,784){\circle*{12}}
\put(605,924){\circle*{12}}
\put(605,876){\circle*{12}}
\put(641,1078){\circle*{12}}
\put(641,722){\circle*{12}}
\put(689,1103){\circle*{12}}
\put(689,697){\circle*{12}}
\put(604,1000){\circle*{12}}
\put(604,800){\circle*{12}}
\put(606,1052){\circle*{12}}
\put(606,748){\circle*{12}}
\put(671,1093){\circle*{12}}
\put(671,707){\circle*{12}}
\put(591,963){\circle*{12}}
\put(591,837){\circle*{12}}
\put(585,1102){\circle*{12}}
\put(585,698){\circle*{12}}
\put(525,962){\circle*{12}}
\put(525,838){\circle*{12}}
\put(523,1030){\circle*{12}}
\put(523,770){\circle*{12}}
\put(550,1090){\circle*{12}}
\put(550,710){\circle*{12}}
\put(600,1153){\circle*{12}}
\put(600,647){\circle*{12}}
\put(784,1107){\circle*{12}}
\put(784,693){\circle*{12}}
\put(867,1274){\circle*{12}}
\put(867,526){\circle*{12}}
\put(937,1336){\circle*{12}}
\put(937,464){\circle*{12}}
\put(966,1333){\circle*{12}}
\put(966,467){\circle*{12}}
\put(1189,1154){\circle*{12}}
\put(1189,646){\circle*{12}}
\put(1166,1254){\circle*{12}}
\put(1166,546){\circle*{12}}
\put(1150,1318){\circle*{12}}
\put(1150,482){\circle*{12}}
\put(1220,1151){\circle*{12}}
\put(1220,649){\circle*{12}}
\put(1210,1215){\circle*{12}}
\put(1210,585){\circle*{12}}
\put(1191,1319){\circle*{12}}
\put(1191,481){\circle*{12}}
\put(1183,1112){\circle*{12}}
\put(1183,688){\circle*{12}}
\put(1179,1038){\circle*{12}}
\put(1179,762){\circle*{12}}
\put(612,900){\circle*{12}}
\put(612,900){\circle*{12}}
\put(474,939){\circle*{12}}
\put(474,861){\circle*{12}}
\put(510,1193){\circle*{12}}
\put(510,607){\circle*{12}}
\put(526,1069){\circle*{12}}
\put(526,731){\circle*{12}}
\put(577,1137){\circle*{12}}
\put(577,663){\circle*{12}}
\put(696,1107){\circle*{12}}
\put(696,693){\circle*{12}}
\put(831,1193){\circle*{12}}
\put(831,607){\circle*{12}}
\put(867,638){\circle*{12}}
\put(867,1162){\circle*{12}}
\put(1229,1183){\circle*{12}}
\put(1229,617){\circle*{12}}
\put(1230,1313){\circle*{12}}
\put(1230,487){\circle*{12}}
\put(1276,1295){\circle*{12}}
\put(1276,505){\circle*{12}}
\put(1259,1370){\circle*{12}}
\put(1259,430){\circle*{12}}
\put(1178,1068){\circle*{12}}
\put(1178,732){\circle*{12}}
\put(1129,1324){\circle*{12}}
\put(1129,476){\circle*{12}}
\put(462,925){\circle*{12}}
\put(462,875){\circle*{12}}
\put(460,1344){\circle*{12}}
\put(460,456){\circle*{12}}
\put(468,1259){\circle*{12}}
\put(468,541){\circle*{12}}
\put(471,1150){\circle*{12}}
\put(471,650){\circle*{12}}
\put(421,1004){\circle*{12}}
\put(421,796){\circle*{12}}
\put(449,1185){\circle*{12}}
\put(449,615){\circle*{12}}
\put(454,1371){\circle*{12}}
\put(454,429){\circle*{12}}
\put(469,778){\circle*{12}}
\put(469,1022){\circle*{12}}
\put(1210,1132){\circle*{12}}
\put(1210,668){\circle*{12}}
\put(1230,1241){\circle*{12}}
\put(1230,559){\circle*{12}}
\put(1376,1296){\circle*{12}}
\put(1376,504){\circle*{12}}
\put(1340,1287){\circle*{12}}
\put(1340,513){\circle*{12}}
\put(1305,1351){\circle*{12}}
\put(1305,449){\circle*{12}}
\put(1266,1343){\circle*{12}}
\put(1266,457){\circle*{12}}
\put(1252,1384){\circle*{12}}
\put(1252,416){\circle*{12}}
\put(983,1406){\circle*{12}}
\put(983,394){\circle*{12}}
\put(931,1430){\circle*{12}}
\put(931,370){\circle*{12}}
\put(919,1441){\circle*{12}}
\put(919,359){\circle*{12}}
\put(799,1145){\circle*{12}}
\put(799,655){\circle*{12}}
\put(886,1471){\circle*{12}}
\put(886,329){\circle*{12}}
\put(920,900){\circle*{12}}
\put(920,900){\circle*{12}}
\put(1180,1020){\circle*{12}}
\put(1180,780){\circle*{12}}
\put(1049,1450){\circle*{12}}
\put(1049,350){\circle*{12}}
\put(1205,1436){\circle*{12}}
\put(1205,364){\circle*{12}}
\put(428,1066){\circle*{12}}
\put(428,734){\circle*{12}}
\put(424,1259){\circle*{12}}
\put(424,541){\circle*{12}}
\put(451,1210){\circle*{12}}
\put(451,590){\circle*{12}}
\put(432,1032){\circle*{12}}
\put(432,768){\circle*{12}}
\put(451,1097){\circle*{12}}
\put(451,703){\circle*{12}}
\put(1052,1468){\circle*{12}}
\put(1052,332){\circle*{12}}
\put(854,1504){\circle*{12}}
\put(854,296){\circle*{12}}
\put(1462,1261){\circle*{12}}
\put(1462,539){\circle*{12}}
\put(1431,1280){\circle*{12}}
\put(1431,520){\circle*{12}}
\put(1415,1286){\circle*{12}}
\put(1415,514){\circle*{12}}
\put(1397,1291){\circle*{12}}
\put(1397,509){\circle*{12}}
\put(1354,1298){\circle*{12}}
\put(1354,502){\circle*{12}}
\put(484,900){\circle*{12}}
\put(484,900){\circle*{12}}
\put(395,1141){\circle*{12}}
\put(395,659){\circle*{12}}
\put(437,1332){\circle*{12}}
\put(437,468){\circle*{12}}
\put(816,1520){\circle*{12}}
\put(816,280){\circle*{12}}
\put(834,1508){\circle*{12}}
\put(834,292){\circle*{12}}
\put(869,1493){\circle*{12}}
\put(869,307){\circle*{12}}
\put(447,900){\circle*{12}}
\put(447,900){\circle*{12}}
\put(214,616){\circle*{12}}
\put(214,1184){\circle*{12}}
\put(379,1271){\circle*{12}}
\put(379,529){\circle*{12}}
\put(1576,1190){\circle*{12}}
\put(1576,610){\circle*{12}}
\put(1517,1235){\circle*{12}}
\put(1517,565){\circle*{12}}
\put(1481,1251){\circle*{12}}
\put(1481,549){\circle*{12}}
\put(1447,1272){\circle*{12}}
\put(1447,528){\circle*{12}}
\put(1602,1164){\circle*{12}}
\put(1602,636){\circle*{12}}
\put(1532,1224){\circle*{12}}
\put(1532,576){\circle*{12}}
\put(443,968){\circle*{12}}
\put(443,832){\circle*{12}}
\put(502,1409){\circle*{12}}
\put(502,391){\circle*{12}}
\put(575,1494){\circle*{12}}
\put(575,306){\circle*{12}}
\put(730,1560){\circle*{12}}
\put(730,240){\circle*{12}}
\put(782,1555){\circle*{12}}
\put(782,245){\circle*{12}}
\put(1088,1455){\circle*{12}}
\put(1088,345){\circle*{12}}
\put(1060,1524){\circle*{12}}
\put(1060,276){\circle*{12}}
\put(981,1517){\circle*{12}}
\put(981,283){\circle*{12}}
\put(478,411){\circle*{12}}
\put(478,1389){\circle*{12}}
\put(409,1112){\circle*{12}}
\put(409,688){\circle*{12}}
\put(423,1302){\circle*{12}}
\put(423,498){\circle*{12}}
\put(1097,1498){\circle*{12}}
\put(1097,302){\circle*{12}}
\put(1184,1484){\circle*{12}}
\put(1184,316){\circle*{12}}
\put(1643,681){\circle*{12}}
\put(1643,1119){\circle*{12}}
\put(1547,1213){\circle*{12}}
\put(1547,587){\circle*{12}}
\put(1500,1244){\circle*{12}}
\put(1500,556){\circle*{12}}
\put(1683,1041){\circle*{12}}
\put(1683,759){\circle*{12}}
\put(1562,1202){\circle*{12}}
\put(1562,598){\circle*{12}}
\put(364,677){\circle*{12}}
\put(364,1123){\circle*{12}}
\put(387,1282){\circle*{12}}
\put(387,518){\circle*{12}}
\put(619,1534){\circle*{12}}
\put(619,266){\circle*{12}}
\put(715,1561){\circle*{12}}
\put(715,239){\circle*{12}}
\put(801,1540){\circle*{12}}
\put(801,260){\circle*{12}}
\put(568,1460){\circle*{12}}
\put(568,340){\circle*{12}}
\put(755,1562){\circle*{12}}
\put(755,238){\circle*{12}}
\put(1673,730){\circle*{12}}
\put(1673,1070){\circle*{12}}
\put(1589,1177){\circle*{12}}
\put(1589,623){\circle*{12}}
\put(265,618){\circle*{12}}
\put(265,1182){\circle*{12}}
\put(353,1201){\circle*{12}}
\put(353,599){\circle*{12}}
\put(1615,1150){\circle*{12}}
\put(1615,650){\circle*{12}}
\put(1051,1533){\circle*{12}}
\put(1051,267){\circle*{12}}
\put(1140,1587){\circle*{12}}
\put(1140,213){\circle*{12}}
\put(1233,1477){\circle*{12}}
\put(1233,323){\circle*{12}}
\put(719,1623){\circle*{12}}
\put(719,177){\circle*{12}}
\put(589,1517){\circle*{12}}
\put(589,283){\circle*{12}}
\put(241,1190){\circle*{12}}
\put(241,610){\circle*{12}}
\put(303,1189){\circle*{12}}
\put(303,611){\circle*{12}}
\put(1693,975){\circle*{12}}
\put(1693,825){\circle*{12}}
\put(1627,1137){\circle*{12}}
\put(1627,663){\circle*{12}}
\put(188,764){\circle*{12}}
\put(188,1036){\circle*{12}}
\put(192,1085){\circle*{12}}
\put(192,715){\circle*{12}}
\put(678,1593){\circle*{12}}
\put(678,207){\circle*{12}}
\put(981,1517){\circle*{12}}
\put(981,283){\circle*{12}}
\put(1250,1462){\circle*{12}}
\put(1250,338){\circle*{12}}
\put(1159,1541){\circle*{12}}
\put(1159,259){\circle*{12}}
\put(201,654){\circle*{12}}
\put(201,1146){\circle*{12}}
\put(536,1429){\circle*{12}}
\put(536,371){\circle*{12}}
\put(1691,790){\circle*{12}}
\put(1691,1010){\circle*{12}}
\put(1660,1097){\circle*{12}}
\put(1660,703){\circle*{12}}
\put(1124,1613){\circle*{12}}
\put(1124,187){\circle*{12}}
\put(642,1567){\circle*{12}}
\put(642,233){\circle*{12}}
\end{picture}

%% file: pfn_mp19/gnuploTeX/q4n8n10l3sq310.tex
\setlength{\unitlength}{0.240900pt}
\ifx\plotpoint\undefined\newsavebox{\plotpoint}\fi
\sbox{\plotpoint}{\rule[-0.200pt]{0.400pt}{0.400pt}}%
\begin{picture}(3000,1800)(0,0)
\font\gnuplot=cmr10 at 10pt
\gnuplot
\sbox{\plotpoint}{\rule[-0.200pt]{0.400pt}{0.400pt}}%
\put(440.0,900.0){\rule[-0.200pt]{115.632pt}{0.400pt}}
\put(920.0,900.0){\rule[-0.200pt]{115.632pt}{0.400pt}}
\put(920.0,420.0){\rule[-0.200pt]{0.400pt}{115.632pt}}
\put(920.0,900.0){\rule[-0.200pt]{0.400pt}{115.632pt}}
\put(60.0,40.0){\rule[-0.200pt]{414.348pt}{0.400pt}}
\put(1780.0,40.0){\rule[-0.200pt]{0.400pt}{414.348pt}}
\put(60.0,1760.0){\rule[-0.200pt]{414.348pt}{0.400pt}}
\put(60.0,40.0){\rule[-0.200pt]{0.400pt}{414.348pt}}
\put(1323,1296){\makebox(0,0){$+$}}
\put(1323,504){\makebox(0,0){$+$}}
\put(1241,1302){\makebox(0,0){$+$}}
\put(1241,498){\makebox(0,0){$+$}}
\put(1065,1131){\makebox(0,0){$+$}}
\put(1065,669){\makebox(0,0){$+$}}
\put(1054,1109){\makebox(0,0){$+$}}
\put(1054,691){\makebox(0,0){$+$}}
\put(1016,1259){\makebox(0,0){$+$}}
\put(1016,541){\makebox(0,0){$+$}}
\put(992,1260){\makebox(0,0){$+$}}
\put(992,540){\makebox(0,0){$+$}}
\put(965,1227){\makebox(0,0){$+$}}
\put(965,573){\makebox(0,0){$+$}}
\put(975,1280){\makebox(0,0){$+$}}
\put(975,520){\makebox(0,0){$+$}}
\put(1069,1047){\makebox(0,0){$+$}}
\put(1069,753){\makebox(0,0){$+$}}
\put(1047,1252){\makebox(0,0){$+$}}
\put(1047,548){\makebox(0,0){$+$}}
\put(1040,1345){\makebox(0,0){$+$}}
\put(1040,455){\makebox(0,0){$+$}}
\put(1065,1068){\makebox(0,0){$+$}}
\put(1065,732){\makebox(0,0){$+$}}
\put(1052,1139){\makebox(0,0){$+$}}
\put(1052,661){\makebox(0,0){$+$}}
\put(1078,1026){\makebox(0,0){$+$}}
\put(1078,774){\makebox(0,0){$+$}}
\put(1076,1376){\makebox(0,0){$+$}}
\put(1076,424){\makebox(0,0){$+$}}
\put(1068,1340){\makebox(0,0){$+$}}
\put(1068,460){\makebox(0,0){$+$}}
\put(1060,1103){\makebox(0,0){$+$}}
\put(1060,697){\makebox(0,0){$+$}}
\put(753,1020){\makebox(0,0){$+$}}
\put(753,780){\makebox(0,0){$+$}}
\put(834,1192){\makebox(0,0){$+$}}
\put(834,608){\makebox(0,0){$+$}}
\put(876,1223){\makebox(0,0){$+$}}
\put(876,577){\makebox(0,0){$+$}}
\put(760,1072){\makebox(0,0){$+$}}
\put(760,728){\makebox(0,0){$+$}}
\put(783,1120){\makebox(0,0){$+$}}
\put(783,680){\makebox(0,0){$+$}}
\put(857,1207){\makebox(0,0){$+$}}
\put(857,593){\makebox(0,0){$+$}}
\put(754,1030){\makebox(0,0){$+$}}
\put(754,770){\makebox(0,0){$+$}}
\put(797,1143){\makebox(0,0){$+$}}
\put(797,657){\makebox(0,0){$+$}}
\put(754,1048){\makebox(0,0){$+$}}
\put(754,752){\makebox(0,0){$+$}}
\put(767,1093){\makebox(0,0){$+$}}
\put(767,707){\makebox(0,0){$+$}}
\put(846,1200){\makebox(0,0){$+$}}
\put(846,600){\makebox(0,0){$+$}}
\put(732,1000){\makebox(0,0){$+$}}
\put(732,800){\makebox(0,0){$+$}}
\put(749,1043){\makebox(0,0){$+$}}
\put(749,757){\makebox(0,0){$+$}}
\put(789,1128){\makebox(0,0){$+$}}
\put(789,672){\makebox(0,0){$+$}}
\put(1092,1026){\makebox(0,0){$+$}}
\put(1092,774){\makebox(0,0){$+$}}
\put(1056,1219){\makebox(0,0){$+$}}
\put(1056,581){\makebox(0,0){$+$}}
\put(1063,1142){\makebox(0,0){$+$}}
\put(1063,658){\makebox(0,0){$+$}}
\put(1062,1157){\makebox(0,0){$+$}}
\put(1062,643){\makebox(0,0){$+$}}
\put(1068,1087){\makebox(0,0){$+$}}
\put(1068,713){\makebox(0,0){$+$}}
\put(1059,1180){\makebox(0,0){$+$}}
\put(1059,620){\makebox(0,0){$+$}}
\put(1089,1034){\makebox(0,0){$+$}}
\put(1089,766){\makebox(0,0){$+$}}
\put(1068,1243){\makebox(0,0){$+$}}
\put(1068,557){\makebox(0,0){$+$}}
\put(1088,1007){\makebox(0,0){$+$}}
\put(1088,793){\makebox(0,0){$+$}}
\put(1055,1200){\makebox(0,0){$+$}}
\put(1055,600){\makebox(0,0){$+$}}
\put(972,1231){\makebox(0,0){$+$}}
\put(972,569){\makebox(0,0){$+$}}
\put(741,1022){\makebox(0,0){$+$}}
\put(741,778){\makebox(0,0){$+$}}
\put(825,1184){\makebox(0,0){$+$}}
\put(825,616){\makebox(0,0){$+$}}
\put(893,1236){\makebox(0,0){$+$}}
\put(893,564){\makebox(0,0){$+$}}
\put(930,1292){\makebox(0,0){$+$}}
\put(930,508){\makebox(0,0){$+$}}
\put(942,1318){\makebox(0,0){$+$}}
\put(942,482){\makebox(0,0){$+$}}
\put(728,977){\makebox(0,0){$+$}}
\put(728,823){\makebox(0,0){$+$}}
\put(812,1163){\makebox(0,0){$+$}}
\put(812,637){\makebox(0,0){$+$}}
\put(884,1230){\makebox(0,0){$+$}}
\put(884,570){\makebox(0,0){$+$}}
\put(941,1282){\makebox(0,0){$+$}}
\put(941,518){\makebox(0,0){$+$}}
\put(1097,985){\makebox(0,0){$+$}}
\put(1097,815){\makebox(0,0){$+$}}
\put(1112,1258){\makebox(0,0){$+$}}
\put(1112,542){\makebox(0,0){$+$}}
\put(1106,1353){\makebox(0,0){$+$}}
\put(1106,447){\makebox(0,0){$+$}}
\put(1109,1312){\makebox(0,0){$+$}}
\put(1109,488){\makebox(0,0){$+$}}
\put(1084,1048){\makebox(0,0){$+$}}
\put(1084,752){\makebox(0,0){$+$}}
\put(1070,1122){\makebox(0,0){$+$}}
\put(1070,678){\makebox(0,0){$+$}}
\put(1103,967){\makebox(0,0){$+$}}
\put(1103,833){\makebox(0,0){$+$}}
\put(1166,1056){\makebox(0,0){$+$}}
\put(1166,744){\makebox(0,0){$+$}}
\put(1162,1311){\makebox(0,0){$+$}}
\put(1162,489){\makebox(0,0){$+$}}
\put(1177,1112){\makebox(0,0){$+$}}
\put(1177,688){\makebox(0,0){$+$}}
\put(1174,1339){\makebox(0,0){$+$}}
\put(1174,461){\makebox(0,0){$+$}}
\put(1184,1145){\makebox(0,0){$+$}}
\put(1184,655){\makebox(0,0){$+$}}
\put(1188,1176){\makebox(0,0){$+$}}
\put(1188,624){\makebox(0,0){$+$}}
\put(1187,1360){\makebox(0,0){$+$}}
\put(1187,440){\makebox(0,0){$+$}}
\put(541,999){\makebox(0,0){$+$}}
\put(541,801){\makebox(0,0){$+$}}
\put(585,1128){\makebox(0,0){$+$}}
\put(585,672){\makebox(0,0){$+$}}
\put(600,1120){\makebox(0,0){$+$}}
\put(600,680){\makebox(0,0){$+$}}
\put(608,1090){\makebox(0,0){$+$}}
\put(608,710){\makebox(0,0){$+$}}
\put(576,1085){\makebox(0,0){$+$}}
\put(576,715){\makebox(0,0){$+$}}
\put(541,1002){\makebox(0,0){$+$}}
\put(541,798){\makebox(0,0){$+$}}
\put(580,1092){\makebox(0,0){$+$}}
\put(580,708){\makebox(0,0){$+$}}
\put(598,1092){\makebox(0,0){$+$}}
\put(598,708){\makebox(0,0){$+$}}
\put(671,993){\makebox(0,0){$+$}}
\put(671,807){\makebox(0,0){$+$}}
\put(771,1103){\makebox(0,0){$+$}}
\put(771,697){\makebox(0,0){$+$}}
\put(761,1084){\makebox(0,0){$+$}}
\put(761,716){\makebox(0,0){$+$}}
\put(793,1137){\makebox(0,0){$+$}}
\put(793,663){\makebox(0,0){$+$}}
\put(866,1215){\makebox(0,0){$+$}}
\put(866,585){\makebox(0,0){$+$}}
\put(929,1265){\makebox(0,0){$+$}}
\put(929,535){\makebox(0,0){$+$}}
\put(1161,962){\makebox(0,0){$+$}}
\put(1161,838){\makebox(0,0){$+$}}
\put(1145,1284){\makebox(0,0){$+$}}
\put(1145,516){\makebox(0,0){$+$}}
\put(1130,1319){\makebox(0,0){$+$}}
\put(1130,481){\makebox(0,0){$+$}}
\put(1171,1082){\makebox(0,0){$+$}}
\put(1171,718){\makebox(0,0){$+$}}
\put(1189,1214){\makebox(0,0){$+$}}
\put(1189,586){\makebox(0,0){$+$}}
\put(1186,1324){\makebox(0,0){$+$}}
\put(1186,476){\makebox(0,0){$+$}}
\put(1197,1158){\makebox(0,0){$+$}}
\put(1197,642){\makebox(0,0){$+$}}
\put(1228,1233){\makebox(0,0){$+$}}
\put(1228,567){\makebox(0,0){$+$}}
\put(1216,1194){\makebox(0,0){$+$}}
\put(1216,606){\makebox(0,0){$+$}}
\put(1203,1293){\makebox(0,0){$+$}}
\put(1203,507){\makebox(0,0){$+$}}
\put(1162,1030){\makebox(0,0){$+$}}
\put(1162,770){\makebox(0,0){$+$}}
\put(729,960){\makebox(0,0){$+$}}
\put(729,840){\makebox(0,0){$+$}}
\put(776,1112){\makebox(0,0){$+$}}
\put(776,688){\makebox(0,0){$+$}}
\put(818,1175){\makebox(0,0){$+$}}
\put(818,625){\makebox(0,0){$+$}}
\put(903,1244){\makebox(0,0){$+$}}
\put(903,556){\makebox(0,0){$+$}}
\put(934,1319){\makebox(0,0){$+$}}
\put(934,481){\makebox(0,0){$+$}}
\put(989,1303){\makebox(0,0){$+$}}
\put(989,497){\makebox(0,0){$+$}}
\put(1003,1322){\makebox(0,0){$+$}}
\put(1003,478){\makebox(0,0){$+$}}
\put(999,1379){\makebox(0,0){$+$}}
\put(999,421){\makebox(0,0){$+$}}
\put(490,934){\makebox(0,0){$+$}}
\put(490,866){\makebox(0,0){$+$}}
\put(505,1175){\makebox(0,0){$+$}}
\put(505,625){\makebox(0,0){$+$}}
\put(550,1068){\makebox(0,0){$+$}}
\put(550,732){\makebox(0,0){$+$}}
\put(555,1099){\makebox(0,0){$+$}}
\put(555,701){\makebox(0,0){$+$}}
\put(565,1177){\makebox(0,0){$+$}}
\put(565,623){\makebox(0,0){$+$}}
\put(536,1056){\makebox(0,0){$+$}}
\put(536,744){\makebox(0,0){$+$}}
\put(541,1162){\makebox(0,0){$+$}}
\put(541,638){\makebox(0,0){$+$}}
\put(525,1039){\makebox(0,0){$+$}}
\put(525,761){\makebox(0,0){$+$}}
\put(513,1013){\makebox(0,0){$+$}}
\put(513,787){\makebox(0,0){$+$}}
\put(504,990){\makebox(0,0){$+$}}
\put(504,810){\makebox(0,0){$+$}}
\put(496,973){\makebox(0,0){$+$}}
\put(496,827){\makebox(0,0){$+$}}
\put(488,911){\makebox(0,0){$+$}}
\put(488,889){\makebox(0,0){$+$}}
\put(441,1012){\makebox(0,0){$+$}}
\put(441,788){\makebox(0,0){$+$}}
\put(460,1356){\makebox(0,0){$+$}}
\put(460,444){\makebox(0,0){$+$}}
\put(468,1367){\makebox(0,0){$+$}}
\put(468,433){\makebox(0,0){$+$}}
\put(492,954){\makebox(0,0){$+$}}
\put(492,846){\makebox(0,0){$+$}}
\put(566,1088){\makebox(0,0){$+$}}
\put(566,712){\makebox(0,0){$+$}}
\put(599,1147){\makebox(0,0){$+$}}
\put(599,653){\makebox(0,0){$+$}}
\put(1109,954){\makebox(0,0){$+$}}
\put(1109,846){\makebox(0,0){$+$}}
\put(1191,1128){\makebox(0,0){$+$}}
\put(1191,672){\makebox(0,0){$+$}}
\put(1211,1263){\makebox(0,0){$+$}}
\put(1211,537){\makebox(0,0){$+$}}
\put(1188,707){\makebox(0,0){$+$}}
\put(1188,1093){\makebox(0,0){$+$}}
\put(1214,1242){\makebox(0,0){$+$}}
\put(1214,558){\makebox(0,0){$+$}}
\put(1255,1323){\makebox(0,0){$+$}}
\put(1255,477){\makebox(0,0){$+$}}
\put(1256,1374){\makebox(0,0){$+$}}
\put(1256,426){\makebox(0,0){$+$}}
\put(1219,1154){\makebox(0,0){$+$}}
\put(1219,646){\makebox(0,0){$+$}}
\put(1224,1199){\makebox(0,0){$+$}}
\put(1224,601){\makebox(0,0){$+$}}
\put(1227,1156){\makebox(0,0){$+$}}
\put(1227,644){\makebox(0,0){$+$}}
\put(1188,1048){\makebox(0,0){$+$}}
\put(1188,752){\makebox(0,0){$+$}}
\put(1076,1066){\makebox(0,0){$+$}}
\put(1076,734){\makebox(0,0){$+$}}
\put(1016,1328){\makebox(0,0){$+$}}
\put(1016,472){\makebox(0,0){$+$}}
\put(970,1423){\makebox(0,0){$+$}}
\put(970,377){\makebox(0,0){$+$}}
\put(755,1063){\makebox(0,0){$+$}}
\put(755,737){\makebox(0,0){$+$}}
\put(804,1151){\makebox(0,0){$+$}}
\put(804,649){\makebox(0,0){$+$}}
\put(915,1253){\makebox(0,0){$+$}}
\put(915,547){\makebox(0,0){$+$}}
\put(1039,1401){\makebox(0,0){$+$}}
\put(1039,399){\makebox(0,0){$+$}}
\put(444,906){\makebox(0,0){$+$}}
\put(444,894){\makebox(0,0){$+$}}
\put(451,1187){\makebox(0,0){$+$}}
\put(451,613){\makebox(0,0){$+$}}
\put(465,1233){\makebox(0,0){$+$}}
\put(465,567){\makebox(0,0){$+$}}
\put(433,1112){\makebox(0,0){$+$}}
\put(433,688){\makebox(0,0){$+$}}
\put(450,1344){\makebox(0,0){$+$}}
\put(450,456){\makebox(0,0){$+$}}
\put(444,1153){\makebox(0,0){$+$}}
\put(444,647){\makebox(0,0){$+$}}
\put(432,1104){\makebox(0,0){$+$}}
\put(432,696){\makebox(0,0){$+$}}
\put(439,1220){\makebox(0,0){$+$}}
\put(439,580){\makebox(0,0){$+$}}
\put(437,1052){\makebox(0,0){$+$}}
\put(437,748){\makebox(0,0){$+$}}
\put(437,985){\makebox(0,0){$+$}}
\put(437,815){\makebox(0,0){$+$}}
\put(425,1034){\makebox(0,0){$+$}}
\put(425,766){\makebox(0,0){$+$}}
\put(425,1184){\makebox(0,0){$+$}}
\put(425,616){\makebox(0,0){$+$}}
\put(457,935){\makebox(0,0){$+$}}
\put(457,865){\makebox(0,0){$+$}}
\put(1159,978){\makebox(0,0){$+$}}
\put(1159,822){\makebox(0,0){$+$}}
\put(1277,1287){\makebox(0,0){$+$}}
\put(1277,513){\makebox(0,0){$+$}}
\put(1252,1341){\makebox(0,0){$+$}}
\put(1252,459){\makebox(0,0){$+$}}
\put(1238,1339){\makebox(0,0){$+$}}
\put(1238,461){\makebox(0,0){$+$}}
\put(1323,1267){\makebox(0,0){$+$}}
\put(1323,533){\makebox(0,0){$+$}}
\put(1304,1312){\makebox(0,0){$+$}}
\put(1304,488){\makebox(0,0){$+$}}
\put(1281,1348){\makebox(0,0){$+$}}
\put(1281,452){\makebox(0,0){$+$}}
\put(1231,1174){\makebox(0,0){$+$}}
\put(1231,626){\makebox(0,0){$+$}}
\put(1188,1071){\makebox(0,0){$+$}}
\put(1188,729){\makebox(0,0){$+$}}
\put(990,1434){\makebox(0,0){$+$}}
\put(990,366){\makebox(0,0){$+$}}
\put(973,1450){\makebox(0,0){$+$}}
\put(973,350){\makebox(0,0){$+$}}
\put(938,1449){\makebox(0,0){$+$}}
\put(938,351){\makebox(0,0){$+$}}
\put(1158,1002){\makebox(0,0){$+$}}
\put(1158,798){\makebox(0,0){$+$}}
\put(1188,1056){\makebox(0,0){$+$}}
\put(1188,744){\makebox(0,0){$+$}}
\put(917,1474){\makebox(0,0){$+$}}
\put(917,326){\makebox(0,0){$+$}}
\put(1050,1409){\makebox(0,0){$+$}}
\put(1050,391){\makebox(0,0){$+$}}
\put(671,992){\makebox(0,0){$+$}}
\put(671,808){\makebox(0,0){$+$}}
\put(444,962){\makebox(0,0){$+$}}
\put(444,838){\makebox(0,0){$+$}}
\put(420,1074){\makebox(0,0){$+$}}
\put(420,726){\makebox(0,0){$+$}}
\put(428,1300){\makebox(0,0){$+$}}
\put(428,500){\makebox(0,0){$+$}}
\put(438,1239){\makebox(0,0){$+$}}
\put(438,561){\makebox(0,0){$+$}}
\put(442,735){\makebox(0,0){$+$}}
\put(442,1065){\makebox(0,0){$+$}}
\put(1399,1282){\makebox(0,0){$+$}}
\put(1399,518){\makebox(0,0){$+$}}
\put(1360,1306){\makebox(0,0){$+$}}
\put(1360,494){\makebox(0,0){$+$}}
\put(1341,1309){\makebox(0,0){$+$}}
\put(1341,491){\makebox(0,0){$+$}}
\put(1327,1311){\makebox(0,0){$+$}}
\put(1327,489){\makebox(0,0){$+$}}
\put(1422,1276){\makebox(0,0){$+$}}
\put(1422,524){\makebox(0,0){$+$}}
\put(1387,1285){\makebox(0,0){$+$}}
\put(1387,515){\makebox(0,0){$+$}}
\put(1352,1311){\makebox(0,0){$+$}}
\put(1352,489){\makebox(0,0){$+$}}
\put(1481,1249){\makebox(0,0){$+$}}
\put(1481,551){\makebox(0,0){$+$}}
\put(1443,1269){\makebox(0,0){$+$}}
\put(1443,531){\makebox(0,0){$+$}}
\put(1411,1279){\makebox(0,0){$+$}}
\put(1411,521){\makebox(0,0){$+$}}
\put(1370,1287){\makebox(0,0){$+$}}
\put(1370,513){\makebox(0,0){$+$}}
\put(906,1485){\makebox(0,0){$+$}}
\put(906,315){\makebox(0,0){$+$}}
\put(886,1493){\makebox(0,0){$+$}}
\put(886,307){\makebox(0,0){$+$}}
\put(867,1503){\makebox(0,0){$+$}}
\put(867,297){\makebox(0,0){$+$}}
\put(624,1554){\makebox(0,0){$+$}}
\put(624,246){\makebox(0,0){$+$}}
\put(421,1150){\makebox(0,0){$+$}}
\put(421,650){\makebox(0,0){$+$}}
\put(435,1334){\makebox(0,0){$+$}}
\put(435,466){\makebox(0,0){$+$}}
\put(405,1111){\makebox(0,0){$+$}}
\put(405,689){\makebox(0,0){$+$}}
\put(437,1342){\makebox(0,0){$+$}}
\put(437,458){\makebox(0,0){$+$}}
\put(505,1385){\makebox(0,0){$+$}}
\put(505,415){\makebox(0,0){$+$}}
\put(535,1433){\makebox(0,0){$+$}}
\put(535,367){\makebox(0,0){$+$}}
\put(1021,1433){\makebox(0,0){$+$}}
\put(1021,367){\makebox(0,0){$+$}}
\put(943,1498){\makebox(0,0){$+$}}
\put(943,302){\makebox(0,0){$+$}}
\put(1182,1422){\makebox(0,0){$+$}}
\put(1182,378){\makebox(0,0){$+$}}
\put(1182,1456){\makebox(0,0){$+$}}
\put(1182,344){\makebox(0,0){$+$}}
\put(449,924){\makebox(0,0){$+$}}
\put(449,876){\makebox(0,0){$+$}}
\put(575,1495){\makebox(0,0){$+$}}
\put(575,305){\makebox(0,0){$+$}}
\put(1574,1185){\makebox(0,0){$+$}}
\put(1574,615){\makebox(0,0){$+$}}
\put(1499,1237){\makebox(0,0){$+$}}
\put(1499,563){\makebox(0,0){$+$}}
\put(1462,1259){\makebox(0,0){$+$}}
\put(1462,541){\makebox(0,0){$+$}}
\put(1433,1273){\makebox(0,0){$+$}}
\put(1433,527){\makebox(0,0){$+$}}
\put(1538,1214){\makebox(0,0){$+$}}
\put(1538,586){\makebox(0,0){$+$}}
\put(1509,1231){\makebox(0,0){$+$}}
\put(1509,569){\makebox(0,0){$+$}}
\put(1583,1178){\makebox(0,0){$+$}}
\put(1583,622){\makebox(0,0){$+$}}
\put(1490,1243){\makebox(0,0){$+$}}
\put(1490,557){\makebox(0,0){$+$}}
\put(1452,1264){\makebox(0,0){$+$}}
\put(1452,536){\makebox(0,0){$+$}}
\put(390,1136){\makebox(0,0){$+$}}
\put(390,664){\makebox(0,0){$+$}}
\put(437,1319){\makebox(0,0){$+$}}
\put(437,481){\makebox(0,0){$+$}}
\put(525,1418){\makebox(0,0){$+$}}
\put(525,382){\makebox(0,0){$+$}}
\put(924,1478){\makebox(0,0){$+$}}
\put(924,322){\makebox(0,0){$+$}}
\put(825,1528){\makebox(0,0){$+$}}
\put(825,272){\makebox(0,0){$+$}}
\put(439,944){\makebox(0,0){$+$}}
\put(439,856){\makebox(0,0){$+$}}
\put(617,1543){\makebox(0,0){$+$}}
\put(617,257){\makebox(0,0){$+$}}
\put(412,1248){\makebox(0,0){$+$}}
\put(412,552){\makebox(0,0){$+$}}
\put(412,1043){\makebox(0,0){$+$}}
\put(412,757){\makebox(0,0){$+$}}
\put(876,1511){\makebox(0,0){$+$}}
\put(876,289){\makebox(0,0){$+$}}
\put(1011,1496){\makebox(0,0){$+$}}
\put(1011,304){\makebox(0,0){$+$}}
\put(398,1266){\makebox(0,0){$+$}}
\put(398,534){\makebox(0,0){$+$}}
\put(1183,1475){\makebox(0,0){$+$}}
\put(1183,325){\makebox(0,0){$+$}}
\put(1211,1410){\makebox(0,0){$+$}}
\put(1211,390){\makebox(0,0){$+$}}
\put(1159,1530){\makebox(0,0){$+$}}
\put(1159,270){\makebox(0,0){$+$}}
\put(823,1529){\makebox(0,0){$+$}}
\put(823,271){\makebox(0,0){$+$}}
\put(492,1396){\makebox(0,0){$+$}}
\put(492,404){\makebox(0,0){$+$}}
\put(484,1380){\makebox(0,0){$+$}}
\put(484,420){\makebox(0,0){$+$}}
\put(966,1490){\makebox(0,0){$+$}}
\put(966,310){\makebox(0,0){$+$}}
\put(1622,1134){\makebox(0,0){$+$}}
\put(1622,666){\makebox(0,0){$+$}}
\put(1529,1220){\makebox(0,0){$+$}}
\put(1529,580){\makebox(0,0){$+$}}
\put(1471,1254){\makebox(0,0){$+$}}
\put(1471,546){\makebox(0,0){$+$}}
\put(1629,1124){\makebox(0,0){$+$}}
\put(1629,676){\makebox(0,0){$+$}}
\put(1547,1207){\makebox(0,0){$+$}}
\put(1547,593){\makebox(0,0){$+$}}
\put(1599,1161){\makebox(0,0){$+$}}
\put(1599,639){\makebox(0,0){$+$}}
\put(1519,1226){\makebox(0,0){$+$}}
\put(1519,574){\makebox(0,0){$+$}}
\put(434,999){\makebox(0,0){$+$}}
\put(434,801){\makebox(0,0){$+$}}
\put(327,1192){\makebox(0,0){$+$}}
\put(327,608){\makebox(0,0){$+$}}
\put(356,1222){\makebox(0,0){$+$}}
\put(356,578){\makebox(0,0){$+$}}
\put(416,1281){\makebox(0,0){$+$}}
\put(416,519){\makebox(0,0){$+$}}
\put(598,1512){\makebox(0,0){$+$}}
\put(598,288){\makebox(0,0){$+$}}
\put(585,1501){\makebox(0,0){$+$}}
\put(585,299){\makebox(0,0){$+$}}
\put(705,1617){\makebox(0,0){$+$}}
\put(705,183){\makebox(0,0){$+$}}
\put(820,1549){\makebox(0,0){$+$}}
\put(820,251){\makebox(0,0){$+$}}
\put(795,1582){\makebox(0,0){$+$}}
\put(795,218){\makebox(0,0){$+$}}
\put(1652,707){\makebox(0,0){$+$}}
\put(1652,1093){\makebox(0,0){$+$}}
\put(1557,1200){\makebox(0,0){$+$}}
\put(1557,600){\makebox(0,0){$+$}}
\put(1678,1044){\makebox(0,0){$+$}}
\put(1678,756){\makebox(0,0){$+$}}
\put(1566,1193){\makebox(0,0){$+$}}
\put(1566,607){\makebox(0,0){$+$}}
\put(375,1132){\makebox(0,0){$+$}}
\put(375,668){\makebox(0,0){$+$}}
\put(276,1176){\makebox(0,0){$+$}}
\put(276,624){\makebox(0,0){$+$}}
\put(367,1236){\makebox(0,0){$+$}}
\put(367,564){\makebox(0,0){$+$}}
\put(797,1581){\makebox(0,0){$+$}}
\put(797,219){\makebox(0,0){$+$}}
\put(850,1536){\makebox(0,0){$+$}}
\put(850,264){\makebox(0,0){$+$}}
\put(1168,1552){\makebox(0,0){$+$}}
\put(1168,248){\makebox(0,0){$+$}}
\put(1176,1497){\makebox(0,0){$+$}}
\put(1176,303){\makebox(0,0){$+$}}
\put(1181,1513){\makebox(0,0){$+$}}
\put(1181,287){\makebox(0,0){$+$}}
\put(348,1120){\makebox(0,0){$+$}}
\put(348,680){\makebox(0,0){$+$}}
\put(258,1169){\makebox(0,0){$+$}}
\put(258,631){\makebox(0,0){$+$}}
\put(385,1253){\makebox(0,0){$+$}}
\put(385,547){\makebox(0,0){$+$}}
\put(558,1465){\makebox(0,0){$+$}}
\put(558,335){\makebox(0,0){$+$}}
\put(647,1579){\makebox(0,0){$+$}}
\put(647,221){\makebox(0,0){$+$}}
\put(673,1600){\makebox(0,0){$+$}}
\put(673,200){\makebox(0,0){$+$}}
\put(1241,1412){\makebox(0,0){$+$}}
\put(1241,388){\makebox(0,0){$+$}}
\put(1257,1440){\makebox(0,0){$+$}}
\put(1257,360){\makebox(0,0){$+$}}
\put(1120,1583){\makebox(0,0){$+$}}
\put(1120,217){\makebox(0,0){$+$}}
\put(1015,1492){\makebox(0,0){$+$}}
\put(1015,308){\makebox(0,0){$+$}}
\put(1643,1105){\makebox(0,0){$+$}}
\put(1643,695){\makebox(0,0){$+$}}
\put(1591,1169){\makebox(0,0){$+$}}
\put(1591,631){\makebox(0,0){$+$}}
\put(1670,1062){\makebox(0,0){$+$}}
\put(1670,738){\makebox(0,0){$+$}}
\put(1615,1143){\makebox(0,0){$+$}}
\put(1615,657){\makebox(0,0){$+$}}
\put(1692,959){\makebox(0,0){$+$}}
\put(1692,841){\makebox(0,0){$+$}}
\put(1684,1025){\makebox(0,0){$+$}}
\put(1684,775){\makebox(0,0){$+$}}
\put(240,1157){\makebox(0,0){$+$}}
\put(240,643){\makebox(0,0){$+$}}
\put(299,1185){\makebox(0,0){$+$}}
\put(299,615){\makebox(0,0){$+$}}
\put(198,1086){\makebox(0,0){$+$}}
\put(198,714){\makebox(0,0){$+$}}
\put(235,1174){\makebox(0,0){$+$}}
\put(235,626){\makebox(0,0){$+$}}
\put(605,1528){\makebox(0,0){$+$}}
\put(605,272){\makebox(0,0){$+$}}
\put(566,1483){\makebox(0,0){$+$}}
\put(566,317){\makebox(0,0){$+$}}
\put(555,1448){\makebox(0,0){$+$}}
\put(555,352){\makebox(0,0){$+$}}
\put(727,1626){\makebox(0,0){$+$}}
\put(727,174){\makebox(0,0){$+$}}
\put(1093,1603){\makebox(0,0){$+$}}
\put(1093,197){\makebox(0,0){$+$}}
\put(1140,1558){\makebox(0,0){$+$}}
\put(1140,242){\makebox(0,0){$+$}}
\put(1155,1591){\makebox(0,0){$+$}}
\put(1155,209){\makebox(0,0){$+$}}
\put(1066,1616){\makebox(0,0){$+$}}
\put(1066,184){\makebox(0,0){$+$}}
\put(192,755){\makebox(0,0){$+$}}
\put(192,1045){\makebox(0,0){$+$}}
\put(209,1158){\makebox(0,0){$+$}}
\put(209,642){\makebox(0,0){$+$}}
\put(688,192){\makebox(0,0){$+$}}
\put(688,1608){\makebox(0,0){$+$}}
\put(635,1568){\makebox(0,0){$+$}}
\put(635,232){\makebox(0,0){$+$}}
\put(1689,1004){\makebox(0,0){$+$}}
\put(1689,796){\makebox(0,0){$+$}}
\put(1607,1152){\makebox(0,0){$+$}}
\put(1607,648){\makebox(0,0){$+$}}
\put(1636,1114){\makebox(0,0){$+$}}
\put(1636,686){\makebox(0,0){$+$}}
\put(1265,1473){\makebox(0,0){$+$}}
\put(1265,327){\makebox(0,0){$+$}}
\put(784,1646){\makebox(0,0){$+$}}
\put(784,154){\makebox(0,0){$+$}}
\put(191,1015){\makebox(0,0){$+$}}
\put(191,785){\makebox(0,0){$+$}}
\put(221,1178){\makebox(0,0){$+$}}
\put(221,622){\makebox(0,0){$+$}}
\put(1662,721){\makebox(0,0){$+$}}
\put(1662,1079){\makebox(0,0){$+$}}
\put(1160,1577){\makebox(0,0){$+$}}
\put(1160,223){\makebox(0,0){$+$}}
\put(1043,1624){\makebox(0,0){$+$}}
\put(1043,176){\makebox(0,0){$+$}}
\put(506,1412){\makebox(0,0){$+$}}
\put(506,388){\makebox(0,0){$+$}}
\put(756,1634){\makebox(0,0){$+$}}
\put(756,166){\makebox(0,0){$+$}}
\put(660,1590){\makebox(0,0){$+$}}
\put(660,210){\makebox(0,0){$+$}}
\put(1692,818){\makebox(0,0){$+$}}
\put(1692,982){\makebox(0,0){$+$}}
\put(1262,1461){\makebox(0,0){$+$}}
\put(1262,339){\makebox(0,0){$+$}}
\put(201,1130){\makebox(0,0){$+$}}
\put(201,670){\makebox(0,0){$+$}}
\put(1337,1317){\circle*{12}}
\put(1337,483){\circle*{12}}
\put(1275,1312){\circle*{12}}
\put(1275,488){\circle*{12}}
\put(1056,1134){\circle*{12}}
\put(1056,666){\circle*{12}}
\put(1065,1085){\circle*{12}}
\put(1065,715){\circle*{12}}
\put(1050,1144){\circle*{12}}
\put(1050,656){\circle*{12}}
\put(1047,1337){\circle*{12}}
\put(1047,463){\circle*{12}}
\put(1020,1234){\circle*{12}}
\put(1020,566){\circle*{12}}
\put(994,1239){\circle*{12}}
\put(994,561){\circle*{12}}
\put(993,1304){\circle*{12}}
\put(993,496){\circle*{12}}
\put(1069,1045){\circle*{12}}
\put(1069,755){\circle*{12}}
\put(1056,1152){\circle*{12}}
\put(1056,648){\circle*{12}}
\put(1042,1294){\circle*{12}}
\put(1042,506){\circle*{12}}
\put(1063,1065){\circle*{12}}
\put(1063,735){\circle*{12}}
\put(1060,1374){\circle*{12}}
\put(1060,426){\circle*{12}}
\put(1060,1094){\circle*{12}}
\put(1060,706){\circle*{12}}
\put(1065,1111){\circle*{12}}
\put(1065,689){\circle*{12}}
\put(1049,1168){\circle*{12}}
\put(1049,632){\circle*{12}}
\put(748,1011){\circle*{12}}
\put(748,789){\circle*{12}}
\put(838,1194){\circle*{12}}
\put(838,606){\circle*{12}}
\put(897,1240){\circle*{12}}
\put(897,560){\circle*{12}}
\put(935,1234){\circle*{12}}
\put(935,566){\circle*{12}}
\put(951,1287){\circle*{12}}
\put(951,513){\circle*{12}}
\put(956,1302){\circle*{12}}
\put(956,498){\circle*{12}}
\put(749,1027){\circle*{12}}
\put(749,773){\circle*{12}}
\put(831,1187){\circle*{12}}
\put(831,613){\circle*{12}}
\put(888,1231){\circle*{12}}
\put(888,569){\circle*{12}}
\put(924,1274){\circle*{12}}
\put(924,526){\circle*{12}}
\put(936,1331){\circle*{12}}
\put(936,469){\circle*{12}}
\put(739,1021){\circle*{12}}
\put(739,779){\circle*{12}}
\put(818,1170){\circle*{12}}
\put(818,630){\circle*{12}}
\put(765,1083){\circle*{12}}
\put(765,717){\circle*{12}}
\put(846,1201){\circle*{12}}
\put(846,599){\circle*{12}}
\put(897,1248){\circle*{12}}
\put(897,552){\circle*{12}}
\put(937,1276){\circle*{12}}
\put(937,524){\circle*{12}}
\put(1079,1023){\circle*{12}}
\put(1079,777){\circle*{12}}
\put(1066,1209){\circle*{12}}
\put(1066,591){\circle*{12}}
\put(1059,1235){\circle*{12}}
\put(1059,565){\circle*{12}}
\put(1064,1167){\circle*{12}}
\put(1064,633){\circle*{12}}
\put(1062,1176){\circle*{12}}
\put(1062,624){\circle*{12}}
\put(1079,1078){\circle*{12}}
\put(1079,722){\circle*{12}}
\put(1083,1242){\circle*{12}}
\put(1083,558){\circle*{12}}
\put(1079,1344){\circle*{12}}
\put(1079,456){\circle*{12}}
\put(1062,1117){\circle*{12}}
\put(1062,683){\circle*{12}}
\put(1046,1246){\circle*{12}}
\put(1046,554){\circle*{12}}
\put(723,955){\circle*{12}}
\put(723,845){\circle*{12}}
\put(802,1152){\circle*{12}}
\put(802,648){\circle*{12}}
\put(862,1214){\circle*{12}}
\put(862,586){\circle*{12}}
\put(754,1056){\circle*{12}}
\put(754,744){\circle*{12}}
\put(798,1145){\circle*{12}}
\put(798,655){\circle*{12}}
\put(854,1208){\circle*{12}}
\put(854,592){\circle*{12}}
\put(748,1001){\circle*{12}}
\put(748,799){\circle*{12}}
\put(812,1163){\circle*{12}}
\put(812,637){\circle*{12}}
\put(870,1219){\circle*{12}}
\put(870,581){\circle*{12}}
\put(600,1063){\circle*{12}}
\put(600,737){\circle*{12}}
\put(723,976){\circle*{12}}
\put(723,824){\circle*{12}}
\put(792,1137){\circle*{12}}
\put(792,663){\circle*{12}}
\put(758,1063){\circle*{12}}
\put(758,737){\circle*{12}}
\put(825,1179){\circle*{12}}
\put(825,621){\circle*{12}}
\put(1096,982){\circle*{12}}
\put(1096,818){\circle*{12}}
\put(1083,1069){\circle*{12}}
\put(1083,731){\circle*{12}}
\put(1083,1274){\circle*{12}}
\put(1083,526){\circle*{12}}
\put(1073,1064){\circle*{12}}
\put(1073,736){\circle*{12}}
\put(1073,1093){\circle*{12}}
\put(1073,707){\circle*{12}}
\put(1083,1044){\circle*{12}}
\put(1083,756){\circle*{12}}
\put(1095,1012){\circle*{12}}
\put(1095,788){\circle*{12}}
\put(1163,1053){\circle*{12}}
\put(1163,747){\circle*{12}}
\put(1156,1312){\circle*{12}}
\put(1156,488){\circle*{12}}
\put(1170,1081){\circle*{12}}
\put(1170,719){\circle*{12}}
\put(1161,1289){\circle*{12}}
\put(1161,511){\circle*{12}}
\put(1099,1004){\circle*{12}}
\put(1099,796){\circle*{12}}
\put(1090,1026){\circle*{12}}
\put(1090,774){\circle*{12}}
\put(676,983){\circle*{12}}
\put(676,817){\circle*{12}}
\put(763,1076){\circle*{12}}
\put(763,724){\circle*{12}}
\put(786,1129){\circle*{12}}
\put(786,671){\circle*{12}}
\put(760,1069){\circle*{12}}
\put(760,731){\circle*{12}}
\put(778,1115){\circle*{12}}
\put(778,685){\circle*{12}}
\put(753,1048){\circle*{12}}
\put(753,752){\circle*{12}}
\put(806,1157){\circle*{12}}
\put(806,643){\circle*{12}}
\put(878,1225){\circle*{12}}
\put(878,575){\circle*{12}}
\put(935,1230){\circle*{12}}
\put(935,570){\circle*{12}}
\put(971,1300){\circle*{12}}
\put(971,500){\circle*{12}}
\put(987,1344){\circle*{12}}
\put(987,456){\circle*{12}}
\put(995,1338){\circle*{12}}
\put(995,462){\circle*{12}}
\put(1158,1028){\circle*{12}}
\put(1158,772){\circle*{12}}
\put(1187,1172){\circle*{12}}
\put(1187,628){\circle*{12}}
\put(1179,1341){\circle*{12}}
\put(1179,459){\circle*{12}}
\put(1190,1209){\circle*{12}}
\put(1190,591){\circle*{12}}
\put(1178,1321){\circle*{12}}
\put(1178,479){\circle*{12}}
\put(1180,1141){\circle*{12}}
\put(1180,659){\circle*{12}}
\put(1175,1109){\circle*{12}}
\put(1175,691){\circle*{12}}
\put(1183,1094){\circle*{12}}
\put(1183,706){\circle*{12}}
\put(1187,1123){\circle*{12}}
\put(1187,677){\circle*{12}}
\put(1086,1002){\circle*{12}}
\put(1086,798){\circle*{12}}
\put(498,961){\circle*{12}}
\put(498,839){\circle*{12}}
\put(566,1090){\circle*{12}}
\put(566,710){\circle*{12}}
\put(581,1111){\circle*{12}}
\put(581,689){\circle*{12}}
\put(589,1138){\circle*{12}}
\put(589,662){\circle*{12}}
\put(602,1147){\circle*{12}}
\put(602,653){\circle*{12}}
\put(598,1120){\circle*{12}}
\put(598,680){\circle*{12}}
\put(593,1087){\circle*{12}}
\put(593,713){\circle*{12}}
\put(507,983){\circle*{12}}
\put(507,817){\circle*{12}}
\put(552,1077){\circle*{12}}
\put(552,723){\circle*{12}}
\put(574,1096){\circle*{12}}
\put(574,704){\circle*{12}}
\put(551,999){\circle*{12}}
\put(551,801){\circle*{12}}
\put(491,926){\circle*{12}}
\put(491,874){\circle*{12}}
\put(539,1159){\circle*{12}}
\put(539,641){\circle*{12}}
\put(499,974){\circle*{12}}
\put(499,826){\circle*{12}}
\put(529,1053){\circle*{12}}
\put(529,747){\circle*{12}}
\put(552,1000){\circle*{12}}
\put(552,800){\circle*{12}}
\put(579,1144){\circle*{12}}
\put(579,656){\circle*{12}}
\put(616,1078){\circle*{12}}
\put(616,722){\circle*{12}}
\put(730,1000){\circle*{12}}
\put(730,800){\circle*{12}}
\put(768,1091){\circle*{12}}
\put(768,709){\circle*{12}}
\put(1105,962){\circle*{12}}
\put(1105,838){\circle*{12}}
\put(1100,1357){\circle*{12}}
\put(1100,443){\circle*{12}}
\put(1121,1269){\circle*{12}}
\put(1121,531){\circle*{12}}
\put(1225,1184){\circle*{12}}
\put(1225,616){\circle*{12}}
\put(1205,1287){\circle*{12}}
\put(1205,513){\circle*{12}}
\put(1185,1309){\circle*{12}}
\put(1185,491){\circle*{12}}
\put(1201,1205){\circle*{12}}
\put(1201,595){\circle*{12}}
\put(1196,1250){\circle*{12}}
\put(1196,550){\circle*{12}}
\put(1199,1162){\circle*{12}}
\put(1199,638){\circle*{12}}
\put(1214,1166){\circle*{12}}
\put(1214,634){\circle*{12}}
\put(1206,1271){\circle*{12}}
\put(1206,529){\circle*{12}}
\put(1179,1031){\circle*{12}}
\put(1179,769){\circle*{12}}
\put(490,909){\circle*{12}}
\put(490,891){\circle*{12}}
\put(511,1210){\circle*{12}}
\put(511,590){\circle*{12}}
\put(541,1102){\circle*{12}}
\put(541,698){\circle*{12}}
\put(558,1172){\circle*{12}}
\put(558,628){\circle*{12}}
\put(506,1006){\circle*{12}}
\put(506,794){\circle*{12}}
\put(493,943){\circle*{12}}
\put(493,857){\circle*{12}}
\put(516,1046){\circle*{12}}
\put(516,754){\circle*{12}}
\put(541,1065){\circle*{12}}
\put(541,735){\circle*{12}}
\put(447,881){\circle*{12}}
\put(447,919){\circle*{12}}
\put(461,1363){\circle*{12}}
\put(461,437){\circle*{12}}
\put(464,1371){\circle*{12}}
\put(464,429){\circle*{12}}
\put(511,1026){\circle*{12}}
\put(511,774){\circle*{12}}
\put(601,1053){\circle*{12}}
\put(601,747){\circle*{12}}
\put(771,1099){\circle*{12}}
\put(771,701){\circle*{12}}
\put(911,1250){\circle*{12}}
\put(911,550){\circle*{12}}
\put(1063,1188){\circle*{12}}
\put(1063,612){\circle*{12}}
\put(1198,1106){\circle*{12}}
\put(1198,694){\circle*{12}}
\put(1216,1238){\circle*{12}}
\put(1216,562){\circle*{12}}
\put(1230,1294){\circle*{12}}
\put(1230,506){\circle*{12}}
\put(1231,1376){\circle*{12}}
\put(1231,424){\circle*{12}}
\put(1218,1212){\circle*{12}}
\put(1218,588){\circle*{12}}
\put(1180,1069){\circle*{12}}
\put(1180,731){\circle*{12}}
\put(1110,1325){\circle*{12}}
\put(1110,475){\circle*{12}}
\put(1112,947){\circle*{12}}
\put(1112,853){\circle*{12}}
\put(985,1241){\circle*{12}}
\put(985,559){\circle*{12}}
\put(1021,1319){\circle*{12}}
\put(1021,481){\circle*{12}}
\put(1154,1001){\circle*{12}}
\put(1154,799){\circle*{12}}
\put(1132,1335){\circle*{12}}
\put(1132,465){\circle*{12}}
\put(1153,973){\circle*{12}}
\put(1153,827){\circle*{12}}
\put(1202,1125){\circle*{12}}
\put(1202,675){\circle*{12}}
\put(1220,1179){\circle*{12}}
\put(1220,621){\circle*{12}}
\put(1259,1298){\circle*{12}}
\put(1259,502){\circle*{12}}
\put(1236,1345){\circle*{12}}
\put(1236,455){\circle*{12}}
\put(677,982){\circle*{12}}
\put(677,818){\circle*{12}}
\put(746,1039){\circle*{12}}
\put(746,761){\circle*{12}}
\put(775,1106){\circle*{12}}
\put(775,694){\circle*{12}}
\put(925,1263){\circle*{12}}
\put(925,537){\circle*{12}}
\put(1000,1415){\circle*{12}}
\put(1000,385){\circle*{12}}
\put(1021,1417){\circle*{12}}
\put(1021,383){\circle*{12}}
\put(1177,1385){\circle*{12}}
\put(1177,415){\circle*{12}}
\put(1179,1022){\circle*{12}}
\put(1179,778){\circle*{12}}
\put(1179,1046){\circle*{12}}
\put(1179,754){\circle*{12}}
\put(460,978){\circle*{12}}
\put(460,822){\circle*{12}}
\put(458,1214){\circle*{12}}
\put(458,586){\circle*{12}}
\put(434,1151){\circle*{12}}
\put(434,649){\circle*{12}}
\put(451,1355){\circle*{12}}
\put(451,445){\circle*{12}}
\put(425,1122){\circle*{12}}
\put(425,678){\circle*{12}}
\put(455,1346){\circle*{12}}
\put(455,454){\circle*{12}}
\put(434,1075){\circle*{12}}
\put(434,725){\circle*{12}}
\put(437,1181){\circle*{12}}
\put(437,619){\circle*{12}}
\put(437,1136){\circle*{12}}
\put(437,664){\circle*{12}}
\put(440,1045){\circle*{12}}
\put(440,755){\circle*{12}}
\put(441,1016){\circle*{12}}
\put(441,784){\circle*{12}}
\put(445,976){\circle*{12}}
\put(445,824){\circle*{12}}
\put(502,1154){\circle*{12}}
\put(502,646){\circle*{12}}
\put(1156,846){\circle*{12}}
\put(1156,954){\circle*{12}}
\put(1289,1286){\circle*{12}}
\put(1289,514){\circle*{12}}
\put(1284,1368){\circle*{12}}
\put(1284,432){\circle*{12}}
\put(1319,1288){\circle*{12}}
\put(1319,512){\circle*{12}}
\put(1301,1310){\circle*{12}}
\put(1301,490){\circle*{12}}
\put(1279,1350){\circle*{12}}
\put(1279,450){\circle*{12}}
\put(1243,1318){\circle*{12}}
\put(1243,482){\circle*{12}}
\put(986,1411){\circle*{12}}
\put(986,389){\circle*{12}}
\put(971,1442){\circle*{12}}
\put(971,358){\circle*{12}}
\put(1040,1393){\circle*{12}}
\put(1040,407){\circle*{12}}
\put(925,1470){\circle*{12}}
\put(925,330){\circle*{12}}
\put(939,1328){\circle*{12}}
\put(939,472){\circle*{12}}
\put(951,1447){\circle*{12}}
\put(951,353){\circle*{12}}
\put(1198,1096){\circle*{12}}
\put(1198,704){\circle*{12}}
\put(1233,1172){\circle*{12}}
\put(1233,628){\circle*{12}}
\put(455,925){\circle*{12}}
\put(455,875){\circle*{12}}
\put(426,1177){\circle*{12}}
\put(426,623){\circle*{12}}
\put(451,1251){\circle*{12}}
\put(451,549){\circle*{12}}
\put(432,1100){\circle*{12}}
\put(432,700){\circle*{12}}
\put(439,1233){\circle*{12}}
\put(439,567){\circle*{12}}
\put(439,1002){\circle*{12}}
\put(439,798){\circle*{12}}
\put(422,1037){\circle*{12}}
\put(422,763){\circle*{12}}
\put(427,1205){\circle*{12}}
\put(427,595){\circle*{12}}
\put(444,905){\circle*{12}}
\put(444,895){\circle*{12}}
\put(781,1122){\circle*{12}}
\put(781,678){\circle*{12}}
\put(966,1466){\circle*{12}}
\put(966,334){\circle*{12}}
\put(1051,1408){\circle*{12}}
\put(1051,392){\circle*{12}}
\put(1078,1384){\circle*{12}}
\put(1078,416){\circle*{12}}
\put(1359,1295){\circle*{12}}
\put(1359,505){\circle*{12}}
\put(1338,1322){\circle*{12}}
\put(1338,478){\circle*{12}}
\put(1320,1311){\circle*{12}}
\put(1320,489){\circle*{12}}
\put(1324,1265){\circle*{12}}
\put(1324,535){\circle*{12}}
\put(1398,1285){\circle*{12}}
\put(1398,515){\circle*{12}}
\put(1380,1290){\circle*{12}}
\put(1380,510){\circle*{12}}
\put(1370,1292){\circle*{12}}
\put(1370,508){\circle*{12}}
\put(1343,1299){\circle*{12}}
\put(1343,501){\circle*{12}}
\put(1234,1204){\circle*{12}}
\put(1234,596){\circle*{12}}
\put(1201,1410){\circle*{12}}
\put(1201,390){\circle*{12}}
\put(1231,1396){\circle*{12}}
\put(1231,404){\circle*{12}}
\put(912,1479){\circle*{12}}
\put(912,321){\circle*{12}}
\put(896,1492){\circle*{12}}
\put(896,308){\circle*{12}}
\put(863,1504){\circle*{12}}
\put(863,296){\circle*{12}}
\put(1180,1465){\circle*{12}}
\put(1180,335){\circle*{12}}
\put(427,1064){\circle*{12}}
\put(427,736){\circle*{12}}
\put(431,1314){\circle*{12}}
\put(431,486){\circle*{12}}
\put(407,1129){\circle*{12}}
\put(407,671){\circle*{12}}
\put(449,1334){\circle*{12}}
\put(449,466){\circle*{12}}
\put(514,1392){\circle*{12}}
\put(514,408){\circle*{12}}
\put(441,951){\circle*{12}}
\put(441,849){\circle*{12}}
\put(566,1470){\circle*{12}}
\put(566,330){\circle*{12}}
\put(1489,1245){\circle*{12}}
\put(1489,555){\circle*{12}}
\put(1440,1270){\circle*{12}}
\put(1440,530){\circle*{12}}
\put(1416,1280){\circle*{12}}
\put(1416,520){\circle*{12}}
\put(1408,1283){\circle*{12}}
\put(1408,517){\circle*{12}}
\put(1389,1288){\circle*{12}}
\put(1389,512){\circle*{12}}
\put(1481,1249){\circle*{12}}
\put(1481,551){\circle*{12}}
\put(1447,1266){\circle*{12}}
\put(1447,534){\circle*{12}}
\put(1424,1277){\circle*{12}}
\put(1424,523){\circle*{12}}
\put(1519,1227){\circle*{12}}
\put(1519,573){\circle*{12}}
\put(1465,1257){\circle*{12}}
\put(1465,543){\circle*{12}}
\put(1432,1274){\circle*{12}}
\put(1432,526){\circle*{12}}
\put(1555,1201){\circle*{12}}
\put(1555,599){\circle*{12}}
\put(1497,1241){\circle*{12}}
\put(1497,559){\circle*{12}}
\put(1456,1261){\circle*{12}}
\put(1456,539){\circle*{12}}
\put(402,1114){\circle*{12}}
\put(402,686){\circle*{12}}
\put(441,1324){\circle*{12}}
\put(441,476){\circle*{12}}
\put(521,1412){\circle*{12}}
\put(521,388){\circle*{12}}
\put(564,1462){\circle*{12}}
\put(564,338){\circle*{12}}
\put(383,1237){\circle*{12}}
\put(383,563){\circle*{12}}
\put(423,1303){\circle*{12}}
\put(423,497){\circle*{12}}
\put(872,1505){\circle*{12}}
\put(872,295){\circle*{12}}
\put(918,1489){\circle*{12}}
\put(918,311){\circle*{12}}
\put(946,1464){\circle*{12}}
\put(946,336){\circle*{12}}
\put(858,1506){\circle*{12}}
\put(858,294){\circle*{12}}
\put(846,1514){\circle*{12}}
\put(846,286){\circle*{12}}
\put(510,1410){\circle*{12}}
\put(510,390){\circle*{12}}
\put(596,1516){\circle*{12}}
\put(596,284){\circle*{12}}
\put(416,1055){\circle*{12}}
\put(416,745){\circle*{12}}
\put(412,1169){\circle*{12}}
\put(412,631){\circle*{12}}
\put(397,1103){\circle*{12}}
\put(397,697){\circle*{12}}
\put(414,1255){\circle*{12}}
\put(414,545){\circle*{12}}
\put(1010,1455){\circle*{12}}
\put(1010,345){\circle*{12}}
\put(1018,1496){\circle*{12}}
\put(1018,304){\circle*{12}}
\put(969,1499){\circle*{12}}
\put(969,301){\circle*{12}}
\put(438,940){\circle*{12}}
\put(438,860){\circle*{12}}
\put(501,1397){\circle*{12}}
\put(501,403){\circle*{12}}
\put(808,1542){\circle*{12}}
\put(808,258){\circle*{12}}
\put(1182,1429){\circle*{12}}
\put(1182,371){\circle*{12}}
\put(1192,1475){\circle*{12}}
\put(1192,325){\circle*{12}}
\put(1170,1499){\circle*{12}}
\put(1170,301){\circle*{12}}
\put(545,1446){\circle*{12}}
\put(545,354){\circle*{12}}
\put(1023,1496){\circle*{12}}
\put(1023,304){\circle*{12}}
\put(1610,1151){\circle*{12}}
\put(1610,649){\circle*{12}}
\put(1534,1217){\circle*{12}}
\put(1534,583){\circle*{12}}
\put(1504,1236){\circle*{12}}
\put(1504,564){\circle*{12}}
\put(1473,1253){\circle*{12}}
\put(1473,547){\circle*{12}}
\put(1582,1178){\circle*{12}}
\put(1582,622){\circle*{12}}
\put(1527,1222){\circle*{12}}
\put(1527,578){\circle*{12}}
\put(1569,1190){\circle*{12}}
\put(1569,610){\circle*{12}}
\put(1512,1232){\circle*{12}}
\put(1512,568){\circle*{12}}
\put(1636,1119){\circle*{12}}
\put(1636,681){\circle*{12}}
\put(1548,1206){\circle*{12}}
\put(1548,594){\circle*{12}}
\put(1642,1110){\circle*{12}}
\put(1642,690){\circle*{12}}
\put(1541,1212){\circle*{12}}
\put(1541,588){\circle*{12}}
\put(419,1013){\circle*{12}}
\put(419,787){\circle*{12}}
\put(403,1254){\circle*{12}}
\put(403,546){\circle*{12}}
\put(415,1281){\circle*{12}}
\put(415,519){\circle*{12}}
\put(609,1540){\circle*{12}}
\put(609,260){\circle*{12}}
\put(618,1549){\circle*{12}}
\put(618,251){\circle*{12}}
\put(628,1557){\circle*{12}}
\put(628,243){\circle*{12}}
\put(668,1593){\circle*{12}}
\put(668,207){\circle*{12}}
\put(657,1584){\circle*{12}}
\put(657,216){\circle*{12}}
\put(1261,338){\circle*{12}}
\put(1261,1462){\circle*{12}}
\put(1214,1447){\circle*{12}}
\put(1214,353){\circle*{12}}
\put(1157,1532){\circle*{12}}
\put(1157,268){\circle*{12}}
\put(807,1544){\circle*{12}}
\put(807,256){\circle*{12}}
\put(361,1119){\circle*{12}}
\put(361,681){\circle*{12}}
\put(407,1301){\circle*{12}}
\put(407,499){\circle*{12}}
\put(303,1194){\circle*{12}}
\put(303,606){\circle*{12}}
\put(346,1209){\circle*{12}}
\put(346,591){\circle*{12}}
\put(403,1303){\circle*{12}}
\put(403,497){\circle*{12}}
\put(581,1500){\circle*{12}}
\put(581,300){\circle*{12}}
\put(823,1557){\circle*{12}}
\put(823,243){\circle*{12}}
\put(887,1511){\circle*{12}}
\put(887,289){\circle*{12}}
\put(734,1631){\circle*{12}}
\put(734,169){\circle*{12}}
\put(431,981){\circle*{12}}
\put(431,819){\circle*{12}}
\put(477,1377){\circle*{12}}
\put(477,423){\circle*{12}}
\put(486,1391){\circle*{12}}
\put(486,409){\circle*{12}}
\put(1173,1519){\circle*{12}}
\put(1173,281){\circle*{12}}
\put(1137,1559){\circle*{12}}
\put(1137,241){\circle*{12}}
\put(1159,1551){\circle*{12}}
\put(1159,249){\circle*{12}}
\put(1258,1438){\circle*{12}}
\put(1258,362){\circle*{12}}
\put(1677,1047){\circle*{12}}
\put(1677,753){\circle*{12}}
\put(1596,1165){\circle*{12}}
\put(1596,635){\circle*{12}}
\put(1562,1195){\circle*{12}}
\put(1562,605){\circle*{12}}
\put(1623,1136){\circle*{12}}
\put(1623,664){\circle*{12}}
\put(1589,1172){\circle*{12}}
\put(1589,628){\circle*{12}}
\put(1653,1091){\circle*{12}}
\put(1653,709){\circle*{12}}
\put(1603,1158){\circle*{12}}
\put(1603,642){\circle*{12}}
\put(1658,1081){\circle*{12}}
\put(1658,719){\circle*{12}}
\put(1669,1061){\circle*{12}}
\put(1669,739){\circle*{12}}
\put(1693,994){\circle*{12}}
\put(1693,806){\circle*{12}}
\put(240,1156){\circle*{12}}
\put(240,644){\circle*{12}}
\put(270,1184){\circle*{12}}
\put(270,616){\circle*{12}}
\put(324,1199){\circle*{12}}
\put(324,601){\circle*{12}}
\put(366,1223){\circle*{12}}
\put(366,577){\circle*{12}}
\put(327,1114){\circle*{12}}
\put(327,686){\circle*{12}}
\put(229,668){\circle*{12}}
\put(229,1132){\circle*{12}}
\put(287,1191){\circle*{12}}
\put(287,609){\circle*{12}}
\put(680,1602){\circle*{12}}
\put(680,198){\circle*{12}}
\put(719,1627){\circle*{12}}
\put(719,173){\circle*{12}}
\put(797,1587){\circle*{12}}
\put(797,213){\circle*{12}}
\put(849,1539){\circle*{12}}
\put(849,261){\circle*{12}}
\put(637,1566){\circle*{12}}
\put(637,234){\circle*{12}}
\put(591,1515){\circle*{12}}
\put(591,285){\circle*{12}}
\put(527,1432){\circle*{12}}
\put(527,368){\circle*{12}}
\put(796,1589){\circle*{12}}
\put(796,211){\circle*{12}}
\put(1145,1579){\circle*{12}}
\put(1145,221){\circle*{12}}
\put(1115,1580){\circle*{12}}
\put(1115,220){\circle*{12}}
\put(1089,1599){\circle*{12}}
\put(1089,201){\circle*{12}}
\put(1189,1527){\circle*{12}}
\put(1189,273){\circle*{12}}
\put(1630,1127){\circle*{12}}
\put(1630,673){\circle*{12}}
\put(1695,949){\circle*{12}}
\put(1695,851){\circle*{12}}
\put(1575,1184){\circle*{12}}
\put(1575,616){\circle*{12}}
\put(221,1120){\circle*{12}}
\put(221,680){\circle*{12}}
\put(255,1177){\circle*{12}}
\put(255,623){\circle*{12}}
\put(192,1120){\circle*{12}}
\put(192,680){\circle*{12}}
\put(217,1179){\circle*{12}}
\put(217,621){\circle*{12}}
\put(1256,1458){\circle*{12}}
\put(1256,342){\circle*{12}}
\put(962,1497){\circle*{12}}
\put(962,303){\circle*{12}}
\put(693,1611){\circle*{12}}
\put(693,189){\circle*{12}}
\put(601,1530){\circle*{12}}
\put(601,270){\circle*{12}}
\put(647,1576){\circle*{12}}
\put(647,224){\circle*{12}}
\put(1132,1599){\circle*{12}}
\put(1132,201){\circle*{12}}
\put(1198,1508){\circle*{12}}
\put(1198,292){\circle*{12}}
\put(1689,787){\circle*{12}}
\put(1689,1013){\circle*{12}}
\put(1617,1144){\circle*{12}}
\put(1617,656){\circle*{12}}
\put(1696,973){\circle*{12}}
\put(1696,827){\circle*{12}}
\put(555,344){\circle*{12}}
\put(555,1456){\circle*{12}}
\put(1031,1618){\circle*{12}}
\put(1031,182){\circle*{12}}
\put(760,1633){\circle*{12}}
\put(760,167){\circle*{12}}
\put(191,1036){\circle*{12}}
\put(191,764){\circle*{12}}
\put(236,1180){\circle*{12}}
\put(236,620){\circle*{12}}
\put(192,1083){\circle*{12}}
\put(192,717){\circle*{12}}
\put(1663,730){\circle*{12}}
\put(1663,1070){\circle*{12}}
\put(1684,1032){\circle*{12}}
\put(1684,768){\circle*{12}}
\put(1059,1611){\circle*{12}}
\put(1059,189){\circle*{12}}
\put(706,1619){\circle*{12}}
\put(706,181){\circle*{12}}
\put(193,1000){\circle*{12}}
\put(193,800){\circle*{12}}
\put(1124,1611){\circle*{12}}
\put(1124,189){\circle*{12}}
\put(1647,1100){\circle*{12}}
\put(1647,700){\circle*{12}}
\put(1252,1425){\circle*{12}}
\put(1252,375){\circle*{12}}
\put(576,1484){\circle*{12}}
\put(576,316){\circle*{12}}
\put(209,645){\circle*{12}}
\put(209,1155){\circle*{12}}
\put(791,1639){\circle*{12}}
\put(791,161){\circle*{12}}
\end{picture}

%% file: martineta.tex.vcleaned.bbl
\providecommand{\bysame}{\leavevmode\hbox to3em{\hrulefill}\thinspace}
\providecommand{\MR}{\relax\ifhmode\unskip\space\fi MR }
\providecommand{\MRhref}[2]{%
  \href{http://www.ams.org/mathscinet-getitem?mr=#1}{#2}
}
\providecommand{\href}[2]{#2}
\begin{thebibliography}{10}

\bibitem{AshkinTeller43}
J~Ashkin and E~Teller, Phys Rev \textbf{64} (1943), 178--184.

\bibitem{Baxter82}
R~J Baxter, \emph{Exactly solved models in statistical mechanics}, Academic
  Press, New York, 1982.

\bibitem{BazhanovBaxter94}
V~V Bazhanov and R~J Baxter, \emph{Star triangle relation for a 3-dimensional
  model}, J Stat Phys \textbf{71} (1994), 839--864.

\bibitem{Brauer37}
R~Brauer, \emph{On algebras which are connected with the semi--simple
  continuous groups}, Annals of Mathematics \textbf{38} (1937), 854--872.

\bibitem{ClineParshallScott88}
E~Cline, B~Parshall, and L~Scott, \emph{Finite-dimensional algebras and highest
  weight categories}, J. reine angew. Math. \textbf{391} (1988), 85--99.

\bibitem{ClineParshallScott99}
\bysame, \emph{Generic and $q$--rational representation theory}, Publ RIMS
  \textbf{35} (1999), 31--90.

\bibitem{DasmahapatraMartin96}
S~Dasmahapatra and P~P Martin, \emph{On the algebraic approach to cubic lattice
  {P}otts models}, Journal of Physics A \textbf{29} (1996), 263--278.

\bibitem{DlabRingel89B}
V.~Dlab and C.M. Ringel, \emph{A construction for quasi-hereditary algebras},
  Compos. Math. \textbf{70} (1989), 155--175.

\bibitem{Donkin93}
S.~Donkin, \emph{On tilting modules for algebraic groups}, Math. Z.
  \textbf{212} (1993), 39--60.

\bibitem{Donkin98agc}
\bysame, \emph{{The $q$-Schur} algebra}, LMS Lecture Notes Series, vol. 253,
  Cambridge University Press, 1998.

\bibitem{Donkin01}
\bysame, \emph{Symmetric and exterior powers, linear source modules and
  representations of {S}chur superalgebras}, Proc LMS \textbf{83} (2001),
  647--680.

\bibitem{EinhornSavitRabinovici}
M~B Einhorn, R~Savit, and E~Rabinovici, \emph{A physical picture for the phase
  transition in ${Z}_{N}$ symmetric models}, Nucl Phys B \textbf{170 [FS1]}
  (1980), 16.

\bibitem{ElitzurPearsonShigemitsu}
S~Elitzur, R~B Pearson, and J~Shigemitsu, \emph{The phase structure of discrete
  abelian spin and gauge systems}, Phys Rev D \textbf{19} (1979), 3698.

\bibitem{Fan72}
C~Fan, Phys Lett \textbf{39A} (1972), 136--136.

\bibitem{GoodmandelaHarpeJones89}
F~M Goodman, P~de~la Harpe, and V~F~R Jones, \emph{Coxeter--{D}ynkin diagrams
  and towers of algebras}, Math Sci Research Inst Publ, Springer--Verlag,
  Berlin, 1989.

\bibitem{Green80}
J~A Green, \emph{Polynomial representations of ${GL}_n$}, Springer-Verlag,
  Berlin, 1980.

\bibitem{GAP99}
The~GAP Group, \emph{{GAP} --- {G}roups, {A}lgorithms, and {P}rogramming,
  version 4.2}, Aachen, St. Andrews, 1999.

\bibitem{HalversonRam95}
T~Halverson and A~Ram, \emph{Characters of algebras containing a {J}ones basic
  construction...}, Adv Math \textbf{116} (1995), 263--321.

\bibitem{HanlonWales94}
P~Hanlon and D~Wales, \emph{A tower construction for the radical in {B}rauer's
  centralizer algebras}, J Algebra \textbf{164} (1994), 773--830.

\bibitem{HenkeKoenig01}
A~Henke and S~Koenig, \emph{Relating polynomial $gl(n)$--representations in
  different degrees}, J reine angew Math (to appear).

\bibitem{Jones94}
V~F~R Jones, in {S}ubfactors (Singapore), World Scientific, 1994.

\bibitem{Mattson93}
H~F~Mattson Jr., \emph{Discrete mathematics with applications}, Wiley,
  Singapore, 1993.

\bibitem{Macdonald95}
I~Macdonald, \emph{Symmetric functions and {H}all polynomials}, 2 ed., Oxford,
  1995.

\bibitem{Martin91}
P~P Martin, \emph{Potts models and related problems in statistical mechanics},
  World Scientific, Singapore, 1991.

\bibitem{Martin94}
\bysame, \emph{Temperley--{L}ieb algebras for non--planar statistical mechanics
  --- the partition algebra construction}, Journal Of Knot Theory and its
  Ramifications \textbf{3} (1994), no.~1, 51--82.

\bibitem{Martin96}
\bysame, \emph{The structure of the partition algebras}, J Algebra \textbf{183}
  (1996), 319--358.

\bibitem{Martin2000}
\bysame, \emph{The partition algebra and the {P}otts model transfer matrix
  spectrum in high dimensions}, J Phys A \textbf{32} (2000), 3669--3695.

\bibitem{MartinRollet98corrigendum}
P~P Martin and G~Rollet, \emph{The {P}otts model representation and a
  {R}obinson--{S}chensted correspondence for the partition algebra}, Compositio
  Math \textbf{112} (1998), 237--254. Corrigendum: multiple substitutions of
  $\implies$ for -- symbol!

\bibitem{MartinSaleur94b}
P~P Martin and H~Saleur, \emph{Algebras in higher dimensional statistical
  mechanics --- the exceptional partition algebras}, Lett. Math. Phys. (1994),
  no.~30, 179--185.

\bibitem{MartinWoodcock98}
P~P Martin and D~Woodcock, \emph{The partition algebras and a new deformation
  of the {S}chur algebras}, J Algebra \textbf{203} (1998), 91--124.

\bibitem{MatveevShrock96b}
V~Matveev and R~Shrock, \emph{Complex-temperature phase diagram of the 1{D}
  ${Z}(6)$ clock model and its connection with higher-dimensional models},
  PHYSICS LETTERS A \textbf{221} (1996), 343--349.

\bibitem{Pathria72}
R~K Pathria, \emph{Statistical {M}echanics}, Pergamon, 1972.

\bibitem{ReggeZecchina2000}
T~Regge and R~Zecchina, \emph{Combinatorial and topological approach to the
  3{D} {I}sing model}, JOURNAL OF PHYSICS A-MATHEMATICAL AND GENERAL
  \textbf{33} (2000), 741--761.

\bibitem{VershikYakubovich98}
A~M Vershik and Yu~V Yakubovich, \emph{On continuous partition lattice},
  preprint (St. Petersburg) (1998).

\bibitem{Weyl46}
H~Weyl, \emph{Classical groups}, Princeton, Princeton, 1946, ch.III-V.

\bibitem{Xi00}
C.~Xi, \emph{Partition algebras are cellular}, Compositio Mathematica
  \textbf{119} (1999), 107--118.

\end{thebibliography}
